# The escaping set in transcendental dynamics

Walter Bergweiler and Lasse Rempe

*To our Ukrainian friends and colleagues*

**Abstract**

The escaping set of an entire function consists of the points in the complex plane that tend to infinity under iteration. This set plays a central role in the dynamics of transcendental entire functions. The goal of this survey is to explain this role, to summarise some of the main results in the area, and to identify a number of open questions.

# Contents

# 1 Introduction

Let $f\colon \mathbb{C} \to \mathbb{C}$ be an *entire function*; that is, a function holomorphic in the entire complex plane. We may think of $f$ as the evolution rule of a system that changes in discrete steps of time; that is, given a state $z \in \mathbb{C}$, the next state is given by $f(z)$. Then the state of the system after $n$ steps, where $n \in \mathbb{N}$, is described by the $n$-fold composition of $f$ with itself, denoted by $f^n$ and called the $n$-th *iterate* of $f$. More formally, $f^0$ is the identity, $f^1 = f$ and $f^n = f \circ f^{n-1}$ for $n \geq 2$.

The subject of this article is the *escaping set*

$$I(f) := \left\{ z \in \mathbb{C} \colon \lim_{n \to \infty} |f^n(z)| = \infty \right\}, \tag{1.1}$$

which plays a crucial role in the dynamics of entire functions. It was introduced by Eremenko [113] in 1989, who gave the first systematic study of it. The name "escaping set" seems to have been used first by Devaney [103, p. 225] for this set. (Devaney's paper appeared in the same conference proceedings where Eremenko's paper appeared, and Devaney refers to Eremenko's paper.)

But points escaping to infinity had been considered much earlier: In 1926, in what can be considered as the birth of transcendental dynamics, Fatou published an article [132] that studied the iteration of transcendental entire functions, and noticed that in certain cases the escaping set contains curves tending to infinity. Today the escaping set, despite or perhaps because of its deceptively simple definition (1.1), plays a major role in transcendental dynamics. The goal of this survey is to explain this role, to summarise some of the main results in the area, and to identify a number (62, to be precise) of open questions.

The main object studied in complex dynamics (the iteration of a function $f$ of one complex variable) is the Julia set $J(f)$. To exclude certain exceptional cases, one assumes in this theory that $f$ is neither constant nor a polynomial of degree 1. We will in fact mainly be concerned with the case where $f$ is transcendental; that is, $f$ is not a polynomial. Intuitively speaking, the Julia set is the set of points at which the eventual behaviour is unstable under small perturbations. More formally, it is the set where the iterates do not form a normal family in the sense of Montel; see §3 below for a brief introduction to complex dynamics.

The sets $I(f)$ and $J(f)$ are closely related, since $J(f) = \partial I(f)$ (Theorem 4.3). In many important cases we also have $I(f) \subset J(f)$ and thus $\overline{I(f)} = J(f)$; see §6.1. Since the definition (1.1) of the escaping set is much simpler and more elementary than that of the Julia set, often results about the Julia set are proved using the escaping set. For example, the fact that $I(f) \neq \emptyset$ (Theorem 4.1) yields a comparatively simple proof of the crucial fact that $J(f) \neq \emptyset$ (Theorem 4.2). Other examples are results about Julia sets of permutable or semiconjugate entire functions (see §5.2 for a discussion) or results about the Lebesgue measure and Hausdorff dimension of Julia sets (see §9).

Another motivation for investigating the escaping set is the study of *wandering domains*: connected components of $\mathbb{C} \setminus J(f)$ that do not eventually become periodic. Wandering domains are a subject of considerable current interest in transcendental dy-namics; see for example [41]. Connected components of the interior of the escaping set are often wandering domains, and vice versa, so there is a close connection between the subjects. A famous theorem of Sullivan [274] states that rational functions do not have wandering domains. Earlier, Baker [13] had already shown that transcendental entire functions may have wandering domains. As we shall see in §2.4 and §6.3, Baker's example does indeed consist of escaping points, and a crucial role in Baker's argument is played by the rate at which these points escape. Functions with similar features play a key role in recent work on transcendental dynamics, for example in Bishop's construction [83] of a transcendental entire function with Julia set of dimension 1, or the work about structures called "spiders' webs" in escaping sets and Julia sets that will be discussed in §7.4.

In the case of a polynomial, the escaping set is an open set. Yoccoz pioneered the use of "puzzles" (see [190]), using curves in the escaping set to cut up the Julia set into "puzzle pieces". This allows a combinatorial study of the dynamics of polynomials and has been a major ingredient in breakthroughs in polynomial dynamics over the past several decades. For a transcendental entire function, the escaping set is never open, but nonetheless frequently contains structures, such as curves to infinity, that allow one to study the global dynamics of the function. Techniques similar to those developed for polynomials, using curves in the escaping set in order to introduce combinatorial partitions, have been employed successfully for certain transcendental entire functions, e.g. by Schleicher [253], Rempe und Schleicher [229] and Benini [40].

While it may seem natural that methods employed for polynomials may be modified to yield results for transcendental entire functions, it is surprising that the connection between polynomial and transcendental dynamics also works the other way round. Indeed, Dudko and Lyubich [112] have used techniques developed for escaping sets of transcendental functions to study polynomial dynamics.

We give a brief overview of the structure of this article. We begin with some motivating examples in §2. This section can be read without any prior knowledge of complex dynamics; we have endeavoured to keep the discussion of the examples as elementary as possible, assuming only basic knowledge of complex analysis and topology. In order to keep the survey as self-contained as possible, we give a brief introduction to complex dynamics in §3, recalling the main results that we will use in the discussions that follows.

In §4, we begin the study of the escaping set $I(f)$ in earnest. In particular, we give two proofs of the fundamental fact that $I(f)$ is non-empty. In §5 we study the rate at which points in $I(f)$ escape to infinity. Connected components of the interior of $I(f)$ will be studied §6 while §7 is devoted to topological properties of $I(f)$. In §8 we will see that, for certain families of entire functions, the set $I(f)$ has certain "rigidity" properties that mean that certain structures are preserved across all maps in a given parameter space. In particular, this gives (in §8.8) explicit examples of the principle, mentioned above, that the dynamics of a function on its Julia set can be understood using the escaping set. The measure and dimension of $I(f)$ are discussed in §9. The final section, §10, briefly discusses the role of the escaping set beyond the setting of transcendental entire functions, such as quasiregular and meromorphic functions.

Throughout the survey, we attempt to keep the discussion as self-contained as possible, but do not hesitate to refer to the relevant literature for background when we feel that this benefits the discussion. Likewise, we include proofs, or ideas of proofs, of key results where this is feasible, but refer to the original articles when the ideas would require too much further discussion. In some cases (particularly in §8) we present a perspective on existing results that is not currently found in the literature. In such cases, we provide the proofs for completeness, using the terminology of the relevant published articles; these proofs are included primarily for the experts.

We conclude this introduction by fixing some terminology used throughout this paper. For $a \in \mathbb{C}$ and $r > 0$, we use the notation $D(a,r) \coloneqq \{z\colon |z - a| < r\}$ for the open disc of radius $r$ around $a$. We also put $\mathbb{D} \coloneqq D(0,1)$. For $0 < r_1 < r_2$, we use the notation $\operatorname{ann}(r_1, r_2) \coloneqq \{z\colon r_1 < |z| < r_2\}$ for the open annulus centred at 0 with radii $r_1$ and $r_2$. The closed disc and the closed annulus are denoted by $\overline{D}(a,r)$ and $\overline{\operatorname{ann}}(r_1, r_2)$. For $t \in \mathbb{R}$, we put $\mathbb{H}_{>t} \coloneqq \{z\colon \operatorname{Re} z > t\}$. The half-planes $\mathbb{H}_{\geq t}$, $\mathbb{H}_{<t}$ and $\mathbb{H}_{\leq t}$ are defined analogously.

The closure of a set $A \subset \mathbb{C}$ is denoted by $\overline{A}$; all closures are taken in $\mathbb{C}$ unless explicitly stated otherwise. The interior of $A$ is denoted by $\operatorname{int}(A)$. A *Jordan curve* is a simple closed curve; i.e., the homeomorphic image of a circle. A *Jordan domain* is a simply connected domain bounded by a Jordan curve. Finally, we use the convention $\mathbb{N} = \{1, 2, 3, \dots\}$ and put $\mathbb{N}_0 \coloneqq \mathbb{N} \cup \{0\}$.

*Acknowledgment.* We thank Anke Pohl for encouraging us to write this article and for valuable comments on a preliminary version of it. We also thank Krzysztof Barański, Weiwei Cui, Alexandre Eremenko, Phil Rippon, Dierk Schleicher, David Sixsmith, Gwyneth Stallard, James Waterman and Yitming Zhang for helpful discussions and comments.

## 2  Motivating examples

We begin by investigating the escaping set of several specific entire functions using elementary means; the reader may wish to keep these in mind throughout the remainder of the article. A good starting point is given by two examples already investigated in Fatou's seminal article [132], whose escaping sets exhibit important properties.

### 2.1  A sine function

Fatou [132, §16] studied the family of entire functions

$$f\colon \mathbb{C} \to \mathbb{C}, \quad f(z) \coloneqq h \sin z + a,$$

where $0 < h < 1$ and $a \in \mathbb{R}$. It turns out that the qualitative behaviour of these functions is independent of the parameters, so for simplicity let us fix one particular choice: $h = 1/2$ and $a = 0$. Then

$$f(z) = \frac{\sin z}{2}. \tag{2.1}$$

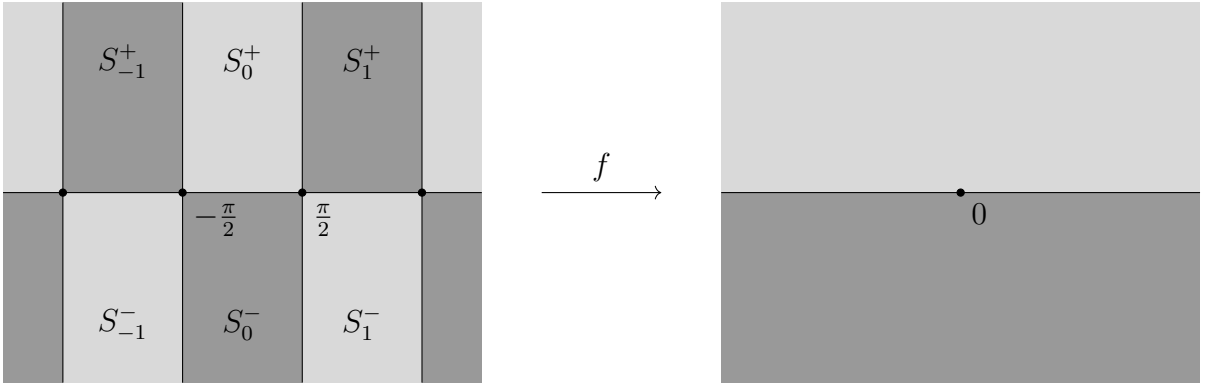

Figure 1: The function $f(z) = \sin(z)/2$ maps each half-strip $S_m^\sigma$ bijectively to the upper or lower half-plane.

We investigate this first example in some detail. Observe that $0$ is a fixed point of $f$; from the fact that $|f(x)| < |x|/2$ for real $x \neq 0$, we deduce that $f^n(x) \to 0$ for all $x \in \mathbb{R}$. Observe that $f^{-1}(\mathbb{R})$ consists of $\mathbb{R}$ together with the vertical lines $\{z \in \mathbb{C} \colon \operatorname{Re} z \in \pi/2 + \pi \mathbb{Z}\}$; hence it cuts the plane into the half-strips

$$S_m^\sigma := \left\{ x + iy \colon \frac{2m-1}{2}\pi < x < \frac{2m+1}{2}\pi \text{ and } \sigma y > 0 \right\},$$

where $m \in \mathbb{Z}$ and $\sigma \in \{+, -\}$. Each of these half-strips is mapped bijectively by $f$ onto either the upper or the lower half-plane; see Figure 1. Real starting values have bounded orbits, so the escaping set lies within the union of these half-strips.

We now consider the shape of the attracting basin

$$U := \{z \in \mathbb{C} \colon f^n(z) \to 0\}$$

in the complex plane (which is disjoint from $I(f)$). Observe that $U$ contains a horizontal strip $S$ centred on the real axis. Indeed, we have

$$|\operatorname{Im} f(z)| = \left| \frac{\cos(\operatorname{Re} z) \sinh(\operatorname{Im} z)}{2} \right| \leq \frac{\sinh(|\operatorname{Im} z|)}{2} < |\operatorname{Im} z|$$

whenever $0 < |\operatorname{Im} z| < t_0$, where $t_0 \approx 2.177$ is the unique positive fixed point of $\sinh/2$. Hence, $|\operatorname{Im} f^n(z)|$ is non-increasing for such $z$. In particular, $|f^n(z)|$ is bounded (by periodicity of $f$), and any limit point $z_0$ of the sequence satisfies $|\operatorname{Im} f(z_0)| = |\operatorname{Im} z_0|$. We see that $z_0 \in \mathbb{R}$, but the only possible limit point there is the fixed point $0$.

In particular

$$S := \{x + iy \colon |y| \leq 2\} \subset U.$$

Any point in $U$ is eventually mapped into $S$ under iteration. We consider the sets $\mathbb{C} \setminus f^{-n}(S)$. For $n = 1$, we obtain a union of simply connected domains, one within each $S_m^\sigma$, and each mapped bijectively to one of the two connected components of $\mathbb{C} \setminus S$. Each of these domains itself contains countably many connected components of $\mathbb{C} \setminus f^{-2}(S)$, and so on. Figure 2 shows the preimages $f^{-k}(S)$ of the horizontal strip $S$ under the map $f(z) = \sin(z)/2$ for $k = 0, \ldots, 3$. The interior of $f^{-k}(S)$ is shown in grey, its complement in white, and its boundary in black. Each white region is mapped univalently to a white region in the previous picture. For $k = 3$, the white regions and the grey channels between them become so thin close to infinity that they are no longer visible at the resolution of the picture.

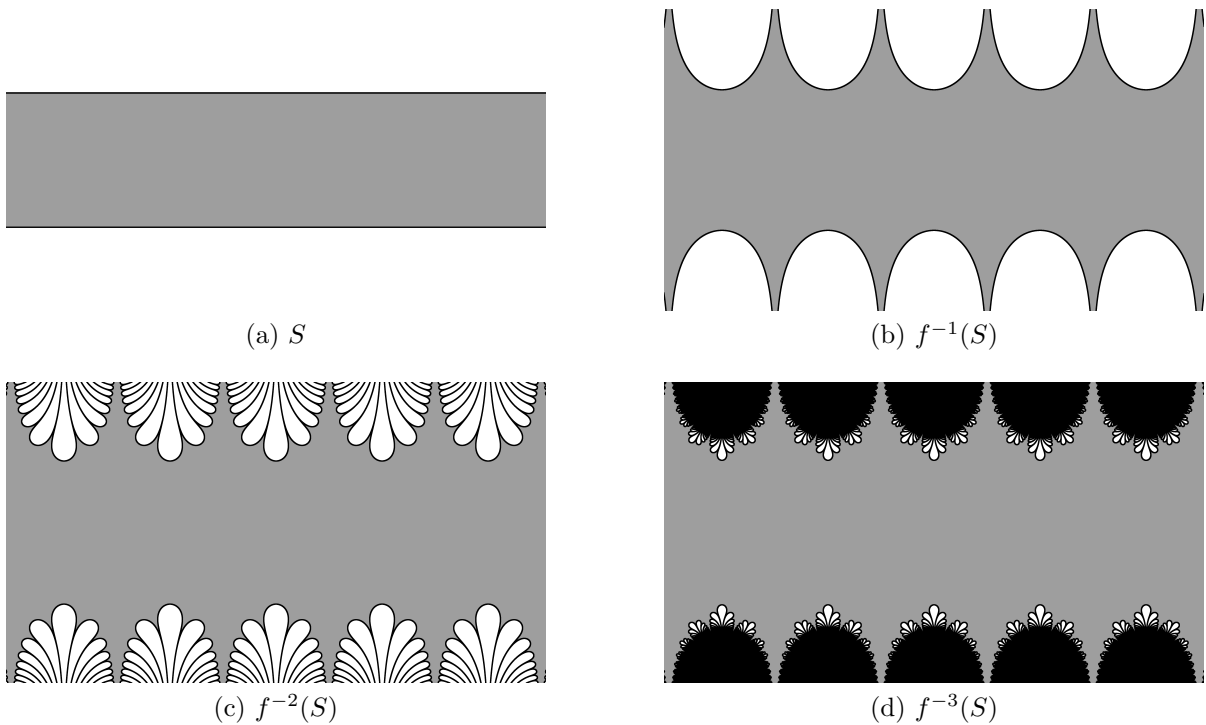

Figure 2: The preimages $f^{-k}(S)$ of the horizontal strip $S$ under the map $f(z) = \sin(z)/2$.

It follows that $U = \bigcup_{n=0}^{\infty} f^{-n}(S)$ is connected. What can be said about its complement? Outside $S$, the map $f$ is expanding:

$$|f'(z)| = \left|\frac{\cos(z)}{2}\right| \geq \frac{\sinh(2)}{2} =: \lambda > 1 \tag{2.2}$$

when $|\operatorname{Im} z| \geq 2$. We claim that, for $n \geq 1$, the set $\mathbb{C} \setminus f^{-n}(S)$ contains no disc of radius greater than $\pi/(2\lambda^{n-1})$.

To see this, let $z_0 \in \mathbb{C} \setminus f^{-n}(S)$. We can connect $z_{n-1} := f^{n-1}(z)$ with a horizontal segment $\gamma$ of length at most $\pi/2$ to a point $w_{n-1}$ whose real part is an odd multiple of $\pi/2$, and which hence belongs to $f^{-1}(S)$. Let $H$ be the connected component of $\mathbb{C} \setminus \mathbb{R}$ containing $\gamma$, and $V$ be the connected component of $f^{-(n-1)}(H)$ containing $z_0$. As noted above, the map $f^{n-1} \colon V \to H$ is bijective and its derivative has modulus at least $\lambda^{n-1}$ everywhere. So there is a connected component of $f^{-(n-1)}(\gamma)$ connecting $z_0$ to a point $w_0 \in f^{-n}(S)$, and of length at most $\pi/(2\lambda^{n-1})$. This proves the claim.

In particular, $U$ is dense in the plane; that is, its complement $J := \mathbb{C} \setminus U$ has no interior. Since the escaping set $I(f)$ is a subset of $J$, it follows that $I(f)$ has no interior.

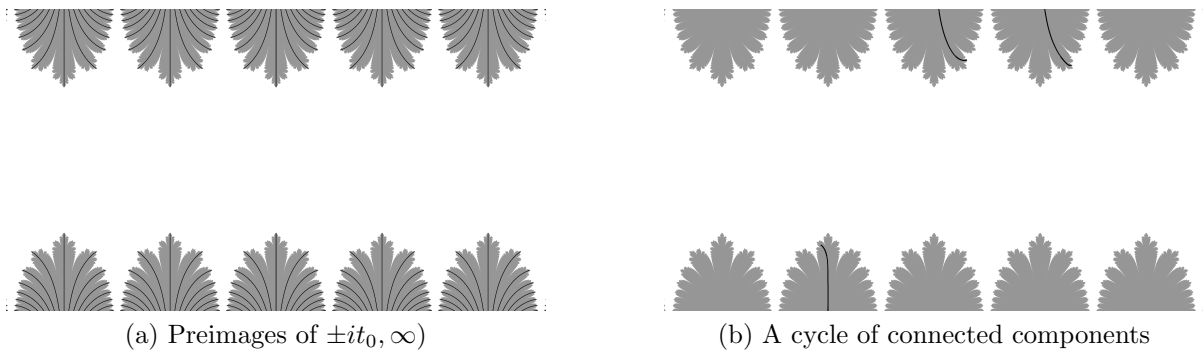

(a) Preimages of $\pm i t_0, \infty)$  (b) A cycle of connected components

Figure 3: Connected components of $J = \mathbb{C} \setminus U$.

What more can we say about $J$? Clearly it contains the half-lines $L^+ := i \cdot [t_0, \infty]$ and $L^- := i \cdot [-t_0, \infty)$, each of which is invariant under $f$. Moreover, since $\sinh(t)/t > t$ for $t > t_0$, all points of $L^\pm$ apart from the endpoints $it_0$ and $-it_0$ tend to infinity under iteration, and thus belong to $I(f)$. The translates of $L^\pm$ by integer multiples of $\pi$ are likewise in $J$. Taking further preimages we obtain a countable collection of arcs to infinity that lie dense in $J$, with all but the finite endpoint of each arc lying in the escaping set. Here and in the following we mean by *arc* a homeomorphic image of an interval. The preimages of $\pm i \cdot [t_0, \infty)$ under $f^2$ are shown in black in Figure 3 (a).

The set $J$ contains many points that are not on preimages of $L^\pm$. For instance, it is easy to show that each half-strip $S_{2m}^\sigma$, where $\sigma \in \{+, -\}$, contains a fixed point of $f$. (First observe that, for sufficiently large $M$, the inverse $(f|_{S_{2m}^\sigma})^{-1}$ maps the rectangle $\{z \in S_{2m}^\sigma : 2 \leq |\operatorname{Im} z| \leq M\}$ inside itself. Then apply the contraction mapping theorem.) Similarly, we can see that there exist many periodic points of higher periods. A cycle of connected components of $J$ containing a periodic point of period 3 is depicted in Figure 3 (b). All of the facts discussed so far can already be found in Fatou's paper [132].

Let $A$ be any connected component of $J$ and consider the set $\bigcap_{k=0}^\infty A_n$, where $A_n$ is the connected component of $\mathbb{C} \setminus f^{-n}(S)$ containing $A$. Each $A_n$ is one of the white domains shown in Figure 2 (for $n \leq 3$), and the pictures suggest that their intersection will again be an arc to infinity (see Figure 3). This is indeed correct:

**Theorem 2.1.** *Every connected component of $J$ is an arc $A$ connecting a finite endpoint to infinity. All points of $A$ except possibly this endpoint belong to the escaping set $I(f)$.*

An analogue of Theorem 2.1 for certain exponential maps (which we discuss in §2.3) was stated by Devaney and Krych [107, p. 50]; see also [101, Theorem on p. 170], [104, §3] and [106, Theorem 3.4]. Devaney [104, §4] as well as Devaney and Tangerman [108, Theorem 4.1] extended the ideas to other entire functions and in particular to the sine family, but the first full proof of Theorem 2.1 of which we are aware is by Aarts and Oversteegen [1, Theorem 5.7] who prove a stronger result that completely describes the topology of $J$.

A natural approach to proving that $A$ is an arc is to construct an explicit parameterisation of $A$ as a curve $\gamma\colon [0,\infty) \to A$ with $\gamma(t) \to \infty$ as $t \to \infty$. This is the approach used in the papers cited above. Another possibility is to use a topological classification of arcs. This is a more flexible approach that was introduced independently by Barański [25] and by Rottenfußer, Rückert, Rempe and Schleicher [251], in order to prove versions of Theorem 2.1 for larger classes of transcendental entire functions (see Theorem 7.4). We follow [251] in using the cut point characterisation of the arc. In order to state this characterisation, we introduce the following terminology. A compact, connected metric space is called a *continuum*. It is called *non-degenerate* if it consists of more than one point. A point $a$ of a continuum $A$ is called a *cut point* if $A \setminus \{a\}$ is disconnected; otherwise $a$ is called a *non-cut point*. The following result can be found in the book by Nadler [195, Theorem 6.17].

**Proposition 2.2.** *A non-degenerate continuum is an arc if and only if it has exactly two non-cut points.*

The two non-cut points in the interval $[0,1]$ are of course the endpoints $0$ and $1$. One way to distinguish these points from the other points in the interval is to use the natural ordering of the real numbers, noting that $0$ and $1$ are the minimum and maximum of the interval $[0,1]$ with respect to this ordering. This indicates that in order to show that a continuum is an arc one may consider an ordering on it. The following result, which combines Proposition 2.2 with [195, Theorem 6.16], makes this idea precise.

**Proposition 2.3.** *Let $A$ be a non-degenerate continuum. Suppose that there is a total ordering $\prec$ on $A$ whose order topology agrees with the topology of $A$. Then $A$ is an arc.*

*Proof of Theorem 2.1.* Let $z_0 \in J$, and consider the set $A$ consisting of all points $z \in J$ such that, for every $n \geq 0$, the points $f^n(z)$ and $f^n(z_0)$ belong to the same connected component of $\mathbb{C} \setminus f^{-1}(\mathbb{R})$. In other words, they lie in the same half-strip $S_m^\sigma$, where $m$ and $\sigma$ depend on $n$. Then

$$\widehat{A} := A \cup \{\infty\} = \bigcap_{n=0}^{\infty} \left(\overline{A_n} \cup \{\infty\}\right),$$

where $A_n$ is the connected component of $\mathbb{C} \setminus f^{-n}(S)$ containing $z_0$. So $\widehat{A}$ is compact and connected as a nested intersection of compact and connected sets. Moreover, every point of $J \setminus A$ is separated from $A$ by $f^{-n}(\mathbb{R})$ for large enough $n$. We will use Proposition 2.3 to show that $\widehat{A}$ is an arc with one endpoint at infinity.

To define a total ordering on $\widehat{A}$, we use the speed of escape of points in $A$. We claim that for all $z, w \in A$ the following hold:

(i) If $|\operatorname{Im} f^n(w)| \geq |\operatorname{Im} f^n(z)| + 2$ for some $n \geq 0$, then $|\operatorname{Im} f^m(w)| > |\operatorname{Im} f^m(z)| + 2$ for all $m > n$.

(ii) $|\operatorname{Im} f^n(z) - \operatorname{Im} f^n(w)| \to \infty$ as $n \to \infty$.

To establish (ii), again connect $f^n(z)$ and $f^n(w)$ by a straight line segment. Then $(f^n|_{A_n})^{-1}$ maps this segment to a curve connecting $z$ and $w$. By (2.2), $|f^n(z)-f^n(w)| \geq \lambda^n|z-w| \to \infty$. But $|\operatorname{Re} f^n(z) - \operatorname{Re} f^n(w)| < \pi$, since the two points belong to a common half-strip $S_m^\sigma$, and (ii) follows.

To prove (i), it is enough to consider the case $n = 0$ and $m = 1$; the general case follows inductively. So suppose that $z, w \in A$ with $|\operatorname{Im} w| \geq |\operatorname{Im} z| + 2$. Since $f$ behaves essentially like an exponential map in the upper and lower half-plane, the image of $w$ will be much larger than that of $z$. Since their real parts differ by at most $\pi$, this must be accounted for by the size of the imaginary parts.

More precisely,
$$|f(w)| \geq \frac{1}{4}\left(e^{|\operatorname{Im} w|} - e^{-|\operatorname{Im} w|}\right) \geq \frac{e^{|\operatorname{Im} w|}}{4} - \frac{1}{4},$$
and similarly, $|f(z)| \leq e^{|\operatorname{Im} z|}/4 + 1/4$. Using the assumption on $z$ and $w$, and noting that $|\operatorname{Im} z| \geq t_0 > 2$ for $z \in A$, a simple calculation shows that
$$|f(w)| \geq e^2 \cdot \frac{e^{|\operatorname{Im} z|}}{4} - \frac{1}{4} > \sqrt{2} \cdot \frac{e^{|\operatorname{Im} z|}}{4} + 10 > \sqrt{2} \cdot |f(z)| + \pi + 2.$$

So indeed
$$|\operatorname{Im} f(w)| \geq |f(w)| - |\operatorname{Re} f(w)| \geq |f(w)| - |\operatorname{Re} f(z)| - \pi$$
$$> \sqrt{2} \cdot |f(z)| - |\operatorname{Re} f(z)| + 2 \geq |\operatorname{Im} f(z)| + 2.$$

We now define $z \prec \infty$ for $z \in A$, and $z \prec w$ for $z, w \in A$ if $|\operatorname{Im} f^n(w)| > |\operatorname{Im} f^n(z)| + 2$ for some $n \geq 0$. It follows easily from (i) and (ii) that $\prec$ is a strict total order on $\widehat{A}$. The topology induced by this order is the same as the topology inherited from the plane; that is, the identity is a homeomorphism from $\widehat{A}$ with the plane topology to $\widehat{A}$ with the order topology. Indeed, since $\widehat{A}$ is compact, one must only check that this map is continuous. This holds because open intervals in $\prec$ are open with respect to the usual topology by definition.

In conclusion, $\widehat{A}$ is a compact and connected ordered metric space. By Proposition 2.3, $\widehat{A}$ is an arc, whose endpoints are the maximal and minimal elements of the order $\prec$. In particular, $\infty$ is one of these endpoints, and $A$ is connected. So $A$ is indeed the connected component of $J$ containing $z_0$.

Let $a$ be the finite endpoint of $\widehat{A}$; that is, the minimal element of the order. If $z \in A \setminus \{a\}$, then for sufficiently large $n$,
$$|\operatorname{Im} f^n(z)| \geq |\operatorname{Im} f^n(z) - \operatorname{Im} f^n(a)| - |\operatorname{Im} f^n(a)|$$
$$\geq |\operatorname{Im} f^n(z) - \operatorname{Im} f^n(a)| - |\operatorname{Im} f^n(z)| + 2$$
by (i) and thus
$$|\operatorname{Im} f^n(z)| \geq \frac{1}{2}|\operatorname{Im} f^n(z) - \operatorname{Im} f^n(a)| + 1.$$
It now follows from (ii) that $|\operatorname{Im} f^n(z)| \to \infty$ as $n \to \infty$. So $z \in I(f)$, which completes the proof. $\square$

We conclude that, for $f(z) = \sin(z)/2$, the escaping set $I(f)$ is nowhere dense, and its connected components are open or half-open arcs to $\infty$.

It is not difficult to see that there are uncountably many such arcs. In fact, a modification of the argument showing that $f$ has periodic points shows that if $(\sigma_k)$ is a sequence in $\{+,-\}$ and if $(n_k)_{k\in\mathbb{N}_0}$ is a bounded sequence in $\mathbb{Z}$, then there exists a point $z \in S_{n_0}^{\sigma_0} \cap J \setminus I(f)$ such that $f^k(z) \in S_{n_k}^{\sigma_k}$ for all $k \in \mathbb{N}$. This yields uncountably many points in $J \setminus I(f)$ and hence uncountably many arcs in $J$ connecting these points with $\infty$.

Loosely speaking, these arcs form a "Cantor set of curves", an expression used by Devaney and Krych [101, 107] to describe the configuration of these arcs. Nowadays this structure is usually called a *Cantor bouquet*, a term first used in papers by Devaney [102, p. 145] and Devaney and Tangerman [108, §1]. For a formal definition of a Cantor bouquet we refer to Barański, Jarque and Rempe [29, Definition 1.1].

## 2.2 Fatou's function

How representative is the structure of the previous example? For example, can the escaping set of a transcendental entire function have interior points? The function

$$f \colon \mathbb{C} \to \mathbb{C}, \quad f(z) := z + 1 + e^{-z}, \tag{2.3}$$

called *Fatou's function* since it was also considered by Fatou [132, §15], answers this question affirmatively. If $\delta > 0$ and $\operatorname{Re} z > \delta$, then $\operatorname{Re} f(z) \geq \operatorname{Re} z + 1 - e^{-\delta} > \operatorname{Re} z$. We conclude that $\operatorname{Re} f^n(z) \geq \operatorname{Re} z + (1 - e^{-\delta})n \to +\infty$. So the interior of $I(f)$ contains the entire right half-plane $\mathbb{H}_{>0}$. Similarly as in §2.1, we may investigate the open set $U$ of all points $z$ for which $\operatorname{Re} f^n(z) \to +\infty$.

The real axis and its translates by integer multiples of $2\pi i$ are contained in $U$. Very similarly to §2.1, within each horizontal strip bounded by these lines, there is a connected component of $f^{-1}(\mathbb{H}_{<0})$ that is mapped conformally to the left half-plane $\mathbb{H}_{<0}$. (See Figure 4 (a). These domains are shown in white, the half-plane $\mathbb{H}_{>0}$ is shown in dark grey, and $f^{-1}(\mathbb{H}_{>0}) \cap \mathbb{H}_{<0}$ is shown in light grey.) Following the same proof as that of Theorem 2.1, one can show that the set $J := \mathbb{C} \setminus U$ (shown in black in Figure 4 (b)) is again a union of curves to infinity, with all points except for some of the finite endpoints escaping (we omit the details here). Moreover, the set $J$ is contained in the *boundary* of the escaping set. Indeed, periodic points of $f$ are dense in $J$. This follows from Theorems 3.1 (d) and Theorem 4.3 below, but here it could also be proved directly, by considering endpoints of periodic connected components of $J$. So we have the following result.

**Theorem 2.4.** *The set $U \subset I(f)$ of points $z$ satisfying $\operatorname{Re} f^n(z) \to +\infty$ is an open and dense subset of the complex plane, in which points tend to infinity at a linear rate.*

*Every connected component of the complement $J := \mathbb{C} \setminus U = \partial U = \partial I(f)$ is an arc connecting a finite endpoint to infinity. Every point of $J$, with the exception of certain endpoints, has an orbit along which the real parts tend to $-\infty$ at a super-exponential rate.*

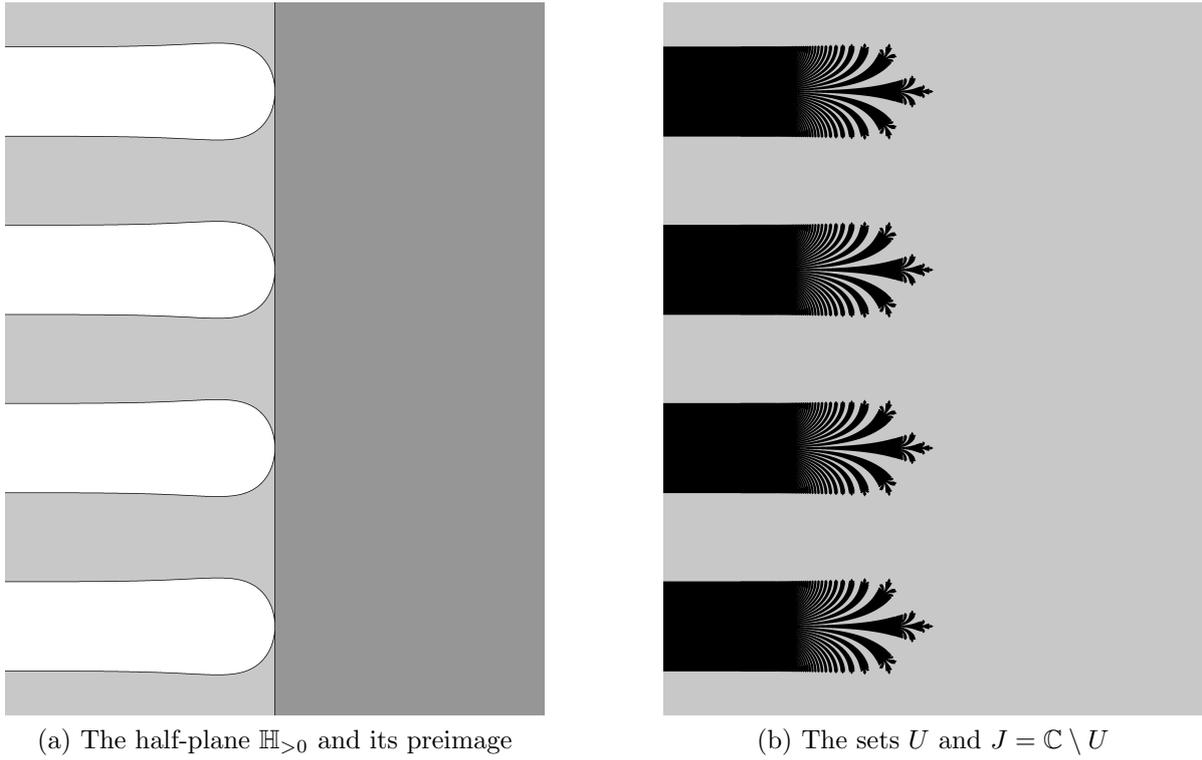

(a) The half-plane $\mathbb{H}_{>0}$ and its preimage

(b) The sets $U$ and $J = \mathbb{C}\setminus U$

Figure 4: Fatou's function.

In particular, $\mathbb{C}\setminus I(f)$ is a totally disconnected set (consisting of certain endpoints of $J$), and $I(f)$ is connected.

As first noted by Kotus and Urbański [163, Proposition 2.5], $J$ is again a Cantor bouquet. For further discussion of the escaping set of Fatou's function we refer to Theorems 7.45–7.47 below.

## 2.3 The exponential family

Next, we consider an entire *family* of functions, given by

$$f_a \colon \mathbb{C} \to \mathbb{C}, \quad f_a(z) \coloneqq ae^z, \tag{2.4}$$

where $a \in \mathbb{C}\setminus\{0\}$ is a parameter. This family, the *exponential family*, is perhaps the most-studied family in all of transcendental dynamics; our goal here is merely to give an impression of what happens to the escaping set of these functions as $a$ varies.

Consider first the case where $a \in (0, 1/e)$. By the intermediate value theorem, $f_a$ has a real fixed point $x_0$ between 0 and 1. Furthermore, the left half-plane $\mathbb{H}_{<1}$ is mapped inside itself:

$$f_a(\mathbb{H}_{<1}) = D(0, a\cdot e)\setminus\{0\} \subset \mathbb{H}_{<1}.$$

In particular, $|f'_a| = |f_a|$ is bounded by $a\cdot e < 1$, and it follows from the contraction mapping theorem that all points in $\mathbb{H}_{<1}$ converge to $x_0$ under iteration. Consider the horizontal lines whose imaginary part is an odd multiple of $\pi$. These are mapped to $(-\infty, a) \subset \mathbb{H}_{<1}$, and cut the plane into strips of height $2\pi$. Similarly as in the previous example, each such strip contains one connected component of $\mathbb{C} \setminus f^{-1}(\mathbb{H}_{\leq 1})$ that is mapped conformally to the right half-plane $\mathbb{H}_{>1}$. Once again, an analogue of Theorem 2.1 holds: the points that never enter $\mathbb{H}_{<1}$ (in particular, all escaping points) are arranged in arcs to infinity. All points of such an arc, with the possible exception of its finite endpoint, belong to the escaping set. Again these arcs form a Cantor bouquet. We omit the details. (As mentioned above, historically Devaney established this result for the maps $f_a$ before it was proved for the function $f$ from (2.1).) The case $a = 1/4$ is illustrated in Figure 5 (a).

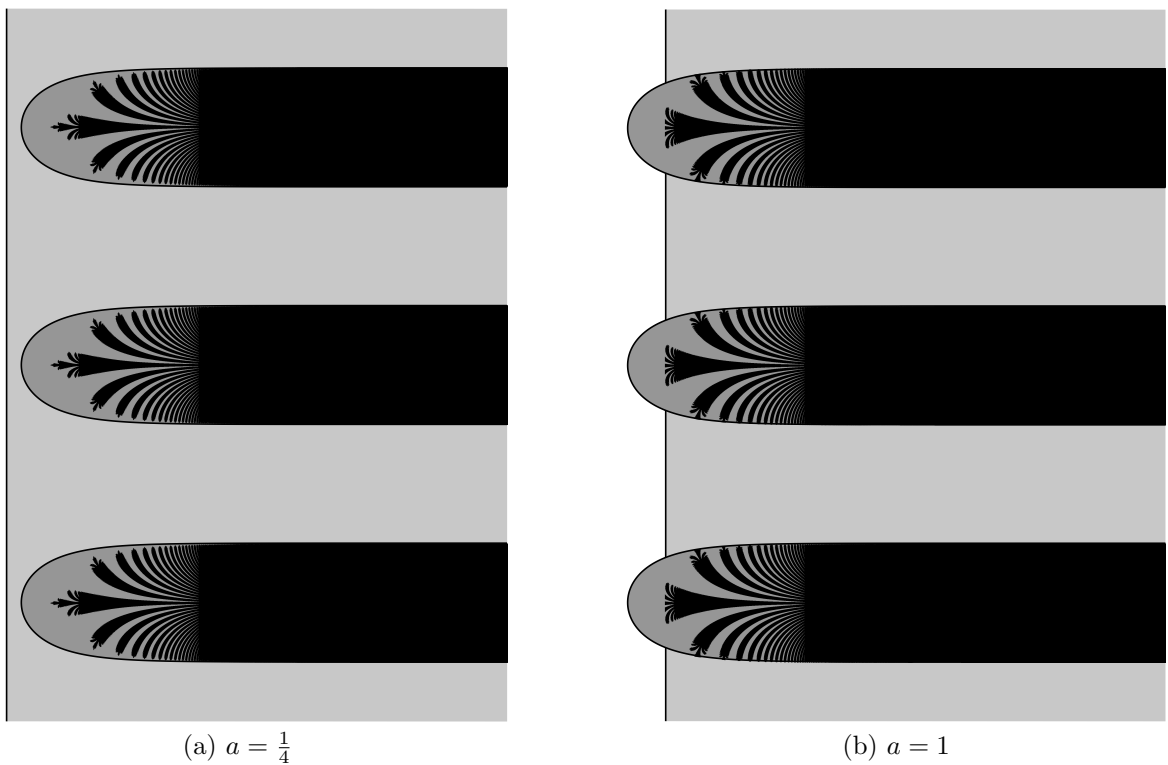

(a) $a = \frac{1}{4}$          (b) $a = 1$

Figure 5: The half-plane $H$ (light grey), its preimage under $f_a$ (dark grey), and the set of points that remain in $H$ under iteration (black)

The dynamics of $f_a$ becomes more complicated for general parameters $a \in \mathbb{C} \setminus \{0\}$. For example, for the exponential map $f_1(z) = e^z$, it is known that the escaping set $I(f_1)$ is dense in the complex plane. This is a consequence of a famous theorem of Misiurewicz [193], which had been conjectured by Fatou in 1926 [132, p. 370]. The same holds for any parameter value for which 0 is preperiodic or in the escaping set; see [23, Corollary 1], [65, §3] or [259, Theorem 4.1] for the result of Misiurewicz and its extensions. Nonetheless, some structure of the escaping set is common to all parameters, as the following result of Schleicher and Zimmer shows [256, Corollary 6.9].

**Theorem 2.5.** *Let $a \in \mathbb{C} \setminus \{0\}$, and let $z \in I(f_a)$. Then there exists an arc $\gamma \subset I(f)$*

*that connects the points $z$ and $\infty$.*

The proof relies on the fact that, for sufficiently large real parts (which is where $|f_a(z)|$ is large), the behaviour of the map is very similar to the case that $a \in (0, 1/e)$. (See Figure 5 (a) for an illustration of this for $a = 1$.) In general, there is no longer an invariant left half-plane, so points with small real parts may eventually map to points with large real parts, and vice versa. However, for those points whose forward orbits *remain* in a right half-plane at all times, we obtain a very similar picture as before. More precisely, the following holds.

**Proposition 2.6.** *Let $a \in \mathbb{C} \setminus \{0\}$. Then there is a number $R > 0$ with the following property. If $z \in \mathbb{C}$ is such that $\operatorname{Re} f^n(z) \geq R$ for all $n \geq 0$, then $z$ can be connected to infinity by a curve consisting of escaping points.*

*Sketch of proof.* We sketch a proof using similar ideas as Theorem 2.1 (the original theorem was, again, proved by constructing a parameterisation of the desired curve). The preimage under $f$ of a right half-plane $\mathbb{H}_{>\rho}$, where $\rho > 0$, consists of an infinite number of domains, one between any two adjacent horizontal lines of the form $\mathbb{R} + i(\pi - \arg a)$. If $\rho$ is chosen sufficiently large, then an analogue of (i) holds: Suppose that $z$ and $w$ belong to a common connected component of $f_a^{-1}(\mathbb{H}_{>\rho})$ and that the imaginary parts of their images $f(z)$ and $f(w)$ differ by at most $2\pi$. If $\operatorname{Re} w \geq \operatorname{Re} z + 2$, then also $\operatorname{Re} f_a(w) > \operatorname{Re} f_a(z) + 2$. (We leave the simple proof as an exercise.)

Now suppose that $R > \rho + 2$, and that $z$ is as in the statement of the proposition. For $n \geq 0$, consider the ray $\alpha_n := [0, \infty) + f_a^n(z)$, and the connected component $\gamma_n$ of $f_a^{-n}(\alpha_n)$ containing $z$. It follows inductively from the above fact that $\operatorname{Re} f_a^j(\zeta) > \operatorname{Re} f_a^j(z) - 2 \geq R - 2$ for $j = 0, \ldots, n$. Taking a Hausdorff limit of the $\gamma_n$, we obtain an unbounded connected set $A$ containing $z$ with $\operatorname{Re} f_a^j(\zeta) \geq R - 2 > \rho$ for all $\zeta \in A$. As in the proof of Theorem 2.1, it follows that $A$ is an arc to infinity, with all points except possibly its finite endpoint escaping. $\square$

*Proof of Theorem 2.5.* Let $R$ be as in Proposition 2.6, and let $z \in I(f_a)$. Since

$$|a| \exp(\operatorname{Re} f_a^n(z)) = |f_a^{n+1}(z)| \to \infty,$$

we see that $\operatorname{Re} f_a^n(z) \to \infty$. Hence there is $n_0$ such that Proposition 2.6 applies to $f_a^{n_0}(z)$, which hence can be connected to infinity by a curve of escaping points.

Observe that, if $\gamma \subset \mathbb{C}$ is an arc to infinity, then every connected component of $f_a^{-1}(\gamma)$ is a curve that tends to infinity in one or both directions. Indeed, if $0 \notin \gamma$, then $f_a^{-1}(\gamma)$ consists of a countable union of arcs to infinity (the images of $\gamma/a$ under the different branches of the logarithm); if $0 \in \gamma$, then $f_a^{-1}(\gamma)$ consists of countably many open arcs tending to $\infty$ in both directions (with the real parts tending to $-\infty$ in one direction and $+\infty$ in the other), together with a countable union of arcs to infinity.

The theorem follows by inductively applying this observation. $\square$

Hence the escaping set still consists of curves for every parameter $a$, as it does in the simplest examples. However, the way in which these curves fit together in the plane varies considerably with the parameter $a$, and it is precisely this fact that makes them a useful tool for studying the dynamics of $f_a$, as mentioned in the introduction.

## 2.4 Baker's example of a wandering domain

While the escaping sets in the previous examples had differing structures, they had one thing in common: the presence of "Cantor bouquets" of curves in $I(f)$. Even though the escaping set of Fatou's function considered in §2.2 has interior points, and the escaping set is dense in the plane, we saw that its boundary still contains such curves. Here we discuss a very different example, given by Baker [11] in 1963. In his example $I(f)$ contains large annuli around 0; in particular, the interior of $I(f)$ separates every point of $\partial I(f)$ from $\infty$.

To define the function, let $C > 0$ and $r_1 > 1$ and define $(r_n)_{n \in \mathbb{N}}$ recursively by

$$r_{n+1} := Cr_n^2 \prod_{k=1}^{n} \left(1 + \frac{r_n}{r_k}\right). \tag{2.5}$$

It is elementary to see that if $C \exp(2/r_1) < 1/4$ and $Cr_1 > 1$, then the sequence $(r_n)$ tends to $\infty$ so rapidly that the infinite product

$$f(z) := Cz^2 \prod_{k=1}^{\infty} \left(1 + \frac{z}{r_k}\right) \tag{2.6}$$

converges and defines an entire function $f$.

Baker showed that the annuli

$$A_n := \text{ann}\left(r_n^2, \sqrt{r_{n+1}}\right) \tag{2.7}$$

satisfy

$$f(A_n) \subset A_{n+1} \tag{2.8}$$

for large $n$.

The idea is that $f(z)$ behaves like a multiple of $z^{n+2}$ in the annulus $A_n$. Indeed, for $z \in A_n$ and $k > n$ the factor $1 + z/r_k$ in the infinite product defining $f$ is very close to 1. Using this and noting that $(r_n)$ increases rapidly it is not difficult to see that

$$f(z) \sim Cz^2 \prod_{k=1}^{n} \frac{z}{r_k} = \varepsilon_n z^{n+2} \tag{2.9}$$

for $z \in A_n$, with $\varepsilon_n := C/\prod_{k=1}^{n} r_k$ The recursion formula (2.5) yields that

$$r_{n+1} = 2Cr_n^2 \prod_{k=1}^{n-1} \left(1 + \frac{r_n}{r_k}\right) \sim 2\varepsilon_{n-1} r_n^{n+1} \tag{2.10}$$

as $n \to \infty$. A simple consequence is that

$$r_{n+1} \leq r_n^{n+1} \tag{2.11}$$

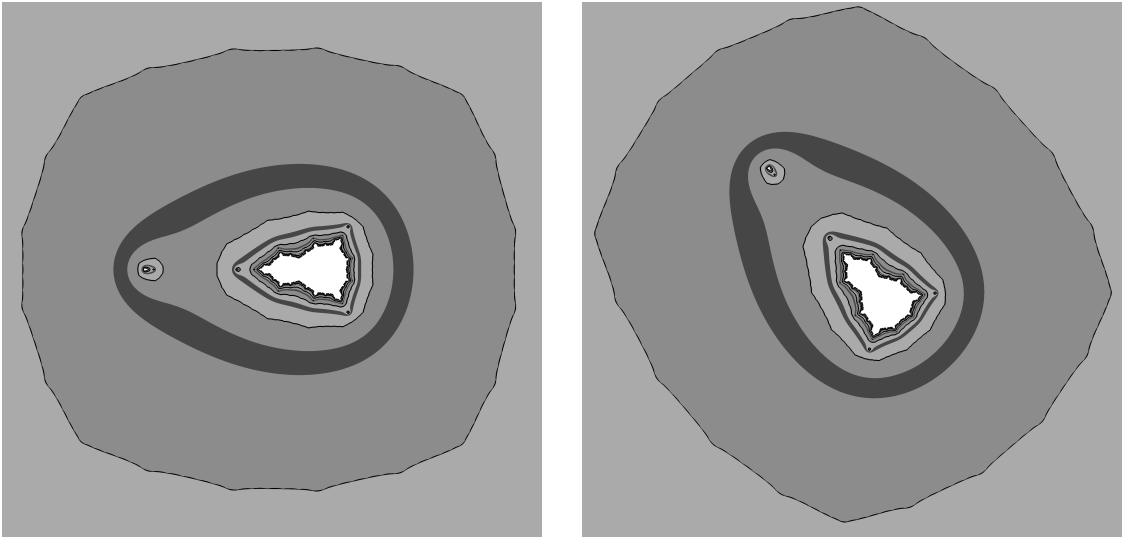

Figure 6: Illustration of Baker's example.

for large $n$. The lower bound $|z| > r_n^2$ for $z \in A_n$ yields together with (2.9), (2.10) and (2.11) that

$$|f(z)| \geq (1-o(1))\varepsilon_n r_n^{2(n+2)} = (1-o(1))\varepsilon_{n-1} r_n^{2(n+2)-1}$$
$$= (1-o(1))\frac{1}{2}r_{n+1}r_n^{n+2} \geq (1-o(1))\frac{1}{2}r_{n+1}^2 r_n > r_{n+1}^2$$

for large $n$. Similarly we deduce from the upper bound $|z| < \sqrt{r_{n+1}}$ and (2.9), using (2.10) and (2.11) with $n$ replaced by $n+1$, that

$$|f(z)| \leq (1+o(1))\varepsilon_n r_{n+1}^{(n+2)/2} = (1+o(1))\frac{r_{n+2}}{2r_{n+1}^{(n+2)/2}} \leq (1+o(1))\frac{\sqrt{r_{n+2}}}{2} < \sqrt{r_{n+2}}$$

for large $n$. The last two inequalities yield (2.8). For more details we refer to Baker's paper.

Property (2.8) implies that $A_n \subset I(f)$ for large $n$. Let $U_n$ be the maximal open subset of $I(f)$ containing $A_n$. Since 0 is a fixed point of $f$, the domains $U_n$ are indeed multiply connected.

Figure 6 illustrates this function for $C = 1/6$ and $r_1 = 12$. Preimages of $A_5$ are shown in dark grey. The maximal open subsets of $I(f)$ containing these preimages are shown in lighter shades of grey and their boundaries in black. Points converging to 0 under iteration are white. The range shown in the left picture is $-85 \leq \operatorname{Re} z \leq 55$ and $|\operatorname{Im} z| \leq 70$.

By a general result in complex dynamics (see Theorem 3.1 $(f)$ below) every point on the boundary of such a grey domain is a limit point of preimages of 0. So there should be white spots near these boundaries. However, these white spots are so small that they are not visible in the picture. For example, $z_0 \approx -59.9095005465 + 39.9450315388i$ satisfies $f^2(z_0) = 0$. The right picture shows the range $-2 \cdot 10^{-8} \leq \operatorname{Re}(z - z_0) \leq 1.5 \cdot 10^{-8}$ and $1.5 \cdot 10^{-8} \leq \operatorname{Im}(z - z_0) \leq 2 \cdot 10^{-8}$.

In [11], Baker left open the question whether the $U_n$ are all distinct or all equal. It was only later [13] that he could prove that the first alternative holds, thereby giving the first example of a *wandering domain*; see §6.3 for a detailed discussion.

## 3 A brief introduction to complex dynamics

In this section, we summarise some concepts and results of complex dynamics. We refer to the books by Beardon [38], Carleson and Gamelin [93], Milnor [192] and Steinmetz [273] for an introduction to the dynamics of rational functions. For transcendental functions, and in particular for the differences that arise in the theory in comparison to rational functions, we refer to articles by Bergweiler [45], Herman [147, Appendix], Eremenko and Lyubich [115, Chapter 4], Kisaka [159], Schleicher [255] and Sixsmith [268] as well as books by Morosawa, Nishimura, Taniguchi and Ueda [194] and Piñeiro [213]. Of these articles and books, [159, Section 7], [194, Section 3.1.3], [213, Section 6.4] and [255, Section 4] contain sections devoted to the escaping set.

A family of functions meromorphic in a domain $D$ is called *normal* in $D$ if every sequence in this family has a subsequence that converges locally uniformly in $D$. Here convergence is understood with respect to the spherical metric. Thus even if all functions in the sequence are holomorphic, we allow the possibility that the limit equals $\infty$. This remark is particularly important for this article, as it is devoted to the set of points that tend to $\infty$ under iteration.

The family is called *normal at a point* $z_0 \in D$, if $z_0$ has a neighbourhood where it is normal. By the Arzelà–Ascoli theorem, a family is normal in a domain if and only if it is equicontinuous at every point of this domain (with respect to the spherical metric).

Let now $f \colon \mathbb{C} \to \mathbb{C}$ be entire. The *Fatou set* $F(f)$ is the set of all points in $\mathbb{C}$ where the family $\{f^n \colon n \in \mathbb{N}\}$ of iterates of $f$ is normal. The complement $J(f) := \mathbb{C} \setminus F(f)$ is called the *Julia set*.

We say that $z_0 \in \mathbb{C}$ is a *periodic point* of $f$ if there exists $p \in \mathbb{N}$ such that $f^p(z_0) = z_0$. The smallest $p$ with this property is called the *period* of $z_0$. For a periodic point $z_0$ of period $p$, we call $(f^p)'(z_0)$ the *multiplier* of $z_0$. A periodic point is called *attracting*, *indifferent*, or *repelling* depending on whether the modulus of its multiplier is less than, equal to, or greater than 1. Periodic points of multiplier 0 are called *superattracting*. The multiplier of an indifferent periodic point is of the form $e^{2\pi i \alpha}$ where $0 \leq \alpha < 1$. We say that $z_0$ is *rationally indifferent* if $\alpha$ is rational and *irrationally indifferent* otherwise. Rationally indifferent periodic points are also called *parabolic*. Also, a point $z_0$ is called *preperiodic* if $f^n(z_0)$ is periodic for some $n \geq 1$. It is called *strictly preperiodic* if it is preperiodic but not periodic. Finally, a periodic point of period 1 is called a *fixed point*.

It is easy to see that attracting periodic points are in the Fatou set while repelling and parabolic periodic points are in the Julia set.

For $z_0 \in \mathbb{C}$, the *(forward) orbit* $O^+(z_0)$ and the *backward orbit* $O^-(z_0)$ of $z_0$ are defined by $O^+(z_0) := \{f^n(z_0) \colon n \geq 0\}$ and
$$O^-(z_0) := \bigcup_{n \geq 0} f^{-n}(z_0),$$
where $f^{-n}(z_0) := \{z \colon f^n(z) = z_0\}$. For $A \subset \mathbb{C}$ we put $O^\pm(A) := \bigcup_{z \in A} O^\pm(z)$. The *exceptional set* $E(f)$ is defined as the set of all points with finite backward orbit. It is easy to see that $E(f)$ contains at most one point. (For transcendental $f$ this is an immediate consequence of Picard's theorem.)

The following theorem summarises some basic results of the theory.

**Theorem 3.1.** *Let $f$ be an entire function, neither constant nor a polynomial of degree $1$. Then:*

(a) *$J(f)$ is a non-empty perfect set.*

(b) *$J(f)$ and $F(f)$ are completely invariant; that is, $z \in J(f) \Leftrightarrow f(z) \in J(f)$ and $z \in F(f) \Leftrightarrow f(z) \in F(f)$.*

(c) *$J(f^n) = J(f)$ and $F(f^n) = F(f)$ for all $n \in \mathbb{N}$.*

(d) *$J(f)$ is the closure of the set of repelling periodic points.*

(e) *If $z_0 \in \mathbb{C} \setminus E(f)$, then $J(f) \subset \overline{O^-(z_0)}$. If $z_0 \in J(f) \setminus E(f)$, then $J(f) = \overline{O^-(z_0)}$.*

(f) *If $K \subset \mathbb{C} \setminus E(f)$ is compact and $D$ is an open set intersecting $J(f)$, then there exists $N \in \mathbb{N}$ such that $f^n(D) \supset K$ for all $n \geq N$.*

(g) *If $A$ is a closed, backward invariant subset of $\mathbb{C}$ containing at least three points, then $J(f) \subset A$.*

Property (f) is sometimes called the *blowing-up property* of the Julia set.

The connected components of $F(f)$ are called *Fatou components*. A Fatou component $U$ is called *periodic* if there exists $p \in \mathbb{N}$ such that $f^p(U) \subset U$. If this holds with $p = 1$ so that $f(U) \subset U$, then $U$ is called *invariant*. We say that $U$ is *preperiodic* if there exists $q \in \mathbb{N}$ such that $f^q(U)$ is contained in a periodic Fatou component. A Fatou component that is not preperiodic is called a *wandering domain*. We will discuss wandering domains in §6.3 and §6.5.

The behaviour of $f^n$ in periodic Fatou components is well understood. Because of Theorem 3.1 (c) it suffices to consider invariant Fatou components. They can be classified as follows; see [45, Theorem 6] and [116, Section 5].

**Theorem 3.2.** *Let $U$ be an invariant Fatou component of the entire function $f$. Then we have one of the following possibilities:*

(a) *$U$ contains an attracting fixed point $z_0$. Then $f^n|_U \to z_0$ as $n \to \infty$. In this case, $U$ is called the immediate attracting basin of $z_0$.*

(b) $\partial U$ contains a fixed point $z_0$ and $f^n|_U \to z_0$ as $n \to \infty$. Then $f'(z_0) = 1$. In this case, $U$ is called an *immediate parabolic basin* of $z_0$.

(c) There exists a biholomorphic map $\varphi\colon U \to \mathbb{D}$ such that $\varphi(f^p(\varphi^{-1}(z))) = e^{2\pi i\alpha}z$ for some $\alpha \in \mathbb{R} \setminus \mathbb{Q}$. In this case, $U$ is called a *Siegel disc*.

(d) $f^n|_U \to \infty$ as $n \to \infty$. In this case, $U$ is called a *Baker domain*.

The *attracting basin* of an attracting fixed point $z_0$ is the defined as the set of all $z \in \mathbb{C}$ for which $f^n(z) \to z_0$ as $n \to \infty$. The immediate attracting basin is the connected component of the attracting basin that contains $z_0$.

There may be more than one immediate parabolic basin associated to a fixed point $z_0$ of $f$ satisfying $f'(z_0) = 1$. In fact, the number $m$ of such basins is given by $m = \min\{k \in \mathbb{N}\colon f^{(k+1)}(z_0) \neq 0\}$. The set of all points that are mapped to an immediate parabolic basin of $z_0$ is called the *parabolic basin* of $z_0$.

We will apply the terminology introduced not only for invariant but also for periodic Fatou components. Periodic Fatou components are closely related to the set $\operatorname{sing}(f^{-1})$ of singularities of the inverse function $f^{-1}$ of $f$. This set consists of the critical and finite asymptotic values of $f$. Here $a \in \mathbb{C}$ is said to be a *critical value* of $f$ if there exists $\xi \in \mathbb{C}$ such that $f(\xi) = a$ and $f'(\xi) = 0$. Then $\xi$ is called a *critical point*. We say that $a$ is an *asymptotic value* of $f$ if there exists a curve $\gamma\colon [0,1) \to \mathbb{C}$ such that $\gamma(t) \to \infty$ and $f(\gamma(t)) \to a$ as $t \to 1$.

The set
$$P(f) := \overline{O^+(\operatorname{sing}(f^{-1}))} \tag{3.1}$$
is called the *postsingular set*. The relation between periodic Fatou components and singularities is described by the following result; see [45, Theorem 7]

**Theorem 3.3.** *Let $f$ be an entire function and let $C = \{U_0, U_1, \ldots, U_{p-1}\}$ be a periodic cycle of Fatou components.*

(a) *If $C$ is a cycle of immediate attracting or parabolic basins, then $U_j \cap \operatorname{sing}(f^{-1}) \neq \emptyset$ for some $j \in \{0, 1, \ldots, p-1\}$.*

(b) *If $C$ is a cycle of Siegel discs, then $\partial U_j \subset P(f)$ for all $j \in \{0, 1, \ldots, p-1\}$.*

With appropriate modifications, the above results also holds for rational functions $f\colon \widehat{\mathbb{C}} \to \widehat{\mathbb{C}}$ and holomorphic functions $f\colon \mathbb{C}\setminus\{0\} \to \mathbb{C}\setminus\{0\}$. In Theorem 3.2 we then have the additional possibility of a Herman ring. This is defined as the Siegel disc, except that $\varphi$ maps $U$ to annulus instead of a disc. Theorem 3.3 (b) holds also for Herman rings.

We mention two classes of entire functions that play an important role in transcendental dynamics, namely the *Speiser class* $\mathcal{S}$ consisting of all transcendental entire functions for which $\operatorname{sing}(f^{-1})$ is finite and the *Eremenko–Lyubich class* $\mathcal{B}$ for which $\operatorname{sing}(f^{-1})$ is bounded. Here we only mention that Eremenko and Lyubich [116, Theorem 3] and, independently, Goldberg and Keen [139, Theorem 4.2] proved that functions in $\mathcal{S}$ do not have wandering domains. Moreover, Eremenko and Lyubich [116, Theorem 1] showed that for functions in $\mathcal{B}$ the escaping set is contained in the Julia set; see Theorem 6.1 below. An excellent survey of the dynamics of functions in these classes has been given by Sixsmith [268].

The exponential function and the trigonometric functions (and thus the examples considered in §2.1 and §2.3) are contained in $\mathcal{S}$. An example of a function $f \in \mathcal{B} \setminus \mathcal{S}$ is given by $f(z) = (\sin z)/z$. Fatou's function and Baker's example, which were considered in §2.2 and §2.4, are not in $\mathcal{B}$. For Fatou's function this follows since it has the critical points $z_k := 2\pi i k$ and the critical values $f(z_k) = z_k + 2$, for $k \in \mathbb{Z}$. For Baker's example it can also be shown by elementary means that the set of critical values is unbounded. But it also follows from general considerations. In fact, the result that $I(f) \subset J(f)$ for $f \in \mathcal{B}$ (Theorem 6.1) that we just mentioned yields that both Fatou's function and Baker's example are not in $\mathcal{B}$.

We conclude this introduction to complex dynamics with the definition of two important subclasses of $\mathcal{B}$.

**Definition 3.4.** A transcendental entire function $f$ is called *hyperbolic* if it is in $\mathcal{B}$ and if every point of $\overline{\text{sing}(f^{-1})}$ belongs to the attracting basin of some attracting periodic point.

It is said to be of *disjoint type* if it belongs to $\mathcal{B}$ and if $\overline{\text{sing}(f^{-1})}$ is contained in the immediate attracting basin of an attracting fixed point.

It is not difficult to see that $f$ is hyperbolic if and only if the postsingular set $P(f)$ defined by (3.1) is a compact subset of $F(f)$, and that for a hyperbolic function $f$ the Fatou set $F(f)$ is the union of finitely many attracting basins; see [62, Proposition 2.1]. Clearly, functions of disjoint type are hyperbolic.

The sine function considered in §2.1 is of disjoint type. As explained at the beginning of §2.3, the exponential functions $f_a$ considered there are of disjoint type if $0 < a < 1/e$. On the other hand, if $a < -e$, then $f_a$ has an attracting periodic point of period 2; see [15, Theorem 4.1]. In this case, $f_a$ is hyperbolic but not of disjoint type.

Functions of disjoint type can be characterised as follows [188, Proposition 2.7].

**Proposition 3.5.** *An entire function $f$ is of disjoint type if and only if there is a bounded Jordan domain $D \supset \overline{\text{sing}(f^{-1})}$ such that $f(\overline{D}) \subset D$.*

# 4 The escaping set is not empty

## 4.1 Basic results about the escaping set

The following theorem of Eremenko [113, Theorem 1] is the keystone of this whole survey.

**Theorem 4.1.** *Let $f$ be a transcendental entire function. Then $I(f) \neq \emptyset$.*

In the following sections, we will give two proofs of this fundamental result: one due to Eremenko [113] and another one by Domínguez [109, Theorem G]. Each proof has certain advantages in the sense that the method yields some additional information and allows some extensions. We will discuss this after the proofs.

A fundamental result in complex dynamics is that the Julia set is non-empty. For transcendental entire functions this result is not that easy to show; cf. Bargmann's paper [32]. One method of doing so is provided by Theorem 4.1.

**Theorem 4.2.** *Let $f$ be a transcendental entire function. Then $J(f) \neq \emptyset$. In fact, $J(f)$ contains infinitely many points.*

*Proof.* An application of Picard's theorem to the function
$$h(z) := \frac{f(f(z)) - z}{f(z) - z}$$
shows that $f^2$ has a fixed point $z_1$. (For the details of this argument, which appears already in Fatou's paper [132, p. 346], see [32, Lemma 1].)

On the other hand, there exists $z_2 \in I(f)$ by Theorem 4.1. Since $z_1$ and $z_2$ cannot be in the same Fatou component of $f$, any curve connecting them will meet $J(f)$, so $J(f) \neq \emptyset$.

To prove that $J(f)$ is infinite we note that $I(f)$ contains the orbit of $z_2$ and hence is infinite. Thus the same argument as above will show that $J(f)$ is infinite once we show that the set of points with bounded orbit is infinite.

Suppose that this is not the case. Then $f^2$ has only finitely many fixed points. Assuming without loss of generality that $z_1 = 0$ so that $f^2(0) = 0$, we find that $g(z) := f^2(z)/z$ defines a transcendental entire function $g$ that takes the value 1 only finitely often. By Picard's theorem, $g$ has infinitely many zeros. But a zero of $g$ is a zero of $f^2$ and thus has bounded orbit with respect to $f$. $\square$

Once it is known that $I(f)$ and $J(f)$ are infinite sets, the following result is fairly easy to prove.

**Theorem 4.3.** *Let $f$ be a transcendental entire function. Then $\partial I(f) = J(f)$.*

*Proof.* If a Fatou component intersects $I(f)$, then – by normality – it is completely contained in $I(f)$. This implies that $\partial I(f) \subset J(f)$.

To prove the opposite inclusion, we note that $I(f)$, $\mathbb{C}\backslash I(f)$ and $\partial I(f)$ are all backward invariant. As shown in the proof of Theorem 4.2, the set of points with bounded orbit and hence $\mathbb{C} \setminus I(f)$ are infinite. This implies that $\partial I(f)$ is infinite. The conclusion now follows since every closed, backward invariant set with at least three points contains the Julia set (Theorem 3.1 $(g)$). $\square$

## 4.2 Eremenko's proof that the escaping set is not empty

The main tool in Eremenko's proof is the Wiman–Valiron theory. For a detailed account of this theory, which describes the behaviour of an entire function near points of maximum modulus, we refer to a paper by Hayman [143]. Essentially it says that near such points the function behaves like a certain power.

For a precise statement, let
$$f(z) = \sum_{n=0}^{\infty} a_n z^n$$
be the Taylor expansion of a transcendental entire function $f$. Then $|a_n|r^n \to 0$ as $n \to \infty$ for each $r > 0$ and thus there exists $\nu(r) \in \mathbb{N}_0$ such that
$$|a_{\nu(r)}|r^{\nu(r)} = \max_{n \in \mathbb{N}_0} |a_n|r^n.$$

If there is more than one integer $\nu(r)$ with this property, let $\nu(r)$ be the largest one of those. Then $\nu(r, f)$ is called the *central index* of $f$. For $r > 0$ there exists $z_r$ such that
$$|z_r| = r \quad \text{and} \quad |f(z_r)| = M(r, f) := \max_{|z|=r} |f(z)|. \tag{4.1}$$

Here $M(r, f)$ is called the *maximum modulus* of $f$, and $z_r$ a *point of maximum modulus*.

We will see that if $r$ is outside a certain exceptional set, then $f$ behaves near $z_r$ like $z \mapsto c_r z^{\nu(r)}$ for some $c_r \in \mathbb{C} \setminus \{0\}$. To formulate the condition for the exceptional set we say that a measurable subset $F$ of $[1, \infty)$ has *finite logarithmic measure* if
$$\int_F \frac{dr}{r} < \infty.$$

The main result of Wiman–Valiron theory is the following; see [143, Section 6.1].

**Theorem 4.4.** *Let $f$ be a transcendental entire function and let $1/2 < \gamma < 1$. Then there exists a subset $F$ of $[1, \infty)$ that has finite logarithmic measure such if $r \in [1, \infty) \setminus F$ and $z_r$ is a point of maximum modulus, then*
$$f(z) \sim \left(\frac{z}{z_r}\right)^{\nu(r)} f(z_r) \quad \text{for } z \in D\left(z_r, \frac{r}{\nu(r)^\gamma}\right) \tag{4.2}$$
*as $r \to \infty$, $r \notin F$.*

Writing $z = z_r e^\tau$ we may write the conclusion as
$$f(z_r e^\tau) \sim e^{\nu(r)\tau} f(z_r) \quad \text{for } |\tau| \leq \frac{1}{\nu(r)^\gamma}. \tag{4.3}$$

If (4.2) holds for every $\gamma > 1/2$ as $r \to \infty$, $r \notin F$, then so does (4.3), and vice versa.

For $w \in \mathbb{C}$ and $\tau = \tau_w := w/\nu(r)$ the right hand side of (4.3) takes the form $e^w f(z_r)$. An application of Rouché's theorem now shows that the left hand side also takes the value $e^w f(z_r)$ in a small disc around $z_r e^\tau$. More precisely, with
$$g(\tau) = f(z_r e^\tau) - e^w f(z_r) \quad \text{and} \quad h(\tau) = e^{\nu(r)\tau} f(z_r) - e^w f(z_r) = \left(e^{\nu(r)(\tau-\tau_w)} - 1\right) e^w f(z_r)$$
we find that if $0 < \delta < 1$ and $|\tau - \tau_w| = \delta/\nu(r)$, then
$$|g(\tau) - h(\tau)| = \left|f(z_r e^\tau) - e^{\nu(r)\tau} f(z_r)\right| = o(|f(z_r)|) < |h(\tau)|$$

for large $r \notin F$. Since $h(\tau_w) = 0$ there exists $\tau \in D(\tau_w, \delta/\nu(r))$ such that $g(\tau) = 0$.
Thus $f(z_r e^\tau) = e^w f(z_r)$. Given $R > 0$ we find that for all $w \in D(0, R)$ there exists
$\tau \in D(0, (R+\delta)/\nu(r))$ with $f(z_r e^\tau) = e^w f(z_r)$, provided $r \notin F$ is large enough.

We will use the abbreviation
$$V(r) := D\left(z_r, \frac{r}{\nu(r)^\gamma}\right).$$

Given $r$ there may be several points $z_r$ satisfying (4.1). We can make the definition of $V(r)$ unique by taking, e.g., the point $z_r$ with the smallest argument, but for our purposes it will not matter which one of the points $z_r$ we choose. Noting that $z_r e^\tau \in V(r)$ for $\tau \in D(0, (R+\delta)/\nu(r))$ and large $r$ we see that $f(V(r)) \supset f(z_r) \exp D(0, R)$. This yields the following result.

**Lemma 4.5.** *Let $f$, $\gamma$, $F$ and $z_r$ be as in Theorem 4.4. Then, given $K > 1$, we have*
$$f(V(r)) \supset \mathrm{ann}\left(\frac{1}{K}M(r, f), KM(r, f)\right) \tag{4.4}$$

*if $r \notin F$ is large enough.*

*Proof of Theorem 4.1 (i.e., $I(f) \neq \emptyset$) following Eremenko.* Fix $\gamma \in (1/2, 1)$ and $K > 4$ and let $F$ be the exceptional set arising in Theorem 4.4. We may choose $R > 1$ such that the following conditions are satisfied:

(i) (4.4) holds for $r \geq R$, $r \notin F$,

(ii) $\int_{F \cap [R, \infty)} dt/t \leq 1$,

(iii) $M(r, f) \geq 2r$ for $r \geq R$,

(iv) $\nu(r)^\gamma \geq 2$ and hence $r/\nu(r)^\gamma \leq r/2$ for $r \geq R$.

Let now $r_0 \geq R$, $r_0 \notin F$. By (i) we have (4.4) for $r = r_0$ so that there exists a compact subset $C_0$ of $V(r_0)$ such that with $M_1 := M(r_0, f)$ we have $f(C_0) = \overline{\mathrm{ann}}(M_1/4, 4M_1)$. By (iii) we have $M_1 \geq 2r_0 \geq 2R > 2$. Hence (ii) yields that there exists $r_1 \in [M_1, 3M_1] \setminus F$. Thus (4.4) also holds for $r = r_1$. Moreover, (iv) guarantees that $V(r_1) \subset \overline{\mathrm{ann}}(M_1/4, 4M_1)$.

The same argument as before yields that there exists a compact subset $C_1$ of $V(r_1)$ such that with $M_2 := M(r_1, f)$ we have $f(C_1) = \overline{\mathrm{ann}}(M_2/4, 4M_2)$. Inductively we obtain a sequence $(r_k)$ satisfying $r_{k+1} \geq M_{k+1} := M(r_k, f) \geq 2r_k$ and a sequence $(C_k)$ of compact sets satisfying $C_k \subset V(r_k) \subset \overline{\mathrm{ann}}(M_k/4, 4M_k)$ such that $f(C_k) = \overline{\mathrm{ann}}(M_{k+1}/4, 4M_{k+1})$.

For $n \in \mathbb{N}$ we now put
$$A_n := \{z \in C_0 \colon f^k(z) \in C_k \text{ for } 0 \leq k \leq n\}.$$

Then $A_n$ is compact and $A_{n+1} \subset A_n$ for all $n \in \mathbb{N}$. Hence
$$A := \bigcap_{n=1}^{\infty} A_n \neq \emptyset.$$
For $z \in A$ we have
$$f^n(z) \in C_n \subset \overline{\operatorname{ann}}\left(\frac{1}{4}M_n, 4M_n\right). \tag{4.5}$$
Since $r_{n+1} \geq 2r_n$ we have $r_n \to \infty$. It follows that $M_n = M(r_{n-1}, f) \to \infty$ so that $z \in I(f)$. □

## 4.3 Further results obtained by Eremenko's method

Eremenko [113, Theorem 2] used his method also to prove the following result.

**Theorem 4.6.** *Let $f$ be a transcendental entire function. Then $I(f) \cap J(f) \neq \emptyset$.*

*Proof of Theorem 4.6.* We distinguish two cases.

*Case 1.* There is no multiply connected Fatou component. We use the notation of the proof of Theorem 4.1 (following Eremenko) in the previous section. Since there is no multiply connected Fatou component, $\operatorname{ann}(M_k/4, 4M_k)$ intersects $J(f)$ for large $k$. The complete invariance of the Julia set (Theorem 3.1 (b)) now implies that $C_k$ intersects $J(f)$ for large $k$. Using the complete invariance of $J(f)$ again we find that $A_n$ also intersect $J(f)$ for large $n$. Hence $A$ intersects $J(f)$. The conclusion follows since $A \subset I(f)$.

*Case 2.* There is a multiply connected Fatou component $U$. We will see in Theorem 6.12 that then $\overline{U} \subset I(f)$. Thus $\partial U \subset I(f) \cap J(f)$. □

Next we will describe a modification of the Wiman–Valiron theory that, together with Eremenko's method to prove that $I(f) \neq \emptyset$, leads to further results about the escaping set.

First we note that by Hadamard's three circles theorem, $\log M(r, f)$ is a convex function of $\log r$. Thus the right derivative $a(r)$ of $\log M(r, f)$ with respect to $\log r$ exists and is increasing. Actually, by a result of Blumenthal (see [281, Section II.3]), $M(r, f)$ is differentiable except for a discrete set of $r$-values. For $r$ outside this discrete set we have
$$a(r) = \frac{d \log M(r, f)}{\log r} = \frac{z_r f'(z_r)}{f(z_r)}.$$
Macintyre [178] developed an alternative approach to Wiman–Valiron theory that, essentially, yields that (4.2) holds with $\nu(r)$ replaced by $a(r)$. Even though the Taylor series expansion is not used, Macintyre's argument still requires that $f$ is entire.

A version of Wiman–Valiron theory that is valid for certain functions holomorphic in subdomains of $\mathbb{C}$ was developed by Bergweiler, Rippon and Stallard [72]. In order to describe their result we will use the following definition, which, with use in later sections in mind, is stated more generally than required for our immediate purposes.

**Definition 4.7.** Let $D$ be an unbounded domain in $\mathbb{C}$ with unbounded complement, and let $f\colon \overline{D} \to \mathbb{C}$ be continuous. Then $D$ is called a *direct tract* of $f$ if $f$ is holomorphic in $D$ and if there exists $R > 0$ such that $|f(z)| = R$ for $z \in \partial D$ and $|f(z)| > R$ for $z \in D$.

If, in addition, $f\colon D \to \mathbb{C} \setminus \overline{D}(0, R)$ is a universal covering, then $D$ is called a *logarithmic tract* of $f$.

If $f$ is a transcendental entire function and $R > 0$, then every connected component of $\{z\colon |f(z)| > R\}$ is a direct tract of $f$. This is not the case for meromorphic functions.

If $D$ is a direct tract of $f$, then $f(z) \to \infty$ as $z \to \infty$ along some curve in $D$. (This follows from the standard proof of Iversen's theorem as given, e.g., in Nevanlinna's book [197, p. 291].) We thus see that if $1/(f - w)$ has a direct (or logarithmic) tract $D$ for some $w \in \mathbb{C}$, then $w$ is an asymptotic value of $f$ and hence a singularity of the inverse of $f$. In this case we call $D$ a *direct (or logarithmic) tract over $w$* and say that $f^{-1}$ has a *direct (or logarithmic) singularity over $w$*.

Let now $D$ be a direct tract of $f$, put

$$M_D(r) := \max_{z \in D, |z|=r} |f(z)|$$

and, as in (4.1), choose $z_r \in D$ with $|z_r| = r$ and $|f(z_r)| = M_D(r)$. Again, $\log M_D(r)$ is a convex function of $\log r$. Denoting the right derivative of $\log M_D(r)$ with respect to $\log r$ by $a_D(r)$, we find as before that $a_D(r)$ is increasing and, except for a discrete set of $r$-values, we have

$$a_D(r) = \frac{z_r f'(z_r)}{f(z_r)}.$$

It is shown in [72, Theorem 2.2] that Theorem 4.4 extends to this setting; that is, the conclusion of this theorem remains valid with $\nu(r)$ replaced by $a_D(r)$. In particular, the disc $D(z_r, r/a_D(r)^\gamma)$ is contained in the direct tract $D$.

Once this extension of Wiman–Valiron theory to direct tracts is known, we obtain the following result [72, Theorem 3.1] with the method of Eremenko described in the previous section.

**Theorem 4.8.** *Let $D$ be a direct tract of $f$. Then there exists $z \in D$ such that $f^n(z) \in D$ for all $n \in \mathbb{N}$ and $f^n(z) \to \infty$ as $n \to \infty$.*

The same method also yields the following result for functions with several tracts.

**Theorem 4.9.** *Let $D_1, \ldots, D_N$ be a direct tracts of $f$ and let $(\sigma(n))_{n \in \mathbb{N}_0}$ be a sequence in $\{1, \ldots, N\}$. Then there exists $z \in D_{\sigma(0)}$ such that $f^n(z) \in D_{\sigma(n)}$ for all $n \in \mathbb{N}$ and $f^n(z) \to \infty$ as $n \to \infty$.*

## 4.4 Domínguez's proof that the escaping set is not empty

One tool used in Domínguez's proof is the following result of Bohr [142, Theorem 6.9].

**Lemma 4.10.** *For every $\rho \in (0,1)$ there exists a positive constant $c$ such that if $h \colon \mathbb{D} \to \mathbb{C}$ is holomorphic with $h(0) = 0$ and*
$$M(\rho, h) \geq 1,$$
*then $\partial D(0, R) \subset h(\mathbb{D})$ for some $R \geq c$.*

*Proof.* Assuming that such $c$ does not exist for some $\rho \in (0, 1)$, there exists a sequence $(h_n)$ of holomorphic functions satisfying the hypotheses of this lemma such that
$$c_n := \sup\{R > 0 \colon \partial D(0, R) \subset h_n(\mathbb{D})\} \to 0$$
as $n \to \infty$. Then there exists $a_n, b_n \in \mathbb{C}$ with $|a_n| = c_n$ and $|b_n| = 2c_n$ such that $h_n(z) \neq a_n$ and $h_n(z) \neq b_n$ for all $z \in \mathbb{D}$. The functions $g_n \colon \mathbb{D} \to \mathbb{C}$ defined by
$$g_n(z) := \frac{h_n(z) - a_n}{b_n - a_n}$$
thus satisfy $g_n(z) \neq 0$ and $g_n(z) \neq 1$ for all $z \in \mathbb{D}$. Hence the $g_n$ form a normal family. But
$$|g_n(0)| = \left|\frac{-a_n}{b_n - a_n}\right| \leq \frac{|a_n|}{|b_n| - |a_n|} = 1$$
while
$$M(\rho, g_n) \geq \frac{M(\rho, h_n) - |a_n|}{|b_n| + |a_n|} \geq \frac{1 - c_n}{3c_n}$$
and thus $M(\rho, g_n) \to \infty$. This contradicts the normality of the sequence $(g_n)$. □

It was shown by Hayman [140, Theorem VIII] that the best possible constant $c$ in Lemma 4.10 is given by
$$c = \frac{(1-\rho)^2}{4\rho},$$
the extremal function being a suitably scaled Koebe function, but we will not need this result.

For a bounded domain $G$ in $\mathbb{C}$ we denote by $U(G)$ the unbounded connected component of $\mathbb{C} \setminus G$. We call $T(G) = \mathbb{C} \setminus U(G)$ the *topological hull* of $G$. Thus $T(G)$ is the union of $G$ and the bounded connected components of its complement. Informally, $T(G)$ is obtained from $G$ by "filling the holes". A (bounded) domain $G$ is simply connected if and only if $T(G) = G$. Often, the notation $\widetilde{G}$ is used instead of $T(G)$.

We also put
$$\mu(G) := \min\{|z| \colon z \in U(G)\}.$$
If $\mu(G) > 0$, then
$$\mu(G) = \sup\{R > 0 \colon D(0, R) \subset T(G)\}.$$
It follows from Lemma 4.10 that if $h$ is as there, then $\mu(h(\mathbb{D})) \geq c$.

Applying this to
$$h(z) = \frac{f(rz) - f(0)}{\max_{|z|=\rho r} |f(z) - f(0)|}$$
with an entire function $f$ leads to the following result.

**Lemma 4.11.** *Let $f$ be entire, $r > 0$ and $0 < \rho < 1$. Let $c > 0$ be as in Lemma 4.10. Then*
$$\mu(f(D(0,r))) \geq cM(\rho r, f) - (1+c)|f(0)|.$$

We will also need the following lemma, whose simple proof we omit.

**Lemma 4.12.** *Let $f$ be entire (and non-constant) and let $G$ be a bounded domain in $\mathbb{C}$. Then*
$$f(T(G)) \subset T(f(G)) \tag{4.6}$$
*and*
$$\partial T(f(G)) \subset f(\partial T(G)). \tag{4.7}$$

*Proof of Theorem 4.1 (i.e., $I(f) \neq \emptyset$) following Domínguez.* We fix $\rho \in (0,1)$. It follows from Lemma 4.11 that if $r$ is large enough, say $r \geq r_0$, then
$$\mu(f(D(0,r))) \geq cM(\rho r, f) - (1+c)|f(0)| \geq \frac{c}{2}M(\rho r, f) \geq 2r. \tag{4.8}$$

By (4.6) and the definition of $\mu(G)$ we have
$$\mu(f(G)) = \mu(T(f(G))) \geq \mu(f(T(G))) \geq \mu(f(D(0, \mu(G)))) \tag{4.9}$$

for every bounded domain $G$. We deduce from (4.8) and (4.9) that that if $\mu(G) \geq r_0$, then
$$\mu(f(G)) \geq 2\mu(G) \tag{4.10}$$

We now fix a bounded domain $G_0$ satisfying $\mu(G_0) \geq r_0$, for example $G_0 = D(0, r_0)$. Putting $G_n := f(G_{n-1}) = f^n(G)$ for $n \in \mathbb{N}$ we conclude that
$$\mu(G_n) \geq 2\mu(G_{n-1}) \geq \cdots \geq 2^n \mu(G_0) = 2^n r_0. \tag{4.11}$$

Moreover, Lemma 4.12 yields that
$$\partial T(G_n) \subset f(\partial T(G_{n-1})) \subset \cdots \subset f^n(\partial G_0) = f^n(\partial D(0, r_0)).$$

Hence
$$A_n := f^{-n}(\partial T(G_n)) \cap \partial D(0, r_0) \neq \emptyset. \tag{4.12}$$

Moreover, $A_n$ is compact and $A_{n+1} \subset A_n$ for all $n \in \mathbb{N}_0$. Thus
$$A := \bigcap_{n=0}^{\infty} A_n \neq \emptyset.$$

For $z \in A$ we have $f^n(z) \in \partial T(G_n)$. Hence $|f^n(z)| \geq 2^n r_0$ by (4.11) so that $z \in I(f)$. □

## 4.5 Further results obtained by Domínguez's method

Let $r_0$ be as in the previous proof and let $\gamma$ be a curve surrounding $D(0, r_0)$. Taking $G_0$ as the interior of $\gamma$ we obtain the following result already noted by Domínguez [109, Theorem G]; see also [236, Theorem 2].

**Theorem 4.13.** *Let $f$ be a transcendental entire function. Then there exists $r_0$ such that every curve surrounding $D(0, r_0)$ intersects $I(f)$.*

More generally, the argument yields the following result.

**Theorem 4.14.** *Let $f$ be a transcendental entire function and let $D$ be a bounded domain intersecting $J(f)$. Then $\partial D \cap I(f) \neq \emptyset$.*

*Proof.* The blowing-up property of the Julia set (Theorem 3.1 (f)) yields that for sufficiently large $r_0$ there exists $N \in \mathbb{N}$ such that $f^N(D) \supset \partial D(0, r_0)$ and hence $T(f^N(D)) \supset D(0, r_0)$. Applying Domínguez's proof with $G_0 := T(f^N(D))$ we see that $\partial T(f^N(D))$ intersect $I(f)$. It now follows from (4.7) that $f^N(\partial T(D))$ and hence $\partial T(D)$ intersects $I(f)$. Since $\partial T(D) \subset \partial D$ the conclusion follows. $\square$

Domínguez's method also gives an alternative proof of Theorem 4.6 in the case that $f$ has no multiply connected wandering domain. We only have to replace $A_n$ in (4.12) by $A_n \cap J(f)$.

Rippon and Stallard [239] observed that the escaping points constructed by Domínguez's method are contained in unbounded connected components of $I(f)$. We will discuss this in detail in §7.3.

Domínguez's proof [109, Theorem G] extends to quasiregular mappings [63]. This and other results about the escaping set of quasiregular mappings will be discussed briefly in §10.1; see also [50, 70] for iteration of quasiregular maps.

# 5 Escape rates

## 5.1 The fast escaping set

Let $M^n(r, f)$ the $n$-th iterate of $M(r, f)$ with respect to the first variable; that is, $M^1(r, f) = M(r, f)$ and $M^{n+1}(r, f) = M(M^n(r, f), f)$ for $n \geq 1$. The maximum modulus principle implies that if $R > 0$ and $z \in \overline{D}(0, R)$, then

$$|f^n(z)| \leq M(R, f^n) \leq M^n(R, f).$$

Thus the iterated maximum modulus gives an upper bound on how fast points in $I(f)$ can tend to infinity.

In turn, we will see that for each $R > 0$ there exists $z \in \mathbb{C}$ such that

$$|f^n(z)| \geq M^n(R, f). \tag{5.1}$$

So in some sense such points escape to infinity "as fast as possible". Of course, in order to guarantee that a point satisfying (5.1) is escaping we choose $R$ such that $M^n(R, f) \to \infty$ as $n \to \infty$. This is satisfied if

$$M(r, f) > r \quad \text{for } r \geq R. \tag{5.2}$$

We call a nonnegative number $R$ satisfying (5.2) an *admissible radius* (for $f$).

**Definition 5.1.** Let $f$ be a transcendental entire function and let $R$ be an admissible radius. Then

$$A(f) := \{z \colon \text{there exists } L \in \mathbb{N}_0 \text{ such that } |f^{n+L}(z)| \geq M^n(R, f) \text{ for } n \in \mathbb{N}_0\} \tag{5.3}$$

is called the *fast escaping set* of $f$.

If $R$ is an admissible radius and $R' > R$, then there exists $m \in \mathbb{N}$ such that $M^m(R, f) \geq M(R', f)$. This implies that the definition of $A(f)$ is independent of the value of $R$, as long as $R$ is admissible. It is also easily seen that $A(f)$ is completely invariant, as is $I(f)$.

As already mentioned, fast escaping points exist so that we have the following result.

**Theorem 5.2.** *Let $f$ be a transcendental entire function. Then $A(f) \neq \emptyset$.*

This result can be proved with the method of Eremenko as well as the one of Domínguez. In both cases we need the following lemma.

**Lemma 5.3.** *Let $f$ be a transcendental entire function and $C > 1$. Then*

$$\lim_{r \to \infty} \frac{M(Cr, f)}{M(r, f)} = \infty.$$

This lemma can easily be deduced from Hadamard's three circles theorem, which says that $\log M(r, f)$ is convex in $\log r$ and a standard generalisation of Liouville's theorem saying that

$$\lim_{r \to \infty} \frac{\log M(r, f)}{\log r} = \infty \tag{5.4}$$

for transcendental entire $f$.

For later use we also note another consequence of Hadamard's three circles theorem observed by Rippon and Stallard [241, Lemma 2.2].

**Lemma 5.4.** *Let $f$ be a transcendental entire function. Then there exists $R > 0$ such that if $r > R$ and $c > 1$, then*

$$M(r^c, f) \geq M(r, f)^c.$$

*Proof of Theorem 5.2 by Eremenko's method.* Eremenko's proof yields that there exists $z \in \mathbb{C}$ satisfying (4.5); that is,
$$|f^n(z)| \geq \frac{1}{4}M_n = \frac{1}{4}M(r_{n-1}, f).$$

Here the sequence $(r_n)$ satisfies $r_n \geq M(r_{n-1}, f)$. Thus $r_{n-1} \geq M^{n-1}(r_0, f)$ so that
$$|f^n(z)| \geq \frac{1}{4}M^n(r_0, f).$$

Choosing $R$ sufficiently large and $r_0 > R$, we can now deduce (5.1) from Lemma 5.3. □

*Proof of Theorem 5.2 by Domínguez's method.* Let $\rho$ and $c$ be as in Lemma 4.10 and put
$$H(r) := \frac{c}{2}M(\rho r, f).$$

The arguments we used to derive (4.10) actually yield that $\mu(f(G)) \geq H(\mu(G))$. Instead of (4.11) we now obtain
$$\mu(G_n) \geq H^n(r_0). \tag{5.5}$$

This implies that there exists $z \in \mathbb{C}$ such that
$$|f^n(z)| \geq H^n(r_0). \tag{5.6}$$

It follows from Lemma 5.3 that if $r$ is sufficiently large, then
$$H(r) \geq K(r) := \frac{1}{\rho^2}M(\rho^2 r, f). \tag{5.7}$$

For large $r_0$ the last two inequalities yield that
$$|f^n(z)| \geq K^n(r_0) = \frac{1}{\rho^2}M^n(\rho^2 r_0, f) \geq M^n(\rho^2 r_0, f).$$

Given a (large) value of $R$ we again obtain (5.1) by choosing $r_0 := R/\rho^2$. □

The same reasoning as in §4.5 shows that the conclusion of Theorem 4.13 holds with $I(f)$ replaced by $A(f)$ so that we have the following result.

**Theorem 5.5.** *Let $f$ be a transcendental entire function. Then there exists $r_0$ such that every curve surrounding $D(0, r_0)$ intersects $A(f)$.*

We will use Theorem 5.5 in §6.3, but note that Theorem 7.24 below is a considerably stronger result.

We also have the following analogue of Theorem 4.3.

**Theorem 5.6.** *Let $f$ be a transcendental entire function. Then $\partial A(f) = J(f)$.*

Here the inclusion that $J(f) \subset \partial A(f)$ follows by the same argument as in the proof of Theorem 4.3. The other inclusion requires the result that if a Fatou component intersects $A(f)$, then it is contained in $A(f)$. This is Theorem 6.5 below.

It follows from (5.4) that if $K > 0$, then $\log M(r, f) \geq K \log r$ for large $r$. For an admissible radius $R$ we thus have $\log M^{n+1}(R, f) \geq K \log M^n(R, f)$ for large $n$, say $n \geq N$. We conclude that $\log M^n(R, f) \geq K^{n-N} \log M^N(R, f)$ for $n \geq N$. As $K$ can be taken arbitrarily large, we obtain the following result.

**Theorem 5.7.** *Let $f$ be a transcendental entire function. If $z \in A(f)$, then*

$$\lim_{n \to \infty} \frac{\log \log |f^n(z)|}{n} = \infty.$$

The definition of the fast escaping set $A(f)$ via (5.3) is due to Rippon and Stallard [244]. However, with a slightly different definition this set had been introduced already before by Bergweiler and Hinkkanen [66]. They worked with the definition

$$A(f) := \{z \colon \text{there exists } L \in \mathbb{N}_0 \text{ such that } |f^{n+L}(z)| \geq M(R, f^n) \text{ for } n \in \mathbb{N}_0\}, \quad (5.8)$$

with $R$ chosen so large that

$$\lim_{n \to \infty} \frac{\log \log M(R, f^n)}{n} = \infty.$$

Rippon and Stallard [244, Corollary 2.7] showed that the definitions via (5.3) and (5.8) are equivalent.

Earlier they had given in [239, Lemma 2.4] the following important characterisation of $A(f)$. Recall here that $T(G)$ denotes the topological hull of a domain $G$.

**Theorem 5.8.** *Let $f$ be a transcendental entire function and let $D$ be an open disc intersecting $J(f)$. Then*

$$A(f) = \{z \colon \text{there exists } L \in \mathbb{N}_0 \text{ such that } f^{n+L}(z) \notin T(f^n(D)) \text{ for } n \in \mathbb{N}_0\}. \quad (5.9)$$

*Sketch of proof.* Let $B(f)$ be the set on the right hand side of (5.9). To prove that $A(f) \subset B(f)$, choose an admissible radius $R$ such that $D \subset D(0, R)$. It then follows that $T(f^n(D)) \subset T(f^n(D(0, R))) \subset D(0, M^n(R, f))$. We deduce that if $|f^n(z)| \geq M^n(R, f)$, then $f^n(z) \notin f^n(D)$. This yields that $A(f) \subset B(f)$.

To prove the opposite inclusion, suppose first that $D \supset D(0, r_0)$, with $r_0$ taken as large as in the proof that $A(f) \neq \emptyset$ by Domínguez's method. It follows from inequalities (5.5), (5.6) and (5.7) in this proof that if $f^n(z) \notin T(f^n(D))$, then

$$|f^n(z)| \geq M^n(\rho^2 r_0, f).$$

This yields that $B(f) \supset A(f)$ if $D \supset D(0, r_0)$. The general case follows from the blowing-up property of the Julia set (Theorem 3.1 (f)). □

In [244], which is the first systematic study of $A(f)$, Rippon and Stallard introduced, for an admissible radius $R$ and $L \in \mathbb{Z}$, the *levels*

$$A_R^L(f) := \{z \in \mathbb{C} \colon |f^n(z)| \geq M^{n+L}(R, f) \text{ if } n, n + L \in \mathbb{N}\}. \tag{5.10}$$

Then

$$A(f) = \bigcup_{L \in \mathbb{N}} A_R^{-L}(f). \tag{5.11}$$

The maximum modulus principle shows that if $0 \leq m < n$ and $m + L \geq 0$, and if $|f^n(z)| \geq M^{n+L}(R, f)$, then $|f^m(z)| \geq M^{m+L}(R, f)$. Thus, as noted by Rippon and Stallard, the definition of the $A_R^L(f)$ remains unchanged if one requires $|f^n(z)| \geq M^{n+L}(R, f)$ not for $n \geq 1$ and $n + L \geq 1$, but for $n \geq 0$ and $n + L \geq 0$. It also suffices to require this only for large $n$.

We will often consider the case $L = 0$ and thus put

$$A_R(f) := A_R^0(f). \tag{5.12}$$

With

$$C_n := \{z \colon |f^n(z)| \geq M^n(R, f)\} \tag{5.13}$$

we then have $C_{n+1} \subset C_n$ and

$$A_R(f) = \bigcap_{n=0}^{\infty} C_n.$$

Figure 7 shows the sets $C_n \setminus C_{n+1}$ in grey, getting darker as $n$ increases. The disc $D(0, R) = \mathbb{C} \setminus C_0$ is shown in white.

The left picture shows these sets for $f_1(z) = e^z/2$ and $R = 1$. The range displayed is $-2 \leq \operatorname{Re} z \leq 7$ and $|\operatorname{Im} z| \leq 9$. The right picture shows these sets for

$$f_2(z) = \frac{1}{2}\left(\cos z^{1/4} + \cosh z^{1/4}\right) = 1 + \frac{z}{4!} + \frac{z^2}{8!} + \ldots. \tag{5.14}$$

This function was considered by Rippon and Stallard ([240, §6] and [244, §8]). Here $R = 14500$. (The radius $R$ is admissible only for $R > R_0 \approx 14456.61$.) The range displayed is $-75000 \leq \operatorname{Re} z \leq 25000$ and $|\operatorname{Im} z| \leq 50000$.

The pictures suggest that $A_R(f_1)$ is disconnected while $A_R(f_2)$ is connected, and separates $0$ from $\infty$. This is indeed the case. For $f_1$ this follows from the fact that, for each $k \in \mathbb{Z}$, the half-line $\{x + 2\pi i k \colon x \geq 1\}$ is contained in $A_1(f_1)$, while the line $\{x + (2k+1)\pi i \colon x \in \mathbb{R}\}$ does not intersect $A_1(f_1)$. The properties claimed about $A_R(f_2)$ follow from Corollary 7.38 below. We refer to §7 for a detailed discussion of the connectivity properties of $I(f)$, $A(f)$ and $A_R(f)$.

Note that the levels $A_R^{-L}(f)$ are closed. This turns out to be advantageous in some arguments. One immediate consequence of (5.11) is thus that $A(f)$ is $F_\sigma$. In contrast, Rempe [225] has shown that $I(f)$ is never $F_\sigma$.

While we argued above that points in $A(f)$ escape to infinity at the maximal rate, Sixsmith [266] considered the *maximally fast escaping set*

$$A'(f) := \{z \colon \text{there exists } N \in \mathbb{N}_0 \text{ such that } |f^{n+1}(z)| = M(|f^n(z)|, f) \text{ for } n \geq N\}.$$

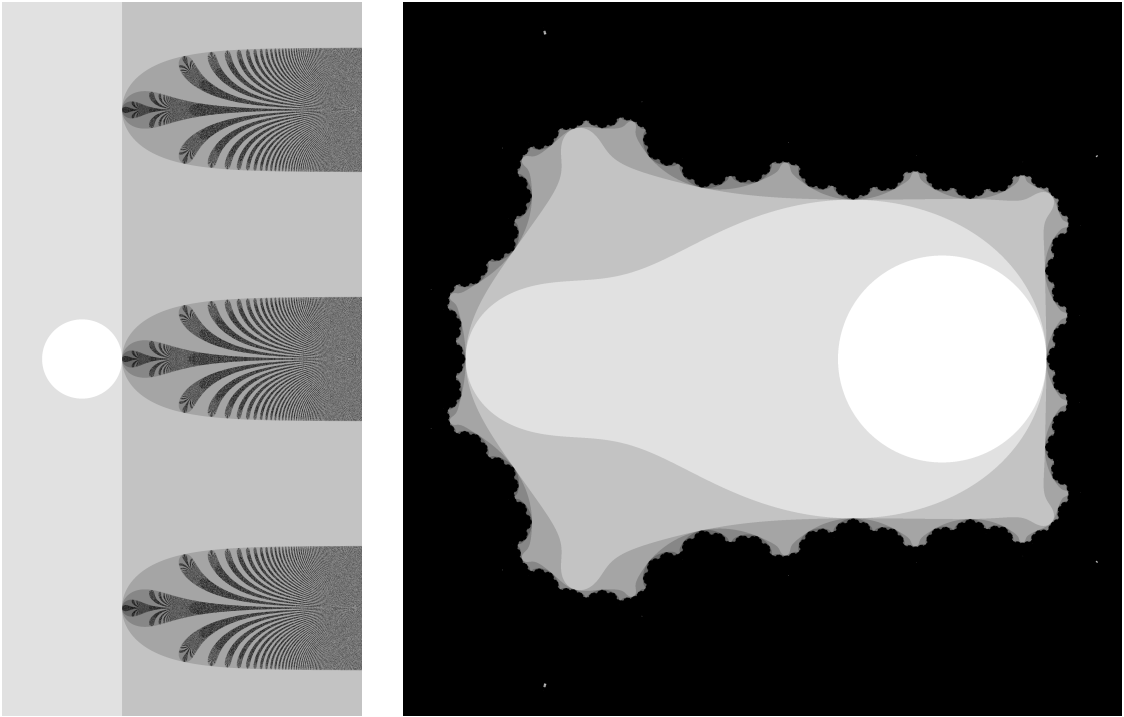

Figure 7: The sets $C_n \setminus C_{n+1}$ for two entire functions.

Clearly, $A'(f) \subset A(f)$. He also considered the *non-maximally fast escaping set* $A''(f) := A(f) \setminus A'(f)$. Sixsmith proved various properties of these two sets. For example, he proved that $A'(f)$ is contained in a countable union of curves each of which is analytic except possibly at its endpoints. On the other hand, he gave an example where $A'(f)$ is uncountable and totally disconnected.

## 5.2 Why do escape rates matter?

Already before the fast escaping set was introduced in [66], and thoroughly studied in [239, 244], points escaping to infinity at a certain rate were considered in several papers – as a tool in the proof of other results. We list some results where escape rates were relevant.

Implicitly, escape rates occurred already in Baker's papers [12, 13] showing that there exists an entire function with a wandering domain and that multiply connected Fatou components of transcendental entire functions are always wandering. Essentially, he showed that points in a multiply connected Fatou component must escape at a certain minimal rate, which is faster than the maximal escape rate in an invariant domain; see Theorems 6.9 and 6.10 below.

McMullen [185] used points escaping at a certain rate to prove the following results.

**Theorem 5.9.** *Let $\alpha, \beta \in \mathbb{C}$, $\alpha \neq 0$. Then $J(\sin(\alpha z + \beta))$ has positive area.*

**Theorem 5.10.** *Let $\lambda \in \mathbb{C} \setminus \{0\}$. Then $J(\lambda e^z)$ has Hausdorff dimension 2.*

In his proof he constructed, for $f(z) = \sin(\alpha z + \beta)$ or $f(z) = \lambda e^z$, a set $E$ of points $z$ satisfying $|f^k(z)| > c \exp^k(1)$ for some $c > 0$ and all $k \in \mathbb{N}$. He then estimated the area (or Hausdorff dimension) of $E$ from below. To obtain the corresponding estimates for $J(f)$ he showed that $E \subset J(f)$. This he deduced from the result [185, Appendix] that, for the above functions $f$,
$$\log \log |f^n(z)| = O(n) \tag{5.15}$$
if $z \in F(f)$.

We note here that for the above functions $f$ we actually have $I(f) \cap F(f) = \emptyset$; see Theorem 6.1 below. Moreover, we will see in Theorem 6.21 that if $f$ is a transcendental entire function and $z$ is in a preperiodic Fatou component of $f$, then (5.15) can be improved to
$$\log |f^n(z)| = O(n).$$
A detailed discussion of escaping points in $F(f)$ will be given in §6

Bergweiler [48] used escape rates to study Julia sets of holomorphic self-maps of the punctured plane $\mathbb{C} \setminus \{0\}$. Many of the results of complex dynamics described in §3 hold, with appropriate modifications, also for such maps. From the large literature on this topic we mention only a paper by Martí-Pete [179], which is concerned with the escaping set of such maps.

For a holomorphic map $g \colon \mathbb{C} \setminus \{0\} \to \mathbb{C} \setminus \{0\}$ there exists an entire function $f$ such that $\exp \circ f = g \circ \exp$. In other words, the diagram

$$\begin{array}{ccc} \mathbb{C} & \xrightarrow{f} & \mathbb{C} \\ {\scriptstyle \exp}\downarrow & & \downarrow{\scriptstyle \exp} \\ \mathbb{C} \setminus \{0\} & \xrightarrow{g} & \mathbb{C} \setminus \{0\} \end{array} \tag{5.16}$$

commutes.

The result of [48] is the following.

**Theorem 5.11.** *If $g \colon \mathbb{C} \setminus \{0\} \to \mathbb{C} \setminus \{0\}$ is holomorphic and if $f$ is a transcendental entire function satisfying $\exp \circ f = g \circ \exp$, then $J(f) = \exp^{-1} J(g)$.*

A crucial role in the proof is played by the set
$$I'(f) := \left\{ z \colon \lim_{n \to \infty} \frac{\log |f^{n+1}(z)|}{\log |f^n(z)|} = \infty \right\}.$$

Points in this set escape at a certain rate, and we briefly indicate why this escape rate matters in the proof. It is easy to see that if $z_0$ is a repelling periodic point of $f$, then $\exp z_0$ is a repelling periodic point of $g$. Noting that the Julia set is the closure of the set of repelling periodic points (Theorem 3.1 (d)), this yields that $\exp^{-1} J(g) \supset J(f)$. To prove the opposite inclusion, it is tempting to start with a repelling periodic point $\exp z_1$ of $g$. However, $z_1$ need not be a periodic point of $f$. But it can be shown that $\log |f^n(z_1)| = O(n)$ as $n \to \infty$. On the other hand, if $\exp z_2 \in I'(g)$, then $z_2 \in I(f)$ and $\log |f^n(z_2)| \neq O(n)$. The different escape rates of $z_1$ and $z_2$ are then used to show that they cannot be in the same Fatou component.

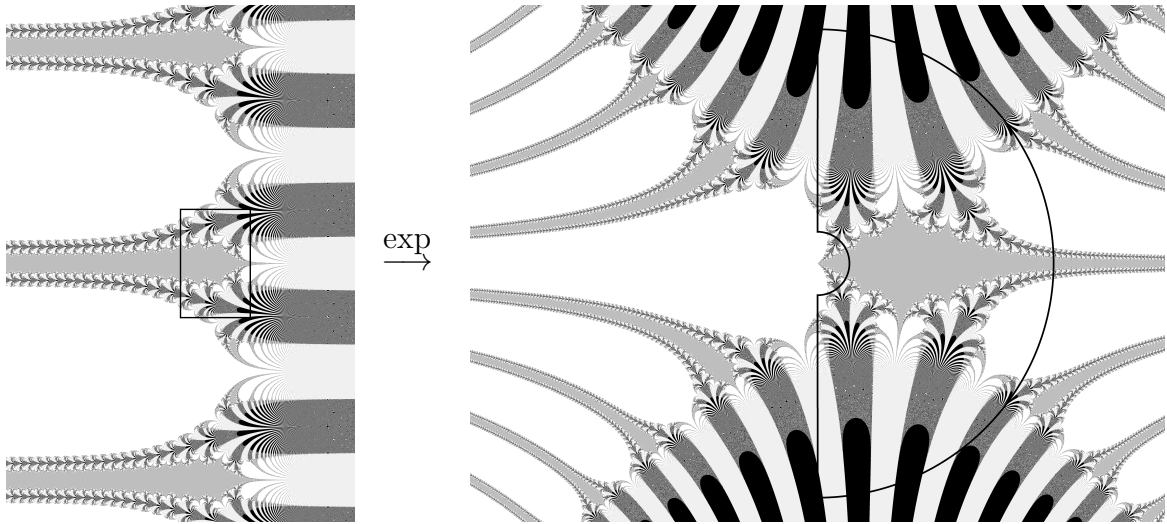

Figure 8: Julia sets of the functions $f$ and $g$ given by (5.17).

We note that (5.16) holds for $g(z) = z^2$ and $f(z) = 2z$, or more generally for $f(z) = 2z + 2\pi i k$ with $k \in \mathbb{Z}$. Then $J(g) = \{z : |z| = 1\}$. Ignoring for a moment that we excluded polynomials of degree 1 from our considerations in complex dynamics, we find that $J(f)$ consists of only one point. Thus the conclusion of Theorem 5.11 does not hold.

Figure 8 shows the Julia sets of the functions

$$f(z) = z - e^{2z} + e^z + \frac{1}{4} \quad \text{and} \quad g(z) = z \exp\left(-z^2 + z + \frac{1}{4}\right), \qquad (5.17)$$

which are related by (5.16). In the left picture for $f$ the range shown is $-6 \leq \operatorname{Re} z \leq 4$ and $|\operatorname{Im} z| \leq 7.5$, in the right one for $g$ it is $|\operatorname{Re} z| \leq 4$ and $|\operatorname{Im} z| \leq 3$. The rectangle given by $|\operatorname{Re} z| \leq 1$ and $|\operatorname{Im} z| \leq \pi/2$ and its image under the exponential function are marked. The function $g$ has the attracting fixed points $(1 \pm \sqrt{2})/2$. The right picture shows the attracting basin of $(1 + \sqrt{2})/2$ in grey and that of $(1 - \sqrt{2})/2$ in white. The same colouring is used in the left picture for the preimages of these basins under the exponential function.

Rippon and Stallard [236] introduced the set

$$Z(f) := \left\{ z : \lim_{n \to \infty} \frac{\log \log |f^n(z)|}{n} = \infty \right\},$$

of points "zipping towards infinity", noting that $I'(f) \subset Z(f) \subset A(f)$. (The second inclusion was already stated as Theorem 5.7.) They established various properties of these sets. For example, they used them to prove that certain invariant curves are contained in the Julia set. As a specific example [236, p. 535], they showed that for $f(z) = z(1 + e^{-z})$ the negative real axis is in $J(f)$.

The reason why the fast escaping set, in the form (5.8), was introduced in [66] was also a problem about Julia sets. This problem concerned Julia sets of permutable entire function and, more generally, of semiconjugate entire functions. Here two entire functions $f$ and $g$ are called *permutable* if $f \circ g = g \circ f$ and they are called semiconjugate if there exists a continuous function $h \colon \mathbb{C} \to \mathbb{C}$ such that $f \circ h = h \circ g$. Note that if the latter equation holds with $f = g$, then $f$ and $h$ are permutable.

Fatou [131] and Julia [153] had used the iteration theory they had developed to study permutable rational functions. In the course of these investigations they had shown in particular that if $f$ and $g$ are permutable rational functions, then $J(f) = J(g)$. It is an open question whether this also holds for permutable transcendental entire functions. A corollary of the results about semiconjugation obtained in [66] says that this is the case if $f$ and $g$ have no wandering domains.

With an additional growth condition this corollary had been obtained by Langley [166]. His proof also involved escaping points. More specifically, he used Eremenko's method to show [166, Lemma 2] that if $f$ is a transcendental entire function with a Baker domain, then there exists $z \in I(f)$ such that

$$\liminf_{n \to \infty} \frac{\log |f^{n+1}(z)|}{|f^n(z)|^{1/2}} > 0.$$

The fast escaping set plays an important role also in the more recent work of Benini, Rippon and Stallard [43]. They show that for permutable entire functions $f$ and $g$ we have $J(f) = J(g)$ if neither $f$ nor $g$ has a fast escaping simply connected wandering domain. Examples of functions with such domains will be discussed in §6.5.

Rippon and Stallard [239] used the fast escaping set $A(f)$ to prove that $I(f)$ has at least one unbounded connected component. This is related to a conjecture of Eremenko that will be discussed in detail in §7. More precisely, they used their characterisation of $A(f)$ given by Theorem 5.8 to strengthen Theorems 4.13 and 5.5 by showing that every connected component of $A(f)$ is unbounded; see Theorem 7.24.

Escape rates also play a role in certain results about measurable dynamics of entire functions. For $z \in \mathbb{C}$ and an entire function $f$, the $\omega$-limit set $\omega(z, f)$ is defined as the set of all limit points of the sequence $(f^n(z))$. Lyubich [174, 176] and Rees [217] proved that for $f(z) = e^z$ we have

$$\omega(z, f) = \{f^n(0) \colon n \geq 0\} \tag{5.18}$$

for almost all $z \in \mathbb{C}$. Note that $\operatorname{sing}(f^{-1}) = \{0\}$ so that the right hand side of (5.18) is equal to the postsingular set $P(f)$. Also, $J(f) = \mathbb{C}$ by a result of Misiurewicz [193].

Urbański and Zdunik [280, Remark 5.5] noted that the proof of (5.18) extends to the functions $f(z) = \lambda e^z$ if $\lambda \in \mathbb{C} \setminus \{0\}$ is such that $0 \in A(f)$. They also obtained this result as a corollary of a more general result [280, Corollary 5.4]. (Urbański and Zdunik used the term "super-growing", but it turns out that this concept is equivalent to "fast escaping"; to verify the equivalence of the two conditions compare, e.g., the arguments in [55, Lemma 2.1].)

Hemke [145] considered functions of the form

$$f(z) = \int_0^z p(t)e^{q(t)}dt + c, \tag{5.19}$$

with polynomials $p$ and $q$ and $c \in \mathbb{C}$. He showed [145, Theorem 1.2] that if all asymptotic values of $f$ are contained in $A(f)$, then $J(f)$ has positive measure and $\omega(z,f) \subset P(f)$ for almost all $z \in J(f)$. Moreover, if all critical points of $f$ are fast escaping, preperiodic, or contained in attracting basins, then $\omega(z,f) = O^+(A)$ for almost all $z \in J(f)$, where $A$ denotes the set of asymptotic values of $f$ [145, Theorem 1.3].

Lyubich [175] proved that the exponential function is not ergodic. Again an important ingredient in the proof is that the singular value 0 is in the fast escaping set. On the other hand, Cui and Wang [100] showed that there exist parameters $\lambda$ such that $f(z) = \lambda e^z$ is ergodic, and $0 \in I(f)$. In their examples, $f^n(0)$ tends to infinity very slowly.

**Question 5.12.** Let $f(z) := \lambda e^z$, with $\lambda \in \mathbb{C} \setminus \{0\}$ such that $0 \in I(f)$. Which growth rate of the sequence $(f^n(0))$ ensures that $f$ is not ergodic? Which growth rate ensures that (5.18) holds?

Similar questions may also be asked for functions of the form (5.19), or more general functions.

## 5.3 Regularity and growth conditions

For a transcendental entire function $f$ and $\alpha, \beta > 0$, let

$$H(r) := \alpha M(\beta r, f).$$

The set $A(f)$ does not change if the condition $|f^{n+L}(z)| \geq M^n(R, f)$ in its definition is replaced by the condition that $|f^{n+L}(z)| \geq H^n(R')$, with $R'$ chosen such that

$$H(r) > r \quad \text{for } r \geq R'. \tag{5.20}$$

This result was already used in the proofs that $A(f) \neq \emptyset$. It seems to have been stated first (for $\beta = 1$) by Rippon and Stallard [244, Theorem 2.9]. An extension to points escaping fast in a tract was given by Xuan and Zheng [288, Theorem 3.9]. Sixsmith [263, Theorem 2] gave a version where the constant $\alpha$ is replaced by a function $\varepsilon(r)$ satisfying certain conditions.

Rippon and Stallard [247] considered, for $0 < \varepsilon < 1$, the function

$$H_\varepsilon(r) := M(r, f)^\varepsilon \tag{5.21}$$

and the resulting analogue

$$Q_\varepsilon(f) := \{z \colon \text{there exists } L \in \mathbb{N} \text{ such that } |f^{n+L}(z)| \geq H_\varepsilon^n(R') \text{ for } n \in \mathbb{N}\}$$

of $A(f)$, with $R'$ such that (5.20) holds with $H$ replaced by $H_\varepsilon$. They called

$$Q(f) := \bigcup_{0<\varepsilon<1} Q_\varepsilon(f). \tag{5.22}$$

the *quite fast escaping set* and gave [247, Theorem 1.2] regularity conditions characterising when $A(f) = Q(f)$ and when $A(f) = Q_\varepsilon(f)$. They used this to prove the following result [247, Theorem 1.2].

**Theorem 5.13.** *If $f \in \mathcal{B}$, then $Q(f) = A(f)$.*

Rippon and Stallard noted that results by Bergweiler, Karpińska and Stallard [68, Theorem 1.1] and by Peter [212, Theorem 1.1] were concerned with lower bounds for the Hausdorff dimension of $I(f)$, but their proofs, although this was not stated explicitly, actually gave lower bounds for the Hausdorff dimension of $Q(f)$. In view of Theorem 5.13 one thus has the same bounds for $A(f)$.

Evdoridou [119] considered, for $m \in \mathbb{N}$ and $0 < \varepsilon < 1$, the function

$$H_{m,\varepsilon}(r) := \exp^{m-1}((\log^{m-1} M(r,f)^\varepsilon)$$

and gave conditions when the resulting set $Q_m(f)$ is equal to $A(f)$. The case $m=2$ had been considered by her previously [117].

## 5.4 Points escaping at a given rate

The previous sections have been concerned with points that escape to infinity fast, in various senses of "fast". It is natural to ask whether there are also points that escape to infinity slowly. The following result of Rippon and Stallard [243, Theorem 1] says that this is indeed the case

**Theorem 5.14.** *Let $f$ be a transcendental entire function and let $(a_n)$ be a sequence of positive numbers that tends to $\infty$. Then there exist $z \in I(f)$ with $|f^n(z)| \leq a_n$ for all large $n$.*

The analogous result for the escape in tracts (cf. Theorem 4.8) is not known.

**Question 5.15.** *Let $D$ be a direct tract of $f$ and let $(a_n)$ be a sequence of positive numbers that tends to $\infty$. Does there exist $z \in D$ with $f^n(z) \in D$ for all $n \in \mathbb{N}$ and $|f^n(z)| \leq a_n$ for all large $n$?*

This question is open also for tracts of entire functions. Some partial results were obtained by Waterman [284, Theorem 1.3].

One may of course be more ambitious than in Theorem 5.14 and ask whether, given a sequence $(a_n)$ of positive numbers tending to infinity, one may achieve not only that $|f^n(z)| \leq a_n$ for all large $n$, but perhaps even that $|f^n(z)|$ is close to $a_n$ in some sense, say $\varepsilon a_n \leq |f^n(z)| \leq a_n$. Since $|f^{n+1}(z)| \leq M(|f^n(z)|, f)$, an obvious condition that has to be made is that $a_{n+1} \leq M(a_n, f)$. Also, if $|f(z)|$ is large in the annulus $\operatorname{ann}(\varepsilon a_n, a_n)$ and if $f^n(z)$ is in this annulus, then $|f^{n+1}(z)|$ will also be large. Thus we cannot achieve $|f^{n+1}(z)| \leq a_{n+1}$ for a prescribed $a_{n+1}$ in this case.

Rippon and Stallard [243, Theorem 2] showed that these are the only restrictions that have to be made.

**Theorem 5.16.** *A transcendental entire function $f$ has the property that for every sequence $(a_n)$ of positive numbers satisfying $a_n \to \infty$ and $a_{n+1} = O(M(a_n, f))$ as $n \to \infty$, there exist $C > 1$ and $z \in I(f)$ with $a_n \leq |f^n(z)| \leq Ca_n$ for all large $n$ if and only if there exist positive constants $c$, $d$ and $r_0$ such that*

$$\min_{r \leq |z| \leq dr} |f(z)| \leq c \tag{5.23}$$

*for $r \geq r_0$.*

The proofs of Theorems 5.14 and 5.16 are based on covering results for annuli that are of independent interest.

If the answer to Question 5.15 is positive, one may also ask for an analogue of Theorem 5.16 for tracts. If we assume $f$ to be entire, then this is of interest only if $f$ has at least two tracts. But then for each tract its boundary contains a curve tending to infinity. Since $f$ is bounded on the boundary of a tract, condition (5.23) is satisfied automatically.

**Question 5.17.** Let $D$ be a direct tract of $f$ such that $\partial D$ contains a curve tending to infinity. Let $(a_n)$ be a sequence of positive numbers that tends to $\infty$. Does there exist $z \in D$ and $C > 0$ such that $f^n(z) \in D$ for all $n \in \mathbb{N}$ and $a_n \leq |f^n(z)| \leq Ca_n$ for all large $n$?

Again this is open also for tracts of entire functions. Waterman [284, Theorem 1.2] showed that for a logarithmic tract one may even achieve that $|f^n(z)| \sim a_n$ as $n \to \infty$.

For admissible $R$ and $n \in \mathbb{N}$, Rippon and Stallard [248] considered the half-open annuli

$$A_n(R) := \{z \colon M^{n-1}(R, f) \leq |z| < M^n(R, f)\}$$

and the disc $A_0(R) := D(0, R)$. For $z \in \mathbb{C}$ the sequence $(s_n)_{n \in \mathbb{N}_0}$ defined by $f^n(z) \in A_{s_n}(R)$ is called the *annular itinerary* of $z$. By the maximum modulus principle, we have $s_{n+1} \leq s_n + 1$ for all $n \in \mathbb{N}_0$. Moreover, $z$ is in the fast escaping set $A(f)$ if and only if $s_{n+1} = s_n + 1$ for all large $n$. Rippon and Stallard give a thorough study of these annular itineraries.

# 6 Escaping Fatou components

## 6.1 No escaping Fatou components in the Eremenko–Lyubich class

This section is concerned with escaping Fatou components; that is, connected components of the Fatou set that are contained in the escaping set. We begin, however, with a result due to Eremenko and Lyubich [116, Theorem 1] that says that functions in the class now named after them do not have such domains.

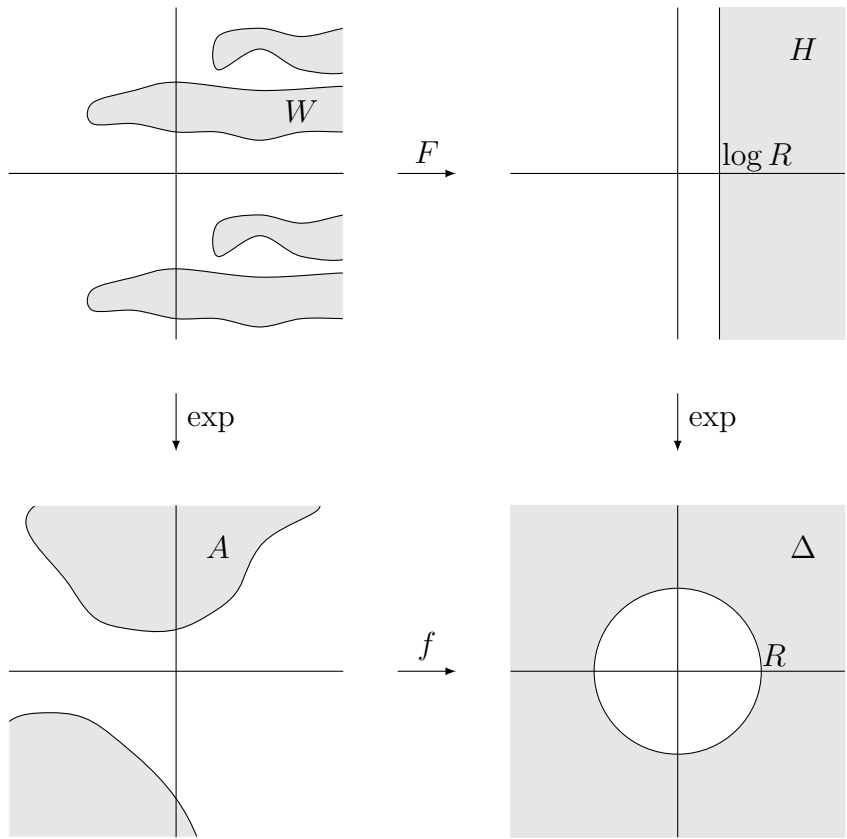

Figure 9: The logarithmic change of variable.

**Theorem 6.1.** *Let $f \in \mathcal{B}$. Then $I(f) \subset J(f)$.*

Together with Theorem 4.3 this yields the following result.

**Theorem 6.2.** *Let $f \in \mathcal{B}$. Then $J(f) = \overline{I(f)}$.*

We sketch the idea behind the proof of Theorem 6.1 only briefly. For the details we refer to the paper by Eremenko and Lyubich [116] as well as the survey by Sixsmith [268].

We choose $R > |f(0)|$ such that $\operatorname{sing}(f^{-1}) \subset D(0, R)$. Let $\Delta := \{z \colon |z| > R\}$, $A := f^{-1}(\Delta)$, $W := \exp^{-1}(A)$ and $H := \mathbb{H}_{>\log R}$. The key result is that the map $f \colon A \to \Delta$ can be lifted to a map $F \colon W \to H$ that has the property that for every connected component $U$ of $W$ the map $F \colon U \to H$ is biholomorphic; see Figure 9. The passage from $f$ to $F$ is called the *logarithmic change of variable*. This change of variable was already made in (5.16), but the difference is that now the functions are considered only on suitable subsets of $\mathbb{C}$.

Let $U$ be a connected component of $W$ and let $G\colon H \to U$ be the inverse of $F\colon U \to H$. For $z \in U$ and $w = f(z) \in H$ we apply Koebe's one quarter theorem to the restriction of $G$ to the disc $D(w, \operatorname{Re} w - \log R)$. This theorem yields that the image of this disc under $G$ contains a disc around $z$ of radius $|G'(w)|(\operatorname{Re} w - \log R)/4$. On the other hand, $U$ does not contain any disc of radius greater than $\pi$. Hence

$$\frac{1}{4}|G'(w)|(\operatorname{Re} w - \log R) \leq \pi.$$

Expressing $G$ and $G'$ in terms of $F$ and $F'$ now yields the following result.

**Theorem 6.3.** *Let $f \in \mathcal{B}$ and let $R, W, H$ and $F$ be as above. Then*

$$|F'(z)| \geq \frac{\operatorname{Re} F(z) - \log R}{4\pi} \quad \text{for } z \in W. \tag{6.1}$$

Rempe [227] showed that (6.1) holds with $4\pi$ replaced by 2, which is best possible.

*Sketch of proof of Theorem 6.1.* Suppose that $z_0 \in I(f) \cap F(f)$. Without loss of generality we may assume that $f^n(z_0) \in A$ for all $n \geq 0$. Let $w_0 \in W$ with $\exp w_0 = z_0$. Then $F^n(w_0) \in W$ for all $n \geq 0$ and $\operatorname{Re} F^n(w_0) \to \infty$ as $n \to \infty$. In fact, there exists $\delta > 0$ such that $F^n(D(w_0, \delta)) \subset W$ for all $n \geq 0$. Koebe's one quarter theorem now yields that $F^n(D(w_0, \delta))$ contains a disc of radius $r_n$ where

$$r_n := \frac{1}{4}|(F^n)'(w_0)|\delta = \frac{1}{4}\delta \prod_{k=0}^{n-1} |F'(F^k(w_0))| \to \infty$$

by (6.1). This is a contradiction since $W$ contains no disc of radius greater than $\pi$. □

## 6.2 A general result on the dynamics in escaping Fatou components

We begin our discussion of escaping Fatou components with a general lemma describing the behaviour of the iterates in such domains. This lemma is taken from Bergweiler's survey [45, Lemma 7], but – as noted there – the results are from papers by Baker [17, Lemmas 1 and 2] and Baker, Kotus and Lü [22, Lemma 4.1]; see also papers by Baker [14, Theorem 6] and McMullen [185, Proposition A.1].

**Lemma 6.4.** *Let $\Omega$ be an unbounded open set in $\mathbb{C}$ with at least two finite boundary points and let $g\colon \Omega \to \mathbb{C}$ be holomorphic. Let $D$ be a domain contained in $\Omega$ and suppose that $g^n(D) \subset \Omega$ for all $n$ and that $g^n|_D \to \infty$ as $n \to \infty$. Then, for any compact subset $K$ of $D$, there exist $C > 0$ and $n_0 \in \mathbb{N}$ such that*

$$|g^n(z_2)| \leq |g^n(z_1)|^C \quad \text{for } z_1, z_2 \in K \text{ and } n \geq n_0. \tag{6.2}$$

*If, in addition, $g(D) \subset D$, then also*

$$\log\log|g^n(z)| = O(n) \quad \text{for } z \in D \tag{6.3}$$

as $n \to \infty$ and there exist $A > 1$ and a curve $\gamma$ in $D$ tending to $\infty$ that satisfies $g(\gamma) \subset \gamma$ such that
$$|z|^{1/A} \leq |g(z)| \leq |z|^A \quad \text{for } z \in \gamma. \tag{6.4}$$
If $\widehat{\mathbb{C}} \setminus \Omega$ contains a connected set $\Gamma$ such that $\{a, \infty\} \subset \Gamma$ for some $a \in \mathbb{C}$, then the estimates in (6.2), (6.3), and (6.4) may be replaced by
$$|g^n(z_2)| \leq C|g^n(z_1)|, \tag{6.5}$$
$$\log|g^n(z)| = O(n), \tag{6.6}$$
and
$$\frac{|z|}{A} \leq |g(z)| \leq A|z|. \tag{6.7}$$
In particular, this is the case if $\Omega$ is simply connected.

We sketch two proofs of Lemma 6.4.

*Sketch of proof of Lemma 6.4 using the hyperbolic metric.* Let $d_\Omega(z_1, z_2)$ denote the hyperbolic distance of $z_1$ and $z_2$ in $\Omega$; see the paper by Beardon and Minda [39] for the definition and basic properties of the hyperbolic metric. Thus
$$d_\Omega(z_1, z_2) = \inf_\gamma \int_\gamma \rho_\Omega(z)|dz|,$$
where $\rho_\Omega$ denotes the density of the hyperbolic metric and the infimum is taken over all smooth paths in $\Omega$ that connect $z_1$ and $z_2$. By the Schwarz–Pick lemma we have
$$d_\Omega(g^n(z_1), g^n(z_2)) \leq d_\Omega(z_1, z_2). \tag{6.8}$$

Standard estimates of the hyperbolic metric show that there exist positive constants $c$ and $R$ such that
$$\rho_\Omega(z) \geq \frac{c}{|z|\log|z|} \quad \text{for } |z| \geq R. \tag{6.9}$$
Assuming that $|g^n(z_2)| > |g^n(z_1)|$ we find that
$$d_\Omega(g^n(z_1), g^n(z_2)) \geq \int_{|g^n(z_1)|}^{|g^n(z_2)|} \frac{c}{|z|\log|z|}|dz| = c\log\left(\frac{\log|g^n(z_2)|}{\log|g^n(z_1)|}\right).$$

Combining this with (6.8) we obtain (6.2).

In order to prove (6.3), let $\sigma$ be a curve in $D$ that connects a point $z_0 \in D$ with $g(z_0)$ and define $\gamma = \bigcup_{n=0}^\infty g^n(\sigma)$. Then (6.3) can be deduced from (6.2) with $K = \sigma \cup g(\sigma)$. Similarly, (6.2) yields for $z_2 = g(z_1)$ that
$$|g^n(z_1)| = |g^{n-1}(z_2)| \leq |g^{n-1}(z_1)|^C = |g^{n-2}(z_2)|^C \leq |g^{n-2}(z_1)|^{C^2} \leq \ldots$$
for large $n$, and (6.4 follows by induction.

If there exists a connected set $\Gamma$ as described in the hypothesis of the lemma, then (6.9) can be improved to
$$\rho_\Omega(z) \geq \frac{c}{|z|} \quad \text{for } |z| \geq R.$$

Instead of (6.2), (6.3) and (6.4) we now obtain (6.5), (6.6) and (6.7). □

*Sketch of proof of Lemma 6.4 using Harnack's inequality.* Harnack's inequality (see, for example, Ransford's book [216, §1.3]) says that if $h$ is a positive harmonic function in a disc $D(a, r)$ and if $0 < \rho < r$, then
$$\frac{r-\rho}{r+\rho} h(a) \leq h(z) \leq \frac{r+\rho}{r-\rho} h(a) \quad \text{for } z \in \overline{D}(a, \rho).$$

Given $K$ and $D$ as in the hypothesis, there exists a connected compact set $K'$ and a domain $G$ with $\overline{G} \subset D$ such that $K \subset K' \subset G$. The set $K'$ can be covered by finitely many discs $D(a_j, r_j/2)$ such that $D(a_j, r_j) \subset G$. It follows that there exists a positive constant $C'$ such that if $h$ is a positive harmonic function in $G$, then
$$h(z_2) \leq C' h(z_1) \quad \text{for } z_1, z_2 \in K'. \tag{6.10}$$

Since $\overline{G} \subset D$ it follows from the hypotheses that $\log|g^n|$ is a positive harmonic function in $G$ for large $n$. Thus (6.2) follows from (6.10), with $C = C'$.

Let $\Gamma$ be as in the hypothesis. Consider an arc that connects some point of $\Gamma$ with $0$ and let $\Gamma_0$ be the union of $\Gamma$ and this arc. The connected components of $\mathbb{C} \setminus \Gamma_0$ are simply connected, so in each such component we can define a branch of the logarithm. As $g^n|_G \to \infty$ uniformly as $n \to \infty$ and since $g^n(G) \subset \Omega \subset \mathbb{C} \setminus \Gamma$ for all $n$, we have $g^n(G) \subset \mathbb{C} \setminus \Gamma_0$ for large $n$. Hence $\log g^n$ is a holomorphic function in $G$ for large $n$ and we may apply (6.10) to $h(z) = \log|\log g^n(z)|$. Here the branch of the logarithm may be chosen such that $|\text{Im}(\log g^n(z_1))| = |\arg g^n(z_1)| \leq \pi$. For large $n$ we thus have
$$|\log g^n(z_1)| \leq |\text{Re}(\log g^n(z_1))| + \pi = \log|g^n(z_1)| + \pi \leq 2\log|g^n(z_1)|.$$

Hence (6.10) yields with $h = \log|\log g^n|$ that
$$\log\log|g^n(z_2)| \leq \log|\log g^n(z_2)| \leq C'\log|\log g^n(z_1)| \leq C'\log\log|g^n(z_1)| + C'\log 2.$$

Thus (6.5) holds for large $n$ if $C > C'$.

As in the previous proof we can now deduce (6.3) and (6.4) from (6.2) and deduce (6.6) and (6.7) from (6.5). □

Inequality (6.2) says that within an escaping Fatou component all points tend to infinity at about the same rate. In particular, it implies the following result; cf. [239, p. 1125, Remark 1].

**Theorem 6.5.** *Let $f$ be a transcendental entire function and let $U$ be a Fatou component of $f$. If $U \cap A(f) \neq \emptyset$, then $U \subset A(f)$.*

In fact, a stronger result holds. Recall that the *levels* $A_R^L(f)$ were defined in (5.10). The following theorem is due to Rippon and Stallard [244, Theorem 1.2].

**Theorem 6.6.** *Let $f$ be a transcendental entire function, $R$ an admissible radius for $f$, and $L \in \mathbb{Z}$. If $U$ is a Fatou component with $U \cap A_R^L \neq \emptyset$, then $\overline{U} \subset A_R^{L-1}(f)$. If $U$ is simply connected, then even $\overline{U} \subset A_R^L(f)$.*

*Proof.* Let $z_2 \in U \cap A_R^L(f)$ and $z_1 \in U$. Then $|f^n(z_2)| \geq M^{n+L}(R,f)$ if $n, n+L \in \mathbb{N}$. We apply Lemma 6.4 with $g = f$ and $\Omega = \bigcup_{n=0}^\infty f^n(U)$ and find that there exists $C > 0$ such that
$$|f^n(z_1)| \geq |f^n(z_2)|^{1/C} \geq M^{n+L}(R,f)^{1/C} \tag{6.11}$$
for large $n$. It follows from (5.4) that $M(r,f) \geq r^C$ for large $r$. We conclude that
$$M^{n+L}(R,f)^{1/C} = M\bigl(M^{n+L-1}(R,f),f\bigr)^{1/C} \geq M^{n+L-1}(R,f)$$
and hence $|f^n(z_1)| \geq M^{n+L-1}(R,f)$ for large $n$. As noted after the definition of the levels $A_R^L(f)$ this yields that $|f^n(z_1)| \geq M^{n+L-1}(R,f)$ whenever $n, n+L-1 \in \mathbb{N}$. Thus $z_1 \in A_R^{L-1}$. Since $A_R^{L-1}$ is closed we conclude that not only $U \subset A_R^{L-1}$, but also $\overline{U} \subset A_R^{L-1}$.

Suppose now that $U$ is simply connected. Then $\widehat{\mathbb{C}} \setminus \Omega$ contains a connected set $\Gamma$ as stated in the lemma and instead of (6.11) we obtain
$$|f^n(z_1)| \geq \frac{1}{C}|f^n(z_2)| \geq \frac{1}{C}M^{n+L}(R,f) \tag{6.12}$$
for large $n$. Suppose that $z_1 \notin A_R^L$. Then there exists $m \in \mathbb{N}$ such that $|f^m(z_1)| < M^{m+L}(R,f)$; in fact, this holds for all large $m$. Fix such a large $m$. Then
$$|f^m(z_1)| = M^{m+L}(R,f)^\alpha$$
for some $\alpha \in (0,1)$. The maximum modulus principle and Lemma 5.4 (applied with $c = 1/\alpha$) imply that
$$|f^{m+1}(z_1)| = |f(f^m(z_1))| \leq M(|f^m(z_1)|,f) = M(M^{m+L}(R,f)^\alpha, f)$$
$$\leq M(M^{m+L}(R,f),f)^\alpha = M^{m+L+1}(R,f)^\alpha$$
if $m$ and hence $M^{m+L}(R,f)$ are large enough. Inductively we find that
$$|f^{m+k}(z_1)| \leq M^{m+L+k}(R,f)^\alpha \tag{6.13}$$
for $k \in \mathbb{N}$. Combining (6.12) and (6.13), with $n = m + k$, yields that
$$\frac{1}{C}M^{n+L}(R,f) \leq M^{n+L}(R,f)^\alpha$$
and thus $M^{n+L}(R,f)^{1-\alpha} \leq C$. This is a contradiction for large $n$ since $M^n(R,f) \to \infty$ as $n \to \infty$. $\square$

An immediate consequence is the following result due to Rippon and Stallard [244, Theorem 1.2, Remark].

**Theorem 6.7.** *Let $f$ be a transcendental entire function and let $U$ be a Fatou component of $f$. If $U \subset A(f)$, then $\overline{U} \subset A(f)$.*

For multiply connected Fatou components Rippon and Stallard proved this result in [239, Theorem 2, (a)], the general case is in [244, Theorem 1.2, Remark].

## 6.3 Multiply connected wandering domains

Already Töpfer [275, p. 67] stated that any multiply connected Fatou component is contained in the escaping set. More precisely, we have the following result of Baker [15, Theorem 3.1]. Here, for a (piecewise smooth) closed curve $\gamma$ and $z \notin \gamma$,

$$n(\gamma, z) := \frac{1}{2\pi i} \int_\gamma \frac{d\zeta}{\zeta - z}$$

denotes the *winding number* of $\gamma$ with respect to $z$.

**Theorem 6.8.** *Let $f$ be a transcendental entire function and let $U$ be a multiply connected Fatou component of $f$. Then $U \subset I(f)$.*

*Moreover, $f^k(U)$ contains a curve $\gamma_k$ such that $\operatorname{dist}(\gamma_k, 0) \to \infty$ as $k \to \infty$ and $n(\gamma_k, 0) \neq 0$ for all large $k$.*

*Proof.* Suppose that $U \not\subset I(f)$. Then there exists a sequence $(k_j)$ tending to $\infty$ such that $f^{k_j}|_U \to \varphi$ for some holomorphic function $\varphi \colon U \to \mathbb{C}$. Since $U$ is multiply connected, there exists a smooth Jordan curve $\gamma$ in $U$ whose interior intersects $J(f)$. By Theorem 3.1 (d), the interior of $\gamma$ contains a repelling periodic point $z_0$. It follows that $|(f^{k_j})'(z_0)| \to \infty$.

Assuming that $\gamma$ is positively oriented, we have

$$f^{k_j}(z) = \frac{1}{2\pi i} \int_\gamma \frac{f^{k_j}(\zeta)}{\zeta - z} d\zeta$$

for $z$ in the interior of $\gamma$ by Cauchy's integral formula. As $f^{k_j} \to \varphi$ uniformly on $\gamma$, we find that $f^{k_j} \to \varphi$ locally uniformly in the interior of $\gamma$. By Weierstraß's theorem, $(f^{k_j})'(z_0) \to \varphi'(z_0)$. This is a contradiction to our previous finding that $|(f^{k_j})'(z_0)| \to \infty$. Hence $U \subset I(f)$.

We now put $\gamma_k := f^k(\gamma)$. Since $U \subset I(f) \cap F(f)$, we see that $f^k(z) \to \infty$ uniformly for $z \in \gamma$. Thus $\operatorname{dist}(\gamma_k, 0) \to \infty$ as $k \to \infty$. Since $n(\gamma, z_0) \neq 0$ we have $n(\gamma_k, f^k(z_0)) \neq 0$. As $(f^k(z_0))$ is bounded and $\operatorname{dist}(\gamma_k, 0) \to \infty$, the points $0$ and $f^k(z_0)$ are in the same connected component of $\mathbb{C} \setminus \gamma_k$ for large $k$. Thus $n(\gamma_k, 0) = n(\gamma_k, f^k(z_0)) \neq 0$ for large $k$. □

The first example of a transcendental entire function with a multiply connected Fatou component is the function $f$ defined by (2.5) and (2.6) that was discussed already in §2.4. As mentioned there, it was given by Baker [11] in 1963. It follows from (2.8) that the annuli $A_n$ defined by (2.7) are contained in $F(f) \cap I(f)$ for large $n$. Since 0 is a fixed point, the Fatou component $U_n$ that contains $A_n$ is multiply connected.

In 1976, Baker [13] showed that these multiply connected Fatou components $U_n$ are distinct. This implies that they are wandering. This was the first example of a wandering domain, predating Sullivan's famous theorem [274] that rational functions do not have wandering domains.

Soon after [13] was written, Baker [12] showed that this is not a particular property of the example (2.6), but that multiply connected Fatou components of transcendental entire functions are always wandering. This paper, although written later, was published earlier than [13].

**Theorem 6.9.** *Let $f$ be a transcendental entire function. Then any multiply connected Fatou component of $f$ is wandering.*

Theorem 6.9 can be proved using the following result of independent interest [244, Theorem 4.4].

**Theorem 6.10.** *Let $f$ be a transcendental entire function and let $U$ be a multiply connected Fatou component of $f$. Then $U \subset A(f)$.*

*Proof.* Let $\gamma_n$ be as in Theorem 6.8. Theorem 5.5 yields that $\gamma_n \cap A(f) \neq \emptyset$, if $n$ is sufficiently large. Since $A(f)$ is completely invariant, we deduce that $\gamma \cap A(f) \neq \emptyset$ and hence $U \cap A(f) \neq \emptyset$. The conclusion now follows from Theorem 6.5. □

*Proof of Theorem 6.9.* Suppose that $F(f)$ has a preperiodic multiply connected Fatou component. Passing to an iterate we may assume that $f$ has an invariant multiply connected Fatou component $U$. For $z \in U$ we then have $\log \log |f^n(z)| = O(n)$ by (6.3), while $z \in A(f)$ by Theorem 6.8. A contradiction is now obtained from Theorem 5.7. □

Theorems 6.8 and 6.9 imply that a multiply connected Fatou component $U$ of $f$ is in a bounded complementary component of $f^n(U)$ for large $n$. In particular, this yields the following result.

**Theorem 6.11.** *Every multiply connected Fatou component of a transcendental entire function is bounded.*

Combining Theorems 6.7 and 6.10 we see that the closure of a multiply connected Fatou component is contained in the fast escaping set.

In particular, we have the following result already used in the proof of Eremenko's Theorem 4.6, which says that $I(f) \cap J(f) \neq \emptyset$. For completeness we include a sketch of its short proof, which does not use Theorem 6.7.

**Theorem 6.12.** *Let $f$ be a transcendental entire function and let $U$ be a multiply connected Fatou component of $f$. Then $\overline{U} \subset I(f)$.*

*Sketch of proof.* Let $\gamma_n$ be as Theorem 6.8. Then, given $k \in \mathbb{N}$, for large $n$ the curve $\gamma_n$ is contained in the unbounded connected component of $\mathbb{C} \setminus \gamma_k$. Since $U$ is a wandering domain by Theorem 6.9, the curves $\gamma_k$ and $\gamma_n$ are in different Fatou components for $k \neq n$. Thus $f^n(U)$ and in fact $f^n(\overline{U})$ are contained in the unbounded connected component of $\mathbb{C} \setminus \gamma_k$ for large $n$. Since $k$ can be taken arbitrarily large, the conclusion follows. □

As already mentioned, inequality (6.2) says that within an escaping Fatou component, all points tend to infinity at about the same rate. The following result by Bergweiler, Rippon and Stallard [73, Theorem 1.1] is a more precise version of this fact for multiply connected Fatou components.

**Theorem 6.13.** *Let $f$ be a transcendental entire function, $U$ a multiply connected Fatou component of $f$, and $z_0 \in U$. Then*

$$h(z) := \lim_{n \to \infty} \frac{\log |f^n(z)|}{\log |f^n(z_0)|}$$

*defines a non-constant positive harmonic function $h \colon U \to \mathbb{R}$, with $h(z_0) = 1$.*

We note that it follows easily from (6.2) and Harnack's inequality that there exists a sequence $(n_k)$ tending to infinity such that

$$h(z) = \lim_{k \to \infty} \frac{\log |f^{n_k}(z)|}{\log |f^{n_k}(z_0)|}$$

exists. However, additional arguments are required to show that $h$ is non-constant and that actually the whole sequence converges, not only a subsequence; see [73, Theorem 1.1].

Figure 10 is an illustration of the function $h$ from Theorem 6.13 for Baker's example (2.6), with $C = 1/5$ and $r_1 = 40$. Dark grey corresponds to large values of $h$ while light grey corresponds to small values. Note that $h$ is large near the outer boundary of the multiply connected domains and small near the other boundary components. One can see that the connectivity of the domains is greater than 2. In fact, as noted by Bergweiler and Zheng [75, Remark after Theorem 1.3], they are infinitely connected. The attracting basin of the superattracting fixed point at 0 is white in this figure. The range shown is $-82 \leq \operatorname{Re} z \leq 42$ and $|\operatorname{Im} z| \leq 62$ in the left picture and $-23 \leq \operatorname{Re} z \leq 15$ and $|\operatorname{Im} z| \leq 19$ in the right picture.

One consequence of Theorem 6.13 obtained in [73, Theorem 1.2] is that, for large $n$, there exist $r_n$ and $R_n$ such that

$$f^n(z_0) \in \operatorname{ann}(r_n, R_n) \subset f^n(U) \tag{6.14}$$

and

$$\liminf_{n \to \infty} \frac{\log R_n}{\log r_n} > 1. \tag{6.15}$$

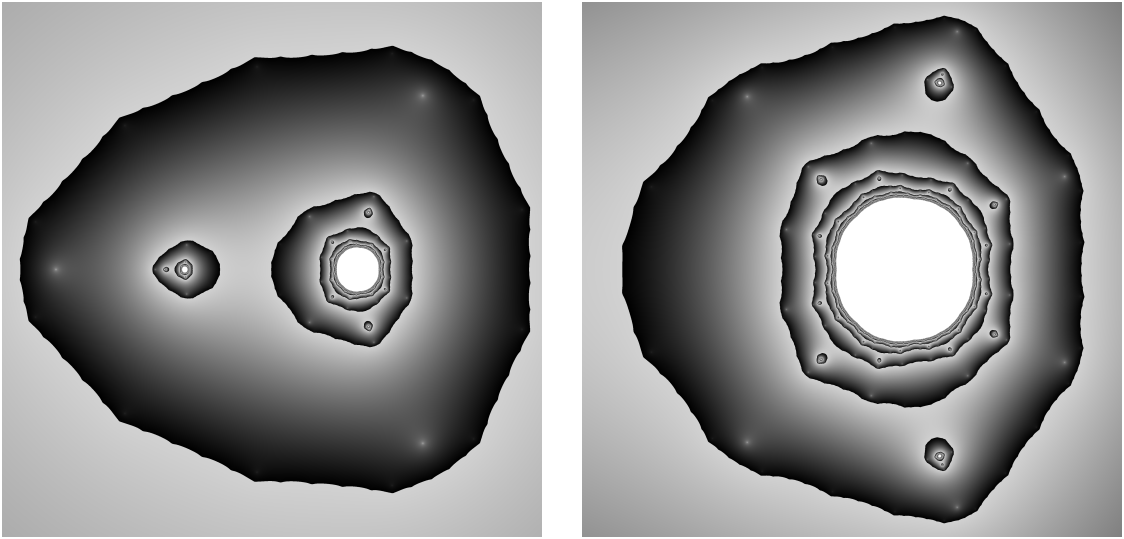

Figure 10: Illustration of the function $h$ from Theorem 6.13.

Previously, Zheng [291] had shown (6.14) holds for large $n$ with $r_n$ and $R_n$ satisfying
$$\lim_{n\to\infty} \frac{R_n}{r_n} = \infty.$$
This result implies that the Julia set of a transcendental entire function with a multiply connected Fatou component is not uniformly perfect in the sense of Pommerenke [214]. This is in contrast with the result [148, 177] that Julia sets of rational functions are uniformly perfect. It is also shown in [73, Section 10] that (6.15) is best possible.

The harmonic function $h$ in Theorem 6.13 can be considered as a measure of how fast points in $U$ escape in comparison to the point $z_0$. It has an interpretation in terms of harmonic measure.

To describe this interpretation, let $U$ be a hyperbolic domain (i.e., a domain $U$ such that $\mathbb{C} \setminus U$ contains at least two points) and let $E$ a Borel subset of $\partial U$. Let $\omega(\cdot, E, U)$ denote the harmonic measure of $E$ in $U$; see, e.g., Ransford's book [216, §4.3] for the definition and properties of harmonic measure. Under the hypotheses of Theorem 6.13, let $r_n := |f^n(z_0)|$ and, for large $n$, let $a_n \in (0,1)$ and $b_n \in (1,\infty)$ be such that
$$B_n := \mathrm{ann}(r_n^{a_n}, r_n^{b_n})$$
is the maximal annulus that contains $f^n(z_0)$ and is contained in $f^n(U)$. It is shown in [73, Theorem 1.5] that the limits $a := \lim_{n\to\infty} a_n$ and $b := \lim_{n\to\infty} b_n$ exist. By (6.15) we have $a < b$. Let $\partial_\infty U := \partial T(U)$ be the "outer boundary" of $U$. It is shown in [73, Theorem 1.6] that if $b < \infty$, then
$$h(z) = a + (b-a)\omega(z, \partial_\infty U, U) \quad \text{for } z \in U,$$
while if $b = \infty$, then
$$\omega(z, \partial_\infty U, U) = 0 \quad \text{for } z \in U.$$

Moreover, it is shown in [73, Theorem 1.3] that the annuli $B_n$ are *absorbing* in the sense that if $K$ is a compact subset of $U$, then $f^n(K) \subset B_n$ for large $n$. In fact, there are absorbing subannuli $C_n$ of $B_n$ that have the additional feature that $f$ and in fact all iterates $f^m$ behave like monomials in $C_n$. The result says that for large $n$ and all $m \in \mathbb{N}$ we have
$$f^m(z) = \varphi_{n,m}(z)^{d_{n,m}} \quad \text{for } z \in C_n \tag{6.16}$$
for some conformal map $\varphi_{n,m}$ and some $d_{m,n} \in \mathbb{N}$; see [73, Theorem 5.2] for the precise statement.

Baker [16] constructed examples of multiply connected wandering domains having various additional properties. Recall here that the *order* $\rho(f)$ of an entire function $f$ is defined by
$$\rho(f) := \limsup_{r \to \infty} \frac{\log \log M(r, f)}{\log r}. \tag{6.17}$$

**Theorem 6.14.** [16, Theorem 1] *Given $\rho$ with $0 \leq \rho \leq \infty$ there exists a transcendental entire function of order $\rho$ that has a multiply connected wandering domain.*

**Theorem 6.15.** [16, Theorem 2] *There exists a transcendental entire function with a multiply connected wandering domain of infinite connectivity.*

Further examples with multiply connected wandering domains of infinite connectivity were given by Bergweiler and Zheng [75]. In particular, as already mentioned, they showed [75, Section 6] that Baker's example (2.6) has this property.

Theorem 6.15 left open the question whether there are multiply connected wandering domain of finite connectivity. An affirmative answer was given by Kisaka and Shishikura [160].

**Theorem 6.16.** *For all $m \in \mathbb{N}$ with $m \geq 2$ there exists a transcendental entire function with a wandering domain of connectivity $m$.*

In their proof they used quasiconformal surgery. A slight modification of their method was used by Bergweiler [52] to construct a transcendental entire function with a multiply connected Fatou component that also has a fast escaping simply connected wandering domain. Burkart and Lazebnik [91, 92] used quasiconformal surgery for a quite general construction of transcendental entire functions with multiply connected Fatou components, and they give several applications of their construction.

All other examples of multiply connected Fatou components that we are aware of have been constructed using infinite products as in (2.6). This applies in particular to Baker's proof [16] of Theorems 6.14 and 6.15 as well as the examples in [73, Section 10] and [75] mentioned above.

We mention some further examples constructed using infinite products. Hinkkanen [149] used this technique to show that a transcendental entire function with a multiply connected wandering domain can have arbitrarily slow growth. (Thus the result that Julia set of polynomials are uniformly perfect does not extend to functions of small growth.)

Baker and Domínguez [20, Theorem G] constructed a transcendental entire function $f$ with a multiply connected Fatou component which has infinitely many repelling fixed points $p$ such that $\{p\}$ is a singleton component of $J(f)$. Moreover, $\{p\}$ is a buried component of $J(f)$, meaning that there exists no Fatou component $V$ with $p \in \partial V$.

Bishop [83] constructed a transcendental entire function (of arbitrarily small growth) with a multiply connected wandering domain for which the Julia set has Hausdorff dimension 1 and in fact packing dimension 1. The boundaries of these wandering domains are $C^1$ curves. Bishop [83, Section 21, (5)] asked the following question.

**Question 6.17** (Bishop, 2018). How smooth can the boundary of a multiply connected Fatou component of a transcendental entire function be? Can it consist of $C^2$ curves? Or $C^\infty$ curves? Or analytic curves?

Even stronger one might ask the following.

**Question 6.18.** Can the boundary of a multiply connected Fatou component of a transcendental entire function consist of circles?

A more ambitious problem is to characterise multiply connected Fatou components.

**Question 6.19.** For which bounded, multiply connected domains $D$ does there exist a transcendental entire function $f$ such that $D$ is a Fatou component of $f$?

Some restrictions on the domain $D$ are necessary. Indeed, by a result of Martí-Pete, Rempe and Waterman [180, Corollary 3.5], there are domains bounded by finitely many Jordan curves that cannot be realised as such a wandering domain.

For simply connected domains $D$ a sufficient condition to be an escaping Fatou component is given in Theorem 6.36 below.

Baumgartner [37] has given general criteria implying that the boundary of a multiply connected Fatou component consists of Jordan curves. These criteria apply, in particular, to Baker's first example (2.6) of a wandering domain as well as his example in Theorem 6.15.

Burkart [89] extended Bishop's construction to show that for every $d \in [1, 2]$ and every $\varepsilon > 0$ there exists a transcendental entire function $f$ such that the packing dimension of $J(f)$ differs from $d$ by less than $\varepsilon$. Zhang [290] showed that one can construct examples $f$ with $\dim_H J(f) = 1$ of any preassigned (including infinite) order.

We conclude the section by briefly discussing the bounded components of the complement of a multiply connected wandering domain $U$, following Rippon and Stallard [250, §7–8]. Let $K$ be a bounded component of $\mathbb{C} \setminus U$. First suppose that $\text{int}(K) \cap J(f) \neq \emptyset$. Then for sufficiently large $n$, the image $K_n \coloneqq f^n(K)$ is the connected component $A_n$ of $\mathbb{C} \setminus f^n(U)$ containing 0. On the other hand, if $\text{int}(K) \subset F(f)$, then $K_n \cap A_n = \emptyset$ for all $n \geq 0$, and in particular $K \subset A(f)$.

It follows from [250, Theorem 8.1] that some of the functions constructed by Bergweiler and Zheng [75] have wandering domains of uncountable connectivity. In particular, most of the bounded complementary components of such a wandering domain have empty interior, and hence are contained in the escaping set. Answering a question of Rippon and Stallard [250, Question 10.4], Burkart and Lazebnik [92] gave an example of a multiply connected wandering domain such that every complementary component $K$ is either a singleton in $I(f)$, or satisfies $\text{int}(K) \cap J(f) \neq \emptyset$. It is unknown whether there can be other types of complementary components.

**Question 6.20.** Does there exist a transcendental entire function, a multiply connected wandering domain $U$ of $f$, and a bounded connected component $K$ of $\mathbb{C} \setminus U$ such that $K$ is not a singleton and $\text{int}(K) \cap J(f) = \emptyset$?

If so, can $K$ be chosen such that $\text{int}(K) \neq \emptyset$?

(The second part of Question 6.20 is [250, Question 10.5].)

## 6.4 Baker domains

By definition (see Theorem 3.2 $(d)$), a Baker domain is contained in the escaping set, while all other types of (pre)periodic Fatou components do not intersect the escaping set. An example of a function with a Baker domain is given by Fatou's function $f(z) = e^{-z} + z + 1$, which was already discussed in §2.2. Here we only review some selected results on Baker domains and refer to the excellent survey by Rippon [235] for a thorough treatment.

By Theorem 6.9, Baker domains are simply connected. The following result of Baker [17, Theorem 1] is thus an immediate consequence of equation (6.6) of his Lemma 6.4.

**Theorem 6.21.** *Let $f$ be a transcendental entire function and let $U$ be a Baker domain of $f$. Then $\log |f^n(z)| = O(n)$.*

We note that for many applications the weaker estimate $\log \log |f^n(z)| = O(n)$ given by (6.3) suffices; e.g., McMullen [185, Appendix] worked with this estimate.

A result of Cowen [97] leads to a classification of Baker domains. We follow König [162] and Bargmann [34] in the description of this classification; see also the book by Abate [2, §4], an article by Fagella and Henriksen [125] and the survey by Rippon [235, §5] already mentioned.

Let $U$ be a domain and let $f \colon U \to U$ be holomorphic. A subdomain $V$ of $U$ is called *absorbing* for $f$, if $V$ is simply connected, $f(V) \subset V$ and for any compact subset $K$ of $U$ there exists $N \in \mathbb{N}$ such that $f^N(K) \subset V$.

**Definition 6.22.** Let $U$ be a domain and let $f \colon U \to U$ be a holomorphic function. Then $(V, \varphi, T, \Omega)$ is called an *eventual conjugacy* of $f$ in $U$, if the following four statements hold:

(i) $V$ is absorbing for $f$;

(ii) $\varphi \colon U \to \Omega \in \{\mathbb{H}_{>0}, \mathbb{C}\}$ is holomorphic and $\varphi$ is univalent in $V$;

(iii) $T$ is a Möbius transformation mapping $\Omega$ onto itself and $\varphi(V)$ is absorbing for $T$;

(iv) $\varphi(f(z)) = T(\varphi(z))$ for $z \in U$.

The terminology is not uniform: König [162] used *conformal* conjugacy instead of *eventual* conjugacy. Cowen used *fundamental* instead of *absorbing*.

In an immediate attracting (but not superattracting) basin or in an immediate parabolic basin, an eventual conjugacy is given by the solution of the Schröder–Kœnigs functional equation or the Abel functional equation, respectively; see [273, §§3.4–3.5]. We are interested in the case where $U$ is a Baker domain. Then $f$ has no fixed point in $U$. Also, as already mentioned, $U$ is simply connected by Theorem 6.9.

Cowen's result says the following.

**Theorem 6.23.** *Let $U$ be a simply connected domain, $U \neq \mathbb{C}$, and let $f\colon U \to U$ be a holomorphic function without fixed point in $U$. Then $f$ has an eventual conjugacy $(V, \varphi, T, \Omega)$. Moreover, we may choose $T$ and $\Omega$ such that we have exactly one of the following cases:*

(a) $\Omega = \mathbb{H}_{>0}$ and $T(z) = \lambda z$, where $\lambda > 1$;

(b) $\Omega = \mathbb{H}_{>0}$ and $T(z) = z + i$ or $T(z) = z - i$;

(c) $\Omega = \mathbb{C}$ and $T(z) = z + 1$.

**Definition 6.24.** Let $f$ be a transcendental entire function and let $U$ be a Baker domain of $f$. Then, according to the three cases arising in Theorem 6.23, the Baker domain $U$ is called *hyperbolic* in case $(a)$, *simply parabolic* in case $(b)$, and *doubly parabolic* in case $(c)$.

This terminology arises from the case where $U = \mathbb{D}$ and $f$ is a Blaschke product of degree at least 2. Since $f$ has no fixed point in $\mathbb{D}$, the Denjoy–Wolff theorem [192, Theorem 5.8] yields that there exists $\xi \in \partial\mathbb{D}$ such that $f^n|_{\mathbb{D}} \to \xi$ as $n \to \infty$. It turns out that we have case $(a)$ when $\xi$ is an attracting fixed point. In cases $(b)$ and $(c)$ the fixed point $\xi$ is parabolic, with one immediate parabolic basin in case $(b)$ and two such basins in case $(c)$. We will use the terminology of Definition 6.24 also for holomorphic self-maps of simply connected domains.

One criterion (due to Bargmann [34, Lemma 2.6]) to determine the type according to the above classification depends on the hyperbolic distance of successive iterates.

**Theorem 6.25.** *Let $U$ be a simply connected domain, $U \neq \mathbb{C}$, and let $f\colon U \to U$ be a holomorphic function without fixed point in $U$. For $z \in U$, put*

$$\rho(z) := \lim_{n \to \infty} d_U(f^{n+1}(z), f^n(z)). \tag{6.18}$$

*Then we have the following:*

(a) *$U$ is hyperbolic $\Leftrightarrow \inf_{z \in U} \rho(z) > 0$;*

(b) *$U$ is simply parabolic $\Leftrightarrow \rho(z) > 0$ for all $z \in U$, but $\inf_{z \in U} \rho(z) = 0$;*

(c) *$U$ is doubly parabolic $\Leftrightarrow \rho(z) = 0$ for all $z \in U$.*

Note that the limit in (6.18) exists since the sequence $(d_U(f^{n+1}(z), f^n(z)))$ is non-increasing by the Schwarz–Pick theorem.

In the case where $U$ is a Baker domain, König [162, Theorem 3] (see also [235, Theorem 5.2]) has given a characterisation of the different types of domain using the *Euclidean* distance of successive iterates, and their distance to $\partial U$.

Besides the above classification, one can draw further conclusions from the sequence $(d_U(f^{n+1}(z), f^n(z)))$. For example, it was shown by Bergweiler [49] that if $U$ is a Baker domain of $f$ and if $\operatorname{sing}(f^{-1}) \cap U$ is bounded, then

$$d_U(f^{n+1}(z), f^n(z)) = \frac{1}{2n} + a\frac{\log n}{n^2} + O\left(\frac{1}{n^2}\right)$$

for some $a \in \mathbb{R}$ if $U$ is doubly parabolic and

$$d_U(f^{n+1}(z), f^n(z)) = \rho(z) + \frac{b}{n^3} + O\left(\frac{1}{n^4}\right),$$

for some $b \in \mathbb{R}$ if $U$ is simply parabolic.

If $f \colon U \to U$ is univalent, then the sequence $(d_U(f^{n+1}(z), f^n(z)))$ is constant. This implies that $f$ cannot be univalent in a doubly parabolic Baker domain. In contrast, there are hyperbolic and simply parabolic Baker domains where $f$ is univalent. In this case we say that the Baker domain is univalent. The first example of a univalent Baker domain seems to be due to Herman [147, p. 609f], who showed that

$$f_1(z) = z + \frac{a}{2\pi} \sin(2\pi z) + b \quad \text{and} \quad f_2(z) = z + 2\pi i \alpha + e^z \qquad (6.19)$$

have univalent Baker domains for suitable $a, b, \alpha \in \mathbb{R}$.

To explain the idea behind these examples, we recall that if $g \colon \mathbb{C} \setminus \{0\} \to \mathbb{C} \setminus \{0\}$ is holomorphic, then there exists an entire function $f$ such that $\exp \circ f = g \circ \exp$; that is, we have the commuting diagram (5.16). Theorem 5.11 says that in this situation we have $J(f) = \exp^{-1} J(g)$, provided $f$ is transcendental.

For $f = f_2$ we have (5.16) with

$$g(z) = g_2(z) := e^{2\pi i \alpha} z e^z. \qquad (6.20)$$

Choosing $\alpha \in \mathbb{R} \setminus \mathbb{Q}$ such that $g_2$ has a Siegel disc $D$ at $0$ yields that $U := \exp^{-1} D$ is a univalent Baker domain of $f_2$.

The right picture of Figure 11 shows this Baker domain of $f_2$ with $\alpha$ chosen as the "golden mean"; that is, $\alpha = (\sqrt{5} - 1)/2$. The range shown is $|\operatorname{Re} z| \leq 5$ and $|\operatorname{Im} z| \leq 8$.

The argument for $f_1$ is similar. Here we replace $e^z$ by $e^{2\pi i z}$ in (5.16). Then $f_1$ corresponds to the self-map $g_1(z) := e^{2\pi i b} \exp\left(\frac{1}{2}a(z - 1/z)\right)$ of $\mathbb{C} \setminus \{0\}$. The map $g_1$ has a Herman ring if $a$ and $b$ are suitably chosen, and this Herman ring lifts to a univalent Baker domain $U$ of $f_1$ under $z \mapsto e^{2\pi i z}$. The Herman ring of $g_1$ and the Baker domain of $f_1$ have been studied in detail; see, e.g., [27, 123, 124, 129, 137]. The papers [27, 123, 124] include pictures of them, and [124, §5] describes the numerical algorithms used to create them.

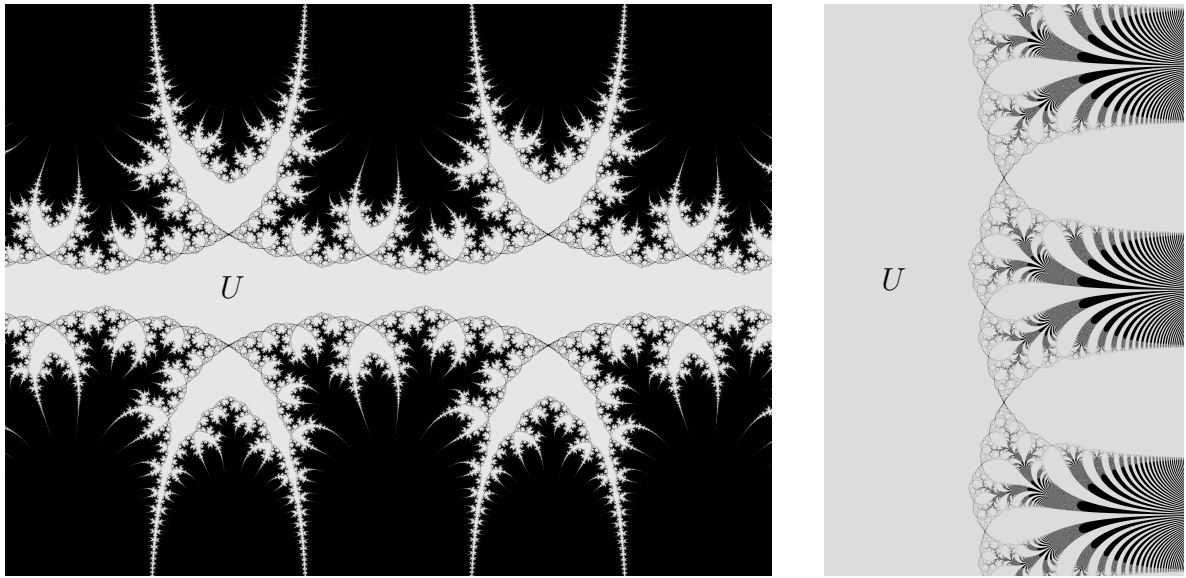

Figure 11: Baker domains of $f_1(z) = z + \frac{a}{2\pi}\sin(2\pi z) + b$ and $f_2(z) = z + 2\pi i\alpha + e^z$.

For $a = 0.6$ and $b \approx 0.61317$ the function $g_1$ has a Herman ring whose rotation number is the golden mean. The corresponding Baker domain $U$ of $f_1$ is shown in the left picture of Figure 11. The range displayed is $|\operatorname{Re} z| \leq 1.2$ and $|\operatorname{Im} z| \leq 0.9$.

Herman did not discuss the types of the Baker domains according to Cowen's classification, but it easily follows from Theorem 6.25 that the Baker domain of $f_1$ is hyperbolic while that of $f_2$ is simply parabolic.

Eremenko and Lyubich [114, Example 3] constructed a univalent Baker domain using approximation theory. For further examples and discussion of univalent Baker domains we refer to [24, Theorem 3], [47], [249] and [261].

Since, as already mentioned, a Baker domain $U$ is simply connected by Theorem 6.9, there exists a biholomorphic map $\varphi\colon \mathbb{D} \to U$. For an invariant Baker domain $U$ thus $g := \varphi^{-1} \circ f \circ \varphi\colon \mathbb{D} \to \mathbb{D}$. It turns out that $g$ is in fact an inner function.

There are a number of papers where Baker domains have been studied using this conformal map $\varphi$ and the inner function $g$. Kisaka [156, 157, 158] considered the set

$$\Theta := \left\{ e^{i\theta} \colon \lim_{r \to 1} \varphi(re^{i\theta}) = \infty \right\}.$$

Since $\infty$ is accessible in a Baker domain, it follows that $\Theta \neq \emptyset$.

The following result was proved by Baker and and Domínguez [19, Theorem 1.2].

**Theorem 6.26.** *Let $f$ be a transcendental entire function with an invariant Baker domain $U$. If $f$ is not univalent in $U$, then $\overline{\Theta}$ contains a perfect subset of $\partial \mathbb{D}$.*

Kisaka had proved this under some additional hypothesis. The hypothesis that $f$ is not univalent in $U$ is essential. In fact, there are examples of univalent Baker domains that are bounded by a Jordan curve (on the Riemann sphere); hence $\Theta$ consists of only one point. The first example of a Baker domain bounded by a Jordan curve is due to

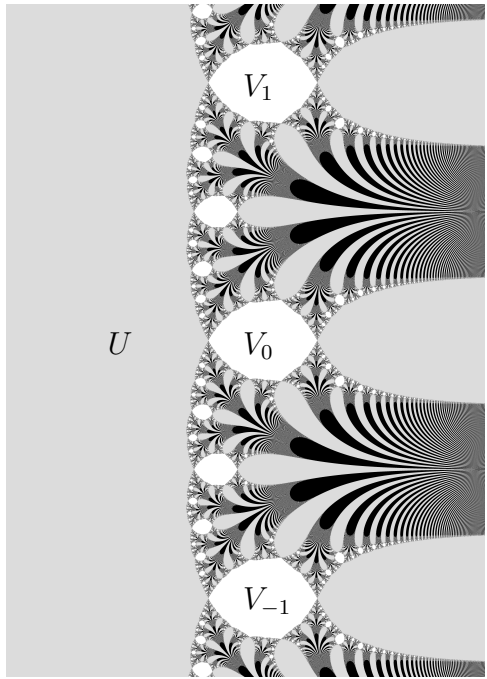

Figure 12: Baker domain of $f_3(z) = 2 - \log 2 + 2z - e^z$.

Baker and Weinreich [24, Theorem 3] who showed that there exists $\alpha \in \mathbb{R} \setminus \mathbb{Q}$ such that the Baker domain of the function $f_2$ in (6.19) has this property. They used deep results of Herman yielding that for suitable $\alpha$ the Siegel disc of the function $g_2$ given by (6.20) is bounded by a Jordan curve. For the example

$$f_3(z) = 2 - \log 2 + 2z - e^z \tag{6.21}$$

given by Bergweiler [47, Theorem 1], the argument is simpler. Here a superattracting basin is lifted via (5.16) to a univalent Baker domain.

Figure 12 shows the Baker domain $U$ of $f_3$. The range displayed is $|\operatorname{Re} z| \leq 5.8$ and $|\operatorname{Im} z| \leq 8.2$. Besides the Baker domain and its preimages, there are additional Fatou components $V_k$ (shown in white) containing the critical points $z_k := \log 2 + 2\pi i k$ of $f_3$. It turns out that $V_0$ is invariant, containing the superattracting fixed point $z_0$, while $V_k$ is wandering if $k \neq 0$, with $f_3(V_k) = V_{2k}$.

Bargmann [34, Theorem 3.1] strengthened the conclusion of Theorem 6.26 for doubly parabolic Baker domains.

**Theorem 6.27.** *Let $f$ be a transcendental entire function with an invariant Baker domain $U$. If $U$ is doubly parabolic, then $\overline{\Theta} = \partial \mathbb{D}$.*

A closely related result had been obtained by Kisaka [157]; see [34, p. 28] for a comparison of these results. On the other hand, Bergweiler and Zheng [76] gave examples of hyperbolic and simply parabolic Baker domains where $f$ is not univalent such that $\overline{\Theta}$ is a proper subset of $\partial \mathbb{D}$.

By Theorem 3.3, a cycle of immediate attracting or parabolic basins contains a singularity of the inverse, while the boundary of a Siegel disc is contained in the postsingular set. It is not clear what the relation between Baker domains and the singularities of the inverse is, if there is any. As already mentioned, Baker domains need not contain singularities of the inverse.

Theorem 6.1 implies in particular that functions in the Eremenko–Lyubich class $\mathcal{B}$ have no Baker domains. Bargmann [33] showed that in order to exclude Baker domains, a weaker hypothesis than the boundedness of $\operatorname{sing}(f^{-1})$ suffices.

**Theorem 6.28.** *Let $f$ be a transcendental entire function with an invariant Baker domain. Then there exists $K > 1$ such that $\operatorname{sing}(f^{-1}) \cap \operatorname{ann}(r, Kr) \neq \emptyset$ for all large $r$.*

Fleischmann [133] gave an example of a transcendental entire function $f$ with an invariant Baker domain for which there exist sequences $(r_n)$ and $(R_n)$ tending to infinity such that $R_n - r_n \to \infty$ and $\operatorname{sing}(f^{-1}) \cap \operatorname{ann}(r_n, R_n) = \emptyset$ for all $n$.

**Question 6.29.** How sparse can the set of singularities be for an entire function with an invariant Baker domain?

More specifically, we may ask the following question.

**Question 6.30** (Fleischmann, 2008)**.** Does there exist a transcendental entire function $f$ with an invariant Baker domain for which there exist $K > 1$ and a sequence $(r_n)$ tending to infinity such that $\operatorname{sing}(f^{-1}) \cap \operatorname{ann}(r_n, Kr_n) = \emptyset$ for all $n$?

Rempe [226, Theorem 1.4] has shown that there exists a *meromorphic* function with this property.

The following result of Bergweiler [47, Theorem 3] indicates some relation between a Baker domain and the postsingular set if the Baker domain contains no singularities of the inverse.

**Theorem 6.31.** *Let $f$ be a transcendental entire function with an invariant Baker domain $U$. If $U \cap \operatorname{sing}(f^{-1}) = \emptyset$, then there exists a sequence $(p_n)$ in $P(f)$ such that $|p_n| \to \infty$, $|p_{n+1}/p_n| \to 1$, and $\operatorname{dist}(p_n, U) = o(|p_n|)$ as $n \to \infty$.*

On the other hand, it was shown in [47, Theorem 1] that the function $f$ given by (6.21) has a Baker domain $U$ such that $\operatorname{dist}(P(f), U) > 0$. This shows that the condition $\operatorname{dist}(p_n, U) = o(|p_n|)$ in Theorem 6.31 cannot be improved to $\operatorname{dist}(p_n, U) = o(1)$. But there is still a gap between this example and Theorem 6.31.

**Question 6.32.** Can the conditions $|p_{n+1}/p_n| \to 1$ and $\operatorname{dist}(p_n, U) = o(|p_n|)$ in Theorem 6.31 be improved?

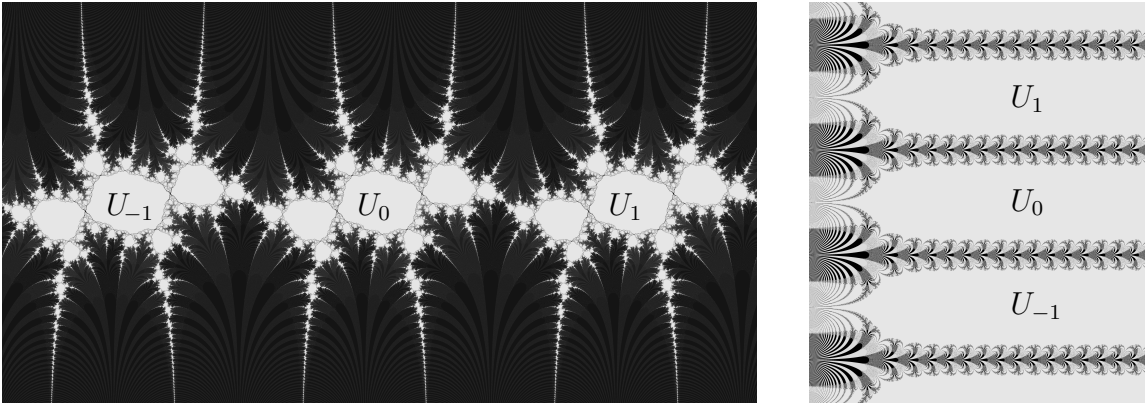

Figure 13: Wandering domains of $f_1(z) = z + \lambda \sin(2\pi z) + 1$ and $f_2(z) = z + 2\pi i - 1 + e^{-z}$.

## 6.5 Simply connected wandering domains

The first examples of simply connected wandering domains are due to Herman who showed that, for suitable $\lambda$, the functions $f_1$ and $f_2$ given by

$$f_1(z) := z + \lambda \sin(2\pi z) + 1 \quad \text{and} \quad f_2(z) := z + 2\pi i - 1 + e^{-z}$$

have simply connected wandering domains; see [146, Section II.11] for $f_1$ and [15, Example 2] and [274, p. 414] for $f_2$.

To explain why these functions have wandering domains, we consider the function $F_2(z) := z - 1 + e^{-z} = f_2(z) - 2\pi i$, which is obtained by applying Newton's method to $h(z) := e^z - 1$. Thus $F_2(z) = z - h(z)/h'(z)$. The functions $f_2$ and $F_2$ are both lifts of $g_2(z) := eze^{-z}$ under $z \mapsto \exp(-z)$, as in (5.16). Using Theorem 5.11 with $\exp(z)$ replaced by $\exp(-z)$ we find that $J(f_2) = J(F_2) = -\exp^{-1} J(g_2)$. (Actually, for these specific examples the equation $J(f_2) = J(F_2)$ is easier to prove, so we do not need to use the general Theorem 5.11.) The zeros of $h$ are given by the points $z_k := 2\pi i k$ with $k \in \mathbb{Z}$. Hence the $z_k$ are superattracting fixed points of $F_2$. Let $U_k$ be the immediate attracting basin of $z_k$ (with respect to $F_2$).

Next we note that $\exp(-z_k) = 1$ and that 1 is a superattracting fixed point of $g_2$. Let $V$ be the immediate attracting basin of 1 (with respect to $g_2$). Then $V = \exp(-U_k)$ for all $k$. Moreover, $F_2(U_k) \subset U_k$ and hence $f_2(U_k) \subset U_{k+1}$ for all $k$. But since $J(f_2) = J(F_2)$ we see that $U_k$ is also a Fatou component of $f_2$. We conclude that the $U_k$ are wandering domains of $f_2$.

The argument for $f_1$ is similar. Here $\lambda$ is chosen such that $1 + 2\pi\lambda = e^{2\pi i \alpha}$, with $\alpha \in \mathbb{R} \setminus \mathbb{Q}$ such that $F_1(z) := z + \lambda \sin(2\pi z)$ has a Siegel disc at 0. Denoting, for $k \in \mathbb{Z}$, by $U_k$ the Fatou component containing $k$ we again find that the $U_k$ are distinct and that $f_1(U_k) \subset U_{k+1}$ such that the $U_k$ are wandering domains of $f_1$.

Figure 13 shows the wandering domains $U_k$ of $f_1$ (left), with $\alpha$ chosen as the golden mean, and $f_2$ (right). The range shown is $|\operatorname{Re} z| \leq 1.5$ and $|\operatorname{Im} z| \leq 0.8$ in the left picture and $-5 \leq \operatorname{Re} z \leq 15$ and $|\operatorname{Im} z| \leq 12$ in the right picture.

Many other examples of wandering domains constructed using this method have been studied, for example $z \mapsto z + \lambda \sin z$ for suitable $\lambda \in \mathbb{R}$ was considered by Devaney and Krych [107, p. 52] and Devaney [103, p. 222] while $z \mapsto z + 2\pi + \sin z$ was studied Fagella and Henriksen [126, Example 2].

The examples constructed this way are all in the escaping set, but not in the fast escaping set. Examples of simply connected wandering domains in the fast escaping set have been given by Bergweiler [52] and Sixsmith [263]. The domains in these examples are bounded. Answering a question of Rippon and Stallard [244, p. 802, Question 1], Evdoridou, Glücksam and Pardo-Simón [120] have recently constructed examples of transcendental entire functions with unbounded, simply connected, fast escaping Fatou components.

A difference between the examples by Bergweiler [52] and Sixsmith [263] mentioned above is that the function constructed in [52] also has a multiply connected wandering domain while the function in [263] does not. Sixsmith notes [263, p. 1033] that the function $f$ that he constructs satisfies $\log M(r, f) = O((\log r)^2)$ as $r \to \infty$. Already before, Baker [18, p. 369] had asked whether there are examples of this type of arbitrarily slow growth.

**Question 6.33** (Baker, 2001)**.** Do there exist transcendental entire functions of arbitrarily slow growth for which all Fatou components are simply connected and at least one Fatou component is in the escaping set?

Baker [18], and independently also Boyd [87], had proved that there are transcendental entire functions of arbitrarily small growth for which the Fatou set consists of a single attracting basin, all components of this basin being bounded and simply connected.

We can strengthen Question 6.33 as follows.

**Question 6.34.** Do there exist transcendental entire functions of arbitrarily slow growth for which all Fatou components are simply connected and escaping?

We may ask which domains $D$ can arise as a simply connected escaping wandering domain. Boc Thaler [85, Theorem 1] proved the following.

**Theorem 6.35.** *Let $D$ be a bounded simply connected domain. Suppose that $\mathbb{C} \setminus \overline{D}$ is connected and that $\partial D = \partial(\mathbb{C} \setminus \overline{D})$. Then there is a transcendental entire function $f$ such that $D$ is an escaping wandering domain of $f$.*

In particular, every Jordan domain can arise as such a wandering domain of a transcendental entire function. Theorem 6.35 was strengthened by Martí-Pete, Rempe and Waterman [181] as follows.

**Theorem 6.36.** *Let $D$ be a bounded simply connected domain, and let $W$ be the unbounded connected component of $\mathbb{C} \setminus \overline{D}$. Suppose that $\partial W = \partial D$. Then there is a transcendental entire function such that $D$ is an escaping wandering domain of $f$.*

Domains satisfying the hypotheses of Theorem 6.36 can be rather complicated. For example, $D$ can be chosen to be a *Lakes of Wada* continuum, for which $\mathbb{C} \setminus \overline{D}$ has two or more (possibly infinitely many) connected components, all of which have boundary equal to $\partial D$. This is the first occurrence of Lakes of Wada continua in holomorphic dynamics (we refer to [181] for a discussion). Moreover, the function $f$ in Theorem 6.36 can either be chosen such that $D \subset A(f)$, or such that $D \subset I(f) \setminus A(f)$. The proof of Theorem 6.36 (as well as that of [85, Theorem 1]) uses approximation theory.

It is natural to ask whether Theorem 6.36 completely describes the bounded simply connected wandering domains of transcendental entire functions [181, Question 1.16].

**Question 6.37** (Martí-Pete, Rempe and Waterman, 2025)**.** Let $D$ be a bounded simply connected wandering domain of a transcendental entire function, and let $W$ be the unbounded connected component of $\mathbb{C} \setminus \overline{D}$. Is it true that $\partial W = \partial D$?

A domain $D$ is called *regular* if $\partial D = \partial \overline{D}$. The blowing-up property of the Julia set (Theorem 3.1 $(f)$) implies that a Fatou component of a transcendental entire function with disconnected Fatou set is regular. See [85, p. 1664] or [180, Lemma 2.2]. In particular, simply connected wandering domains are regular. So any domain for which the answer to Question 6.37 is negative would have to be a regular domain $D$ for which some point of $\partial D$ is separated from infinity by a proper subset of $\partial D$.

Simply connected wandering domains need not be in the escaping set. Eremenko and Lyubich [114, Example 1] constructed a transcendental entire function with a wandering domain where some subsequences of the iterates tend to finite limits; all domains in Theorem 6.36 can also be realised as such "oscillating" wandering domains. Bishop [80, Theorem 17.1] showed that such examples also exist in the Eremenko–Lyubich class $\mathcal{B}$. Recall that by Theorem 6.1 there are no escaping Fatou components in this class. Other examples of non-escaping wandering domains in class $\mathcal{B}$ are given in [127, 167, 182].

The behaviour of the iterates in Baker domains is described by Cowen's Theorem 6.21 while the behaviour in multiply connected Fatou components is described by the result of Bergweiler, Rippon and Stallard [73, Theorem 5.2] already mentioned, a simplified version of which is (6.16). There is no such description of the iterative behaviour in simply connected wandering domains. However, Benini, Evdoridou, Fagella, Rippon and Stallard [41] have given a classification of the dynamics in simply connected wandering domains. This classification takes into account the hyperbolic distances between iterates as well as the behaviour of orbits in relation to the boundaries of the wandering domains. Overall this leads to nine different types of behaviour. It is shown that all nine cases actually occur. We refer to [41] for more details.

## 6.6 Dynamics on the boundary of an escaping Fatou component

We begin this section with the following question asked in [242, p. 2808], which is still open.

**Question 6.38** (Rippon and Stallard, 2011)**.** Let $f$ be a transcendental entire function and let $U$ be an escaping Fatou component of $f$. Does $\partial U$ contain escaping points?

We will discuss some results related to this question. First we recall that Theorem 6.7 says that if $U$ is a Fatou component contained in $A(f)$, then even $\overline{U} \subset A(f)$. By Theorem 6.10 this applies in particular to multiply connected Fatou components.

The answer to Question 6.38 is also positive for simply connected wandering domains and hence for all wandering domains, as shown by the following result due to Rippon and Stallard [242, Theorem 1.1]. Recall here that $\omega(\cdot, E, U)$ denotes the harmonic measure of $E$ in $U$.

**Theorem 6.39.** *Let $f$ be a transcendental entire function and let $U$ be a wandering domain of $f$. If $U \subset I(f)$, then $\partial U \cap I(f) \neq \emptyset$. In fact, $\omega(\cdot, \partial U \cap I(f), U) \equiv 1$.*

Since the only escaping Fatou components that are not wandering are Baker domains or preimages thereof, Question 6.38 reduces to the following question asked by Rippon and Stallard [249, p. 802].

**Question 6.40** (Rippon and Stallard, 2018)**.** Let $f$ be a transcendental entire function and let $U$ be a Baker domain of $f$. Do we have $\partial U \cap I(f) \neq \emptyset$?

For a Baker domain $U$ of $f$ we need not have $\overline{U} \subset I(f)$. An example is given by Fatou's function $f(z) = z + 1 + e^{-z}$ considered in §2.2. It follows from the discussion there that $f$ has a Baker domain $U$ for which $\partial U = J(f)$. Thus periodic points are dense in $\partial U$. For this example we have $\partial U \cap I(f) \neq \emptyset$, but note that if $z \in U$, then $\operatorname{Re} f^n(z) \to \infty$ as $n \to \infty$, while $\operatorname{Re} f^n(z) \to -\infty$ for $z \in \partial U \cap I(f)$.

For univalent Baker domains, the answer to Question 6.40 is positive. In fact, Rippon and Stallard [249, Theorem 1.1] proved an analogue of Theorem 6.39 for such domains.

**Theorem 6.41.** *Let $f$ be a transcendental entire function and let $U$ be a univalent Baker domain of $f$. Then $\omega(\cdot, \partial U \cap I(f), U) \equiv 1$.*

Recall that univalent Baker domains are hyperbolic or simply parabolic according to the classification in Definition 6.24. Thus the following result of Barański, Fagella, Jarque and Karpińska [28, Theorem A] is a generalisation of this result.

**Theorem 6.42.** *Let $f$ be a transcendental entire function and let $U$ be a Baker domain of $f$. Suppose that $U$ is hyperbolic or simply parabolic and that the map $f|_U \colon U \to U$ has finite degree. Then $\omega(\cdot, \partial U \cap I(f), U) \equiv 1$.*

**Question 6.43.** Is the hypothesis that $f|_U \colon U \to U$ has finite degree necessary in Theorem 6.42?

For further results concerning the dynamics on the boundary of a hyperbolic or simply parabolic Baker domain we refer to a recent paper by Jové [152].

The conclusion of Theorem 6.42 result does not hold for doubly parabolic domains. In fact, for such domains a very different conclusion holds [28, Theorem B].

**Theorem 6.44.** *Let $f$ be a transcendental entire function and let $U$ be a Baker domain of $f$. Suppose that $U$ is doubly parabolic and that the map $f|_U \colon U \to U$ has finite degree. Let $D$ be the set of points in $\partial U$ whose forward orbit is dense in $\partial U$. Then $\omega(\cdot, D, U) \equiv 1$. In particular, $\omega(\cdot, \partial U \cap I(f), U) \equiv 0$.*

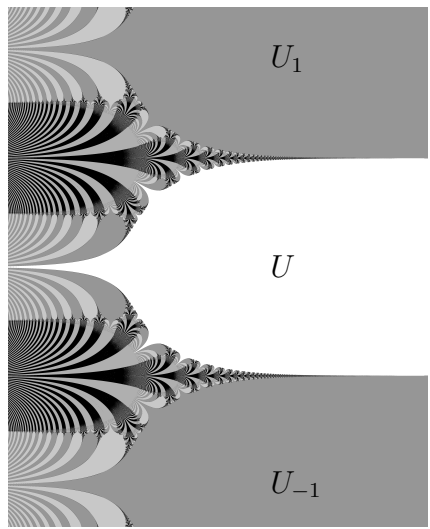

Figure 14: Baker domains of $f(z) = z + e^{-z}$.

An example is given by $f(z) := e^{-z} + z$. This function has a doubly parabolic Baker domain $U$ contained in $\{z\colon |\operatorname{Im} z| \leq \pi\}$ such that $f|_U\colon U \to U$ has degree 2; see the papers by Baker and Domínguez [19, §5] as well as Fagella and Henriksen [126, Example 3]. The dynamics of $f$ on $\partial U$ have been studied in detail by Fagella and Jové [128]. They show that $\partial U$ contains escaping points, but these points are not accessible from $U$ [128, Theorem C]. Periodic points are dense in $\partial U$ [128, Theorem D].

Actually, $f$ has infinitely many Baker domains, since for each $k \in \mathbb{Z}$ the domain $U_k := U + 2\pi i k$ is also a doubly parabolic invariant Baker domain. Figure 14 shows the Baker domain $U$ in white. Preimages of the $U_k$ with even $k$ (and thus in particular the preimages of $U = U_0$) are shown in light grey, while the Baker domains $U_k$ with odd $k$ and their preimages are in darker grey. The range displayed is $-5 \leq \operatorname{Re} z \leq 7$ and $|\operatorname{Im} z| \leq 7.5$.

The conclusion of Theorem 6.44 need not hold if $U$ is doubly parabolic and $f|_U\colon U \to U$ has infinite degree [28, Example 1.3].

The following result of Rippon and Stallard [242, Theorem 1.2] is a counterpart of Theorem 6.39.

**Theorem 6.45.** *Let $f$ be a transcendental entire function and let $U$ be a Fatou component of $f$.*

*(a) If $\omega(\cdot, \partial U \cap I(f), U) \not\equiv 0$, then $U \subset I(f)$.*

*(b) If $\omega(\cdot, \partial U \cap A(f), U) \not\equiv 0$, then $U \subset A(f)$.*

Every example of a simply connected wandering domain $D \subset I(f)$ that we have so far encountered in fact satisfies $\overline{D} \subset I(f)$. This prompted Rippon to ask the following question [144, Problem 2.94]: If $U$ is a bounded escaping wandering domain of a transcendental entire function $f$, is $\partial U \subset I(f)$? This question was answered in the negative by Martí-Pete, Rempe and Waterman [181, Theorem 1.6]. In fact, in [181, Theorem 1.7], they show that each of the domains whose existence is asserted in Theorem 6.36 may be chosen to have such points on the boundary. In particular, the result implies the following.

**Theorem 6.46.** *Let $D$ be a bounded simply connected domain, and let $W$ be the unbounded connected component of $\mathbb{C} \setminus \overline{D}$. Suppose that $\partial W = \partial D$, and let $Z \subset \partial W$ be countable. Then there is a transcendental entire function such that $D$ is an escaping wandering domain of $f$ and $Z \cap I(f) = \emptyset$.*

We note that it is also possible to realise $D$ as an oscillating wandering domain for which all points of $Z$ are escaping.

In view of Theorem 6.46, it is natural to ask how large the set of non-escaping points on the boundary of a simply connected escaping wandering domain may be. Theorem 6.39 says that it has harmonic measure zero, while Theorem 6.46 implies that it may be dense in the boundary.

To state a specific question we note that a theorem of Beurling [215, §9.5] says that if $\varphi \colon \mathbb{D} \to U$ is biholomorphic, then there exists a subset $E$ of $\partial \mathbb{D}$ of capacity zero such that the radial limit
$$\varphi(\xi) := \lim_{r \to 1} \varphi(r\xi)$$
exists for all $\xi \in \partial \mathbb{D} \setminus E$. For the definition and properties of capacity we refer to the books by Pommerenke [215, §9] and Ransford [216, §5]. Here we only note that having capacity zero is a much stronger restriction than having Lebesgue measure zero. In fact, sets of capacity zero have Hausdorff dimension zero.

Theorem 6.39 says that the set of all $\xi \in \partial \mathbb{D}$ for which $\varphi(\xi) \notin I(f)$ has measure zero (with respect to the 1-dimensional Lebesgue measure on $\partial \mathbb{D}$). The following question was asked by Bishop [79, p. 133]. See also [181, Question 1.12].

**Question 6.47** (Bishop, 2014)**.** Let $U$ be a simply connected escaping wandering domain of a transcendental entire function. Let $\varphi \colon \partial \mathbb{D} \to U$ be biholomorphic. Does the set of all $\xi \in \partial \mathbb{D}$ for which $\varphi(\xi) \notin I(f)$ have capacity zero?

This would be best possible since it was shown in [181, Theorem 1.13] that if $E$ is a compact subset of $\partial \mathbb{D}$ of capacity zero, then there exists a transcendental entire function $f$ for which $\mathbb{D}$ is a wandering domain and $E \subset J(f) \setminus I(f)$.

# 7 Topological properties of the escaping set

In this section, we turn our attention to the topological properties of the escaping set. As we saw in our initial discussion of examples, the escaping set often contains or even consists of "Cantor bouquets" of curves tending to infinity. A very different structure arises for functions with multiply connected wandering domains: these domains contain larger and larger annuli surrounding zero. Such structures are extremely useful for studying the wider dynamical behaviour of the functions under consideration. At the same time, topological properties of the escaping set can be surprisingly subtle. One reason for this is that the escaping set is not a topologically simple set. By definition, it is a countable intersection of countable unions of closed sets, or an $F_{\sigma\delta}$:

$$I(f) = \bigcap_{j\geq 0} \bigcup_{k\geq 0} \{z \in \mathbb{C} \colon |f^n(z)| \geq j \text{ for all } n \geq k\}.$$

As mentioned in §5, Rempe [225] proved that, on the other hand, $I(f)$ is never a countable union of closed sets (an $F_\sigma$). This may help to explain some of the topological subtleties that arise in its study. Lipham [172] showed that for the exponential family given by (2.4), the escaping set is not even $G_{\delta\sigma}$. It is not known whether this fact holds in general.

**Question 7.1.** Is there a transcendental entire function $f$ such that $I(f)$ is a $G_{\delta\sigma}$ set? If yes, can $f$ be chosen such that $f \in \mathcal{B}$, or even $f \in \mathcal{S}$?

As discussed in §2, Fatou [132, p. 369] noticed the existence of curves to infinity in the escaping sets of certain examples. He then asked whether such curves exist for "much more general" transcendental entire functions: "Il serait intéressant de rechercher si cette proprieté n'appartiendrait pas à des substitutions beaucoup plus générales."

Motivated by examples such as the ones discussed in §2, Eremenko [113, p. 334] strengthened Fatou's question (and made it more precise) as follows.

**Conjecture 7.2** (Eremenko, 1989). Let $f$ be a transcendental entire function. Then every point of $I(f)$ can be connected to $\infty$ by a curve in $I(f)$.

He also proposed the following weaker version [113, p. 333].

**Conjecture 7.3** (Eremenko, 1989). Let $f$ be a transcendental entire function. Then the set $I(f)$ has no bounded connected components.

Conjectures 7.3 and 7.2 have become known as *Eremenko's conjecture* and the *strong Eremenko conjecture*. They have motivated a substantial amount of research in transcendental dynamics over the last several decades. While both conjectures turn out to be false in general, there are large classes of functions where they hold, and where the corresponding structure in the escaping set has been used successfully to study the wider dynamical behaviour.

There is another type of structure often encountered in the escaping set, different from the case of a "Cantor bouquet": that of a "spider's web". It arises when $I(f)$ is connected and contains the boundaries of a sequence of nested simply connected domains that exhaust the plane; see Definition 7.29 below for the formal definition. For example, a spider's web occurs when $f$ has a multiply connected Fatou component (see Theorem 7.28), but it may appear also when no such components are present. When it does, this has important consequences for the structure of the dynamical plane. We discuss this property in §7.4.

The two structures mentioned are not mutually exclusive. For example, in the case of Fatou's function (2.3) we saw in §2.2 that the escaping set contains a Cantor bouquet, but, as shown by Evdoridou [118], it is also a spider's web; see Theorem 7.47 below. Further examples of entire functions where both structures coexist are discussed in Section 7.4. On the other hand, there also exist entire functions whose escaping sets contain neither Cantor bouquets nor spiders' webs; see Remark 7.50.

## 7.1 "Hairs" in the escaping set

As we saw in Theorem 2.1, every escaping point of the function $f\colon \mathbb{C} \to \mathbb{C}, z \mapsto \sin(z)/2$ studied by Fatou can be connected to infinity by a curve $\gamma$ consisting entirely of escaping points (in fact, the iterates of $f$ tend to infinity uniformly on $\gamma$). Thus $f$ satisfies the conclusion of the strong Eremenko conjecture (Conjecture 7.2). The same is true for the functions $z \mapsto ae^z$ with $a \in \mathbb{C} \setminus \{0\}$, as shown in Theorem 2.5.

These curves are called "hairs" or "Devaney hairs," in honour of Devaney who pioneered their study. They are also sometimes referred to as "external rays" or "dynamic rays", to stress the analogy to external rays of polynomials. For a polynomial of degree at least 2, such rays foliate the escaping set, which is the basin of attraction of the superattracting fixed point at $\infty$. Often, several such rays land together at a common non-escaping point, and the resulting structure plays a crucial role in the study of polynomial dynamics; see [190, 191]. Curves (and even more general unbounded connected subsets) in the escaping can be used in a similar way to study transcendental entire functions; see, e.g., papers by Rempe and Schleicher [229] and Benini [40] for examples of this approach in the context of the exponential family.

The strongest currently known positive result concerning the strong Eremenko conjecture is the following result of Rottenfußer, Rückert, Rempe and Schleicher [251, Theorem 1.2]. Recall that the order of an entire function was defined in (6.17).

**Theorem 7.4.** *Suppose that $f \in \mathcal{B}$ has finite order, or that, more generally, $f$ is a finite composition of finite-order functions in $\mathcal{B}$. Then every point of $I(f)$ can be connected to $\infty$ by an arc $\gamma$ in $I(f)$ such that $f^n|_\gamma \to \infty$ uniformly.*

Theorem 7.4 was proved independently by Barański [25] for the important special case of a finite-order function of disjoint type; see Definition 3.4. The paper [25] appeared before [251], but (as acknowledged in both papers) the results were proved independently and announced around the same time.

The proof of Theorem 7.4 uses a similar idea as the proof we gave of Theorem 2.1. (Indeed, our proof of Theorem 2.1 was adapted from the general proof given in [251].) The underlying idea is as follows. Apply the logarithmic change of variable to $f$, as described in §6.1, and denote the resulting function by $F\colon W \to H$. By choosing the number $R$ used in the definition of $F$ sufficiently large, we may assume that $\operatorname{Re} z > 1$ for all $z \in W$. The condition that $f$ has finite order implies that the connected components of $W$ have two particular geometric properties: *bounded slope* (the arguments within each component are bounded away from $\pi/2$ and $-\pi/2$) and *bounded wiggling* (any point $z \in W$ can be connected to infinity by a curve in $W$ consisting of points with real part at least $c \cdot \operatorname{Re} z$, where $c$ is a constant depending on $f$). These properties in turn imply analogues of (i) and (ii) in the proof of Theorem 2.1. More precisely, there exists a constant $A > 1$ with the following properties.

(i) Let $U, V$ be connected components of $W$, and let $z, w \in U$ with $F(z), F(w) \in V$. If $\operatorname{Re} w \geq A \cdot \operatorname{Re} z$, then also $\operatorname{Re} F(w) > A \cdot \operatorname{Re} F(z)$.

(ii) Let $z, w \in W$ be such that, for every $n \geq 0$, the points $F^n(z)$ and $F^n(w)$ are defined and belong to the same connected component of $W$. Then there is some $n$ such that either $\operatorname{Re} F^n(w) > A \cdot \operatorname{Re} F^n(z)$ or $\operatorname{Re} F^n(z) > A \cdot \operatorname{Re} F^n(w)$.

In other words, if $w$ has a sufficient "head start" over $z$ (i.e., $w$ starts sufficiently far to the right of $z$), then this remains true of their images. Furthermore, of any two points $z$ and $w$ whose orbits remain in the same connected component of $W$, one will eventually gain a head start over the other. This type of condition is called a "head-start condition"; see [251, Definition 4.1]. Theorem 7.4 follows from this condition in the same manner as Theorem 2.5. We refer to [251] for the details of the proof, and to Sixsmith's survey [268, §8.2] for further discussion.

Theorem 7.4 shows that there is a large and natural class of transcendental entire functions for which the strong Eremenko conjecture holds. In some sense, it thus gives a positive answer to Fatou's question from [132] that we discussed at the beginning of the section. On the other hand, without the assumption of finite order, the conjecture is false in general, as shown in [251, Theorem 1.1].

**Theorem 7.5.** *There exists a function $f \in \mathcal{B}$ such that $J(f)$, and hence $I(f)$, contains no curve to $\infty$.*

The counterexample $f$ is constructed to be of disjoint type. Recall that, in the proof of Theorem 7.4, the key fact was that the connected components of $\{z \colon |f(z)| > R\}$ are geometrically tame, which in turn implies that the function satisfies a "head-start condition". The proof of Theorem 7.5 turns this idea around. The authors first construct a "model" for the function: a simply connected domain $T \subset \mathbb{C}$ and a universal covering $g \colon T \to \mathbb{C} \setminus \overline{\mathbb{D}}$, such that for the function $g$, the head-start condition fails: There are many points $z$ and $w$ such that there are infinitely many $n$ such that $|g^n(z)| > |g^n(w)|$ but $|g^{n+1}(w)| > |g^{n+1}(z)|$. This is achieved by drawing the domain $T$ in such a way that it has large "wiggles" (where the real parts along the domain first increase, then decrease, and then increase again).

This model is then realised by a transcendental entire function $f$ using an approximation result using Cauchy integrals; the construction is sufficiently robust that the same properties as for the model hold for the resulting entire function. We refer to [251] for the details of the proof, and to §8 and, in particular, Theorem 8.17 for a more general discussion concerning the realisation of model functions.

The paper [251] left open the question whether the strong Eremenko conjecture holds in the Speiser class $\mathcal{S}$. This was addressed by Bishop [80, §18], who used his technique of *quasiconformal folding* to show that the model from [251] can be realised, in a certain sense, by a Speiser class function. The same is true for all variations on Theorem 7.5 discussed below.

Following [221, Appendix A], let us say that a function has the *strong Eremenko property* if it satisfies the conclusion of Conjecture 7.2. In view of Theorems 7.4 and 7.5, it is natural to ask what can be said about the growth of a function for which the strong Eremenko property fails. The following is proved in [251, Theorem 8.1].

**Theorem 7.6.** *There exists a constant $M > 1$ such that the function in Theorem 7.5 may be chosen such that*
$$|f(z)| \leq \exp(\exp(|z|^M)) \qquad (7.1)$$
*whenever $|z|$ is sufficiently large.*

In [251, §8], the authors also sketch how to modify the construction so that (7.1) holds for *every* $M > 1$. In other words,
$$\frac{\log\log\log|f(z)|}{\log\log|z|} \to 1 \qquad (7.2)$$
as $|f(z)| \to \infty$. (More precisely, in [251, §8], the authors sketch how to construct a model, in the sense discussed above, for which (7.2) holds and such that there is no curve to infinity in the analogue of the escaping set. That this model can be realised by a transcendental entire function $f \in \mathcal{B}$ was first shown by Rempe [232] and also follows from the work of Bishop [81]; see Theorem 8.17.)

Recently, Brown [88] has obtained significant improvements to the estimates (7.1) and (7.2). His growth condition is somewhat technical, but implies in particular that $f$ can be constructed such that, for all $\alpha > 0$,
$$\frac{\log\log\log|f(z)|}{\log\log|z|} < 1 + \frac{1}{(\log\log\log|z|)^\alpha}$$
when $|f(z)|$ is sufficiently large. Brown's analysis also suggests that the condition of finite order in Theorem 7.4 may not be optimal, i.e. that the strong Eremenko conjecture holds for all functions $f \in \mathcal{B}$ of sufficiently small infinite order of growth.

**Question 7.7.** What growth conditions on $f \in \mathcal{B}$ ensure that $f$ has the strong Eremenko property (that is, every point of $I(f)$ can be connected to infinity by a curve in $I(f)$)?

While this question asks to what extent the hypothesis that $f$ has finite order is sharp in Theorem 7.4, one may also ask whether the hypothesis that $f$ is in $\mathcal{B}$ is essential.

**Question 7.8.** Does every transcendental entire function of finite order have the strong Eremenko property?

In [251, Theorem 7.5], the authors sketch a proof of the following stronger version of Theorem 7.5.

**Theorem 7.9.** *There exists a function $f \in \mathcal{B}$ such that $J(f)$, and hence $I(f)$, does not contain any non-degenerate curve.*

This result was subsequently strengthened further (see [224, Theorem 1.5] or [44, Theorem 1.4]) to show that $f$ can be chosen such that every continuum (i.e., a compact and connected set) contained in $J(f)$ is indecomposable (i.e., cannot be written as the union of two proper subcontinua). The constructions are inductive and hence do not immediately yield estimates on the growth of the function, but it is clear from the proof that the resulting maps grow much faster than (7.1).

**Question 7.10.** What growth conditions on $f \in \mathcal{B}$ ensure that $I(f)$ contains an arc?

Benitez and Rempe [44, Theorem 1.5] have shown that the *lower order* (defined by replacing $\limsup$ by $\liminf$ in (6.17)) of a function as in Question 7.10 may be finite, and indeed equal to $1/2$. A classical result of Wiman implies that this is the minimum lower order achievable for an entire function bounded on a curve to infinity, and hence for any function in class $\mathcal{B}$; see [60, Proof of Corollary 2] or [165, p. 1788]. Compare also [237, Lemma 3.5] and [227, Corollary 1.2].

In view of Theorems 7.4 and 7.5, it is natural to study classes of functions for which hairs do exist. The following definition is due to Benini and Rempe [42, Definition 1.2].

**Definition 7.11.** An entire function $f$ is called *criniferous* if, for every $z \in I(f)$, there exist $n_0$ and a sequence $(\gamma_n)_{n=n_0}^{\infty}$ of arcs with the following properties.

- For all $n \geq n_0$, the arc $\gamma_n$ connects $f^n(z)$ to $\infty$.

- For all $n \geq n_0$, the map $f$ maps $\gamma_n$ homeomorphically onto $\gamma_{n+1}$.

- The arcs $\gamma_n$ tend uniformly to infinity as $n \to \infty$.

By Theorem 7.4, all functions $f \in \mathcal{B}$ of finite order, and their compositions, are criniferous. The argument sketched at the end of the proof of Theorem 7.4 shows that all criniferous functions have the strong Eremenko property. It is natural to ask whether the converse holds, at least in the class $\mathcal{B}$.

**Question 7.12.** Let $f \in \mathcal{B}$ have the strong Eremenko property. Is $f$ necessarily criniferous?

If $\gamma$ is an arc in the escaping set $I(f_a)$ of an exponential map (2.4), then, as shown by Förster, Rempe and Schleicher [136, Theorem 4.3], $f_a^n|_\gamma \to \infty$ uniformly. This is not true for general $f \in \mathcal{B}$, even those that are of disjoint type, by a result of Rempe [224, Theorem 1.6]. Nonetheless, Pardo-Simón and Rempe [211, Theorem 1.2] show that the answer to Question 7.12 is positive in this setting.

**Theorem 7.13.** *Suppose that $f \in \mathcal{B}$ is of disjoint type and has the strong Eremenko property. Then $f$ is criniferous.*

They also ask [211, Question 1.3] whether the converse holds.

**Question 7.14** (Pardo-Simón and Rempe, 2023)**.** Is there a disjoint-type entire function $f$ such that every point in $I(f)$ is contained in an unbounded connected set on which $f^n \to \infty$ uniformly, but such that $f$ is not criniferous?

For exponential maps, [136, Theorem 4.3] implies, in particular, that all path-connected components of $I(f)$ are curves (each obtained as a union of connected components of iterated preimages of arcs as constructed in Proposition 2.6). For a general criniferous function, the path-connected components of $I(f)$ may be more complicated. Indeed, if $f$ has an escaping critical point $c$ that is eventually mapped to the interior of one of the arcs from Definition 7.11, then $c$ will be a branch point at which more than two different arcs in the escaping set meet. Hence a path-connected component of $I(f)$ may contain the injective continuous image of a (possibly infinite) tree; compare [244, Example after Theorem 5.1]. It is natural to ask whether this is the only complication that can occur.

**Question 7.15.** Suppose that $f \in \mathcal{B}$ is criniferous and that $I(f)$ contains no critical points. Let $C$ be a path-connected component of $I(f)$. Is $C$ the continuous injective image of an open or half-open interval?

As already mentioned, [136, Theorem 4.3] gives a positive answer to this question for exponential maps; in fact, it does more than this, providing a complete description of the path-connected components of the escaping set of an exponential map. To explain this and strengthen Question 7.15, let us begin with a general observation. This concerns arcs in $I(f)$ on which the iterates tend to infinity uniformly, when $f \in \mathcal{B}$.

**Proposition 7.16.** *Let $f \in \mathcal{B}$ and suppose that $I(f)$ contains no critical points of $f$. Let $z \in I(f)$, and let $C$ be the set of all points $w \in I(f)$ that can be connected to $z$ using an arc on which $f^n$ tends to $\infty$ uniformly. If $C$ is non-degenerate, then $C$ is the continuous and injective image of an open, half-open or closed interval. Moreover, $f^n$ is injective on $C$ for all $n \geq 0$, and if $\alpha \subset C$ is a closed arc, then $f^n|_\alpha \to \infty$ uniformly.*

*Proof.* In the case of exponential maps, this follows from the proof given by Schleicher and Zimmer [256, Corollary 6.9]. In our more general setting, the proof can be deduced from a result of Rempe [224] using similar ideas, but we are not aware of the result having been previously stated in this form. We therefore sketch the argument, using results and ideas from [42], [224] and [221]. These are not employed elsewhere in this survey, so readers unfamiliar with the techniques may choose to skip the proof.

We first show the following. If $\alpha$ and $\beta$ are arcs on each of which the iterates of $f$ tend to infinity uniformly and such that $\alpha \cap \beta \neq \emptyset$, then $\alpha \cup \beta$ is also an arc. If $f$ is of disjoint type, this follows directly from [224, Theorem 1.4]. Indeed, by this result, $\alpha \cup \beta$ is contained in a continuum of *span zero*, and hence itself has span zero. See [224, Definition 2.1] for the definition of span zero continua; what is important for us is that a span zero continuum contains neither a simple closed curve nor a simple triod (the union of three arcs that meet at a single common point). Hence any locally connected non-degenerate span zero continuum is an arc, which yields the desired result when $f$ is of disjoint type.

If $f$ is not of disjoint type, then using the same method of proof as in [224, Theorem 1.4], one may deduce that $A_n := f^n(\alpha) \cup f^n(\beta)$ is of span zero, and hence an arc, for sufficiently large $n$. (Alternatively, by [221, Theorem 1.1], for large $n$ the set $A_n$ is homeomorphic to a subset of the escaping set of a disjoint-type function; see §8. Hence the claim also follows directly from [224, Theorem 1.4].) Since $f$ has no critical points, and $\alpha \cup \beta$ is bounded, it follows that there is a branch of $f^{-n}$ that maps $A_n$ to $\alpha \cup \beta$, which is thus an arc, as claimed.

Now let $C$ be as in the statement of the proposition. Using the claim, it is straightforward to deduce that $C$ can be written as a countable increasing union of sub-arcs, on each of which the iterates tend to infinity uniformly. This implies the first claim of the proposition. That is, there is a bijective and continuous function $\gamma\colon I \to C$, where $I = (0,1)$, $I = (0,1]$ or $I = [0,1]$. (Compare [42, Proposition 5.3].) Moreover, if $J \subset I$ is a compact sub-interval, then $f^n$ tends to $\infty$ uniformly on $\gamma(J)$.

It remains to show that every sub-arc of $C$ is of the form $\gamma(J)$ with $J$ is above. This means that we must show that $\gamma(t)$ cannot tend to a limit on $C$ as $t \to 0 \notin I$ or as $t \to 1 \notin I$. This follows from the fact that there are curves disjoint from $C$ that accumulate on $C$ from both sides by [224, Proposition 8.1], again via [221, Theorem 1.1]; see [42, Proof of Corollary 6.7] for a similar argument. We omit the details. □

For exponential maps, the sets $C$ appearing in Proposition 7.16 are exactly the path-connected components of $I(f)$. As mentioned above, this is not true in general, but is it the case when (in addition to the hypothesis of Proposition 7.16) $f$ is criniferous? In view of the final statement of Proposition 7.16, this is equivalent to the following, which in particular would imply a positive answer to Question 7.15.

**Question 7.17.** Suppose that $f \in \mathcal{B}$ is criniferous, and that $I(f)$ contains no critical points. If $\gamma \subset I(f)$ is an arc, does $f^n|_\gamma \to \infty$ uniformly?

The proof of this fact for exponential maps relies on a result of Schleicher and Zimmer [256, Corollary 6.9], which uses the periodicity of the exponential map in an essential manner. The same idea gives a positive answer to Questions 7.15 and 7.17 for functions in the family $z \mapsto a\sin(z) + b\cos(z)$; see the paper by Rottenfußer and Schleicher [252]. The question is open (and interesting) even for general functions $f \in \mathcal{S}$ of finite order.

Barański, Jarque and Rempe [29] showed that the functions covered by Theorem 7.4 have a stronger property than being criniferous: there is a certain continuity in the dependence of the curves $\gamma_n$ in Definition 7.11 on the point $z$. More precisely, suppose that $K \subset \mathbb{C}$ is a compact set such that $f^n|_K \to \infty$ uniformly. Then there is $n_0$ such that the curves $\gamma_n = \gamma_n(z)$ in Definition 7.11 are defined for all $z \in K$ and all $n \geq n_0$, the curve $\gamma_{n_0}(z)$ depends continuously on $f^{n_0}(z)$ with respect to the Hausdorff metric, and $\gamma_n(z) \to \infty$ uniformly in $z$ as $n \to \infty$. Pardo-Simón and Rempe [211, Theorem 1.5] show that this property characterises a certain subclass $\mathcal{CB}$ of $\mathcal{B}$ that had been previously introduced by Pardo-Simón [208]. In [211, Theorem 1.1], it is also shown that there exist criniferous functions $f \in \mathcal{B} \setminus \mathcal{CB}$. We refer to these articles for further details.

So far, except for Fatou's function studied in §2.2, we have discussed the existence of hairs only for functions in the Eremenko–Lyubich class $\mathcal{B}$, which is a natural setting due to their strong expansion properties near infinity. The methods and results used to study these functions naturally apply to some larger classes of functions that share similar behaviour. For example, suppose that $f$ is a transcendental entire function with a logarithmic tract $T$. If the order of $f$ restricted to $T$ is finite, then it follows from [251] that those points that escape to infinity through $T$ are organised in hairs. Similar results hold for functions such as Fatou's function (2.3), which are related to class $\mathcal{B}$ functions by a logarithmic change of variable as in (5.16).

Dourekas [111] constructed hairs in the escaping set for a class of functions that have a single, non-logarithmic, direct tract. Although these functions have no logarithmic tracts, they are obtained as sums of exponentials, and with a careful analysis Dourekas is able to apply similar methods as those used to study exponential maps in certain parts of the complex plane.

A setting that is markedly different to all of the ones mentioned so far is given by functions with multiply connected wandering domains. We are not aware of any work studying the strong Eremenko property in this context. By Theorem 6.8, any point of $\mathbb{C}$ is separated from $\infty$ by sufficiently high iterates of this wandering domain, so any curve to $\infty$ in $I(f)$ would have to pass through infinitely many Fatou components. This is in contrast to the hairs mentioned above, which are contained in $J(f)$.

It is plausible that at least Kisaka and Shishikura's example from Theorem 6.16 and Bishop's function from [83], both of which are well-controlled as a result of their construction, satisfy the strong Eremenko property. The following is a natural question.

**Question 7.18.** Does there exist a transcendental entire function $f$ with a multiply connected Fatou component such that $I(f)$ does not contain a curve to $\infty$?

As we shall see below, different multiply connected Fatou components of $f$ are always connected by continua in $I(f)$ (and even $A(f)$), so the question is whether these continua are curves.

We conclude this subsection on hairs by recalling Theorem 2.5, which says that for maps in the exponential family given by (2.4) the escaping set consists of hairs. Viana [282] proved that these hairs are $C^\infty$-differentiable. He asked the following question that is still open.

**Question 7.19** (Viana, 1988)**.** Are the hairs in the escaping sets of the exponential family analytic?

Kisaka and Shishikura [161] showed that Viana's result holds for more general functions. In particular, they extended it to functions of the form $pe^q$ with polynomials $p$ and $q$. On the other hand, Comdühr [95] showed that their exist functions in $\mathcal{B}$ of finite order for which the hairs in the escaping set (which exist by Theorem 7.4) are nowhere differentiable. It would be interesting to find further conditions under which the hairs are (or are not) smooth, say $C^1$, $C^2$, or $C^\infty$. It appears likely that there are finite-order functions with hairs that are $C^n$ but not $C^{n+1}$. We also note that it should be possible to adapt Comdühr's argument to construct hairs that are smooth at some points but not at others.

The preceding discussion concerns differentiability of the hairs at interior points; when a hair has an escaping endpoint, we may also ask about the differentiability in this endpoint. For the exponential family, a necessary and sufficient condition for differentiability in an escaping endpoint was given by Rempe [218, Theorem 6.2]; this condition depends on the combinatorics of the hair in question, but is independent of the specific choice of exponential map. We refer to [218] for further details.

## 7.2 Topology of escaping compacta

We mentioned in Theorem 6.36 that certain bounded simply connected domains can be exactly realised as escaping Fatou components of transcendental entire functions. In fact, Theorem 6.36 is a special case of the following more general result by Martí-Pete, Rempe and Waterman [181, Theorem 1.7], concerning wandering compact subsets of $I(f)$.

**Theorem 7.20.** *Suppose that $K \subset \mathbb{C}$ is compact and that $\mathbb{C} \setminus K$ is connected. Then there exists a transcendental entire function $f \colon \mathbb{C} \to \mathbb{C}$ such that $f^n|_K \to \infty$ uniformly, such that $\partial K \subset J(f)$ and such that $f^n(K) \cap f^m(K) = \emptyset$ for $n \neq m$.*

This immediately implies Theorem 6.36, setting $K \coloneqq \mathbb{C} \setminus W$. As mentioned there, the proof of Theorem 7.20 uses approximation theory and builds on previous work by Boc Thaler [85].

The proof shows that the function $f$ may be chosen such that $f^n \to \infty$ uniformly on $K$, and such that either $K \subset A(f)$ or $K \subset I(f) \setminus A(f)$. Moreover, the construction ensures that there is no path in $I(f)$ connecting a point of $K$ to a point of $I(f) \setminus K$. This yields an alternative, and more elementary, proof that there are entire functions without the strong Eremenko property.

**Theorem 7.21.** *Let $K \subset \mathbb{C}$ be compact and path-connected, and suppose that $\mathbb{C} \setminus K$ is connected. Then there exists a transcendental entire function $f$ such that $K$ is a path-connected component of $I(f)$ which is contained in $A(f)$. In particular, no point of $K$ can be connected to infinity by a curve in $I(f)$.*

Since the function in this theorem is constructed using approximation theory, there is no information about the set of singular values of $f$, or its rate of growth, in contrast to Theorem 7.5.

Theorem 7.20 shows that many plane continua can be realised as wandering subsets of $I(f) \cap J(f)$. Following Blokh, Oversteegen and Timorin [84], a continuum $X \subset \mathbb{C}$ is called *unshielded* if it agrees with the boundary of the unbounded connected component of $\mathbb{C} \setminus X$. Then $X$ is unshielded if and only if it is the boundary of a continuum $K$ as in Theorem 7.20. So every unshielded plane continuum is a wandering subset of some transcendental entire function. There are, however, many continua that are not unshielded, such as the union of two circles that are tangent to each other.

**Question 7.22.** *Let $K \subset \mathbb{C}$ be a planar continuum. When does there exist a transcendental entire function $f$ such that $K \subset I(f) \cap J(f)$, and such that $f^n(K) \cap f^m(K) = \emptyset$ for $n \neq m$?*

## 7.3 Positive results on Eremenko's conjecture

Recall that Eremenko's conjecture asks whether $I(f)$ only has unbounded connected components, for every transcendental entire function $f$. A number of weaker statements are known to be true. The first of these was shown by Eremenko himself [113, Theorem 3], and partly motivated his conjecture.

**Theorem 7.23.** *Let $f$ be a transcendental entire function. Then every connected component of $\overline{I(f)}$ is unbounded.*

*Proof.* Let $C$ be a bounded component of $\overline{I(f)}$. Then there exists a component $V$ of $\mathbb{C} \setminus \overline{I(f)}$ that separates $C$ from $\infty$. Since $J(f) = \partial I(f) \subset \overline{I(f)}$ by Theorem 4.3, we see that $V \subset F(f)$. Thus $V$ is contained in a non-escaping Fatou component $U$. By Theorem 6.12, $U$ is simply connected. This implies that $C \subset U \subset \mathbb{C} \setminus \overline{I(f)}$, a contradiction. □

In Theorem 5.5, we saw that the fast escaping set $A(f)$ intersects every curve surrounding a sufficiently large disc. Using a more careful argument, Rippon and Stallard [239, Theorem 1] proved the stronger result that all connected components of $A(f)$ are unbounded. In other words, Eremenko's conjecture holds when $I(f)$ is replaced by $A(f)$. Later, Rippon and Stallard [244, Theorem 1.1] strengthened the result even further by showing that the same statement is true even for the levels of $A(f)$ defined in (5.10).

**Theorem 7.24.** *Let $f$ be a transcendental entire function, $R$ an admissible radius for $f$, and $L \in \mathbb{Z}$. Then every connected component of $A_R^L(f)$ is unbounded.*

*In particular, every connected component of $A(f)$ is unbounded.*

*Proof.* By the "boundary bumping theorem" (see Nadler's book [195, Theorem 5.6]), a closed subset $A$ of $\mathbb{C}$ has only unbounded connected components if and only if $A \cup \{\infty\}$ is connected. For $n \geq \max\{0, -L\}$, consider the set

$$E_n := f^{-n}(\mathbb{C} \setminus D(0, M^{n+L}(R))) \cup \{\infty\}.$$

Note that for $L = 0$ we have $E_n = C_n \cup \{\infty\}$ where $C_n$ is the set already considered in (5.13).

The preimage $f^{-1}(A)$ of an unbounded closed connected set $A$ by an open mapping $f \colon \mathbb{C} \to \mathbb{C}$ has only unbounded connected components. (Indeed, otherwise there is a bounded open set $U \subset \mathbb{C}$ with $U \cap f^{-1}(A) \neq \emptyset$ and $\partial U \cap f^{-1}(A) = \emptyset$. But then $\partial f(U) \subset \mathbb{C} \setminus A$ separates some point of $A$ from infinity.)

So all $E_n$ are connected. We have $E_{n+1} \subset E_n$ by the maximum modulus principle, so

$$A_R^L(f) \cup \{\infty\} = \bigcap_{n=\max\{0,-L\}}^{\infty} E_n$$

is connected as a nested intersection of compact connected sets. This proves the first claim; the second follows because $A(f) = \bigcup_{L=-\infty}^{-1} A_R^L$ by (5.11). □

Since $A(f) \neq \emptyset$ by Theorem 5.2, this strengthens Theorem 4.14 as follows.

**Corollary 7.25.** *Let $f$ be a transcendental entire function. Then $I(f)$ has at least one unbounded connected component.*

Since Eremenko's conjecture asks whether *every* connected component of $I(f)$ is unbounded, Corollary 7.25 can be considered as establishing a weak version of Eremenko's conjecture. Using Theorem 6.39, Rippon and Stallard [242, Theorem 4.1] proved another weakened version of Eremenko's conjecture, as follows.

**Theorem 7.26.** *Let $f$ be a transcendental entire function. Then $I(f) \cup \{\infty\}$ is connected.*

*Proof.* The claim is equivalent to the following statement: If $U$ is a bounded domain that intersects $I(f)$, then $\partial U \cap I(f) \neq \emptyset$. If $U$ intersects $J(f)$, this follows from Theorem 4.14. Otherwise, $U$ is an escaping wandering Fatou component, and $\partial U \cap I(f) \neq \emptyset$ by Theorem 6.39. □

At first glance, the statement of Theorem 7.26 appears to be very similar to Eremenko's conjecture. Indeed, if $A \subset \mathbb{C}$ is either open or closed, then $A \cup \{\infty\}$ is connected if and only if every connected component of $A$ is unbounded. However, in general this is not the case. Indeed, consider a set $A \subset \mathbb{C}$ that consists of a sequence of horizontal lines $\{z \colon \operatorname{Im} z = \varepsilon_n\}$, with $\varepsilon_n \to 0$, and the point $\{0\}$. Then $\{0\}$ is a connected component of $A$, but $A \cup \{\infty\}$ is connected. It is even possible that a set $A$ is totally disconnected, but $A \cup \{\infty\}$ is connected. Indeed, let $0 < a < 1/e$ and $f_a(z) = ae^z$ as in §2.3, and let $A$ be the set of the endpoints of the arcs that form the Cantor bouquet $J(f_a)$. Mayer [183] showed that $A$ is totally disconnected, but $A \cup \{\infty\}$ is connected.

By Theorem 7.24, for a transcendental entire function $f$ every connected component of $I(f)$ that intersects $A(f)$ is unbounded, but from Theorem 7.26 we cannot conclude anything about the unboundedness of other connected components of $I(f)$.

Theorem 7.24, Corollary 7.25 and Theorem 7.26 establish properties similar to, but weaker than, Eremenko's property for all transcendental entire functions. There are also several known sub-classes of functions where Eremenko's property itself is known to hold. We already encountered one such class in Theorem 7.4: all functions $f \in \mathcal{B}$ of finite order, and compositions of such. Observe that both hypotheses – that $f \in \mathcal{B}$ and that $\rho(f) < \infty$ – are *non-dynamical*: they concern the mapping behaviour of the function $f$, but not of its iterates. In contrast, the following result by Rempe [219] applies to functions $f \in \mathcal{B}$ of arbitrarily fast growth, but imposes a dynamical restriction on these. We say that a transcendental entire function $f$ is *postsingularly bounded* if the postsingular set $P(f)$ defined in (3.1) is bounded.

**Theorem 7.27.** *Let $f$ be a postsingularly bounded transcendental entire function. Then every connected component of $I(f)$ is unbounded.*

The assumption implies, in particular, that $f \in \mathcal{B}$. The class of postsingularly bounded functions includes many important cases; in particular that of hyperbolic functions $f \in \mathcal{B}$. Recall that the example from Theorem 7.5 is of disjoint type and thus, in particular, hyperbolic. So it is an example of a function having the Eremenko property, but not the strong Eremenko property.

The proof of Theorem 7.27 associates to each escaping point of the postsingularly bounded function $f$ an *external address*, which is analogous to the sequence $S_m^\sigma$ of half-strips considered in the proof of Theorem 2.1. Following Benini and Rempe [42], this address can be described as follows. Let $R > 0$ be so large that $P(f)$ is contained in the disc $D := D(0, R)$, and consider the set $W := \mathbb{C} \setminus \overline{D}$. We say that two escaping points *have the same address* if, for all sufficiently large $n$, they belong to the same connected component of $f^{-n}(W)$.

Let $z_0 \in I(f)$, and let $\Gamma = \Gamma(z_0) \subset I(f)$ be the set of points having the same address as $z_0$. It is shown in [219] that $\Gamma$ is the closure of an increasing union of unbounded closed and connected subsets on which the iterates of $f$ tend to infinity uniformly. In particular, $\Gamma$ itself is unbounded, establishing Eremenko's property.

As already discussed, one motivation for studying curves to infinity in the escaping sets is that, on the one hand, these curves can be indexed using combinatorial information, and furthermore often several of these curves end at common non-escaping points (often, repelling periodic points and their iterated preimages). These facts allow one to study the dynamics from a combinatorial point of view. As we have just seen, for postsingularly bounded entire functions, even though there may be no curves to infinity, the escaping set is still organised in unbounded connected sets indexed by "addresses" (also called "itineraries"). These sets are called *filaments*, and it is natural to ask whether they can serve the same function as hairs in cases where the latter do not exist. This is indeed the case. For example, in [42], it is shown that the *Douady–Hubbard landing theorem* for external rays of polynomials, a cornerstone of polynomial dynamics, has an analogue for postsingularly bounded transcendental entire functions and their filaments.

## 7.4 Spiders' webs

We have seen in §2.1 that for the example considered there the escaping set has uncountably many connected components. On the other hand, Fatou's function considered in §2.2 is an example of a function for which the escaping set is connected. Indeed, we saw that for the function $f$ given by (2.3) the Fatou set $F(f)$ consists of a single connected component $U$, which is a Baker domain. So $U \subset I(f) \subset \overline{U} = \mathbb{C}$ and we conclude that $I(f)$ is connected. Thus Eremenko's conjecture is trivially true for this function.

The preceding argument applies to any entire function $f$ for which $F(f)$ consists of a single completely invariant Baker domain, but this is a rather exceptional occurrence. Can the escaping set be connected in other cases? Let us postpone a discussion of this question for the exponential family considered in §2.3 until the next section, and first consider Baker's example discussed in §2.4. Here the escaping set contains a sequence $(A_n)$ of nested annuli tending to infinity, which are contained in the fast escaping set $A(f)$ by Theorem 6.10. Every connected component of $A(f)$ is unbounded by Theorem 7.24, and hence contains $A_n$ for sufficiently large $n$. It follows that $A(f)$ is connected.

We will show that $I(f)$ is also connected. To this end, let $X$ be the connected component of $I(f)$ containing $A(f)$. We want to show that $X = I(f)$. Since $J(f) \subset \overline{A(f)}$

by Theorem 5.6, we have

$$A(f) \cup (I(f) \cap J(f)) \subset \overline{A(f)}.$$

Thus $A(f) \cup (I(f) \cap J(f))$ is connected and hence $I(f) \cap J(f) \subset X$. Let now $U$ be a connected component of $(I(f) \setminus A(f)) \cap F(f)$. By Theorem 6.10, $U$ is a simply connected escaping wandering domain. For every $m$ and sufficiently large $n$, the set $f^n(\overline{U})$ belongs to the unbounded connected component of $\mathbb{C} \setminus A_m$. Hence $f^n|_{\overline{U}} \to \infty$ uniformly, and $\overline{U} \subset I(f)$; in particular, $\overline{U} \cap X \neq \emptyset$. We conclude that $U \subset X$. Overall we see that $X = I(f)$ is connected. (That $\overline{U} \cap X \neq \emptyset$ also follows from Theorem 6.39.)

This argument is due to Rippon and Stallard [239, Theorem 2], and applies to any transcendental entire function with a multiply connected Fatou component.

**Theorem 7.28.** *Suppose that $f$ is a transcendental entire function with a multiply connected wandering domain. Then $A(f)$ and $I(f)$ are connected.*

Observe that the key property in the above argument is that $A(f)$ contains a sequence of continua (here, essential curves in the annuli $A_n$), tending to infinity and each of which separates all preceding ones from infinity. By Theorem 7.24, each of these continua is connected to the next by a connected subset of $A(f)$, reminiscent of the structure that one encounters in physical spiders' webs. Rippon and Stallard [244] introduced the following definition.

**Definition 7.29.** A set $E \subset \mathbb{C}$ is a *spider's web* if $E$ is connected and there exists a sequence $(G_n)_{n=0}^\infty$ of bounded simply connected domains such that $\partial G_n \subset E$ and $G_n \subset G_{n+1}$ for all $n \geq 0$, and such that $\bigcup_{n=0}^\infty G_n = \mathbb{C}$.

The proof of Theorem 7.28 shows that $I(f)$, $A(f)$ and even the "levels" $A_R(f)$ defined by (5.10) and (5.12) are spiders' webs, with $R$ an admissible radius of course. (Everything we say below concerning these levels holds more generally for $A_R^L(f)$, where $L \in \mathbb{Z}$, but we restrict to the case $L = 0$ for simplicity.)

It was observed by Evdoridou [121, Theorem 1.5] that a connected subset $E$ of $\mathbb{C}$ is a spider's web if and only if, for every $z \in \mathbb{C}$, there exists a compact subset of $E$ that separates $z$ from $\infty$. The notion of a spider's web is purely topological and makes no reference to dynamics, but we are usually interested in the case where $E$ is a dynamically defined subset, such as $E = I(f)$, $E = A(f)$ or $E = A_R(f)$. In these cases, the condition from [121, Theorem 1.5] can be further weakened, as noted by Rippon and Stallard for $I(f)$ [240, Theorem 2] and for $A_R(f)$ [244, Lemma 7.1(b)] and by Sixsmith [269, Theorem 1.5] for $A(f)$. The following is a reformulation of this fact and shows that, for the sets in question, a set $E$ being a spider's web is in fact a topological property of its complement $\widehat{\mathbb{C}} \setminus E$ in the Riemann sphere.

**Theorem 7.30.** *Let $f$ be a transcendental entire function, and let $E$ be one of the sets $I(f)$, $A(f)$ or $A_R(f)$, where $R$ is an admissible radius. Then $E$ is a spider's web if and only if $\widehat{\mathbb{C}} \setminus E$ is disconnected.*

*Proof.* First suppose that $E$ is a spider's web, and let $G_n$ be a corresponding sequence of simply connected domains. Since $E \subset I(f) \neq \mathbb{C}$, the set $\mathbb{C} \setminus E$ is non-empty. If $z \in \mathbb{C} \setminus E$, then $z$ is contained in $G_n$ for sufficiently large $n$, and hence separated from $\infty$ by $\partial G_n \subset E$. In particular, $\widehat{\mathbb{C}} \setminus E$ is disconnected.

Now suppose that $\widehat{\mathbb{C}} \setminus E$ is disconnected. Then there is an open subset $U$ of $\widehat{\mathbb{C}}$ such that $\partial U \subset E$ (where we understand the boundary to be taken in $\widehat{\mathbb{C}}$) and $\widehat{\mathbb{C}} \setminus E$ intersects both $U$ and $\widehat{\mathbb{C}} \setminus U$. Without loss of generality, we may suppose that $\infty \in \widehat{\mathbb{C}} \setminus U$ (otherwise, replace $U$ by the interior of its complement). Then $U$ is a bounded domain. Replacing $U$ by its topological hull $T(U)$, we may suppose that $U$ is simply connected.

We claim that $U$ must intersect $J(f)$. Indeed, otherwise, $U$ is a Fatou component of $F$. If $E = I(f)$ or $E = A(f)$, then $U \subset E$ by Theorem 6.45, which contradicts the fact that $U$ intersects $\mathbb{C} \setminus E$. So suppose that $E = A_R(f)$. Since $\partial U \subset A_R(f)$ but $U \not\subset A_R(f)$, the minimum modulus theorem implies that there is $n \geq 0$ with $0 \in f^n(U)$. As $U$ is bounded and $E$ is forward invariant, we have $\partial f^n(U) \subset f^n(\partial U) \subset E$. So $f^n(U)$ contains the disc $D(0, R)$. This disc intersects the Julia set, and hence $U$ does also, as claimed.

Similarly as in Domínguez's proof of Theorem 4.1 in §4.4, define $G_0 := U$ and inductively $G_n := T(f(G_{n-1}))$. As in the proof of Theorem 4.14, the blowing-up property of the Julia set (Theorem 3.1 (f)) implies that there is a subsequence $(G_{n_k})$ such that $G_{n_k} \subset G_{n_{k+1}}$ and $\bigcup_{k=0}^\infty G_{n_k} = \mathbb{C}$. Hence $E$ satisfies the second part of the definition of a spider's web.

It remains to show that $E$ is connected. If $E = A(f)$ or $E = A_R(f)$, then all connected components of $E$ are unbounded by Theorem 7.24, and hence each such connected component contains $\partial G_{n_k}$ for all sufficiently large $k$. It follows that $E$ is connected.

If $E = I(f)$, then, for the same reason, $I(f)$ has a connected component $C$ that includes $A(f)$, all $\partial G_{n_k}$, and all unbounded connected components of $F(f) \cap I(f)$ for sufficiently large $k$. Since $J(f) \subset \overline{A(f)}$, the component $C$ also contains all points of $J(f) \cap I(f)$. Finally, if $U$ is a bounded connected component of $F(f) \cap I(f)$, then $U$ is a wandering domain, and $\partial U \cap I(f) \neq \emptyset$ by Theorem 6.39. Since $\partial U \cap I(f) \subset J(f) \cap I(f) \subset C$, we find that $U \subset C$. It follows that $I(f) = C$, as required. $\square$

**Corollary 7.31.** *If $A_R(f)$ is a spider's web for some admissible radius $R$, then so is $A(f)$, and if $A(f)$ is a spider's web, then so is $I(f)$.*

*Proof.* For the sets in question, the first condition in the definition of a spider's web (that $E$ be connected) follows from the second by Theorem 7.30. Clearly if the second condition holds for $E$, it also holds for any larger set. $\square$

We also note that the condition that $A_R(f)$ is a spider's web does not depend on the choice of $R$; see [244, Lemma 7.1].

**Theorem 7.32.** *Let $f$ be a transcendental entire function, and suppose that there is an admissible radius $R$ such that $\mathbb{C} \setminus A_R(f)$ has a bounded complementary component. Then $A_{\tilde{R}}(f)$ is a spider's web for every admissible radius $\tilde{R}$.*

*Proof.* Suppose that $A_R(f)$ has a bounded complementary component $U$, and let $(G_n)$ be the sequence of domains obtained in the proof of Theorem 7.30. Let $\tilde{R}$ be an admissible radius and let $n_0$ be such that $M^n(R,f) \geq \tilde{R}$ for $n \geq n_0$. Then $\partial G_n \subset f^n(\partial G_0) \subset A_{\tilde{R}}(f)$ for $n \geq n_0$. By Theorem 7.30, $A_{\tilde{R}}(f)$ is a spider's web, as claimed. $\square$

The property that $A_R(f)$ is a spider's web can be characterised as follows (see [244, Theorem 8.1]).

**Theorem 7.33.** *Let $f$ be a transcendental entire function and let $R$ be an admissible radius for $f$. Then $A_R(f)$ is a spider's web if and only if the following holds: There is a simply connected domain $D$ with $D \cap J(f) \neq \emptyset$ and a sequence $(G_n)_{n=0}^\infty$ of bounded simply connected domains such that, for all $n \geq 0$,*

$$f^n(D) \subset G_n \tag{7.3}$$

*and $f(\partial G_n) \cap G_{n+1} = \emptyset$.*

*Proof.* If $A_R(f)$ is a spider's web, let $D = D(0, R)$, and for each $n \geq 0$, let $G_n$ be the connected component of the complement of $A_{M^n(R,f)}(f) = f^n(A(R,f))$ that contains $f^n(D)$. Then the domains $G_n$ are simply connected by Theorem 7.24 and bounded by Theorem 7.32. By definition, $\partial G_n \subset f^n(A(R,f))$, and hence $f(\partial G_n) \cap G_{n+1} = \emptyset$.

Now suppose that domains $(G_n)$ with the stated properties exist. For $n, k \geq 0$, let $G_n^k$ denote the connected component of $f^{-k}(G_{n+k})$ containing $f^n(D)$. Then $G_n^k$ is simply connected. We claim that $G_n^{k+1} \subset G_n^k$, and $f^j \colon G_n^k \to G_{n+j}^{k-j}$ is a proper map for $j \leq k+1$. Indeed, for $k = 0$ this follows from the fact that $G_n^1$ contains $f^n(D)$ and is disjoint from $\partial G_n$. For $k > 0$, the claim follows inductively by applying the case $k = 0$ to the domains $(G_n^{k-1})_{n=0}^\infty$, which also satisfy the hypotheses of the theorem.

Now consider the intersection $G_n^\infty := \bigcap_{k=0}^\infty G_n^k$. This is a simply connected domain containing $f^n(D)$ for all $n$, and $f \colon G_n^\infty \to G_{n+1}^\infty$ is a proper map. In particular, $f^n(\partial G_0) = \partial G_n$ does not intersect $T(f^n(D))$. For fixed $R$ and sufficiently large $n$, we thus have $\partial G_n \subset A_R(f)$, and it follows that $A_R(f)$ is a spider's web. $\square$

The following result [244, Theorem 1.8] is a simple consequence of Theorem 7.33.

**Corollary 7.34.** *Let $f$ be a transcendental entire function and let $R$ be an admissible radius for $f$. If $f$ is bounded on a connected unbounded set, then $A_R(f)$ is not a spider's web.*

*In particular, $A_R(f)$ is not a spider's web if $f \in \mathcal{B}$.*

*Proof.* Suppose that $A_R(f)$ is a spider's web and let $D$ and $G_n$ be as in Theorem 7.33. It follows from (7.3) and the blowing-up property of the Julia set that $\min_{z \in \partial G_n} |z| \to \infty$ as $n \to \infty$. Since $f(\partial G_n) \cap G_{n+1} = \emptyset$ we have $\min_{z \in \partial G_n} |f(z)| \geq \min_{z \in \partial G_{n+1}} |z|$. Since a connected unbounded set must intersect $\partial G_n$ for all large $n$, we find that $f$ cannot be bounded on such a set.

The second claim follows since if $f \in \mathcal{B}$, say $\text{sing}(f^{-1}) \subset D(0, r)$, then $\{z \colon |f(z)| = r\}$ contains a curve tending to infinity. $\square$

There are large classes of functions for which it is known that $A_R(f)$, and hence by Corollary 7.31 also $A(f)$ and $I(f)$, are spiders' webs. This structure is particularly prevalent for functions of order less than $1/2$. *Baker's conjecture*, named for a question asked by Baker [14, p. 484] in 1981, states that such functions (and also functions of order $1/2$ that have *minimal type*, a condition we do not define here) do not have unbounded Fatou components. Any unbounded Fatou component of such a function would necessarily be an unbounded simply connected wandering domain; see Hinkkanen's survey [150, Theorem 6.2]. The following result by Rippon and Stallard [244, Theorem 1.5(b)] makes a connection between the structure of $A_R(f)$ and this conjecture.

**Theorem 7.35.** *Let $f$ be a transcendental entire function and suppose that $A_R(f)$ is a spider's web for some admissible radius $R$. Then all Fatou components of $f$ are bounded.*

*Proof.* By Theorem 6.11, all multiply connected Fatou components of an entire function are bounded. If $U$ was an unbounded simply connected Fatou component, then it would follow from Theorem 7.32 that $U$ intersects $A_R(f)$ for all $R$. However, by Theorem 6.6, any Fatou component that intersects $A_R(f)$ is entirely contained in $A_R(f)$. This is impossible since $\bigcap_R A_R(f) = \emptyset$. □

Baker's conjecture remains open, but there is a long history of results that establish it for functions of order less than $1/2$ whose growth, additionally, satisfies certain regularity conditions; see [150]. As observed by Rippon and Stallard [244], these proofs implicitly prove that $A_R(f)$ is a spider's web, although this terminology had not been established at the time. From Theorem 7.33, we obtain the following condition [244, Corollary 8.2], which involves the *minimum modulus*

$$m(r, f) := \min_{|z|=r} |f(z)|.$$

**Corollary 7.36.** *Let $f$ be a transcendental entire function and let $R$ be an admissible radius for $f$. Suppose that there is a sequence $(\rho_n)_{n=0}^\infty$ such that*

$$m(\rho_n, f) \geq \rho_{n+1}$$

*and*

$$\rho_n \geq M^n(R, f)$$

*for all $n \geq 0$. Then $A_R(f)$ is a spider's web.*

*Proof.* Set $G_n := D(0, \rho_n)$ and $D := D(0, R)$, and apply Theorem 7.33. □

Observe that the hypotheses of Corollary 7.36 are satisfied, in particular, if there is $\rho > 0$ such that the iterated minimum modulus $m^n(\rho, f)$ tends to infinity faster than $M^n(R, f)$. Osborne, Rippon and Stallard [206, Theorem 1.3] showed that this result is equivalent to the hypotheses of Corollary 7.36. The papers [201, 206] also study, more generally, functions for which $m^n(\rho, f) \to \infty$ for some $\rho$, and the set of points that tend to infinity at least as fast as this sequence.

The condition in Corollary 7.36 has been shown to hold for many classes of transcendental entire functions (not only of order less than 1/2). We refer to [244, Theorem 1.9] for a list of several cases in which Corollary 7.36 applies. Here we only note that this list includes certain gap series and random entire functions in the sense of Littlewood and Offord [173].

To give a flavour of the type of condition for which it can be shown that $A_R(f)$ is a spider's web, we mention the following result [244, Theorem 1.9(c)]. (As noted in [244], this result can also be deduced from the authors' previous results in [240].)

**Theorem 7.37.** *Suppose that $f$ is a transcendental entire function with $\rho(f) < 1/2$. Assume additionally that there is $m > 1$, an admissible radius $R$ and a sequence $(r_n)_{n=0}^\infty$ of positive real numbers such that*

$$r_n \geq M^n(R) \quad and \quad M(r_n) \geq r_{n+1}^m, \tag{7.4}$$

*for $n \geq 0$. Then $A_R(f)$ is a spider's web.*

*Proof.* If $f$ has order less than 1/2, then by a version of the "$\cos \pi\rho$-theorem," for all sufficiently large $r$, there is some $\rho < r^m$ such that $m(\rho, f) \geq M(r, f)$, where $m(\rho, f)$ is again the minimum modulus. (See [36] and also [10, Satz 1].) If $(r_n)_{n=0}^\infty$ is a sequence as in (7.4), where all $r_n$ are sufficiently large, we thus find $\rho_n$ such that $\rho_n \leq r_n^m$ and

$$m(\rho_n, f) \geq M(r_n, f) \geq r_{n+1}^m \geq \rho_{n+1}.$$

By assumption, the sequence $(r_n)_{n=0}^\infty$ grows at least as fast as the maximum modulus. □

As an application of Theorem 7.37 we mention the following result.

**Corollary 7.38.** *Let $f$ be a transcendental entire function and let $R$ be an admissible radius for $f$. Suppose that $f$ has order less than 1/2 and positive lower order. Then $A_R(f)$ is a spider's web.*

Indeed, if $f$ satisfies the hypotheses of this corollary, then (7.4) is satisfied (see [244, §8] and [262, Lemma 5.3]).

In particular, for such $f$ all Fatou components are bounded (that is, Baker's conjecture holds). This was shown already by Anderson and Hinkkanen; see [6, Theorem 2] and the remark following [150, Theorem 8.3].

As an explicit example where Corollary 7.38 applies we mention the function (5.14) considered in §5.1. For this function the order and lower order are equal to 1/4.

Let us also mention a rather explicit class of examples for which it is known that $A_R(f)$ is a spider's web. Let $p$ be a polynomial with a repelling fixed point $z_0$. It is well-known [192, Theorem 8.2] that such fixed points are *linearisable*; that is, there is a conformal map $f$ defined near 0 such that $f(0) = z_0$ and

$$f(\lambda z) = p(f(z)), \tag{7.5}$$

where $\lambda = p'(z_0) \in \mathbb{C} \setminus \overline{\mathbb{D}}$ is the multiplier of $z_0$. Using the functional relation (7.5), $f$ extends to a transcendental entire function $f\colon \mathbb{C} \to \mathbb{C}$, called the *Poincaré function* of $p$ at $z_0$; see [192, Corollary 8.12]. Mihaljević-Brandt and Peter [189] showed the following.

**Theorem 7.39.** *Let $f$ be the Poincaré function of a polynomial $p$ at a repelling fixed point $z_0$, and let $R$ be an admissible radius for $f$. Then $A_R(f)$ is a spider's web if and only if $z_0$ is a singleton component of the Julia set $J(p)$.*

*Sketch of proof.* Let $C$ be the connected component of $f^{-1}(J(p))$ containing 0. This component is invariant under the map $z \mapsto \lambda z$. If the connected component of $J(p)$ containing $z_0$ is non-trivial, then the same is true for $C$, and therefore $C$ is unbounded. Since $f(C) \subset J(p)$ and $J(p)$ is compact, the map $f$ is bounded on $C$. It now follows from Corollary 7.34 that $A_R(f)$ is not a spider's web.

On the other hand, suppose that $\{z_0\}$ is a connected component of $J(p)$. We use *equipotential lines* of the Green's function associated to the polynomial $p$ (see [192, Definition 9.6]). There exists a small equipotential line $\gamma$ that is a Jordan curve and separates $z_0$ from $\infty$. Let $U$ be the interior of $\gamma$. Then $p^k(U)$ is bounded by some other equipotential curve, for all $k \geq 0$, and the map $p^k\colon U \to p^k(U)$ is proper. For sufficiently large $k$, the curve $p^k(\gamma) = \partial p^k(U)$ is approximately a large round circle whose radius grows like $c^{d^k}$ for some $c > 1$, where $d$ is the degree of $p$. (Recall that $p$ is conformally conjugate to $z \mapsto z^d$ near $\infty$, and equipotential lines of $p$ correspond to round circles under this conjugacy.)

On the other hand, if $\gamma$ is small enough, then there is a Jordan domain $V$ containing 0 that is mapped homeomorphically to $U$ by $f$. We define a sequence $(G_n)_{n=n_0}^\infty$ as $G_n := \lambda^{k_n} \cdot V$, where $(k_n)_{n=0}^\infty$ is a sequence of nonnegative integers defined inductively as follows. Let $k_0$ be sufficiently large to ensure that $\lambda^k \cdot V \subset p^k(U)$ for $k \geq k_0$.. If $k_n$ and $G_n$ have been chosen, choose $k_{n+1} \geq k_n$ maximal such that $G_{n+1} \subset p^{k_n}(U) = f(G_n)$. Then $f(\partial G_n) = \partial f(G_n) \cap G_{n+1} = \emptyset$. Moreover, $k_{n+1}$ is exponential in $k_n$. It is known that the function $f$ has finite order $\log d / \log|\lambda|$. It can be deduced that, for any given $R > 0$, if $k_0$ was sufficiently large, the domain $G_n$ contains $D(0, M^n(R, f))$ for all $n \geq 0$. The claim then follows from Theorem 7.33. We refer to [189] for the details of the argument. □

Figure 15 shows the Fatou and Julia sets of the Poincaré functions of $q_\lambda(z) := \lambda z + z^2$ at the repelling fixed point 0, for the parameter values $\lambda = 2$ (left), $\lambda = 2 + i$ (middle) and $\lambda = 2 + 1.4i$ (right). The Poincaré functions $f_\lambda$ are normalised to satisfy $f'_\lambda(0) = 1$. Thus 0 is a parabolic fixed point of $f_\lambda$. The parabolic basins of 0 are shown in light grey in the pictures. Note that the parabolic basins do not intersect $I(f_\lambda)$. In particular, they do not intersect $A_R(f_\lambda)$. The ranges shown are $-3 \leq \operatorname{Re} z \leq 6$ and $|\operatorname{Im} z| \leq 9$ in the left picture, and $|\operatorname{Re} z| \leq 50$ and $|\operatorname{Im} z| \leq 50$ in the two other pictures. For $\lambda = 2$ and $\lambda = 2 + i$ the Julia set of $q_\lambda$ is connected and thus $A_R(f_\lambda)$ is not a spider's web, while for $\lambda = 2 + 1.4i$ the Julia set of $q_\lambda$ is a Cantor set and thus $A_R(f_\lambda)$ is a spider's web. Note that $q_\lambda$ is conjugate to $p(z) := z^2 + c$ for $c = \lambda/2 - \lambda^2/4$. For $\lambda = 2$ and $\lambda = 2 + i$ this value of $c$ is in the Mandelbrot set, while for $\lambda = 2 + 1.4i$ it is not.

Sixsmith [262] gives further classes of functions for which $A_R(f)$ is a spider's web, and exhibits several examples. In particular, he shows [262, Example 4] that for $f(z) := \cos z + \cosh z$, the set $A_R(f)$ is a spider's web.

In view of Baker's conjecture, it is natural to ask whether every $A_R(f)$ is a spider's web for all functions of order less than $1/2$, but this is not the case. Indeed, Rippon and

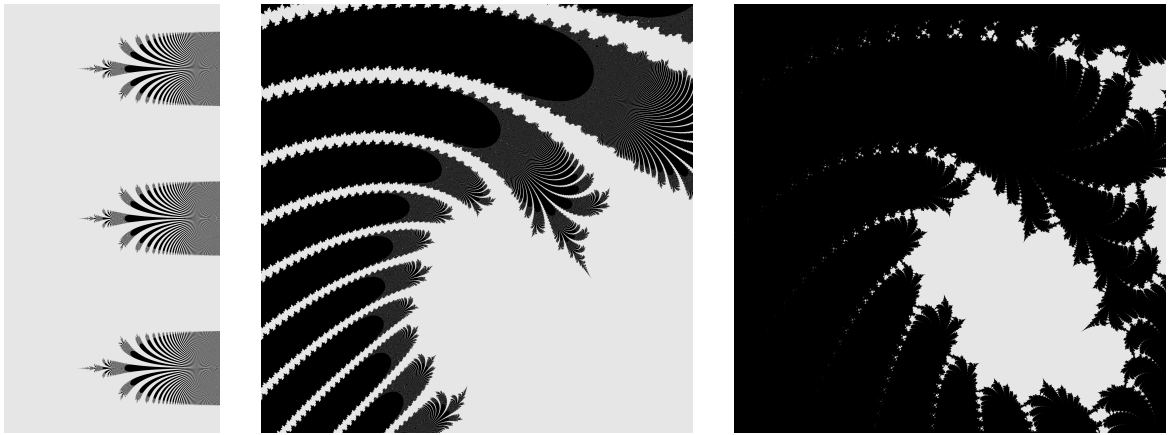

Figure 15: Julia sets of some Poincaré functions.

Stallard [246] gave a growth condition on an entire function $f$ that implies that $A_R(f)$ is a spider's web, and showed that this condition is optimal in a certain sense.

**Theorem 7.40.** *Let $f$ be a transcendental entire function and let $R$ be an admissible radius for $f$. Put, for $n \in \mathbb{N}$,*

$$R_n := M^n(R, f) \quad and \quad \varepsilon_n := \max_{R_n \leq r \leq R_{n+1}} \frac{\log \log M(r)}{\log r}.$$

*If*

$$\sum_{n=1}^{\infty} \varepsilon_n < \infty,$$

*then $A_R(f)$ is a spider's web.*

*On the other hand, for every $c > 1$, there is an entire function $f$ such that the above sequence satisfies*

$$\sum_{n=1}^{\infty} \varepsilon_n^c < \infty, \tag{7.6}$$

*but $A(f)$ is not a spider's web.*

The first statement is [246, Theorem 1.1]. The second is a consequence of [246, Theorem 1.2], which gives a stronger but more technical result.

Observe that any function satisfying (7.6) has order 0, so Theorem 7.40 shows that there are functions of order 0 for which $A_R(f)$ and $A(f)$ are *not* spiders' webs.

Although $A(f)$ is not a spider's web for the examples from Theorem 7.40, the *escaping set* $I(f)$ is still a spider's web, and $f$ has no unbounded Fatou components. Indeed, the function is real with only negative roots, and for such functions of order less than $1/2$ these facts were proved in [245]. Later, this was extended to functions with real zeros by Nicks, Rippon and Stallard; see [199, Corollary 1.2] and [200, Theorems 1.1 and 1.2].

**Theorem 7.41.** *Suppose that $f$ is a real transcendental entire function of order less than $1/2$ with only real roots. Then $f$ has no unbounded Fatou components and $I(f)$ is a spider's web.*

The proof that $I(f)$ is a spider's web considers a certain closed subset of $I(f)$ on which points escape with a certain rate (i.e., a set defined similarly to $A_R(f)$, but with a smaller minimum rate of growth). Similarly as in Theorem 7.30, it is enough to show that the complement of this set has no unbounded connected components. The authors thus assume that this complement contains an unbounded curve and reach a contradiction by carefully considering the way that iterated forward images of this curve must either stretch forward to infinity or wind around the origin. We will not attempt to give further details here and refer to the original articles for the proofs, and for more general classes of functions to which the theorem applies.

In view of Baker's conjecture and the results above, it is natural to ask the following question.

**Question 7.42.** *Let $f$ be a transcendental entire function of order $\rho(f) < 1/2$.*

(a) *Is $I(f)$ necessarily connected?*

(b) *Is $I(f)$ necessarily a spider's web?*

(c) *Is there necessarily a bounded simply connected domain $U$ with $U \cap J(f) \neq \emptyset$ and such that $f^n|_{\partial U} \to \infty$ uniformly?*

*Remark* 7.43. For any transcendental entire function, condition (c) implies (with the same proof as Theorem 7.30) that $I(f)$ contains a spider's web on which the iterates tend to infinity uniformly. This condition does not appear to have been studied explicitly, but it is satisfied whenever $A_R(f)$ is a spider's web, as well as by the functions in Theorem 7.41. Uniformly escaping structures tend to be particularly useful for the further study of the functions under consideration, so the class of functions satisfying (c) warrants further study.

*Remark* 7.44. Nicks, Rippon and Stallard [200, Conjecture 3.3] conjecture that the answer to (b), and hence (a) is positive, even with $I(f)$ replaced by the quite fast escaping set $Q(f)$ defined in (5.22).

We have encountered two prevalent types of structures in the escaping set: "Cantor bouquets" of curves, which occur for many (though not all) functions in the class $\mathcal{B}$, and spiders' webs. As already mentioned in the introduction to this section, these two are not mutually exclusive. Indeed, let us return to Fatou's function studied in §2.2. Let $E$ be the set of all "endpoints" of the arcs that make up the Julia set $J(f)$, as described in Theorem 2.4. It follows easily from the description of $J(f)$ in §2.2 that $E$ is totally disconnected. Evdoridou [118, Theorem 5.1] obtained the following analogue of the result of Mayer [183] for the exponential family that was mentioned after Theorem 7.26.

**Theorem 7.45.** *Let $f$ be Fatou's function and let $E$ be the set of endpoints of $J(f)$. Then $E \cup \{\infty\}$ is connected.*

As we saw in §2.2, the Fatou set $F(f)$ of Fatou's function consists of a single Baker domain and is thus contained in $I(f)$. Except possibly for their endpoints, the arcs that make up $J(f)$ are also contained in $I(f)$. Thus $\mathbb{C} \setminus I(f) \subset E$. Evdoridou [118, Theorem 5.2] showed that the conclusion of Theorem 7.45 changes if we consider only the non-escaping endpoints.

**Theorem 7.46.** *Let $f$ be Fatou's function and let $E$ be the set of endpoints of $J(f)$. Then $E \setminus I(f) \cup \{\infty\} = \widehat{\mathbb{C}} \setminus I(f)$ is totally disconnected.*

Theorem 7.30 yields the following corollary [118, Theorem 1.1].

**Theorem 7.47.** *Let $f(z) = e^{-z} + z + 1$ be Fatou's function. Then $I(f)$ is a spider's web.*

The proof of the theorem shows, in fact, that $I(f)$ contains a spider's web on which the iterates tend to infinity uniformly; that is, $f$ belongs to the class of functions discussed in Remark 7.43. Evdoridou and Rempe [121, Theorem 5.1] showed that $I(f)$ can be replaced by $A(f) \cup F(f)$ in Theorem 7.47.

Further examples of entire functions where both structures coexist are given by Dourekas [111]. Cantor bouquets are most strongly associated to functions in the Eremenko–Lyubich class $\mathcal{B}$, but neither Fatou's function nor Dourekas's examples (the exponential sums mentioned previously in §7.1) belong to this class. It is natural to ask whether spiders' webs can also occur in the Eremenko–Lyubich class. By Corollary 7.34, $A_R(f)$ is never a spider's web when $f \in \mathcal{B}$. Furthermore, Sixsmith [269, Theorem 1.4] proved the following.

**Theorem 7.48.** *Let $f \in \mathcal{B}$ and suppose that $f^{-1}$ has a logarithmic singularity over some finite asymptotic value. Then $I(f)$ is not a spider's web.*

Sixsmith [269] conjectured that $I(f)$ is never a spider's web when $f \in \mathcal{B}$, but this was disproved by Rempe [228].

**Theorem 7.49.** *For $f = \cosh$, the sets $A(f)$ and $I(f)$ are spiders' webs.*

The proof makes use of a model due to Pardo-Simón [209, 210], which describes the dynamics of $f(z) = \cosh(z)$ in terms of the function $g(z) = \cosh(z)/2$; see also §8.8. (The function $g$ behaves essentially in the same way as the function $\sin(z)/2$, studied in §2.1.) In particular, there is a continuous function from $\widehat{\mathbb{C}} \setminus A(f)$ to the set $(J(g) \setminus A(g)) \setminus \{\infty\}$ of "meandering endpoints" of $g$. Generalising the results of [121] mentioned above, it was shown by Evdoridou and Sixsmith [122, Theorem 1.5] that $A(g) \cup F(g)$ is a spider's web, and in fact that the set $\widehat{\mathbb{C}} \setminus (A(g) \cup F(g)) = (J(g) \setminus A(g)) \cup \{\infty\}$ is totally disconnected. In particular, $\widehat{\mathbb{C}} \setminus A(f)$ is disconnected, and $A(f)$ is a spider's web.

Theorem 7.49 also gives the first example of a function for which $A(f)$ is a spider's web, but $A_R(f)$ is not, answering a question of Rippon and Stallard [244, Question 2]. (Note that $A_R(\cosh)$ is not a spider's web by Corollary 7.34.) As far as we are aware, it is also the first known example where $I(f)$ is a spider's web, but $I(f)$ contains no spider's web on which the iterates tend to infinity uniformly, as discussed in Remark 7.43

To conclude our discussion of examples where spiders' webs and hairs coexist, we observe that there are also transcendental entire functions where neither occurs.

*Remark* 7.50. The function $f$ from Theorem 7.5 is of disjoint type, and hence $F(f)$ is connected and consists of a single attracting basin. So this is an example of a function for which $I(f)$ contains neither a spider's web nor a curve to infinity.

Recall the distinction between Eremenko's property and the *strong* Eremenko property, which asks about the existence of *curves* to infinity. Similarly, when discussing spiders' webs, let us say that a spider's web is *Jordan* if the domains $G_n$ in Definition 7.29 may be chosen to be Jordan. If $f$ has a multiply connected Fatou component, then this is clearly true. We are not aware of any other known example of a Jordan spider's web escaping set. However, it is plausible that this property holds at least for sufficiently well-controlled examples.

**Question 7.51.** Is $I(f)$ a Jordan spider's web when

(a) $f$ is Fatou's function (2.3)?

(b) $f = \cosh$?

(c) $f = \cos + \cosh$?

(d) $f$ is a Poincaré function with spider's web escaping set, as in Theorem 7.39?

Part (b) of Question 7.51 was also asked in [228, Question 4.1].

## 7.5 Connected components of the escaping set

In the previous section, we encountered many functions for which $I(f)$ is a spider's web, and hence connected. It is natural to ask whether there are cases where $I(f)$ is connected, but not a spider's web. The first such example was described by Rempe [222].

**Theorem 7.52.** *Let $f(z) = ae^z$, where $a > 1/e$. Then $I(f)$ is connected, but not a spider's web.*

*Sketch of proof.* The proof that $I(f)$ is connected uses a construction of Devaney [105]. Let $X$ denote the set of points $z \in \mathbb{C}$ with $\operatorname{Im} f^n(z) \in [-\pi, \pi]$ for all $n \geq 0$ and $f^n(z) \in \mathbb{R}$ for sufficiently large $n$. The argument given in [105, p. 631] implies that $X$ is a connected subset of $\mathbb{C}$. Observe that $X$ contains the horizontal lines at imaginary parts $-\pi$, $0$ and $\pi$. The union of the $2\pi i \mathbb{Z}$-translates of $X$ thus forms a connected subset $X_0$ that is invariant under translation by $2\pi i$. For the same reason, the sets $X^+ := \bigcup_{j=1}^\infty X + 2\pi i j$ and $X^- := \bigcup_{j=1}^\infty X - 2\pi i j$ are connected; observe that $X_0 = X \cup X^+ \cup X^-$. Also observe that $X_0 = f^{-1}(X)$.

Now consider $X_1 := f^{-1}(X_0) = X_0 \cup f^{-1}(X^+) \cup f^{-1}(X^-)$. The function $f$ maps each strip $\{z \colon (2k-1)\pi < \operatorname{Im} z < (2k+1)\pi\}$ conformally to $\mathbb{C} \setminus (-\infty, 0]$. Hence every such strip contains one connected component of $f^{-1}(X^+)$ and one connected component of $f^{-1}(X^-)$, each containing a connected component of $f^{-1}(\{z\colon |\operatorname{Im} z| = \pi\}) \subset X_0$. We conclude that $X_1$ is connected.

Applying the same argument inductively, it follows that $X_n := f^{-n}(X_0)$ is connected for all $n \geq 0$. So $\bigcup_{n=0}^{\infty} f^{-n}(X) = \bigcup_{n=0}^{\infty} f^{-n}(\mathbb{R})$ is connected. This set is dense in $I(f)$ by the blowing-up property of the Julia set, and hence the escaping set is connected.

On the other hand, $\operatorname{sing}(f^{-1}) = \{0\}$, and $0$ is a finite asymptotic value. So by Theorem 7.48, the set $I(f)$ is not a spider's web, for any $a \in \mathbb{C}$. For $a > 1/e$, it can also be deduced directly from the work of Devaney [105] that the set $\mathbb{C} \setminus I(f)$ contains an unbounded connected subset; in particular, $I(f)$ is not a spider's web. See also a paper by Osborne and Sixsmith [207, Example 2], where it is even shown that the set of points with unbounded but non-escaping orbits has uncountably many pairwise disjoint unbounded connected subsets. $\square$

Jarque [151] and Rempe [223] obtained some generalisations of Theorem 7.52. In particular, the following holds; see [151, Theorem A] and [223, Theorem 1.2, (a)].

**Theorem 7.53.** *Let $f(z) = ae^z$, where $a \in \mathbb{C} \setminus \{0\}$. If $0 \in I(f)$ or if $0$ is preperiodic, then $I(f)$ is connected.*

The result in [223, Theorem 1.2] in fact applies to a much larger set of parameters $a$, but the condition is somewhat technical to state. The article also asks the following question [223, p. 73].

**Question 7.54** (Rempe, 2011). *Let $a \in \mathbb{C}$ be such that all periodic cycles of the exponential map $f(z) = ae^z$ are repelling. Is $I(f)$ connected?*

As we have seen, it may happen that $I(f)$ is connected and $\mathbb{C} \setminus I(f)$ contains an unbounded connected set. Rippon and Stallard [250, Question 10.2] raised the following question.

**Question 7.55** (Rippon and Stallard, 2019). *Does there exist a transcendental entire function $f$ such that $I(f)$ is connected and $\mathbb{C} \setminus I(f)$ contains an unbounded closed connected set?*

It can be deduced from [223, Proposition 3.2] that such sets do not exist for the exponential maps in Theorem 7.53, and in particular for $f = \exp$. (For this function the question was raised in [250, p. 3108].)

What can be said about the structure of the escaping set, and its connected components, when $I(f)$ is *not* connected? In §2.1 (and §2.3, in the case where $f_a$ is of disjoint type) we encountered examples where the escaping set has infinitely many and in fact uncountably many connected components. This extends to all disjoint-type entire functions.

**Proposition 7.56.** *Let $f$ be an entire function of disjoint type and let $R$ be an admissible radius for $f$. Then $J(f)$ has uncountably many connected components that intersect $A_R(f)$.*

*In particular, $J(f)$, $I(f)$, $A(f)$ and $A_R(f)$ all have uncountably many connected components.*

*Proof.* The claim follows from well-known properties of disjoint-type entire functions (compare the discussion before Theorem 1.2 of [228]). The fact that $J(f)$ and $I(f)$ have uncountably many connected components was noted in Theorem 3.10 and Corollary 3.11 of [224]. However, a proof for $A(f)$ and $A_R(f)$ does not seem to appear in the literature, so we sketch the argument here.

Since $f$ is of disjoint type, we have $f \in \mathcal{B}$. Recall that $I(f) \subset J(f)$ by Theorem 6.1. Therefore the first claim in the proposition does indeed imply the final sentence.

Let $\tilde{R} > R$ be sufficiently large (see below). Choose some $z_0 \in A_{\tilde{R}}(f) \neq \emptyset$, and let $C_0$ be the connected component of $J(f)$ containing $z_0$. It follows from Proposition 8.1 and Corollary 8.2 of [224] that there exist $\rho > 1$ and $M > 0$ (depending only on $f$) and an uncountable collection $\mathcal{C}$ of different connected components $C$ of $J(f)$ with the following property. For each $C \in \mathcal{C}$ there is a homeomorphism

$$\theta_C \colon C_0 \to C$$

with

$$|f^n(z)|^{1/\rho} \leq |f^n(\theta_C(z))| \leq |f^n(z)|^\rho$$

whenever $|f^n(z)| \geq M$. Assuming that $\tilde{R}$ was chosen sufficiently large, we have

$$|f^n(\theta_C(z_0))| \geq M^n(\tilde{R}, f)^{1/\rho} \geq M^n(\tilde{R}^{1/\rho}, f) \geq M(R, f)$$

by Lemma 5.4, for all $C \in \mathcal{C}$ and all $n \geq 0$. Hence each $C \in \mathcal{C}$ intersects $M(R, f)$, and the proof is complete. □

For general entire functions we have the following theorem, which is a consequence of a result of Rippon and Stallard [242, Theorem 5.2].

**Theorem 7.57.** *Let $f$ be a transcendental entire function, and let $E \subset I(f)$ be completely invariant under $f$ such that $J(f) = \overline{E \cap J(f)}$. Then either $E$ is connected or it has infinitely many connected components.*

Since the hypothesis holds for $E = I(f)$ and $E = A(f)$, either set is either connected or has infinitely many connected components. The question whether "infinitely many" can always be replaced by "uncountably many" remains open [250, Question 10.1].

**Question 7.58** (Rippon and Stallard, 2019)**.** Let $f$ be a transcendental entire function. For each of the sets $I(f)$ and $A(f)$, is it the case that it is either connected or it has uncountably many connected components?

The following result by Rippon and Stallard [250, Theorem 1.2] gives a partial result towards Question 7.58 for $I(f)$.

**Theorem 7.59.** *Let $f$ be a transcendental entire function such that $I(f)$ is disconnected. If $D$ is a disc intersecting $J(f)$, then $I(f) \setminus D$ has uncountably many connected components that meet $\partial D$.*

Concerning unbounded subsets of the levels $A_R(f)$, Rippon and Stallard [250, Theorem 1.3] have shown the following.

**Theorem 7.60.** *Let $f$ be a transcendental entire function and let $R_0$ be an admissible radius for $f$. Suppose there exists $R \geq R_0$ such that $A_R(f)$ is not a spider's web. Then there is a dense subset of $R \in [R_0, \infty)$ for which $A_R(f)$ contains uncountably many pairwise disjoint unbounded connected $F_\sigma$ sets.*

In [250, Question 10.3], the authors ask a number of questions related to this theorem.

**Question 7.61** (Rippon and Stallard, 2019)**.** In Theorem 7.60, can we replace the "unbounded connected $F_\sigma$ sets" in $A_R(f)$ by "unbounded connected closed sets" in $A_R(f)$ or even by "connected components" of $A_R(f)$?

**Question 7.62** (Rippon and Stallard, 2019)**.** In Theorem 7.60, can we replace the "dense set" of values of $R \in [R_0, \infty)$ by "all" $R \in [R_0, \infty)$?

**Question 7.63** (Rippon and Stallard, 2019)**.** Does there exist a transcendental entire function $f$ and an admissible radius $R$ for $f$ such that $A_R(f)$ is connected, but $A_R(f)$ is not a spider's web?

We recall that, if $f \in \mathcal{B}$, then $A_R(f)$ is not a spider's web by Corollary 7.34. It can be shown in this case that $A_R(f)$ does indeed have uncountably many connected components, and so the answers to Questions 7.61 to 7.63 are positive. This can be shown similarly as in Proposition 7.56, by adapting [224, Theorem 8.1] to the case of general functions in the Eremenko–Lyubich class. Instead, we will deduce this fact directly from Proposition 7.56 in Corollary 8.10 below.

Recall from Theorem 7.53 that the escaping set of an exponential map $f(z) = ae^z$ is connected when 0 is preperiodic. Of course, the reason that the orbit of 0 plays an important role is because $\operatorname{sing}(f^{-1}) = \{0\}$. It is interesting to note that the fact that 0 is an asymptotic value is crucial here. Indeed, Theorem 8.27 below (from [5]) describes a large class of functions in the Eremenko–Lyubich class $\mathcal{B}$ (namely the class of *strongly geometrically finite functions* with no asymptotic values in the Julia set and disconnected escaping sets such that the answer to Question 7.58 is positive. This includes all cosine maps $f(z) = ae^z + be^{-z}$ where $a, b \in \mathbb{C} \setminus \{0\}$ and both critical values of $f$ are strictly preperiodic, providing a counterpart to the second case in Theorem 7.53.

## 7.6 Counterexamples to Eremenko's conjecture

As we have seen, a number of weaker statements than Eremenko's conjecture hold for all transcendental entire functions, while a number of stronger conjectures were shown to fail some time ago. We have also discussed classes of entire functions for which Eremenko's conjecture is known to hold. Very recently Martí-Pete, Rempe and Waterman [181] showed that Eremenko's conjecture is false in general.

**Theorem 7.64.** *There exists a transcendental entire function $f$ such that $\{0\}$ is a connected component of $I(f)$.*

The construction is closely related to the ideas used to prove Theorem 6.46. Recall that this theorem shows that certain *bounded* simply connected domains can be realised as escaping wandering domains with non-escaping points on the boundary. As we already mentioned, the same method of proof also produces oscillating wandering domains with escaping points on the boundary. A key observation for Theorem 7.64 is that, with extra care, the same ideas can also be applied to certain simple *unbounded* domains. Hence it is possible to construct a transcendental entire function such that

(a) The half-strip $U := \{a + ib \colon a > 0 \text{ and } |b| < 1\}$ is an oscillating wandering domain of $f$, and

(b) $0 \in I(f)$.

The construction ensures that there is a sequence of topological strips, each tending to infinity in both directions and surrounding $U$ and approximating $\partial U$ more and more closely, each of which is a subset of an attracting basin of $f$. This not only means that $\partial U$ is a subset of the Julia set (a key fact for establishing (a)), but also that the connected component of $I(f)$ containing 0 is contained in $\overline{U}$. Recall that $U \cap I(f) = \emptyset$, and hence by Theorem 6.45, $\partial U \cap I(f)$ has zero harmonic measure as seen from $U$. Thus $\overline{U} \cap I(f) = \partial U \cap I(f)$ is totally disconnected. We conclude that $\{0\}$ is a connected component of $I(f)$, as claimed. (See [181, Remark 7.4].)

In fact, it is possible to apply the same construction to obtain an example for which 0 is not on the boundary, by replacing the closure $\overline{U}$ by a straight ray (which one might consider a degenerate half-strip of height zero). More precisely, the following holds [181, Theorem 7.1].

**Theorem 7.65.** *There exists a transcendental entire function $f$ with the following properties.*

(a) *The interval $[0, \infty)$ is a connected component of $J(f) \cup I(f)$.*

(b) *$0 \in I(f)$, but $t \notin I(f)$ for $t > 0$.*

This construction is modified in [181, Theorem 1.2] to also construct non-singleton bounded connected components of $I(f)$.

**Theorem 7.66.** *Let $K \subset \mathbb{C}$ be a compact set such that $\mathbb{C} \setminus K$ is connected. Then there exists a transcendental entire function $f$ such that $K$ is a connected component of $I(f)$.*

As noted in [181, Proposition 7.6], a bounded connected component of $I(f)$ need not be compact, but if $K$ is a compact component of $I(f)$, then $\mathbb{C} \setminus K$ is necessarily connected. We refer to [181] for further discussion of the results stated in this section, and their proofs.

Recall from Theorem 7.4 that Eremenko's conjecture holds – even in its strong form – when $f$ belongs to the Eremenko–Lyubich class $\mathcal{B}$ and also has finite order. In contrast, it is clear from the construction in [181] that the functions in the preceding theorems have both infinite order of growth and do not belong to the class $\mathcal{B}$ [181, Remark 7.5]. It is natural to ask whether or not this is necessary [181, Question 1.15].

**Question 7.67** (Martí-Pete, Rempe and Waterman, 2025). Does Eremenko's conjecture hold for all functions in $\mathcal{B}$?

If not, we may ask the same question for the smaller class $\mathcal{S}$.

**Question 7.68.** Does Eremenko's conjecture hold for all functions in $\mathcal{S}$?

**Question 7.69** (Martí-Pete, Rempe and Waterman, 2025). Is Eremenko's conjecture true for functions of finite order?

Note that this is weaker than Question 7.8 which asks whether such functions even have the strong Eremenko property.

# 8 Rigidity of escaping dynamics

Suppose that $f$ and $g$ are polynomials of the same degree. Then, by Böttcher's theorem (see, e.g., Milnor's book [192, Theorem 9.1]), $f$ and $g$ are conformally conjugate in a neighbourhood of infinity; that is, there exists a function $\theta$ conformal in a neighbourhood $U$ of $\infty$ such that
$$\theta(g(z)) = f(\theta(z)) \tag{8.1}$$
for $z \in U$. This is one of the most important foundational results in the study of the dynamics of complex polynomials. It means that the behaviour near infinity is determined only by the degree of the map, and often this universal structure can be used to obtain a good understanding of the function also away from infinity.

More precisely, suppose that $f$ is a polynomial of degree $d \geq 2$, and that the Julia set $J(f)$ is connected (which is equivalent to requiring that the postsingular set $P(f)$ be bounded [192, Theorem 9.5]). In this case the escaping set $I(f)$, which is the basin of the superattracting fixed point at $\infty$ of $f$, is simply connected (considered as a subset of the Riemann sphere). Let us choose $g$ in (8.1) to be the monomial $z \mapsto z^d$. Then the local conjugacy $\theta$ between $g$ and $f$ extends to a biholomorphic map between $I(g) = \mathbb{C} \setminus \overline{\mathbb{D}}$ and $I(f)$, and (8.1) holds whenever $|z| > 1$. If, moreover, $J(f)$ is *locally connected* (i.e., every point has arbitrarily small connected neighbourhoods), then by the Carathéodory–Torhorst theorem [192, Theorem 17.14], the map $\theta$ extends continuously to $\partial \mathbb{D}$, giving rise to a surjection $\theta \colon \partial \mathbb{D} \to J(f)$. In other words, the dynamics of $f$ on its Julia set can be completely described as a quotient of the dynamics of $g$ on the unit circle, which corresponds to the $d$-tupling of angles. For this reason, the local connectivity of the Julia set is an extremely important property in polynomial dynamics. We refer to [192, §18] for further discussion of this construction.

In contrast, even within simple families of entire functions, different maps may have drastically different escaping dynamics. For example, consider the functions
$$f_a \colon \mathbb{C} \to \mathbb{C}, \quad z \mapsto z + a + e^{-z},$$
with $a \in \mathbb{C}$. For $a = 1$, this is Fatou's function (2.3), which has a completely invariant Baker domain, in which $\operatorname{Re} f_a^n(z) \to +\infty$ locally uniformly; the same is true whenever

Re $a > 0$. But for other values of $a$ the escaping set $I(f_a)$ can have a very different structure.

Indeed, for $a = 0$ the map $f_a$ has infinitely many Baker domains; see the discussion after Theorem 6.44. As explained at the beginning of §6.5, the function $f_a$ has wandering domains if $a = 2\pi i - 1$ and infinitely many attracting basins if $a = -1$. In the latter case, $I(f_a)$ has empty interior.

Even within the family of exponential maps (2.4), $f_a(z) = ae^z$, we saw that the *global* structure of the escaping set varies dramatically. For instance, if $0 < a \leq 1/e$, then the set $I(f_a)$ has uncountably many connected components (see §2.3), while Theorem 7.52 says that if $a > 1/e$, then $I(f_a)$ is connected. If, however, we restrict $f_a$ to the set of points *whose orbits under iteration stay sufficiently large*, the situation becomes quite different: In Proposition 2.6 we observed that the structure of this set is essentially independent of the parameter $a$. It was shown by Rempe [218, Theorem 1.1] that there is indeed an analogue of Böttcher's theorem that connects the dynamics of any two exponential maps on sets of this form. In [221], he extended this to all natural parameter spaces (in a sense to be defined) within the Eremenko–Lyubich class $\mathcal{B}$. The goal of this section is to discuss these results and their implications. In particular, we discuss the notion of "docile functions" introduced by Alhamed, Rempe and Sixsmith [5] as an analogue of the local connectivity of polynomial Julia sets.

There have been a number of developments since the original publication of [221] that allow us to simplify and clarify the presentation of the theory. Where the statements we make cannot be found in the existing literature, we provide their proofs, using results and techniques from the relevant literature. These arguments are provided for completeness, and are not required for understanding the other parts of the section.

## 8.1 Basic properties of quasiconformal mappings

Throughout the section, we will use terminology from the theory of *quasiconformal mappings*. However, for most of our discussion, except §8.6, detailed knowledge of the theory of quasiconformal mappings is not required. Informally, quasiconformal mappings are homeomorphisms that do not distort too much. The formal definition is that a homeomorphism $\varphi$ between domains in $\mathbb{C}$ is called *quasiconformal* if it is absolutely continuous on almost all lines parallel to the coordinate axes and if there is a constant $k < 1$ such that

$$\left|\frac{\overline{\partial}\varphi}{\partial\varphi}\right| \leq k \qquad (8.2)$$

almost everywhere. (Here $\overline{\partial}\varphi$ and $\partial\varphi$ denote the Wirtinger derivatives; that is, $\partial\varphi = (\partial\varphi/\partial x - i\partial\varphi/\partial y)/2$ and $\overline{\partial}\varphi = (\partial\varphi/\partial x + i\partial\varphi/\partial y)/2$.) The condition (8.2) means that $\varphi$ preserves angles (almost everywhere) up to a constant factor; the quantity $\overline{\partial}\varphi/\partial\varphi$ is called the *complex dilatation* of $\varphi$. One may extend the definition to the case where $\varphi$ is a homeomorphism between domains on the Riemann sphere.

As mentioned, we will not require detailed knowledge about quasiconformal mappings, and readers unfamiliar with the theory should be able to follow most of the discussion if they keep the following key properties in mind.

**Proposition 8.1.**

(a) *Compositions and inverses of quasiconformal mappings are quasiconformal.*

(b) *Non-constant complex affine maps from $\mathbb{C}$ to $\mathbb{C}$ are quasiconformal.*

(c) *A quasiconformal map $\psi$ is Hölder continuous at every point of the Riemann sphere; in particular there is $K > 1$ such that*
$$|z|^{1/K} \leq |\psi(z)| \leq |z|^K$$
*for sufficiently large $z$.*

Here by a complex affine map we mean a map of the form $z \mapsto az + b$ with $a, b \in \mathbb{C}$. A real affine map can be written in the form $z \mapsto az + b\overline{z} + c$ with $a, b, c \in \mathbb{C}$. We note that such a map is quasiconformal if $|b| < |a|$, with complex dilatation equal to $b/a$.

For a detailed introduction to quasiconformal mappings (including a proof of the properties stated in Proposition 8.1) we refer to the books by Ahlfors [3] and Lehto and Virtanen [168].

## 8.2 Eventual quasiconformal equivalence

In order to formulate an analogue of Böttcher's theorem for the behaviour near infinity of transcendental entire functions, we must first decide what should take the place of the degree of the polynomial. It turns out that the correct notion is that of *eventual quasiconformal equivalence*, introduced in [221]. (In [221], the term "quasiconformal equivalence near $\infty$ is used; the more concise term "eventual quasiconformal equivalence" is from [231].)

**Definition 8.2.** Let $f, g \in \mathcal{B}$. Then we say that $f$ and $g$ are *eventually quasiconformally equivalent*, and write $f \underset{\text{qc}}{\overset{\infty}{\sim}} g$, if there exist quasiconformal homeomorphisms $\varphi, \psi \colon \mathbb{C} \to \mathbb{C}$ and $R > 0$ such that
$$\psi(g(z)) = f(\varphi(z)) \tag{8.3}$$
for all $z \in \mathbb{C}$ with $\max\{|\psi(g(z))|, |f(\varphi(z))|\} \geq R$.

*Remark* 8.3. (a) The definition also makes sense when $f$ and $g$ are polynomials: these are eventually quasiconformally equivalent if and only if they have the same degree.

(b) Clearly pre- and post-composition of $g \in \mathcal{B}$ by affine maps yields a map $f \underset{\text{qc}}{\overset{\infty}{\sim}} g$. In particular, any two maps in the exponential family (2.4) are eventually quasiconformally equivalent. The reader may wish to keep this example in mind in the following.

(c) Definition 8.2 is motivated by work of Eremenko and Lyubich [116, §3]. For $g \in \mathcal{S}$, they consider the set $M_g$ of entire functions $f$ *topologically equivalent* to $g$, i.e. those for which there exist homeomorphisms $\varphi$ and $\psi$ such that (8.3) holds for *all* $z \in \mathbb{C}$. Eremenko and Lyubich show that $M_g$ is a finite-dimensional analytic manifold, which may be considered the natural "parameter space" of $g$. If $f \in M_g$, then the maps $\varphi$ and $\psi$ can be chosen to be quasiconformal, and hence any two maps in the same Eremenko–Lyubich parameter space are eventually quasiconformally equivalent.

(d) Any map of the form $f(z) = p(z)e^z$, where $p$ is a polynomial, is eventually quasi-conformally equivalent to $g(z) = e^z$ [221, Observation 2.7], but if $p$ is non-constant, the two maps are not (globally) topologically equivalent in the sense of (c), since $f$ has critical points but $g$ does not.

Suppose that $f \underset{\mathrm{qc}}{\overset{\infty}{\sim}} g$. It is often useful to be more explicit about the domain where (8.3) holds.

**Definition 8.4.** Let $f \in \mathcal{B}$. We say that an unbounded domain $W_f \subset \mathbb{C}$ is an *initial domain* for $f$ if

(a) $W_f$ contains a punctured neighbourhood of $\infty$ and its boundary (in $\mathbb{C}$) is a Jordan curve, and

(b) $\mathrm{dist}(\mathrm{sing}(f^{-1}), W_f) > 0$.

Let $f, g \in \mathcal{B}$ and let $W_f, W_g \subset \mathbb{C}$ be initial domains for $f$ and $g$, respectively. Suppose that there exist quasiconformal maps $\psi \colon \mathbb{C} \to \mathbb{C}$ and $\varphi \colon \mathbb{C} \to \mathbb{C}$ with $W_f = \psi(W_g)$ and $f^{-1}(W_f) = \varphi(g^{-1}(W_g))$, and such that $\psi \circ g = f \circ \varphi$ on $g^{-1}(W_g)$. Then we say that the restrictions $f \colon f^{-1}(W_f) \to W_f$ and $g \colon g^{-1}(W_g) \to W_g$ are quasiconformally equivalent (via $\psi$ and $\varphi$), and write $(f, W_f) \underset{\mathrm{qc}}{\sim} (g, W_g)$.

Observe that $f \underset{\mathrm{qc}}{\overset{\infty}{\sim}} g$ if and only if there exist initial domains $W_f$ and $W_g$ such that $(f, W_f) \underset{\mathrm{qc}}{\sim} (g, W_g)$. The following result shows that the "only if" direction can be strengthened.

**Proposition 8.5.** *Suppose that $f, g \in \mathcal{B}$ are such that $f \underset{\mathrm{qc}}{\overset{\infty}{\sim}} g$. Let $W_f$ and $W_g$ be initial domains for $f$ and $g$, and suppose that there exists a quasiconformal homeomorphism $\psi \colon \mathbb{C} \to \mathbb{C}$ with $\psi(W_g) = W_f$.*

*Then there exists a quasiconformal map $\varphi \colon \mathbb{C} \to \mathbb{C}$ such that $(f, W_f) \underset{\mathrm{qc}}{\sim} (g, W_g)$ via $\psi$ and $\varphi$.*

*Proof.* This is a result of Rempe and Waterman from a forthcoming article concerning the notion of eventual quasiconformal equivalence [231]. As that article is still in preparation, we briefly outline the argument.

First, suppose that we already know that $(f, W_f) \underset{\mathrm{qc}}{\sim} (g, W_g)$ via quasiconformal homeomorphisms $\widetilde{\psi}$ and $\widetilde{\varphi}$, where $\widetilde{\psi} = \psi$ on $\partial W_g$. We wish to show that there exists a quasiconformal homeomorphism $\varphi$ such that also $(f, W_f) \underset{\mathrm{qc}}{\sim} (g, W_g)$ via $\psi$ and $\varphi$.

(This observation already appears in [221, §3].) By assumption (b), for any connected component $V$ of $g^{-1}(W_g)$, the restriction $g|_V \colon V \to W_g$ is a universal covering. The same is true for $f$ and $W_f$. It follows that $\psi$ lifts to a homeomorphism $\varphi \colon V \to \varphi(V)$ with $\psi \circ g = f \circ \varphi$. Since $g \colon \partial V \to \partial W_g$ is also a covering map, this lift $\varphi$ extends continuously to $\partial V$; the lift may be chosen so that this extension agrees with $\widetilde{\varphi}$ on $\partial V$. This defines $\varphi$ on $g^{-1}(W_g)$. We define $\varphi = \widetilde{\varphi}$ on $g^{-1}(\mathbb{C} \setminus W_g)$, and obtain a homeomorphism $\varphi \colon \mathbb{C} \to \mathbb{C}$. This map is quasiconformal by a classical glueing lemma due to Royden [77, Lemma 2].

Next we prove the proposition in the case where $f = g$ and $\overline{W_f} \subset W_g$. Choose two Jordan domains $D_0$ and $D_1$, bounded by analytic curves, such that $\overline{\operatorname{sing}(f^{-1})} \subset D_0$ and $\overline{D_0} \subset D_1 \subset \mathbb{C} \setminus \overline{W_g}$. If $\partial D_1$ is chosen sufficiently close to $\partial W_g$, then $\psi(\partial D_1)$ is close to $\psi(\partial W_g) = \partial W_f \subset W_g$. In particular, we can choose $D_1$ such that $\psi(\partial D_1) \subset W_g$, and in particular $\overline{D_0} \subset \psi(D_1)$. Hence $\psi(D_1) \setminus D_0$ is an annulus, whose outer boundary is the quasicircle $\psi(\partial D_1)$ and whose inner boundary is the analytic curve $\partial D_0$. Then there is a quasiconformal map $\widetilde{\psi}$ such that $\widetilde{\psi} = \operatorname{id}$ on $D_0$ and $\widetilde{\psi} = \psi$ on $\mathbb{C} \setminus D_1$. (See [169].) Now apply our first observation to the pairs $(f, \mathbb{C} \setminus \overline{D_0}) = (g, \mathbb{C} \setminus \overline{D_0})$, which are trivially quasiconformally equivalent via the identity. Since $\widetilde{\psi} = \operatorname{id}$ on $\partial D_0$, we obtain a quasiconformal map $\varphi$ that satisfies $\widetilde{\psi} \circ g = f \circ \varphi$ on $\mathbb{C} \setminus g^{-1}(D_0)$, and hence $\psi \circ g = \widetilde{\psi} \circ g = f \circ \varphi$ on $g^{-1}(W_g)$.

Finally, let $f$, $g$, $W_f$, $W_g$ and $\psi$ be as in the statement of the proposition. We may assume without loss of generality that $0 \notin \overline{W_g}$. Let $\psi_0$ and $\varphi_0$ be the maps from the definition of eventual quasiconformal equivalence and, for some $R \gg 1$, set $U_g := R \cdot W_g$ and $U_f := \psi_0(U_g)$. Define $\psi_g(z) := R \cdot z$ and $\psi_f := \psi_0 \circ \psi_g \circ \psi^{-1}$; then $\psi_g(W_g) = U_g$ and $\psi_f(W_f) = U_f$. If $R > 0$ was chosen sufficiently large, then $\overline{U_g} \subset W_g$ and $\overline{U_f} \subset W_f$, so we may twice apply the special case we just established and obtain $\varphi_f$ and $\varphi_g$ such that $\psi_f \circ f = f \circ \varphi_f$ on $f^{-1}(W_f)$ and $\psi_g \circ g = g \circ \varphi_g$ on $\mathbb{C} \setminus g^{-1}(W_g)$.

By choice of $\psi_0$ and $\varphi_0$, we have $\psi_0 \circ g = f \circ \varphi_0$ on $g^{-1}(U_g)$, provided that $R$ was chosen sufficiently large. Setting $\varphi := \varphi_f^{-1} \circ \varphi_0 \circ \varphi_g$, the proof is complete. $\square$

In the definition of quasiconformal eqivalence, it is not essential that the functions $f$ and $g$ are globally defined. More precisely, suppose that $W = W_f$ is an initial domain of $f \in \mathcal{B}$, and consider $\Omega := f^{-1}(W)$. Then $\Omega$ has the following properties:

(1) $\Omega$ is a disjoint union $\Omega = \bigcup_j \Omega_j$ of unbounded simply connected domains $\Omega_j$;

(2) any sequence of distinct $\Omega_j$ accumulates only at infinity;

(3) $\partial \Omega_j$ is connected for each $j$.

Let $h \colon W \to \mathbb{C} \setminus \overline{\mathbb{D}}$ be biholomorphic. For every $\Omega_j$, we can define a branch $\tau \colon \Omega_j \to \mathbb{H}_{>0}$ of $\log h \circ f$. This defines $\tau \colon \Omega \to \mathbb{H}_{>0}$ such that

(4) $\tau \colon \Omega_j \to \mathbb{H}_{>0}$ is biholomorphic for each $j$;

(5) if $(z_n)$ is a sequence in $\Omega$ such that $\tau(z_n) \to \infty$, then $z_n \to \infty$.

Conversely, let $W$ be any domain that is the outside of a Jordan curve, and let $\Omega \neq \emptyset$ and $\tau\colon \Omega \to \mathbb{H}_{>0}$ have properties (1) to (5). Following Bishop [81], we say that $f \coloneqq h^{-1} \circ \exp \circ \tau$ is a *model function*. (We remark that Bishop restricts to the case where $W = \mathbb{C} \setminus \overline{\mathbb{D}}$ and $h = \mathrm{id}$, which can always be achieved by restricting $f$ and conjugating with a complex affine map.)

The notion of eventual quasiconformal equivalence extends verbatim to model functions, and it is natural to ask when a model is realised by an entire function in $\mathcal{B}$ up to such equivalence. Following partial results by Rottenfußer, Rückert, Rempe and Schleicher [251, §7] and Rempe [232, Theorems 1.7 and 1.8], Bishop [81] showed that this is always the case.

**Theorem 8.6.** *Suppose that $f$ is a model function. Then there exists $g \in \mathcal{B}$ such that $f \underset{\mathrm{qc}}{\overset{\infty}{\sim}} g$.*

In contrast, Bishop [82, Theorem 1.4] showed that Theorem 8.6 does not hold when the class $\mathcal{B}$ is replaced by the Speiser class $\mathcal{S}$.

The results discussed in §8.3–8.6 all have appropriate versions for model functions, with the same proofs, although for simplicity's sake we will state them only for entire functions. By Theorem 8.6, the possible dynamics near infinity is the same for functions in $\mathcal{B}$ as it is for general model functions. Hence, to construct an entire function $f \in \mathcal{B}$ with specific dynamical behaviour near infinity, it suffices to construct a model with the same properties. See Theorem 8.17 below.

For a further discussion of eventual quasiconformal equivalence, we refer to the forthcoming paper [231].

## 8.3 Rigidity of escaping dynamics

Given a function $f \in \mathcal{B}$ we are interested in studying those points whose orbit under $f$ remains sufficiently large. Therefore we introduce the following notation.

**Definition 8.7.** Let $f \in \mathcal{B}$, and let $\rho > 0$. We define

$$J_{\geq \rho}(f) \coloneqq \{z \in J(f) \colon |f^n(z)| \geq \rho \text{ for all } n \geq 0\}. \tag{8.4}$$

The argument used to prove that $I(f) \subset J(f)$ for $f \in \mathcal{B}$ (Theorem 6.1) shows that, if $\rho$ is sufficiently large and $z \in \mathbb{C}$ satisfies $|f^n(z)| \geq \rho$ for all $n \geq 0$, then $z \in J(f)$. Thus the condition $z \in J(f)$ in (8.4) can be replaced by $z \in \mathbb{C}$ when $\rho$ is sufficiently large (which is the only case of interest to us here).

Now the main theorem of [221] can be stated as follows.

**Theorem 8.8.** *Suppose that $f, g \in \mathcal{B}$ with $f \underset{\mathrm{qc}}{\overset{\infty}{\sim}} g$. Then there exist $\rho > 0$ and $\eta > 0$ and a quasiconformal homeomorphism $\theta \colon \mathbb{C} \to \mathbb{C}$ such that*

$$\theta(g(z)) = f(\theta(z)) \tag{8.5}$$

*whenever $z \in J_{\geq \rho}(g) \cup \theta^{-1}(J_{\geq \eta}(f))$.*

Theorem 8.8 may be considered an analogue of Böttcher's theorem for functions in $\mathcal{B}$. Since every escaping point $z \in I(g)$ eventually enters $J_{\geq \rho}(g)$, for every $\rho > 0$, and similarly for $f$, Theorem 8.8 suggests that there is a close relationship between the escaping dynamics of $f$ and $g$ when $f \underset{\mathrm{qc}}{\overset{\infty}{\approx}} g$. In particular, we have the following.

**Proposition 8.9.** *Let $f$, $g$, $\rho$ and $\theta$ be as in Theorem 8.8, and let $z \in J_{\geq \rho}(g)$. Then $\theta(z) \in I(f)$ if and only if $z \in I(g)$, and $\theta(z) \in A(f)$ if and only if $z \in A(g)$.*

*Moreover, if $R \geq \rho$ is an admissible radius for $g$, then there exists an admissible radius $S \geq \eta$ for $f$ such that $A_S(f) \subset \theta(A_R(g))$. Here $A_S(f)$ and $A_R(g)$ are the levels of $A(f)$ and $A(g)$ defined in (5.12).*

*Proof.* We have $f^n(\theta(z)) = \theta(g^n(z))$, so $f^n(\theta(z)) \to \infty$ if and only if $g^n(z) \to \infty$. This proves the first claim.

We next prove the statement about the fast escaping sets. Let $\psi$ and $\varphi$ be as in Definition 8.2. By Proposition 8.1(c), the maps $\varphi$ and $\varphi^{-1}$ are both Hölder continuous. So there is $K > 1$ such that $|z|^{1/K} \leq |\varphi(z)| \leq |z|^K$ for all sufficiently large $z$, and similarly for $\psi$ and $\theta$. It follows that the orders of growth of $f$ and $g$ are related: When $r$ is large enough,
$$M(r, f) \geq M(r^{1/K}, g)^{1/K} = H_\varepsilon(r^{1/K})^K,$$
where $\varepsilon := 1/K^2$ and $H_\varepsilon(r) = M(r, g)^\varepsilon$ as defined in (5.21). Let $R > 0$ be large and put $S := R^K$. Let $z \in \mathbb{C}$ be such that $\theta(z) \in A(f)$. Then there exists $L \in \mathbb{N}$ such that
$$|f^{n+L}(\theta(z))| \geq M^n(S, f) \geq \bigl(H_\varepsilon^n\bigl(S^{1/K}\bigr)\bigr)^K = H_\varepsilon^n(R)^K$$
and hence
$$|g^{n+L}(z)| = |\theta^{-1}(f^{n+L}(\theta(z)))| \geq |f^{n+L}(\theta(z))|^{1/K} \geq H_\varepsilon^n(R)$$
for all $n \in \mathbb{N}$. Thus $z \in Q(g)$, where $Q(g)$ is the quite fast escaping set defined in (5.22). By Theorem 5.13, and since $g \in \mathcal{B}$, we also have $z \in A(g)$. So $z \in A(g)$ if $\theta(z) \in A(f)$. The converse implication follows by exchanging the roles of $f$ and $g$.

The result about the levels $A_R(g)$ and $A_S(f)$ follows by a similar reasoning from the proof of Theorem 5.13 by Rippon and Stallard [247] (more precisely from Theorem 1.3 and formula (1.3) in that paper). We omit the details. $\square$

As mentioned in §7.5, Theorem 8.8 can be used to give a positive answer to Questions 7.61, 7.62 and 7.63 when $f \in \mathcal{B}$.

**Corollary 8.10.** *Let $f \in \mathcal{B}$ and let $R$ be an admissible radius for $f$. Then $A_R(f)$ has uncountably many connected components.*

*Proof.* Let $g \underset{\mathrm{qc}}{\overset{\infty}{\approx}} f$ be of disjoint type. (For example, we may take $g = \lambda f$, where $|\lambda|$ is sufficiently small.) Let $\theta$ be as in Theorem 8.8. Applying Proposition 8.9 twice (once to $f$ and $g$ and once to $g$ and $f$), for any sufficiently large $S > 0$, there are $R > R' \geq \rho$ with $\theta(A_{R'}(g)) \supset A_S(f) \supset \theta(A_R(g))$. (Here $\rho$ is as in Theorem 8.8.)

By Proposition 7.56, there are uncountably many connected components of $J(g)$ that intersect $A_R(g)$. It follows that $A_S(f)$ also has uncountably many connected components.

This proves the claim of the corollary when $R$ is sufficiently large. For an arbitrary admissible radius $R$, let $n$ be so large that $A_{M^n(R,f)}(f)$ has uncountably many connected components. But $f^n(A_R(f)) = A_{M^n(R,f)}(f)$. Since $f^n$ is an open mapping and $A_R(f)$ is closed, it follows that $A_R(f)$ also has uncountably many connected components, as claimed. $\square$

In order to be able to discuss further properties of the map $\theta$ from Theorem 8.8, we sketch a proof of the theorem. In [221], the theorem is proved by applying a logarithmic change of variable to $f$ and $g$, as discussed in §6.1.

Here, we instead describe how to construct the function $\theta$ directly, without passing to logarithmic coordinates, making use of Proposition 8.5. The underlying idea is similar to the classical proof of Böttcher's theorem: Given $z \in J_{\geq R}(f)$, iterate forward under $g$, cross to the dynamical plane of $f$ using the quasiconformal map $\psi$, and then iterate backwards under $f$, taking care to use the "correct" branches of $f^{-1}$. As the number of iterates taken grows, the value obtained converges to the desired value $\theta(z)$.

To make this precise, we need to clarify how to identify inverse branches of $f$ near infinity with those of $g$. We do so using a well-established way of identifying different branches of $f^{-1}$ (and $g^{-1}$) near infinity, using "fundamental domains". (Compare [42, §2].)

Suppose that $W_f$ is an initial domain for $f$, with the additional property that

$$f^{-1}(W_f) \not\supset \mathbb{C} \setminus W_f, \tag{8.6}$$

which can be ensured for example by requiring that $0, f(0) \notin W_f$. Since every component of $\partial f^{-1}(W_f)$ is an arc tending to infinity in both directions, it follows that there exists an arc $\delta \subset W_f \setminus f^{-1}(W_f)$ that connects a point of $\partial W_f$ to $\infty$. Indeed, if $\partial f^{-1}(W_f) \cap \partial W_f \neq \emptyset$, then we can take $\delta \subset W_f \cap \partial f^{-1}(W_f)$. Otherwise, $\partial W_f \cap f^{-1}(W_f) = \emptyset$ by (8.6), and hence there is a curve connecting $\partial W_f$ to some point of $\partial f^{-1}(W_f)$ within $W_f \setminus f^{-1}(W_f)$, and the latter point can be connected to $\infty$ in $\partial f^{-1}(W_f)$.

Set $W_f^0 := W_f \setminus \delta$. Recall that $f\colon f^{-1}(W_f) \to W_f$ is a covering map. Since $W_f^0$ is simply connected, every connected component $F$ of $f^{-1}(W_f^0)$ is simply connected and $f\colon F \to W_f^0$ is biholomorphic. These domains $F$ are called *fundamental domains of $f$* (defined by $W_f$ and $\delta$). If $z \in \mathbb{C}$ is such that $f^n(z) \in W_f$ for all $n \geq 0$, then $f^n(z) \in f^{-1}(W_f^0)$ for all $n \geq 0$. (This is why it is important to choose $\delta \subset f^{-1}(\mathbb{C} \setminus W_f)$.)

We can carry out the same operation for $g$, obtaining an initial domain $W_g$, a curve $\gamma \subset W_g \setminus f^{-1}(W_g)$ and corresponding fundamental domains of $g$. If $(f, W_f) \underset{\text{qc}}{\sim} (g, W_g)$ via $\psi$ and $\varphi$, where $\psi$ is such that $\psi(\gamma) = \delta$, then the map $\varphi$ induces a correspondence between the fundamental domains of $g$ and $f$, and thus between the branches of $g^{-1}$ and $f^{-1}$. The next proposition shows that this can always be arranged.

**Proposition 8.11.** *Let $f \underset{\text{qc}}{\overset{\infty}{\sim}} g$. Then there exist initial domains $W_f, W_g$ and quasiconformal maps $\psi, \varphi$ such that $(f, W_f) \underset{\text{qc}}{\sim} (g, W_g)$ via $\psi$ and $\varphi$, and satisfying the following additional property. There exists an arc $\gamma \subset W_g \setminus g^{-1}(W_g)$ connecting a point of $\partial W_g$ to $\infty$ with $\psi(\gamma) \cap f^{-1}(W_f) = \emptyset$.*

*Proof.* By assumption, there are initial domains $U_f$ and $U_g$ such that $(f, U_f)$ and $(g, U_g)$ are quasiconformally equivalent via maps $\psi_0$ and $\varphi_0$. Let $\gamma_0$ be a piece of $\partial g^{-1}(U_g)$ connecting some finite point to $\infty$.

Set $\psi \coloneqq \varphi_0$; then
$$\psi(\gamma_0) = \varphi_0(\gamma_0) \subset \varphi_0(\mathbb{C} \setminus g^{-1}(U_g)) = \mathbb{C} \setminus f^{-1}(U_f).$$

Set $W_g \coloneqq \mathbb{C} \setminus \overline{D}(0, \rho)$, where $\rho$ is sufficiently large to ensure that $W_g \subset U_g$, $W_f \coloneqq \psi(W_g) \subset U_f$ and $\gamma_0 \cap \partial W_g \neq \emptyset$. By Proposition 8.5, there exists $\varphi$ such that $(f, U_f)$ and $(g, U_g)$ are quasiconformally equivalent via $\psi$ and $\varphi$. Let $\gamma$ be the piece of $\gamma_0$ that connects $\partial W_g$ to infinity in $W_g$; then $\gamma$ has the desired properties. $\square$

Now we can sketch the proof of Theorem 8.8. Let $f \underset{\text{qc}}{\approx} g$, and apply Proposition 8.11. Define $W_g^0 \coloneqq W_g \setminus \gamma$ and $W_f^0 \coloneqq W_f \setminus \psi(\gamma)$; then the connected components of $g^{-1}(W_g^0)$ and $f^{-1}(W_f^0)$ are fundamental domains of $g$ and $f$, respectively. Since $\psi(W_g^0) = W_f^0$, the map $\varphi$ induces a bijection between the set of fundamental domains of $g$ and the set of fundamental domains of $f$.

To construct the conjugacy $\theta$, define $\theta_0 \coloneqq \psi$, and inductively define
$$\theta_{n+1}(z) = f|_{\varphi(F_z)}^{-1}(\theta_n(g(z))) \tag{8.7}$$

where $F_z$ is the fundamental domain of $g$ containing $z$. If $\theta_n(z)$ is defined, then $g^k(z) \in W_g^0$ for $k = 1, \ldots, n$. Given such $z$, let $F_0, F_1, \ldots, F_{n-1}$ be the fundamental domains of $g$ containing $z, g(z), \ldots, g^{n-1}(z)$. Then $\theta_n(z)$ is the element of $f^{-n}(\psi(g^n(z)))$ that passes through the fundamental domains $\varphi(F_0), \varphi(F_1), \ldots, \varphi(F_{n-1})$, if such a point exists.

It can be shown that, if $z \in J_{\geq \rho}(g)$, where $\rho$ is sufficiently large, then $\theta_n(z)$ is defined for all $n \geq 0$. Moreover, the hyperbolic distance in $W_f$ between $\theta_n(z)$ and $\theta_{n+1}(z)$ tends to zero exponentially fast, uniformly on $J_{\geq \rho}(g)$. (The proof uses the expansion properties of $f$ and $g$, and fundamental properties of quasiconformal mappings; see [221, Lemma 2.6] for details.) Hence the maps $(\theta_n)$ form a Cauchy sequence and converge uniformly to a continuous function $\theta$ on $J_{\geq \rho}(g)$ satisfying (8.5).

If $z \in J_{\geq \rho}(g)$, then the point $w \coloneqq \theta(z)$ has the property that $f^n(w)$ and $\varphi(g^n(z))$ belong to the same fundamental domain of $f$, for all $n \geq 0$, and the hyperbolic distance between the two points (in $f^{-1}(W_f)$) is uniformly bounded. These two properties determine the point $w$ uniquely, and we see that $\theta$ is injective. When we exchange the roles of $f$ and $g$ in our construction, we obtain a map on $J_{\geq \eta}(f)$ for sufficiently large $\eta$. It is not difficult to see that this map agrees with $\theta^{-1}$ when $\eta$ is sufficiently large, and thus (8.5) also holds on $\theta^{-1}(J_{\geq \eta}(f))$.

This establishes the theorem except for the claim that $\theta$ extends to a quasiconformal homeomorphism of $\mathbb{C}$. The latter is deduced in an indirect fashion: the maps $f$ and $g$ are embedded in a holomorphic family of (not necessarily globally defined) model functions in the sense of §8.2. Fixing $g$, the function $\theta$ depends holomorphically on $f$ (see [221, Proposition 3.6]), and hence provides a *holomorphic motion* of the set $J_{\geq R}(f)$. A theorem of Bers and Royden [78, Theorem 1] then implies that $\theta$ extends to a global quasiconformal map, as claimed. We refer to [221] for details.

If $\theta$ is a map as constructed in the proof of Theorem 8.8, we shall say that $\theta$ is a *generalised Böttcher map* for $f$ and $g$. (That is, $\theta\colon \mathbb{C} \to \mathbb{C}$ is quasiconformal, and for some $\rho > 0$ the map $\theta|_{J_{\geq \rho}(g)}$ is the limit of the sequence $(\theta_n)$ defined via (8.7), for some initial configuration as in Proposition 8.11.)

We note that the above proof of the existence of a global generalised Böttcher map is highly non-constructive and non-canonical. In particular, the map $\theta$ does not relate the dynamics of $g$ to that of $f$ outside the nowhere dense set $J_{\geq \rho}(g)$. It is an interesting question whether the extension can also be constructed directly, and whether it can be chosen to conjugate $f$ and $g$ on a larger set. (We will see in §8.5 below that the answer is positive when both functions are of disjoint type.)

**Question 8.12.** Let $f, g \in \mathcal{B}$ with $f \overset{\infty}{\underset{\text{qc}}{\sim}} g$. Is it possible to construct a global generalised Böttcher map $\theta\colon \mathbb{C} \to \mathbb{C}$ explicitly, without using holomorphic motions or parameter arguments?

**Question 8.13.** Let $f, g \in \mathcal{B}$ with $f \overset{\infty}{\underset{\text{qc}}{\sim}} g$. Does there always exist a generalised Böttcher map $\theta$ for $f$ and $g$ so that the conjugacy relation (8.5) holds also on an open neighbourhood of $J_{\geq \rho}(g)$?

## 8.4 Non-uniqueness of the conjugacy

Let $f$ and $g$ be polynomials of degree $d$. The conformal conjugacy $\theta$ near $\infty$ between $f$ and $g$ given by (8.1) is not unique when $d \geq 2$, but depends on a choice of normalisation. More precisely, if $a$ and $b$ are the leading coefficients of $f$ and $g$, then $\theta(z) \sim cz$ as $z \to \infty$, where $c^{d-1} = b/a$. The map for a given $c$ is unique, but each $(d-1)$-th root of $b/a$ can be realised. (Usually the books, e.g. [192, Theorem 9.1], treat only the case $g(z) = z^d$, but the general case follows easily from this.)

A similar phenomenon occurs with the generalised Böttcher maps constructed in the previous subsection. Let us illustrate this with an example, where $f(z) := \sin(z)$ and $g(z) := \cos(z)$. Set $\psi := \operatorname{id}$ and consider the maps $\varphi_1(z) := z + \pi/2$, $\varphi_2(z) := z - 3\pi/2$ and $\varphi_3(z) := -z + \pi/2$. Then $\psi \circ g = f \circ \varphi_j$ for each $j \in \{1, 2, 3\}$. Let $W_f := W_g := g(\{z \in \mathbb{C}\colon |\operatorname{Im} z| > 1\})$. Then $\partial W_f$ is an analytic Jordan curve (in fact, the ellipse with focal points 1 and $-1$ that passes through the point $\cosh(1)$). Set $\gamma := [\cosh(1), \infty)$. Then the preimage of $W_g^0 := W_g \setminus \gamma$ is

$$\{a + ib\colon |b| > 1 \text{ and } a \notin 2\pi\mathbb{Z}\}.$$

So the fundamental domains of $g$ corresponding to $W_g$ and $\gamma$ are vertical half-strips in the upper and lower half-planes. Similarly, the fundamental domains of $f$ corresponding to $W_f$ and $\psi(\gamma) = \gamma$ are the connected components of

$$\{a + ib\colon |b| > 1 \text{ and } a \notin \pi/2 + 2\pi\mathbb{Z}\}.$$

Clearly each of the maps $\varphi_j$ induces a different bijection between the fundamental domains of $f$ and $g$, and hence all three give rise to different generalised Böttcher maps.

(Recall that, for large $R$, the function $\theta$ maps points of $J_R(f)$ belonging to a fundamental domain $F$ to points in the fundamental domain $\varphi(F)$.)

On the other hand, this is the only ambiguity that is possible for the values of a generalised Böttcher map on $J_{\geq\rho}(g)$. To make this statement more precise, we continue the notation from the proof of Theorem 8.8. The following uniqueness statement is a consequence of [221, Corollary 4.2].

**Theorem 8.14.** *Suppose that $f \underset{qc}{\overset{\infty}{\sim}} g$, that $\theta$ is a generalised Böttcher map for $f$ and $g$, and that $\rho$ is as in Theorem 8.8. Then for every $Q > \rho$ there is $Q' > 0$ with the following property.*

*Suppose that $\tilde{\theta}\colon J_{\geq Q}(g) \to J(f)$ is a conjugacy between $f$ and $g$ such that $\tilde{\theta}(z) \to \infty$ as $z \to \infty$, and such that $\theta(z)$ and $\tilde{\theta}(z)$ belong to the same fundamental domain of $f$ (as defined in the proof of Theorem 8.8) for all $z \in J_{\geq Q}(g)$. Then $\tilde{\theta} = \theta$ on $J_{\geq Q'}(g)$.*

In the case where $\theta$ and $\tilde{\theta}$ both are generalised Böttcher maps, this leads to a uniqueness statement that can be phrased as follows. A sequence $\underline{s} = s_0 s_1 s_2 \ldots$ of fundamental domains $s_k$ of $f$ is called an *external address* of $f$. A point $z$ has external address $\underline{s}$ for $f$ if $f^k(z) \in s_k$ for $k \geq 0$. If $\sigma > 0$ is sufficiently large, then every point $z \in J_{\geq \sigma}(f)$ has an external address. (We previously discussed external addresses in connection with Theorem 7.27. If $f$ is postsingularly bounded, and if $\sigma$ is sufficiently large, then two points *have the same address* in the sense discussed there if and only if they visit the same sequence of fundamental domains. See also §8.8.)

We remark that the definition of external addresses depends on the initial configuration given by $W_f$ and the curve $\delta$. However, this dependence is not essential: given two different such configurations, there is a natural bijection between the two corresponding sets of external addresses (see the discussion preceding [42, Observation 4.12]).

Hence it makes sense to speak about external addresses of $g$ and external addresses of $f$. Any generalised Böttcher map $\theta$ induces a natural bijection between the addresses of $g$ and those of $f$. It is this bijection that essentially determines $\theta$: If $\tilde{\theta}$ is another generalised Böttcher map that induces the same bijection, then $\theta$ and $\tilde{\theta}$ agree on $J_{\geq \rho}(g)$ for sufficiently large $\rho$.

## 8.5 The case of disjoint-type maps

We now return to Theorem 8.8, in the special case where both $f$ and $g$ are of disjoint type. It turns out that in this case, the answer to Questions 8.12 and 8.13 is positive; that is, we are able to construct a generalised Böttcher map for $f$ and $g$ more directly. Recall from Proposition 3.5 that $f$ is of disjoint type if and only if $W_f$ can be chosen such that $\overline{f^{-1}(W_f)} \subset W_f$.

**Theorem 8.15.** *Let $f, g \in \mathcal{B}$ be of disjoint type, and let $W_f$ and $W_g$ be initial domains for $f$ and $g$, with the additional property that $f^{-1}(\overline{W_f}) \subset W_f$ and $g^{-1}(\overline{W_g}) \subset W_g$.*

*Assume that $(f, W_f) \underset{qc}{\sim} (g, W_g)$ via quasiconformal maps $\psi$ and $\varphi$ with $\varphi = \psi$ on $\mathbb{C} \setminus W_g$. Then there exists a quasiconformal homeomorphism $\theta \colon \mathbb{C} \to \mathbb{C}$ such that*

(a) $f \circ \theta = \theta \circ g$ on $g^{-1}(W_g)$ and

(b) $\theta = \varphi$ on $\mathbb{C} \setminus g^{-1}(W_g)$.

*The map $\theta$ is unique and depends only on the values of $\varphi$ on $\mathbb{C} \setminus g^{-1}(W_g)$, but not on the values of $\psi$ on $W_g$. Moreover, $\theta$ is a generalised Böttcher map.*

*Proof.* The existence of $\theta$ is proved using a pull-back argument. Using the fact that $f \circ \varphi = \varphi \circ g$ on $\partial g^{-1}(W_g)$, it is easy to construct a sequence $(\theta_n)$ of quasiconformal maps with $\theta_0 = \varphi$, $\theta_n \circ g = f \circ \theta_{n+1}$ on $g^{-1}(W_g)$, and such that $\theta_n = \varphi$ on $\mathbb{C} \setminus g^{-1}(W_g)$. All the $(\theta_n)$ have the same maximal dilation and the sequence $(\theta_n(z))$ is eventually constant for every $z$ in the dense open set $\mathbb{C} \setminus J(g)$. It follows that $\theta_n$ converges to the desired function $\theta$; we refer to the proof of [221, Theorem 3.1] for details.

The conditions (a) and (b) mean that the values of $\varphi$ on $\mathbb{C} \setminus g^{-1}(W_g)$ determine $\theta$ uniquely on the open and dense set $\mathbb{C} \setminus J(g)$, and hence on all of $\mathbb{C}$. This proves the second claim.

Finally, set $\psi_0 := \varphi$ and $\varphi_0 := \theta_1$, where $\theta_1$ is as defined in the construction of $\theta$. Then $(f, W_f) \underset{\mathrm{qc}}{\sim} (g, W_g)$ via $\psi_0$ and $\varphi_0$. Moreover, these maps satisfy the additional property from Proposition 8.11, for any curve $\gamma$ as in the statement of the proposition. (Indeed, $\varphi(W_g \setminus g^{-1}(W_g)) = W_f \setminus f^{-1}(W_f)$ by assumption.) If we follow the proof of Theorem 8.8 starting with this configuration, then the resulting sequence of maps $\theta_n$ agrees with the sequence of quasiconformal maps constructed at the beginning of this proof, which converge to $\theta$. Thus $\theta$ is a generalised Böttcher map, as claimed. □

**Proposition 8.16.** *Let $f$ and $g$ be of disjoint type with $f \underset{\mathrm{qc}}{\overset{\infty}{\approx}} g$. Let $\tilde{\theta}$ be any generalised Böttcher map for $f$ and $g$. Then there exists $\rho \geq 0$ and a generalised Böttcher map $\theta$ with $\theta = \tilde{\theta}$ on $J_{\geq \rho}(g)$ and such that $\theta$ restricts to a conjugacy between $g|_{J(g)}$ and $f|_{J(f)}$.*

*(In fact, $\theta$ has the properties stated in Theorem 8.15, for a suitable choice of $W_f$, $W_g$, $\psi$ and $\varphi$.)*

*Proof.* By Proposition 3.5, we may choose initial domains $W_f$ and $W_g$ with $f^{-1}(\overline{W_f}) \subset W_f$, and similarly for $W_g$. Moreover, $\partial W_f$ and $\partial W_g$ can be chosen to be analytic Jordan domains.

Now suppose that we are given a configuration as in Proposition 8.11; that is, $U_f$ and $U_g$ are initial domains for $f$ and $g$, we have $(f, U_f) \underset{\mathrm{qc}}{\sim} (g, U_g)$ via maps $\psi_0$ and $\varphi_0$, and $\gamma_0 \subset U_g \setminus g^{-1}(U_g)$ is an arc connecting a point of $\partial U_g$ to $\infty$, such that $\delta_0 := \psi_0(\gamma_0) \subset U_f \setminus f^{-1}(U_f)$. We wish to prove the proposition for the generalised Böttcher map $\tilde{\theta}$ associated to this configuration.

We may choose $U_g$ sufficiently small to ensure that $\overline{U_g} \subset W_g$ and $\overline{U_f} \subset W_f$; this does not affect the values of $\tilde{\theta}$ on $J_{\geq \rho}(g)$ for sufficiently large $\rho$. Similarly, we assume that $U_g \cup g^{-1}(W_g)$ does not separate the finite endpoint of $\gamma_0$ from $\partial W_g$, and that the same holds for $f$ and $\delta_0$. (To achieve this, let $U_g$ be a small punctured neighbourhood of $\infty$, together with a thin simply connected neighbourhood of $\gamma_0$.)

Now we may choose a piecewise analytic curve $\gamma_1 \subset W_g \setminus (U_g \cup g^{-1}(W_g))$ that connects $\partial W_g$ to the finite endpoint of $\gamma_0$, and a similar curve $\delta_1$ for $f$. We can then extend $\psi_0$ to a quasiconformal map $\psi \colon W_g \to W_f$ with $\psi(\gamma_1) = \delta_1$.

We claim that there exists a quasiconformal map $\varphi \colon \mathbb{C} \to \mathbb{C}$ such that

(a) $(f, W_f) \underset{\mathrm{qc}}{\sim} (g, W_g)$ via $\psi$ and $\varphi$;

(b) $\varphi = \varphi_0$ on $g^{-1}(U_g)$;

(c) $\varphi = \psi$ on $\mathbb{C} \setminus W_g$.

(Observe that Proposition 8.5 shows the existence of $\varphi$ satisfying (a), but the map constructed there will not satisfy (b).) Observe that there is a unique continuous extension $\varphi \colon g^{-1}(W_g) \to f^{-1}(W_f)$ of $\varphi_0$ with $\psi \circ g = f \circ \varphi$; we claim that this map has a global quasiconformal extension to the complex plane.

This follows again from [231]; we sketch the proof here. Consider the simply connected domains $\Delta_g := \mathbb{C} \setminus g^{-1}(U_g)$, $Y_g := \mathbb{C} \setminus g^{-1}(W_g)$ and $D_g := \mathbb{C} \setminus W_g$. Then $\overline{Y_g} \subset \Delta_g$ and $\overline{D_g} \subset Y_g$. Let $h_g \colon \Delta_g \to \mathbb{D}$ be a Riemann map; then it can be shown that $h_g(Y_g)$ is a quasidisc (i.e., the image of a Euclidean disc under a quasiconformal homeomorphism of the plane). The same is true for $\Delta_f, Y_f, D_f$ and $h_f$ defined analogously.

The map $\underline{\varphi}_0 := h_f \circ \varphi_0 \circ h_g^{-1}$ is a quasiconformal self-map of the unit disc, and hence has a quasiconformal extension to the entire Riemann sphere, which we also denote by $\underline{\varphi}_0$. Define $\underline{\varphi}$ on $\widehat{\mathbb{C}} \setminus \overline{h_g(Y_g)}$ by

$$\underline{\varphi}(z) := \begin{cases} h_f(\varphi(h_g^{-1}(z))) & \text{if } z \in h_g(\Delta_g \setminus \overline{Y_g}) \\ \underline{\varphi}_0(z) & \text{otherwise.} \end{cases}$$

Now
$$\underline{\varphi} \colon \widehat{\mathbb{C}} \setminus \overline{h_g(Y_g)} \to \widehat{\mathbb{C}} \setminus \overline{h_f(Y_f)}$$
is a quasiconformal homeomorphism between two quasidiscs, and again has a quasiconformal extension to all of $\widehat{\mathbb{C}}$. Since $\overline{h_g(D_g)}$ is a compact subset of $h_g(Y_g)$, and similarly for $f$, we may choose the extension in such a way that $\underline{\varphi} = h_f \circ \psi \circ h_g^{-1}$ on $h_g(D_g)$. Defining $\varphi := h_f^{-1} \circ \underline{\varphi} \circ h_g$ on $\overline{\Delta}_g$, we obtain the desired global extension of $\varphi$.

Now we are in a position to apply Theorem 8.15 to obtain the desired map $\theta$. As in the proof of the final claim in Theorem 8.15, it follows directly from the construction that $\theta = \tilde{\theta}$ on $J_{\geq \rho}(g)$ when $\rho$ is sufficiently large. □

Suppose that $f \colon \Omega \to W$ is a model function as discussed in §8.2, and that additionally $\overline{\Omega} \subset W$. Then $f$ is said to be of *disjoint type*, and the conjugacy results we discussed in this section hold also for such models. Together with Theorem 8.6, this implies the following (see [81, Theorem 1.2]).

**Theorem 8.17.** *Suppose that $f$ is a disjoint-type model function. Then there exists a disjoint-type function $g \in \mathcal{B}$ and a quasiconformal homeomorphism $\theta$ that conjugates $f$ and $g$ on their respective Julia sets.*

Here the Julia set $J(f)$ of a disjoint-type model function consists of those points that do not leave $\Omega$ under iteration.

## 8.6 Absence of line fields on the escaping set

Another result proved in [221] is the *absence of invariant line fields* on the escaping set of a function $f \in \mathcal{B}$ [221, Theorem 1.2]. Here an *invariant line field* for an entire function $f$ is a measurable choice of tangent line, defined on a set that is forward-invariant under $f$, and such that the field of lines is invariant under $f$ almost everywhere.

**Theorem 8.18.** *Let $f \in \mathcal{B}$, and suppose that there is an $f$-invariant line field defined on a completely invariant set $X \subset I(f)$. Then $X$ has zero two-dimensional Lebesgue measure.*

Observe that the conclusion Theorem 8.18 is non-trivial only if $I(f)$ has positive Lebesgue measure. This may or may not hold; for example, $I(\sin)$ has positive Lebesgue measure by McMullen's Theorem 5.9, while $I(\exp)$ does not by Theorem 9.24 of Eremenko and Lyubich, discussed below.

Once again, the assumption that $f \in \mathcal{B}$ is crucial: An entire function $f \notin \mathcal{B}$ may support an invariant line field on $J(f) \cap I(f)$. This was first shown by Eremenko and Lyubich [114, Example 5], using approximation theory. This is also a direct corollary of Theorem 7.20: Let $K \subset \mathbb{C}$ be a full compact set of positive area and with empty interior (for example, a Jordan arc of positive area). By Theorem 7.20, there exists an entire function $f$ such that $K \subset J(f)$, and $f^n(K) \cap f^m(K) = \emptyset$ for $n \neq m$. By pullbacks, any measurable line field on $K$ can be extended to an invariant line field on the grand orbit of $K$.

The non-existence of invariant line fields supported on the Julia set of an entire function $f$ is an important question because it is connected to the existence of non-trivial quasiconformal deformations of $f$. For the remainder of the subsection, we assume more familiarity with background and terminology from the theory of quasiconformal mappings than in the rest of the section; we again refer to [3, 168].

Every measurable invariant line field gives rise to a non-trivial measurable invariant Beltrami differential $\mu(z)\frac{\partial \bar{z}}{\partial z}$, and vice versa. By the measurable Riemann mapping theorem, such a differential in turn gives rise to a non-affine conjugacy between $f$ and another entire function $g$. Hence Theorem 8.18 can also be stated as follows.

**Theorem 8.19.** *Let $f, g \in \mathcal{B}$, and suppose that there is a quasiconformal homeomorphism $\varphi \colon \mathbb{C} \to \mathbb{C}$ such that $f \circ \varphi = \varphi \circ g$. Then the complex dilatation of $\varphi$ vanishes almost everywhere on $I(g)$.*

It is conjectured that if $f$ is a polynomial, then $J(f)$ does not support any invariant line fields. This would mean that any quasiconformal deformation of $f$ arises in a well-understood way from invariant Beltrami differentials supported on $F(f)$. (Compare [186, Chapter 1].) If true, this would imply the famous conjecture that *hyperbolicity is dense* in the space of polynomials of any given degree.

It is also conjectured that a function $f \in \mathcal{S}$ does not support invariant line fields on its Julia set. Theorem 8.18 states that such a line field cannot be supported on the escaping set. Rempe and van Strien [233] combined this result with techniques from the study of real polynomials to prove density of hyperbolicity in certain classes of real functions $f \in \mathcal{S}$.

In view of its reformulation in Theorem 8.19, it is perhaps not surprising that Theorem 8.18 is closely related to the study of the dilatation of a generalised Böttcher map $\theta$.

**Theorem 8.20.** *Let $f, g \in \mathcal{B}$ with $f \overset{\infty}{\underset{\mathrm{qc}}{\sim}} g$, and let $\theta$ be a generalised Böttcher map for $f$ and $g$. Then the complex dilatation of $\theta$ vanishes on $I(g) \cap J_{\geq \rho}(g)$, for sufficiently large $\rho$.*

Theorem 8.20 follows from Theorem 8.18. Indeed, if the dilatation did not vanish, it would give rise to an invariant line field on the escaping set. The converse also holds: In [221], Theorem 8.20 is established first. Then Theorem 8.19 (and hence Theorem 8.18) are deduced using the results on the uniqueness of the conjugacy from §8.4.

The proof of Theorem 8.20 in [221] itself is also rather indirect. A key idea is to split the proof into two cases: one where $f$ and $g$ are both of disjoint type, and one where $f$ and $g$ are not only quasiconformally equivalent, but related by *affine* maps. It seems interesting to ask whether there is a more direct argument. To make this question more precise, it is useful to outline the proof from [221], which we present here in a slightly more streamlined version. This approach also shows how to prove Theorem 8.18 without appealing to the results from §8.4. We will not, however, attempt to provide technical details, and refer to [221] for these.

As mentioned, the first step is to consider the case where $f$ and $g$ are both of disjoint type. Recall from Proposition 8.16 that every generalised Böttcher map arises from the construction in Theorem 8.15. In that theorem, we may modify $\psi$ on $W_g$ so that $\psi$ is conformal in a neighbourhood of infinity that includes $g^{-1}(\overline{W_g})$. This does not affect the values of $\varphi$ outside $\mathbb{C} \setminus g^{-1}(W_g)$, and hence, it does not change the map $\theta$. Each of the maps $\theta_n$ constructed in the proof of Theorem 8.15 is conformal on a neighbourhood of $J(g)$, and it follows that the dilatation of the limit function $\theta$ itself vanishes on this set. We again refer to the proof of [221, Theorem 3.1] for details. This proves Theorem 8.20 when $f$ and $g$ are of disjoint type.

This argument also proves Theorem 8.18 when $f$ is of disjoint type, via Theorem 8.19. Indeed, if $f$ and $g$ are of disjoint type and quasiconformally conjugate via $\theta$, then we can take $\psi = \varphi$ in Theorem 8.15, and get $\theta = \varphi$. As we have just shown, the complex dilatation of this map vanishes on $J(g)$.

Now suppose that $f$ and $g$ are *affine equivalent*; that is, $\psi \circ f = g \circ \varphi$ for two complex affine maps $\psi, \varphi \colon \mathbb{C} \to \mathbb{C}$. Recall that, in the proof of Theorem 8.8, the map $\theta$ was first constructed on some $J_{\geq \rho}(g)$, and then extended to a quasiconformal map using the Bers–Royden version of the "$\lambda$-lemma", a tool from the theory of holomorphic motions. In the case of affine equivalence, it is possible to show that the complex dilatation of this extension is small when $\rho$ is large. Since points in $I(g)$ are eventually mapped to $J_{\geq \rho}(g)$ for every $\rho$, this leads to the following result. (See [221, Theorem 3.4].)

**Theorem 8.21.** *Suppose that $f, g \in \mathcal{B}$ are affine equivalent, let $\theta$ be a generalised Böttcher map, and let $\mu = \overline{\partial}\theta/\partial\theta$ denote the complex dilatation of $\theta$. Then for every $\varepsilon > 0$, there is $\rho > 0$ such that $|\mu(z)| \leq \varepsilon$ for almost every $z \in J_{\geq \rho}(g)$.*

*In particular, the complex dilatation of $\theta$ vanishes on $J_{\geq\rho}(g) \cap I(g)$ for sufficiently large $\rho$.*

*Proof of Theorem 8.18.* Suppose, by contradiction, that $I(f)$ supports an invariant line field. Let $\rho \geq 0$ be arbitrary. Since $I(f) \subset \bigcup_{n=0}^{\infty} f^{-n}(J_{\geq\rho}(f))$, it follows that $J_{\geq\rho}(f) \cap I(f)$ supports a forward-invariant line field.

Now let $g$ be a disjoint-type function affine equivalent to $f$ (such as $g \coloneqq \varepsilon \cdot f$, for sufficiently small $\varepsilon$), and let $\theta$ be a generalised Böttcher map for $f$ and $g$. By Theorem 8.21, the dilatation of $\theta$ vanishes on $J_{\geq\rho}(f) \cap I(f)$, for sufficiently large $\rho$. But then we could use $\theta$ to transfer the line field that is invariant under $f$ to a forward-invariant line field of $g$, which in turn can easily be extended to a fully invariant line field supported on $I(g)$. However, we proved above that disjoint-type functions support no invariant line fields on their escaping sets, a contradiction. □

*Proof of Theorem 8.20.* As already noted, Theorem 8.20 follows from Theorem 8.18. Alternatively, let $\tilde{f}$ and $\tilde{g}$ be disjoint-type functions affine equivalent to $f$ and $g$, respectively. Any generalised Böttcher map between $f$ and $g$ can be written (on $J_{\geq\rho}(g)$ for sufficiently large $\rho$) as a composition of generalised Böttcher maps for the pairs $f \underset{\text{qc}}{\overset{\infty}{\sim}} \tilde{f}$, $\tilde{f} \underset{\text{qc}}{\overset{\infty}{\sim}} \tilde{g}$ and $\tilde{g} \underset{\text{qc}}{\sim} g$. As we already saw, the dilatation of each of these vanishes on the corresponding escaping set, proving the theorem. Compare [221, Corollary 3.5]. □

As already noted, Theorems 8.18 and 8.20 trivially hold whenever the area of $I(f)$ is zero. McMullen's Theorem 5.9 says that for maps in the *sine family* $z \mapsto a\sin(z) + b$, the escaping set has positive measure. For maps in this family, a direct proof of Theorem 8.18 was given by Chen, Jiang and Zhang [94], using the geometric mapping properties of the sine family, and without appealing to parameter arguments or the $\lambda$-lemma. It is an interesting question whether it is possible to prove the general form of Theorem 8.18 in a similar manner. (This problem was raised in [221, p. 255].)

**Question 8.22** (Rempe, 2009)**.** Let $f \in \mathcal{B}$. Is there a proof of the absence of $f$-invariant line fields on $I(f)$ that uses only dynamical properties of the function $f$, and avoids arguments involving holomorphic motions in parameter space?

## 8.7 Natural bijections between escaping sets

Recall that, if $f$ is a polynomial of degree $d$ with bounded postsingular set, then the Böttcher map $\theta$ between $g(z) = z^d$ and $f$ extends to a conformal conjugacy between the two maps on their escaping sets.

A similar result holds for an entire function $f \in \mathcal{B}$ with bounded postsingular set, if we replace the model map $z \mapsto z^d$ by any disjoint-type function $g$ with $g \underset{\text{qc}}{\overset{\infty}{\sim}} f$. In contrast to the polynomial case, there is no obvious canonical candidate to use in a given eventual quasiconformal equivalence class. However, as we saw in Proposition 8.16, any two different choices of $g$ are quasiconformally conjugate on their Julia sets, which means that the specific choice of $g$ is not essential. Also recall that such $g$ always exists, since $g \coloneqq \lambda f$ is of disjoint type when $|\lambda|$ is sufficiently small.

The extension of the map $\theta$ can be described as follows.

**Theorem 8.23.** *Let $f, g \in \mathcal{B}$ with $f \underset{qc}{\overset{\infty}{\sim}} g$, where $g$ is of disjoint type and the postsingular set $P(f)$ is bounded. Let $\theta$ be a generalised Böttcher map between $f$ and $g$. If $\rho > 0$ is sufficiently large, then the restriction of $\theta$ to $X := I(g) \cap J_{\geq \rho}(g)$ has a unique extension $\tilde{\theta} \colon I(g) \to I(f)$ such that*

(a) *$\tilde{\theta}$ is a bijection;*

(b) *$\tilde{\theta} \circ g = f \circ \tilde{\theta}$ on $I(g)$;*

(c) *if $A \subset \mathbb{C}$ is such that $\inf_{z \in A} |g^n(z)| \to \infty$ as $n \to \infty$, then $\tilde{\theta} \colon A \to \tilde{\theta}(A)$ is a homeomorphism.*

*Proof.* The theorem is closely related to a result by Alhamed, Rempe and Sixsmith [5, Theorem 2.2]; many of the underlying ideas can already be found in [219]. However, the theorem does not appear in the literature in the form in which we have stated it here, so we outline the proof, using the results and techniques of [5].

The existence of the extension is proved in [5, Theorem 2.2], for the case where $g$ is of the form $g(z) = f(\lambda z)$, with $\lambda > 0$ sufficiently small, and for a specific choice of $\theta$. We remark that [5, Theorem 2.2] does not explicitly state that the extension has the property we have formulated here; the extension is characterised in another way. However, our characterisation follows immediately from the proof in [5]. Indeed, if $z \in I(f)$, then $\tilde{\theta}(z)$ is defined by iterating $f$ a finite number $n$ of times, until $f^n(z) \in X$, and then applying a suitable backward branch of $g^{-1}$. (The key fact is that the value of $\tilde{\theta}(z)$ is then independent of $n$.) If $A$ is as in the statement of the theorem, then continuity of $\tilde{\theta}$ follows because the inverse branches in question are continuous.

Now suppose that $g_1$ is another disjoint-type map with $f \underset{qc}{\overset{\infty}{\sim}} g_1$ and $\theta_1$ is a generalised Böttcher map for $f$ and $g_1$. Then $\theta_1^{-1} \circ \theta$ is a generalised Böttcher map for $g$ and $g_1$. By Proposition 8.16, its restriction to $J_{\geq \rho}(g)$, for large $\rho$, extends to a quasiconformal conjugacy $\pi \colon J(g) \to J(g_1)$. Then $\tilde{\theta}_1 := \tilde{\theta} \circ \pi^{-1}$ satisfies the conditions in the statement of the proposition, for $f$, $g_1$ and $\theta_1$. (Alternatively, the proof of [5, Theorem 2.2] can easily be extended to cover the general case.)

Finally, we must show that $\tilde{\theta}$ is uniquely determined by the condition in the theorem. Let $\eta$ be as in the statement of Theorem 8.8; we may suppose that $\eta$ is so large that $J_{\geq \eta}(f) \subset \theta(J_{\geq \rho}(g))$. Choose $\rho' > \rho$ so large that $D(0, \rho')$ contains $P(g)$, and $\eta' \geq \eta$ so large that $D(0, \eta')$ contains $P(f)$ and furthermore $\theta(D(0, \rho')) \supset D(0, \eta' + 1)$.

First suppose that $A \subset I(g)$ is an unbounded closed connected set and such that $\inf_{z \in A} |g^n(z)| \to \infty$. There is some $n$ such that $g^n(A) \subset J_{\geq \rho'+1}(g)$. Then $g^n(A)$ is contained in a component $U$ of $\{z \colon |g(z)| > \rho'\}$. Since $U$ is simply connected and disjoint from $P(g)$, there is a branch $\varphi$ of $g^{-n}$ defined on $U$ that maps $g^n(A)$ to $A$. In particular, $g^n(A)$ is also unbounded.

Similarly, $B_n := \theta(g^n(A))$ is an unbounded closed connected set with $B_n \subset J_{\geq \eta'+1}(f)$. So $B_n$ is contained in a component $V$ of $\{z \colon |f(z)| > \eta'\}$. If $\tilde{\theta}$ is an extension as in the statement of the theorem, then $B := \tilde{\theta}(A)$ is connected and hence contained in a single connected component $W$ of $f^{-n}(V)$. Let $\alpha$ be the branch of $f^{-n}$ defined on $V$ that maps $B_n$ to $B$; by the functional equation we have $\tilde{\theta}|_A = \alpha \circ \theta \circ g^n$.

We claim that $W$, and hence $\alpha$, is determined already by the original map $\theta$. If $z \in A$ is sufficiently large, then (by the continuity of the branch of $\varphi$ at $\infty$), $z \in J_{\geq \rho}(g)$. Since $\theta(z) = \tilde{\theta}(z) \subset B \subset W$, we see that $W$ is the unique connected component of $g^{-n}(V)$ that contains $\theta(z)$ for sufficiently large $z \in A$. We conclude that any two choices of $\tilde{\theta}$ must coincide on the set $A$.

Now let $z \in I(g)$ be arbitrary. It follows from [219, Proposition 3.2] that there exists a sequence $(z_k)_{k=0}^{\infty}$ of points in $I(g)$ such that

(i) $z_k \to z$ as $k \to \infty$;

(ii) $|g^n| \to \infty$ uniformly on $\{z\} \cup \{z_k\}_{k=0}^{\infty}$;

(iii) each $z_k$ is contained in an unbounded connected closed set $A_k$ on which $g^n$ tends to infinity uniformly.

As we just saw, $\tilde{\theta}(z_k)$ is uniquely determined by the values of $\theta$ on $J_{\geq \rho}(f)$. Since $\tilde{\theta}$ is continuous on $\{z\} \cup \{z_k\}_{k=0}^{\infty}$ by assumption, we see that the value $\tilde{\theta}(z)$ is also uniquely determined, completing the proof. □

**Definition 8.24.** Let $f \in \mathcal{B}$ be postsingularly bounded, and let $g \underset{\text{qc}}{\overset{\infty}{\approx}} f$ be of disjoint type. Any bijection $\tilde{\theta} \colon I(g) \to I(f)$ that satisfies $(b)$ and $(c)$ and that agrees with a generalised Böttcher map on $J_{\geq \rho}(g)$ for sufficiently large $\rho$ is called a *natural bijection* (for $f$ and $g$).

## 8.8 Docile functions

Recall that the Julia set of a polynomial $f$ of degree $d$ is locally connected if and only if the Böttcher map $\theta \colon I(g) \to I(f)$, where $g(z) = z^d$ extends continuously to $J(g) = \partial \mathbb{D}$. The following notion, introduced in [5], is a natural generalisation of this concept to the setting of transcendental entire functions.

**Definition 8.25.** Let $f \in \mathcal{B}$ be postsingularly bounded, let $g \underset{\text{qc}}{\overset{\infty}{\approx}} f$ be of disjoint type, and let $\theta \colon I(g) \to I(f)$ be a natural bijection for $f$ and $g$.
We say that $f$ is *docile* if $\theta$ extends continuously to a map $J(g) \cup \{\infty\} \to J(f) \cup \{\infty\}$.

We saw that any two different natural bijections, say $\theta \colon I(g) \to I(f)$ and $\theta_1 \colon I(g_1) \to I(f)$ differ by pre-composition with a quasiconformal conjugacy between $g$ and $g_1$ on their respective Julia sets. Hence docility is independent of the choice of $\theta$.

Similar as for polynomials, if $f$ is docile then its dynamics can be completely described in terms of the disjoint-type function $g$ and the semiconjugacy $\theta$. The following result, which is [5, Proposition 2.3], gives a number of elementary properties of docile functions.

**Proposition 8.26.** *Suppose that $f$ is docile. Let $g \underset{\text{qc}}{\overset{\infty}{\approx}} f$ be of disjoint type, and let $\theta \colon J(g) \cup \{\infty\} \to J(f) \cup \{\infty\}$ be the continuous extension of a natural bijection for $f$ and $g$. Then the following hold.*

*(i) $\theta(\infty) = \infty$.*

(ii) $\theta$ is surjective onto $J(f) \cup \{\infty\}$.

(iii) $\theta(J(g)) = J(f)$ and $f \circ \theta = \theta \circ g$ on $J(g)$.

(iv) $\theta \colon I(g) \to I(f)$ is a homeomorphism and $\theta^{-1}(I(f)) = I(g)$. In particular, $I(f)$ and $A(f)$ each have uncountably many connected components.

(v) For each component $C$ of $J(g)$ the map $\theta \colon C \cup \{\infty\} \to \theta(C) \cup \{\infty\}$ is a homeomorphism.

The statements about $I(f)$ and $A(f)$ in (iv) follow from from Propositions 7.56 and 8.9.

It is known that many polynomials have locally connected Julia sets, and likewise we expect that there are large classes of docile transcendental entire functions. The following is proved in [5, Theorem 1.3].

**Theorem 8.27.** *Let $f \in \mathcal{B}$ and suppose that the following hold:*

(a) *there is $D \geq 2$ such that all critical points of $f$ in $J(f)$ have degree at most $D$;*

(b) *$J(f)$ contains no asymptotic values of $f$;*

(c) *the set $P(f) \cap J(f)$ is finite;*

(d) *the set $\overline{\operatorname{sing}(f^{-1})} \cap F(f)$ is compact.*

*Then $f$ is docile.*

A function $f \in \mathcal{B}$ satisfying properties (c) and (d) is called *geometrically finite*; see the paper by Mihaljević [187]. If it also satisfies (a) and (b), it is called *strongly geometrically finite*. Although docility was defined in [5], Theorem 8.27 had previously been established implicitly in some special cases. For hyperbolic entire functions (which are exactly the geometrically finite functions for which $P(f) \cap J(f) = \emptyset$), this was done in [221, §5]. For *strongly subhyperbolic* entire functions (which are those strongly geometrically finite functions for which $P(f) \cap F(f)$ is compact), the theorem was proved by Mihaljević [188]. The latter class includes cosine maps where both critical values of $f$ are strictly preperiodic. In this case, the disjoint-type function $g$ can be chosen as the function (2.1) from §2.1. In particular, as for $g$, the Julia set $J(f) = \mathbb{C}$ can be expressed as a union of hairs and their endpoints, a fact that had previously been shown by Schleicher [253].

By Proposition 8.26 (iv), we obtain the following conclusion, which was mentioned in §7.5 in connection with Question 7.58.

**Corollary 8.28.** *Let $f$ be strongly geometrically finite. Then $I(f)$ and $A(f)$ have uncountably many connected components.*

Theorem 8.27 and its predecessors are proved by constructing a suitably expanding metric on a neighbourhood of $J(f)$. For polynomials, local connectivity of the Julia set is often proved by using the "Yoccoz puzzle"; see, e.g., Milnor's article [190]. It seems likely that docility could be established for a wider class of transcendental entire functions by adapting this technique; this is a direction that deserves further research.

If $f \in \mathcal{B}$ is a postsingularly bounded transcendental entire function for which the inverse has a logarithmic singularity over some point in the Julia set, then $f$ cannot be docile. In particular, the only docile maps in the exponential family (2.4) are those for which there is an attracting or parabolic orbit. (See also [220, Theorem 1.1].) It would be interesting to develop a weaker notion of docility that also applies to dynamically tame functions with asymptotic values in the Julia set.

Probably condition $(a)$ is essential in Theorem 8.27, but we are not aware of any examples that show this. It appears likely that, when this condition is dropped but the other conditions in the theorem are satisfied, docility of the function $f$ will depend on how quickly the degree of the critical points is allowed to grow.

**Question 8.29.** Is there a geometrically finite entire function $f$ such that $J(f)$ contains no asymptotic values of $f$, but $f$ is not docile?

Is there an example of such a function for which the degree of critical points in $J(f)$ is unbounded, but $f$ is docile?

So far we have only defined docility for functions with bounded postsingular set. Let $f \in \mathcal{B}$ and suppose that $P(f)$ is unbounded. Let $g \underset{\mathrm{qc}}{\overset{\infty}{\approx}} f$ be of disjoint type as before. A priori, there is no natural bijection between $I(g)$ and $I(f)$. However, when $f$ does not have any escaping critical or asymptotic values, we may still ask whether there exists a map $\theta$ that satisfies the conclusions of Proposition 8.26, and agrees with a generalised Böttcher map on $J_{\geq \rho}(g)$ for some $\rho$. In particular, we may ask the following.

**Question 8.30.** Are there a transcendental entire function $f \in \mathcal{B}$ and a disjoint-type function $g \underset{\mathrm{qc}}{\overset{\infty}{\approx}} f$ for which there exists a continuous surjection $\theta \colon J(g) \cup \{\infty\} \to J(f) \cup \{\infty\}$ such that

(a) $P(f) = \mathbb{C}$;

(b) $\theta(\infty) = \infty$ and $\theta^{-1}(\infty) = \{\infty\}$;

(c) $f \circ \theta = \theta \circ g$ on $J(g)$;

(d) $\theta \colon I(g) \to I(f)$ is a homeomorphism;

(e) for each component $C$ of $J(g)$, the map $\theta \colon C \cup \{\infty\} \to \theta(C) \cup \{\infty\}$ is a homeomorphism;

(f) there is $\rho > 0$ such that $\theta$ agrees with a generalised Böttcher map on $J_{\geq \rho}(g)$?

If $I(f)$ does contain a critical value, or an asymptotic value corresponding to a logarithmic singularity, then the global structure of the escaping set of $f$ and that of a disjoint-type function $g \underset{\mathrm{qc}}{\approx} f$ will be different from each other. For example, $I(f)$ may be connected (recall Theorem 7.52 and 7.49). Hence we cannot expect to obtain a semiconjugacy between the two function as in the case of docility. Nonetheless, as previously mentioned in connection with Theorem 7.49, Pardo-Simón [209] was able to describe the global dynamics of a large class of functions $f$ with escaping critical orbits using the dynamics of a disjoint-type function $g \underset{\mathrm{qc}}{\approx} f$.

# 9 Measure and dimension of the escaping set

The literature on the measure and dimension of escaping sets is very large. Here we can discuss only a small part of it. We note, however, that the articles of Kotus and Urbański [164] and Stallard [272] survey many results on the dimensions of Julia sets, escaping sets and other dynamically relevant sets up to 2008. Moreover, results concerning the measure or dimension of Julia sets often yield corresponding results for escaping sets. We have already mentioned this for the results of McMullen [185].

The Lebesgue measure, Hausdorff dimension and packing dimension of a subset $E$ of $\mathbb{C}$ (or $\mathbb{R}^n$) are denoted by $\operatorname{meas} E$, $\dim_\mathrm{H} E$ and $\dim_\mathrm{P} E$, respectively; see Falconer's book [130] for a thorough treatment of these concepts of dimension. Here we only note that we always have $\dim_\mathrm{H} E \leq \dim_\mathrm{P} E$.

## 9.1 The dimension of the escaping set

The following result is an immediate consequence of Theorem 4.13 or Corollary 7.25.

**Theorem 9.1.** *Let $f$ be a transcendental entire function. Then $\dim_\mathrm{H} I(f) \geq 1$.*

In fact, Theorem 7.24 and also Theorem 5.5 imply that $\dim_\mathrm{H} A(f) \geq 1$.

As a converse to Theorem 9.1, we have the following result by Rempe and Stallard [230, Corollary 1.2].

**Theorem 9.2.** *For each $d \in [1,2]$ there exists a transcendental entire function $f$ such that $\dim_\mathrm{H} I(f) = d$.*

The main contribution of [230] is the case $d = 1$, for which the authors employ Theorem 8.21. The case $1 < d < 2$ is based on a function constructed by Stallard [271] who used it to prove that for all such $d$ there exists a transcendental entire function $f$ with $\dim_\mathrm{H} J(f) = d$. The case $d = 2$ follows from McMullen's Theorems 5.9 or 5.10.

The functions constructed in the proof of Theorem 9.2 are in $\mathcal{B}$.

**Question 9.3** (Albrecht and Bishop, 2020)**.** Given $d \in [1,2]$, does there exist $f \in \mathcal{S}$ such that $\dim_\mathrm{H} I(f) = d$?

Albrecht and Bishop [4] have shown that for every $\varepsilon > 0$ there exists $f \in \mathcal{S}$ such that $\dim_\mathrm{H} J(f) < 1 + \varepsilon$. By Theorem 6.1 we also have $\dim_\mathrm{H} I(f) < 1 + \varepsilon$ for such $f$. The most interesting case in Question 9.3 is the case $d = 1$.

The lower bound 1 for the Hausdorff dimension of $I(f)$ can be improved if $f$ is in $\mathcal{B}$ and satisfies some growth restriction. This is shown by the following result of Bergweiler, Karpińska and Stallard [68, Theorem 1.2].

**Theorem 9.4.** *Let $f \in \mathcal{B}$ and $q \geq 1$. Suppose that, for each $\varepsilon > 0$, there exists $r_\varepsilon > 0$ such that*
$$|f(z)| \leq \exp\bigl(\exp\bigl((\log|z|)^{q+\varepsilon}\bigr)\bigr) \quad \text{for } |z| \geq r_\varepsilon. \tag{9.1}$$
*Then*
$$\dim_H I(f) \geq 1 + \frac{1}{q}. \tag{9.2}$$

The examples in the paper by Stallard [271] already quoted show that the estimate (9.2) is best possible; that is, she constructed examples of functions $f \in \mathcal{B}$ for which we have equality in (9.2) while (9.1) holds for every $\varepsilon > 0$.

**Question 9.5.** *Is the estimate (9.2) best possible also for $f \in \mathcal{S}$?*

A corollary of Theorem 9.4 is the following result obtained previously by Barański [25, Theorem A] and Schubert [257].

**Theorem 9.6.** *Let $f \in \mathcal{B}$ be of finite order. Then $\dim_H I(f) = 2$.*

Barański and Schubert only stated that $\dim_H J(f) = 2$, but their proofs yield the stronger result that $\dim_H I(f) = 2$. (Recall that $I(f) \subset J(f)$ for $f \in \mathcal{B}$ by Theorem 6.1.) Theorem 9.6 applies in particular to the exponential family, so it generalises McMullen's Theorem 5.10.

In the proofs of the above and many other estimates of the Hausdorff dimension of the escaping set, the idea is to construct a subset of the escaping set that can be written as an intersection of (dynamically defined) nested sets. We describe this only briefly. Consider, for $k \in \mathbb{N}$, a collection $\mathcal{E}_k$ of pairwise disjoint compact subsets of $\mathbb{R}^n$ such that every element of $\mathcal{E}_{k+1}$ is contained in a unique element of $\mathcal{E}_k$, while every element of $\mathcal{E}_k$ contains at least one element of $\mathcal{E}_{k+1}$. Let $E_k$ be the union of all elements of $\mathcal{E}_k$. We want to estimate the dimension of $E := \bigcap_{k=1}^{\infty} E_k$.

It seems plausible that $E$ has large dimension if $E_{k+1}$ has a large density in the elements of $\mathcal{E}_k$. Here the *density* $\operatorname{dens}(X, Y)$ of a measurable set $X$ in a measurable set $Y$ of positive measure is defined by
$$\operatorname{dens}(X, Y) := \frac{\operatorname{meas}(X \cap Y)}{\operatorname{meas}(Y)}.$$

However, a large density alone is not enough to achieve large dimension. The following result is due to McMullen [185, Proposition 2.2]. Here $\operatorname{diam} A$ denotes the (Euclidean) diameter of a set $A$.

**Lemma 9.7.** *Let $\mathcal{E}_k$, $E_k$ and $E$ be as above. Suppose that $(\Delta_k)$ and $(d_k)$ are sequences of positive real numbers such that if $A \in \mathcal{E}_k$, then*
$$\operatorname{dens}(E_{k+1}, A) \geq \Delta_k \tag{9.3}$$

*and*
$$\operatorname{diam} A \leq d_k.$$

*Then*
$$\dim_H E \geq n - \limsup_{k\to\infty} \frac{\sum_{j=1}^{k}|\log \Delta_j|}{|\log d_k|}.$$

*Proof of Theorem 9.6.* Let $F$ be the function obtained from the logarithmic change of variable as described in §6.1. With $W$ and $R$ as there we put $x_0 := \inf\{\operatorname{Re} z : z \in W\}$ and define $h \colon (x_0, \infty) \to (\log R, \infty)$,

$$h(x) := \max_{y \in \mathbb{R}} \operatorname{Re} F(x + iy).$$

Then, by the maximum modulus principle and in analogy with the Hadamard three circles theorem, $h$ is increasing and convex. Thus $h$ is differentiable except possibly on a countable set. For $x > x_0$ we choose $z_x \in W$ such that $\operatorname{Re} z_x = x$ and $h(x) = \operatorname{Re} F(z_x)$. (If there is more than one such point $z_x$, we can take any of them.) It is not difficult to see that if $h$ is differentiable at $x$, then

$$h'(x) = F'(z_x).$$

Given $w \in \mathbb{C}$ with $\operatorname{Re} w > 0$, let $S(w)$ be the square with centre $w$ and sidelength $\operatorname{Re} w$. The sets $\mathcal{E}_k$ will be defined recursively in such a way that if $A \in \mathcal{E}_k$, then $F^k(A) = S(F^k(z_x))$ for some $x > x_0$. Let $Q(x)$ be the component of $F^{-1}(S(F(z_x)))$ that contains $z_x$. Note that $F$ is $2\pi i$-periodic, so not only $Q(x)$ but also all translates of $Q(x)$ by multiples of $2\pi i$ are mapped to $S(F(z_x))$.

We start by putting $\mathcal{E}_0 = \{S(w_0)\}$ for some $w_0$ with large real part. We then choose certain points $x_1, \ldots, x_m \in (\operatorname{Re} w_0/2, 3\operatorname{Re} w_0/2)$ and define $\mathcal{E}_1$ as the set of all translates of $Q(x_1), \ldots, Q(x_m)$ by multiples of $2\pi i$ that are contained in the square $S(w_0)$. We will explain later how $x_1, \ldots, x_m$ are chosen. Suppose now that $\mathcal{E}_k$ has been defined and let $A \in \mathcal{E}_k$. Then $F^k(A) = S(F^k(z_a))$ for some $a > x_0$. Again we choose points $x_1, \ldots, x_m \in (\operatorname{Re} F^k(z_a)/2, 3\operatorname{Re} F^k(z_a)/2)$ and consider the collection $\mathcal{X}$ of all translates of the $Q(x_j)$ by multiples of $2\pi i$ that are contained in $S(F^k(z_a))$. We then put

$$\mathcal{E}_{k+1}(A) := \{F^{-k}(B) \cap A : B \in \mathcal{X}\} \quad \text{and} \quad \mathcal{E}_{k+1} := \bigcup_{A \in \mathcal{E}_k} \mathcal{E}_{k+1}(A).$$

To apply Lemma 9.7 we have to estimate the density of $E_{k+1}$ in an element $A$ of $\mathcal{E}_k$. The Koebe distortion theorem yields that this density differs only by a bounded factor from the density of $F^k(E_{k+1} \cap A) = \bigcup_{B \in \mathcal{X}} B$ in $F^k(A) = S(F^k(z_a))$. To estimate the latter density, we have to estimate the size of the elements of $\mathcal{X}$; that is, we have to estimate the size of the $Q(x_j)$.

In order to do so we note that since $f$ has finite order, there exists $\mu > 0$ such that $|f(z)| \leq \exp(|z|^\mu)$ for large $|z|$. In terms of $F$ this means that $\operatorname{Re} F(z) \leq \exp(\mu \operatorname{Re} z)$ if $\operatorname{Re} z$ is large; that is, $\log h(x) \leq \mu x$ for large $x$. An elementary growth lemma (see, e.g., [35, Lemma 3]) now yields that given $\delta \in (0,1)$ we have $h'(x)/h(x) = (\log h)'(x) \leq \mu/\delta$ for $x$ in some subset $T$ of $(x_0, \infty)$ of lower density at least $1 - \delta$; that is,

$$\liminf_{x \to \infty} \frac{\operatorname{meas}(T \cap (x_0, x])}{x} \geq 1 - \delta.$$

Suppose now that $x \in T$. Since $D(F(z_x), \operatorname{Re} F(z_x)/2) \subset S(F(z_x))$, Koebe's one quarter theorem yields that

$$\begin{aligned}
Q(x) &\supset D\left(z_x, \frac{1}{4}|(F^{-1})'(F(z_x))| \cdot \frac{1}{2}\operatorname{Re} F(z_x)\right) \\
&= D\left(z_x, \frac{\operatorname{Re} F(z_x)}{8|F'(z_x)|}\right) = D\left(z_x, \frac{h(x)}{8h'(x)}\right) \supset D\left(z_x, \frac{\delta}{8\mu}\right).
\end{aligned} \tag{9.4}$$

So $Q(x)$ contains a disc of definite size. Similarly, the Koebe distortion theorem shows that there exists $C > 0$ such that

$$Q(x) \subset D(z_x, C). \tag{9.5}$$

We now explain how the points $x_1, \ldots, x_m \in (\operatorname{Re} w/2, 3\operatorname{Re} w/2)$ are chosen, where $w = w_0$ or $w = F^k(z_a)$. We may suppose that $\operatorname{Re} w$ is large. We choose points $x_1, \ldots, x_m \in [\operatorname{Re} w/2 + C, 3\operatorname{Re} w/2 - C] \cap T$ in such a way that $x_j + 2C < x_{j+1}$ for $1 \leq j \leq m-1$. It follows from (9.5) that if $j \neq k$, then the translates of $Q(x_j)$ by multiples of $2\pi i$ are disjoint from the translates of $Q(x_k)$ by such multiples.

Since $T$ has lower density at least $1 - \delta$, where we may take for example $\delta = 1/4$, we can achieve that the number $m$ of these points $x_j$ satisfies $m \geq \eta \operatorname{Re} w$ for some $\eta > 0$. For each $j \in \{1, \ldots, m\}$ there are at least $\operatorname{Re} w/(4\pi)$ of the $2\pi i$-translates of $Q(x_j)$ contained in $S(w)$. Overall there are at least $\eta(\operatorname{Re} w)^2/(4\pi)$ translates of the $Q(x_j)$ in $S(w)$. Thus by (9.4) they form a set of density at least $\eta\delta/(32\mu\pi)$ in $S(w)$. As explained above, the Koebe distortion theorem now yields that there exists $\Delta > 0$ such that $\operatorname{dens}(E_{k+1}, A) \geq \Delta$ for all $A \in \mathcal{E}_k$. Thus (9.3) holds with $\Delta_k = \Delta$.

It is not difficult to see that the diameters of the sets in $\mathcal{E}_k$ decrease so rapidly with $k$ that Lemma 9.7 yields that $\dim_H E = 2$. The construction shows that $\operatorname{Re} F^k(z) \to \infty$ for $z \in E$. Thus $\exp E \subset I(f)$. Since we also have $\dim_H(\exp E) = 2$, the conclusion follows. $\square$

The proof of Theorem 9.4 uses similar ideas. Here $h'(x)/h(x)$ may tend to $\infty$ as $x \to \infty$, but we still get an upper bound on a large set. Specifically, we have $h'(x)/h(x) \leq x^{q-1+\varepsilon}$ on a set of density 1. For details we refer to [63].

As the proofs sketched above use the logarithmic change of variables, they work only for functions in class $\mathcal{B}$.

**Question 9.8.** Is the hypothesis that $f \in \mathcal{B}$ necessary in Theorem 9.4? In other words, does (9.1) imply (9.2) for all transcendental entire functions $f$?

A weaker form of this question is the following.

**Question 9.9.** Is the hypothesis that $f \in \mathcal{B}$ necessary in Theorem 9.6? In other words, do we have $\dim_H I(f) = 2$ for every transcendental entire function $f$ of finite order?

Bergweiler and Karpińska [67] showed that the answer to Question 9.9 is positive if the growth of $f$ is sufficiently regular. Sixsmith [264] showed that this is the case for functions of genus 0 for which the zeros are contained in a sector of opening angle less than $\pi$. If the answer to Question 9.9 is negative, one may still ask the following still weaker question.

**Question 9.10.** Is there a growth rate such that $\dim_H I(f) = 2$, or at least $\dim_H I(f) > 1$, for every transcendental entire function $f$ growing slower than this rate?

We note that Bishop [83] has constructed transcendental entire functions of arbitrarily slow growth for which the Julia set has Hausdorff dimension 1. These functions have multiply connected Fatou components. Thus, by Theorem 6.8, their escaping sets have positive measure and hence in particular Hausdorff dimension 2.

There are various papers where the Hausdorff dimension of certain subsets of $I(f)$ is considered. Here we mention [71, 267] about slow escaping points, [30, 155] about points escaping at a certain rate, and [31, 58] about escaping points in the boundary of attracting periodic basins

If $0 < \lambda < 1/e$, then – as explained at the beginning of §2.3 – the Fatou set $F(\lambda e^z)$ is connected and consists of the basin of an attracting fixed point. The Julia set $J(\lambda e^z)$ is a union of hairs, which except possibly for their endpoints are contained in $I(\lambda e^z)$. (In fact, the hairs without endpoints are in $A(\lambda e^z)$.) Let $C_\lambda$ be the set of endpoints of the hairs. Karpińska [154] proved the following surprising result.

**Theorem 9.11.** *Let $0 < \lambda < 1/e$. Then $\dim_H(J(\lambda e^z) \setminus C_\lambda) = 1$.*

Note that $\dim_H J(\lambda e^z) = 2$ by Theorem 5.10. Thus we have the seemingly paradoxical situation that the set of hairs has dimension 2, but the set of hairs without endpoints has dimension 1.

Barański [26] showed that this result holds much more generally: The function $\lambda e^z$ may be replaced by an entire function in $\mathcal{B}$ that has finite order and is of disjoint type.

Rippon and Stallard [237, Theorem 1.1] showed that the Julia set of an entire function in $\mathcal{B}$ has packing dimension 2. The following theorem by Bergweiler [53, Theorem 1.1] generalises this result.

**Theorem 9.12.** *Let $f$ be a transcendental entire function satisfying*

$$\liminf_{r \to \infty} \frac{\log \log M(r, f)}{\log \log r} = \infty. \tag{9.6}$$

*If $f$ has no multiply connected Fatou component, then $\dim_P(I(f) \cap J(f)) = 2$.*

To see that the hypotheses of this theorem are satisfied if $f \in \mathcal{B}$, we recall that, as already noted after Question 7.10, a function in $\mathcal{B}$ has lower order at least $1/2$. In particular, (9.6) holds if $f \in \mathcal{B}$. Moreover, if $f \in \mathcal{B}$, then $f$ has no multiply connected Fatou component, since such a component would be in $I(f)$ by Theorem 6.8, but $I(f) \subset J(f)$ by Theorem 6.1.

**Corollary 9.13.** *Let $f$ be a transcendental entire function. If $f$ satisfies (9.6), then $\dim_P I(f) = 2$.*

It is not known whether the hypothesis (9.6) is essential here.

**Question 9.14.** Do we have $\dim_P I(f) = 2$ for every transcendental entire function $f$?

We note that the examples of Bishop [83] already quoted also satisfy $\dim_P J(f) = 1$.

## 9.2 The measure of the escaping set

Intersections of nested sets as described in the previous section have been used not only to estimate the Hausdorff dimension of escaping sets from below, but also to prove that certain escaping sets have positive measure. McMullen used the following simple lemma [185, Proposition 2.1] in his proof that the escaping set of a trigonometric function has positive measure (Theorem 5.9).

**Lemma 9.15.** *Let $\mathcal{E}_k$, $E_k$ and $E$ be as before Lemma 9.7 and let $(\Delta_k)$ be a sequence of positive numbers such that (9.3) holds. If*

$$\prod_{k=1}^{\infty} \Delta_k > 0, \tag{9.7}$$

*then $E$ has positive measure.*

Note that (9.7) is equivalent to

$$\sum_{k=1}^{\infty} (1 - \Delta_k) < \infty.$$

There are a number of generalisations of McMullen's Theorem 5.9. We will discuss some of them. Most of the proofs consist of constructing nested subsets of the escaping set to which (at least implicitly) Lemma 9.15 is then applied.

The following result is due to Sixsmith [265, Corollary 1.1].

**Theorem 9.16.** *Let $q \geq 3$, $\omega := \exp(2\pi i/q)$, $a_1, \ldots, a_q \in \mathbb{C} \setminus \{0\}$ and*

$$f(z) := \sum_{j=1}^{q} a_j \exp(\omega^j z).$$

*Then $A(f)$, $I(f)$ and $J(f)$ are spiders' webs of positive measure.*

Sixsmith [265, Remark 1] remarked that, more generally, the same conclusion holds if $f$ has the form

$$f(z) = \sum_{j=1}^{q} a_j \exp(b_j z), \tag{9.8}$$

where $a_j, b_j \in \mathbb{C} \setminus \{0\}$ for $1 \leq j \leq q$,

$$\arg b_j < \arg b_{j+1} < \arg b_j + \pi \tag{9.9}$$

for $1 \leq j \leq q-1$ and

$$\arg b_1 < \arg b_q - \pi. \tag{9.10}$$

Here we choose the value of the argument to lie in $[0, 2\pi)$. The conditions (9.9) and (9.10) ensure that there exists a finite subset $X$ of $[0, 2\pi)$ such that if $\theta \in [0, 2\pi) \setminus X$, then $|f(re^{i\theta})| \to \infty$ as $r \to \infty$.

On the other hand, if $\arg b_1 < \arg b_2 < \cdots < \arg b_q$ and if $\arg b_{k+1} > \arg b_k + \pi$ for some $k$ or $\arg b_1 > \arg b_q - \pi$, then there exists an open subinterval $Y$ of $[0, 2\pi)$ such that $|f(re^{i\theta})| \to 0$ as $r \to \infty$ uniformly for $\theta \in Y$.

**Question 9.17.** Let $f$ be of the form (9.8) and suppose that there is a sector $S$ such that $|f(z)| \to 0$ as $|z| \to \infty$, $z \in S$. Does this imply that $\operatorname{meas} I(f) = 0$? Or at least that $\operatorname{meas} A(f) = 0$?

Schubert [258] generalised McMullen's Theorem 5.9 in a different direction. He showed that for the sine function the complement of the escaping set in a vertical strip of finite width has finite area. The result extends to the functions $\sin(\alpha z + \beta)$ and $P(e^z)/e^z$ with a polynomial $P$, and it also holds for the complement of the fast escaping set; we refer to papers by Schleicher [254, Theorem 4] and by Zhang and Yang [289, Theorem 1.1].

Hemke [145, Theorem 5.1] showed that if $Q_1, Q_2$ and $P$ are polynomials with $\deg P \geq 3$, and

$$f(z) := Q_1(z) \exp(P(z)) + Q_2(z) \exp(-P(z)),$$

then $\operatorname{meas}(\mathbb{C} \setminus I(f)) < \infty$. For example, this applies to $f(z) = \sin(z^3)$. Figure 16 shows the Julia sets of $\sin z$, $\sin z^2$ and $\sin z^3$. (Recall that these functions are in the Speiser class so that for each of these functions the Julia set is the closure of the escaping set by Theorem 6.2. Since the pictures do not distinguish between a set and its closure, we can also consider them as pictures of the escaping sets.) The range shown is $|\operatorname{Re} z| \leq 4$ and $|\operatorname{Im} z| \leq 4$.

Wolff [285, Theorem 1.2] gave very precise conditions implying that $\operatorname{meas}(\mathbb{C} \setminus I(f)) < \infty$ for a function $f$ of the form

$$f(z) = \sum_{j=1}^{q} Q_j(z) \exp\left(b_j z^d + P_j(z)\right),$$

where $Q_j$ and $P_j$ are polynomials, $d \geq 3$, $\deg P_j < d$ and $b_j \in \mathbb{C} \setminus \{0\}$.

Functions like $f(z) = \sin(z^3)$ provide examples where $0 < \operatorname{meas}(\mathbb{C} \setminus I(f)) < \infty$ and $0 < \operatorname{meas} F(f) = \operatorname{meas}(\mathbb{C} \setminus J(f)) < \infty$ Wolff [287] also constructed an example $f$ where $0 < \operatorname{meas} I(f) < \infty$ and $0 < \operatorname{meas} J(f) < \infty$. (She did not state the claim about $I(f)$, but it follows from her construction.)

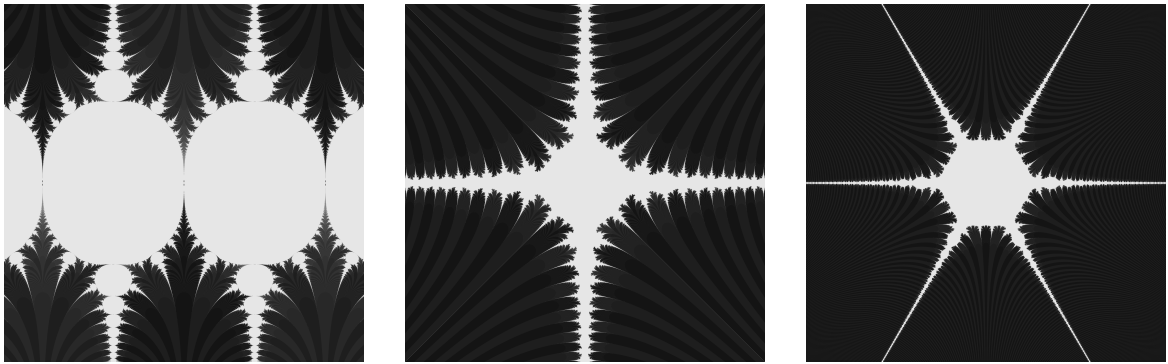

Figure 16: The Julia sets of $\sin z$, $\sin z^2$ and $\sin z^3$.

Another condition guaranteeing that $I(f)$ has positive measure was given by Aspenberg and Bergweiler [7]. It was sharpened by Cui [98], who then proved that this sharpened condition is best possible. To formulate the condition we note first that the Denjoy–Carleman–Ahlfors theorem [138, Chapter 5, Theorem 1.4] implies that if $f \in \mathcal{B}$ and if $\{z\colon |f(z)| > R\}$ has $N$ connected components, with $N \in \mathbb{N}$, then
$$\log \log M(r,f) \geq N \log r - O(1)$$
as $r \to \infty$. (If $N \geq 2$, then this also holds if $f \notin \mathcal{B}$.) The question considered in [7, 98] is whether, for suitable $\varepsilon(r)$ tending to 0, the condition
$$\log \log M(r,f) \leq (N + \varepsilon(r)) \log r + O(1)$$
implies that $I(f)$ has positive measure. Essentially, it is shown that under suitable regularity conditions on $\varepsilon(r)$ this is the case if
$$\sum_{k=1}^{\infty} \varepsilon(\exp^k 0) < \infty,$$
while $I(f)$ may have measure zero if this sum diverges.

Bergweiler and Chyzhykov [56] considered entire functions of completely regular growth in the sense of Levin and Pfluger [170, Chapter III]. Under certain additional hypothesis (see [56] for the precise statement) they show that if the indicator function (see Levin's book [170, §I.15] for the definition) of $f$ is positive, then $I(f) \cap J(f)$ has positive measure. Since functions of order less than $1/2$ that are of completely regular growth have a positive indicator, this leads to the following question.

**Question 9.18.** Do we have $\operatorname{meas} I(f) > 0$ for every transcendental entire function $f$ of order less than $1/2$?

If the answer to this question is negative, we can ask the following sharpening of Question 9.10.

**Question 9.19.** Is there a growth rate such that $\operatorname{meas} I(f) > 0$ for every transcendental entire function $f$ growing slower than this rate?

The idea behind Questions 9.10 and 9.19 is that if $f$ is an entire function of small growth, then $|f(z)|$ is large outside small neighbourhoods of the zeros. As an example of a result of this type we mention a theorem of Hayman [141, Theorem 1] who showed that if $\log M(r,f) = O((\log r)^2)$ as $r \to \infty$, then $\log |f(z)| \sim \log M(|z|,f)$ as $|z| \to \infty$ outside a set of discs subtending angles at the origin whose sum is finite.

Functions that are large outside small neighbourhoods of the zeros are also said to have the *pits effect*. The paper by Littlewood and Offord quoted in §7.4 says that certain random entire functions have the pits effect. Offord [205] showed that certain gap series also have the pits effect.

**Question 9.20.** Let $(\lambda_n)$ be an increasing sequence of nonnegative integers and let

$$f(z) = \sum_{n=0}^{\infty} a_n z^{\lambda_n}$$

be an entire function. Are there are conditions on the $\lambda_n$ that ensure that $\dim_H I(f) = 2$ or even $\operatorname{meas} I(f) > 0$?

**Question 9.21.** Let $f$ be a random entire function in the sense of Littlewood and Offord. Do we have $\dim_H I(f) = 2$ or even $\operatorname{meas} I(f) > 0$ almost surely?

One may restrict to special classes of functions in the last two questions, say functions of finite order (as in the paper by Littlewood and Offord), or functions displaying some regularity in their growth. The random entire functions considered by Littlewood and Offord are so-called Rademacher entire functions, but the question is of interest also for other random entire functions such as Gaussian or Steinhaus entire functions; see, e.g., [171, 196] for recent results on random entire functions.

Questions 9.10 and 9.18–9.21 may also be asked with $I(f)$ replaced by $A(f)$, or by other subsets of $I(f)$ such as the quite fast escaping set $Q(f)$ defined in (5.22).

Further conditions implying that the escaping or the fast escaping set have positive measure were obtained by Bergweiler [55]. For example, the results given there imply that if $f$ has the form (9.8) with polynomials $a_j$ and with $b_j$ satisfying (9.9) and (9.10), then $\operatorname{meas} A(f) > 0$. Moreover, this paper gives criteria for the escaping set of a Poincaré function to have positive measure. For example, this applies to the function displayed in the right picture of Figure 15.

Having discussed a number of conditions implying that the escaping set has positive measure, we now turn to conditions yielding that it has measure zero. The first result of this type is from the paper of McMullen already mentioned. He showed [185, Theorem 1.3] that if $E_\lambda(z) = \lambda e^z$ has an attracting periodic point, then $J(E_\lambda)$ and hence $I(E_\lambda)$ have measure zero.

In his proof, McMullen introduced the notion of thinness at infinity: A measurable subset $A$ of $\mathbb{C}$ is called *thin at infinity* if there exists $R > 0$ and $\varepsilon > 0$ such that the density of $A$ in every disc of radius $R$ is less than $1 - \varepsilon$.

McMullen [185, Proposition 7.3] proved the following result.

**Theorem 9.22.** *Let $f$ be a hyperbolic transcendental entire function and let $E$ be a measurable completely invariant subset of $J(f)$. If $E$ is thin at infinity, then* $\operatorname{meas} E = 0$.

As hyperbolic entire functions are in $\mathcal{B}$, we can apply this with $E = I(f)$ by Theorem 6.1.

Stallard [270, Theorem B] weakened the hypothesis of Theorem 9.22. Instead of requiring that $f$ is hyperbolic she assumed only that $\operatorname{dist}(P(f), J(f)) > 0$ and obtained the same conclusion. As an example where her result applies Stallard considers the function $f(z) = z + 1 - e^z$. Note that via $z \mapsto -z$ this function $f$ is conjugate to the function $F_2(z) = z - 1 + e^{-z}$ discussed at the beginning of §6.5. It follows that the Julia set of the function $f_2$ shown in Figure 13 has measure zero.

A further weakening of the hypothesis of Theorem 9.22 is due to Wolff [286]. She introduced the additional concept of *uniform thinness of a set $A$ at a set $B$*, meaning that there exist $\varepsilon > 0$ and $\delta > 0$ such that $\operatorname{dens}(A, D(z, |z - z_0|)) < 1 - \varepsilon$ for all $z_0 \in B$ and all $z \in D(z_0, \delta)$. A simplified version of her result [286, Theorem 1.3] is as follows.

**Theorem 9.23.** *Let $f$ be a transcendental meromorphic function. If $J(f)$ is thin at infinity and uniformly thin at $P(f)$, then $J(f)$ has measure zero.*

As a specific example where her result applies Wolff considers, for $c \in \mathbb{C}$, the function

$$f_c(z) := z - e^{z^2} \left( \int_0^z \exp(-t^2) \, dt + c \right),$$

which arises from the application of Newton's method to

$$g_c(z) := \int_0^z \exp(-t^2) \, dt + c.$$

She shows that if 0 is contained in a cycle of attracting basins, then $f_c$ satisfies the hypotheses of Theorem 9.23.

Figure 17 shows the Julia sets of $f_c$ for $c = 0$ (left), $c \approx 1.60336 + 2.26711i$ (middle) and $c \approx -1.21138 - 0.99433i$ (right). Clearly $g_0(0) = 0$, and hence 0 is a superattracting fixed point of $f_0$. For $c \approx 1.60336 + 2.26711i$, the point 0 is mapped to a zero of $g_c$ by $f_c$, and for $c \approx -1.21138 - 0.99433i$ (right), it is a superattracting periodic point of $f_c$ of period 2.

In all three pictures, the Fatou component containing 0 and its images and preimages are shown in white, while attracting basins that do not contain 0 are shown in various shades of grey. In the left picture, the range shown is $|\operatorname{Re} z| \leq 4$ and $|\operatorname{Im} z| \leq 5$. In the two other pictures, it is $|\operatorname{Re} z| \leq 3.2$ and $|\operatorname{Im} z| \leq 4$.

It follows from results of Bergweiler [46] that if 0 is contained in an attracting basin of $f_c$, then $F(f_c)$ consists of attracting basins only. Thus $I(f_c) \subset J(f_c)$. Hence in the examples considered not only $J(f_c)$ but also $I(f_c)$ has measure zero. We note that [46, 286] actually cover, more generally, Newton's method for functions $g$ of the form $g(z) = \int_0^z p(t) \exp q(t) dt + c$ where $p$ and $q$ are polynomials and $c \in \mathbb{C}$. (So this also covers $g(z) := e^z - 1$ and its Newton function $F_2(z) = z - 1 + e^{-z}$ discussed above.)

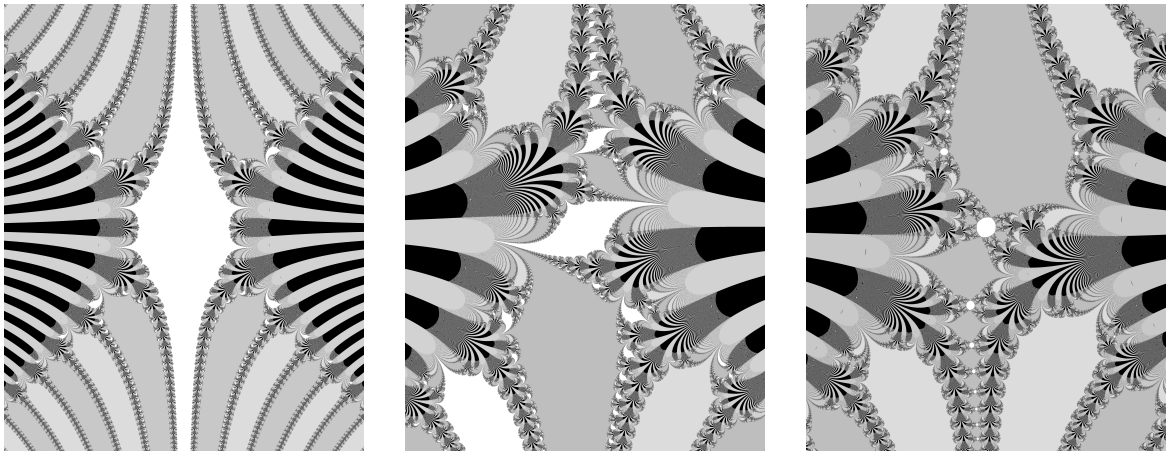

Figure 17: The Julia sets of the Newton functions $f_c$.

Eremenko and Lyubich [116, §7] gave another criterion implying that the escaping set has measure zero. Let $R > 0$ and let
$$\Theta_R(r, f) := \operatorname{meas}\{t \in [0, 2\pi] \colon |f(re^{it})| < R\}. \tag{9.11}$$

Their result [116, Theorem 7] is as follows.

**Theorem 9.24.** *Let $f \in \mathcal{B}$ and $R > 0$. Suppose that*
$$\liminf_{r \to \infty} \frac{1}{\log r} \int_1^r \Theta_R(t, f) \frac{dt}{t} > 0. \tag{9.12}$$

*Then* $\operatorname{meas} I(f) = 0$.

The condition (9.12) is clearly satisfied if $\Theta_R(r, f)$ is bounded from below by a positive constant. For example, this is the case for sufficiently large $R$ if $f$ is bounded in a sector. This does not answer Question 9.17, however, since the functions considered there need not be in $\mathcal{B}$.

The integral in (9.12) is related to the Denjoy–Carleman–Ahlfors theorem; cf. the books by Goldberg and Ostrovskii [138, p. 176] and by Tsuji [279, p. 117, p. 236]. Eremenko and Lyubich [116, Proposition 4] note that it follows from this theorem that (9.12) holds if $f$ has finite order and its inverse $f^{-1}$ has a logarithmic singularity over some finite value. Cui [98, Theorem 1.4] showed that in this situation the hypothesis of finite order can be relaxed slightly.

**Theorem 9.25.** *Let $f \in \mathcal{B}$ and suppose that the inverse $f^{-1}$ has a logarithmic singularity over some finite value. Suppose also that $A \colon [1, \infty) \to \mathbb{R}$ is an increasing function such that $A(r) \leq \log r$ and*
$$\log \log M(r, f) \leq A(r) \log r$$
*for large $r$ and*
$$\sum_{k=1}^{\infty} \frac{1}{A(\exp^k 0)} = \infty. \tag{9.13}$$

*Then* meas $I(f) = 0$.

Cui deduces this theorem from another result showing that meas $I(f) = 0$ under a suitable hypothesis on the function $\Theta_R(r, f)$ defined by (9.11). In his proof, he first considers functions of disjoint type, and then applies Theorem 8.8.

**Question 9.26.** Is condition (9.13) best possible in Theorem 9.25?

# 10 Beyond entire functions

## 10.1 Quasiregular maps

Quasiregular maps are a natural generalisation of holomorphic functions to higher dimensions. Informally, these are functions that map infinitesimal balls to infinitesimal ellipses with bounded eccentricity. For the formal definition, let $G$ be a domain in $\mathbb{R}^d$ with $d \geq 2$ and let $f \colon G \to \mathbb{R}^d$ be continuous. Then $f$ is called *quasiregular* if $f$ is absolutely continuous on almost all lines parallel to the coordinate axes, with partial derivatives locally in $L^d$, and if there exists $K \geq 1$ such that $|Df(x)|^d \leq KJ_f(x)$ almost everywhere. Here $Df$ is the derivative, $|Df|$ its norm and $J_f$ the Jacobian determinant. For more details and a thorough study of quasiregular maps we refer to the book by Rickman [234]. Some general results on the dynamics of quasiregular maps can be found in [50, 70]

Let $|x|$ denote the (Euclidean) norm of a point $x \in \mathbb{R}^d$. A quasiregular map $f \colon \mathbb{R}^d \to \mathbb{R}^d$ is said to be of *polynomial type* if

$$\lim_{|x| \to \infty} |f(x)| = \infty$$

and it is said to be of *transcendental type* otherwise. Quasiregular maps of transcendental type are the analogues of transcendental entire functions.

The definition of $I(f)$ via (1.1) extends to quasiregular maps literally. And the analogue of the basic Theorem 4.1 also holds: *If $f$ is a quasiregular mappings of transcendental type, then $I(f) \neq \emptyset$.* This was proved by Bergweiler, Fletcher, Langley and Meyer [63] using the method of Domínguez. It is not known whether there is an analogue of Theorem 4.8 that was proved by Eremenko's method.

**Question 10.1.** Let $f$ be a quasiregular map of transcendental type. Let $R > 0$ and let $D$ be a connected component of $\{x \colon |f(x)| > R\}$. Does there exist $x \in D \cap I(f)$ such that $f^n(x) \in D$ for all $n \in \mathbb{N}$?

It is also shown in [63] that $I(f)$ has at least one unbounded connected component, as in Corollary 7.25 by Rippon and Stallard. On the other hand, in contrast to Theorem 7.23, $\overline{I(f)}$ may have bounded connected components for quasiregular maps of transcendental type [63, §6].

The definition of the Julia set via non-normality is not appropriate for quasiregular maps. Instead, a definition via the so-called *blowing-up property* (cf. Theorem 3.1 (f)) has been advanced in [54, 70]. However, in contrast to Theorem 4.3 we may have $J(f) \neq \partial I(f)$. On the other hand, we always have $J(f) \subset \partial I(f)$; see [70, Theorem 1.3 and Examples 7.3 and 7.4].

The fast escaping $A(f)$ can be defined for quasiregular maps of transcendental type in the same way as for entire functions. In fact, Bergweiler, Drasin and Fletcher [59, Theorem 1.3] showed that the definitions (5.3), (5.8) and (5.9) are again equivalent. It was shown by Bergweiler, Fletcher and Nicks [64] that $J(f) \subset \partial A(f)$, with equality if

$$\lim_{r \to \infty} \frac{\log \log M(r, f)}{\log \log r} = \infty.$$

**Question 10.2** (Bergweiler, Fletcher and Nicks, 2014)**.** Do we have $J(f) = \partial A(f)$ for every quasiregular map $f$ of transcendental type?

As for entire functions, escape rates are of interest not only in their own right, but they have also applications. Thus Tsantaris [276] used $A(f)$ to study permutable quasiregular maps, as was done in [66] for entire maps.

Nicks [198] showed that for quasiregular maps of transcendental type there exists arbitrarily slow escaping points. In other words, Theorem 5.14 extends literally to this setting. He also extended Theorem 5.16 to quasiregular maps.

A domain in $\mathbb{R}^d$ is called *full* if it does not have a bounded complementary component. Otherwise it is called *hollow*. Nicks and Sixsmith [204, Theorem 1] showed that if $x$ is in a full invariant connected component of the *quasi-Fatou set* $\mathbb{R}^d \setminus J(f)$, then $\log \log |f^n(x)| = O(n)$, and this estimate is best possible [204, §8]. Note that for entire functions we have the stronger estimate $\log |f^n(x)| = O(n)$, by Theorem 6.21.

Nicks and Sixsmith [202, 203] also studied hollow quasi-Fatou components. They showed that they share many properties of multiply connected Fatou components as discussed in §6.3. It remains open, however, whether there is an analogue of Theorem 6.11.

**Question 10.3.** Is a hollow quasi-Fatou component of a quasiregular map of transcendental type necessarily bounded?

Examples of bounded hollow quasi-Fatou components were constructed by Burkart, Fletcher and Nicks [90].

Nicks and Sixsmith [202, Theorem 1.6] showed that a positive answer to Question 10.3 would imply a a positive answer to Question 10.2.

Spiders' webs are defined in $\mathbb{R}^d$ as in §7.4, even though for $d \geq 3$ they do not match the picture of a spider's web that we have in mind. Analogous to Theorem 7.28 and the remarks after Definition 7.29, Nicks and Sixsmith [202, Theorem 1.3] showed that if a quasiregular map $f$ of transcendental type has a bounded hollow quasi-Fatou component, then $A_R(f)$, $A(f)$ and $I(f)$ are spiders' webs. Bergweiler, Drasin and Fletcher [59, Theorem 1.6] showed that if the minimum modulus is large in certain intervals, then the fast escaping set is a spider's web. Recall that Theorem 7.39 by Mihaljević-Brandt and Peter [189] gives a condition when the fast escaping set of a Poincaré function is a spider's web. Corresponding results for quasiregular maps were given by Fletcher [134] and Fletcher and Stoertz [135].

One important example of a quasiregular map is the so-called Zorich map. It is a quasiregular analogue of the exponential map. Bergweiler [51] showed that (for certain parameters) the escaping set of the Zorich maps also consists of hairs, and for the dimension of the hairs and their endpoints there are analogues of Theorem 5.10 and Theorem 9.11. Comdühr [96] showed that under suitable additional hypotheses the hairs are differentiable, in analogy with Viana's [282] result mentioned at the end of §7.1. Tsantaris [277, 278] studied Zorich maps in detail. Here we only mention that he proved [278, Theorem 1] that under suitable hypotheses on the parameters the escaping set of a Zorich map is connected; see Theorems 7.52 and 7.53 for such results for exponential functions.

It follows from Theorem 9.24 that meas $I(\lambda e^z) = 0$ for all $\lambda \in \mathbb{C} \setminus \{0\}$. The analogue for Zorich maps is not known.

**Question 10.4.** Let $f$ be a Zorich map and $\lambda \in \mathbb{R} \setminus \{0\}$. Do we have meas $I(\lambda f) = 0$?

The dynamics of quasiregular analogues of the trigonometric functions were studied by Bergweiler and Eremenko [61]. Vogel [283] extended results about the measure of escaping sets of trigonometric functions to such maps.

Since $I(f)$ has an unbounded connected component, Theorem 9.1 extends literally to the quasiregular setting; that is, $\dim_H I(f) \geq 1$ for every quasiregular map $f$ of transcendental type. It is not known whether this is sharp. Possibly we have $\dim_H I(f) \geq d-1$ for every quasiregular map $f \colon \mathbb{R}^d \to \mathbb{R}^d$ of transcendental type; see [74, Remark 4.7] for a discussion.

**Question 10.5** (Bergweiler and Tsantaris, 2025)**.** What is the sharp lower bound for the Hausdorff dimension of the escaping set for quasiregular maps of transcendental type?

Analogous to Theorem 7.24, it was shown by Bergweiler, Drasin and Fletcher [59, Theorem 1.2] that every connected component of the fast escaping set of a quasiregular maps of transcendental type is unbounded. For such a map $f$ we thus have $\dim_H A(f) \geq 1$. This bound is sharp. In fact, Bergweiler and Tsantaris [74, Theorem 1.1] showed that for every $d \in [1,3]$ there exists a quasiregular map $f \colon \mathbb{R}^3 \to \mathbb{R}^3$ with $\dim_H A(f) = d$.

The following question is an analogue of Questions 9.10 and 9.19.

**Question 10.6.** Is there a growth rate such that for every quasiregular map $f \colon \mathbb{R}^d \to \mathbb{R}^d$ of transcendental type growing slower than this rate we have $\dim_H I(f) = d$, or even meas $I(f) > 0$?

Again we may also ask this with $I(f)$ replaced by $A(f)$.

We conclude this section by mentioning that without quasiregularity much less can be said. This follows from results of Short and Sixsmith [260] who studied the escaping set for functions $f \colon \mathbb{R}^d \to \mathbb{R}^d$ that are merely continuous.

## 10.2 Meromorphic functions

Let $f\colon \mathbb{C} \to \widehat{\mathbb{C}}$ be a transcendental meromorphic function. Here $\widehat{\mathbb{C}} := \mathbb{C} \cup \{\infty\}$ is the Riemann sphere. Then for $z \in f^{-n}(\infty)$ the iterate $f^{n+1}(z)$ is not defined. Nevertheless an iteration theory very similar to the one described in §3 can be developed for such functions; see Bergweiler's survey [45]. Here the Fatou set is defined as the set of all points that have a neighbourhood where the iterates *are defined and form a normal family*.

The escaping set $I(f)$ consists of all points $z$ for which $f^n(z)$ is defined for all $n \in \mathbb{N}$ and $|f^n(z)| \to \infty$ as $n \to \infty$. Thus points in $O^-(\infty)$ are not in $I(f)$. Theorem 4.1 holds literally:

**Theorem 10.7.** *Let $f$ be a transcendental meromorphic function. Then $I(f) \neq \emptyset$.*

This result is due to Domínguez [109, Theorems G and H]. In the case that $f$ has only finitely many poles she used her method described in §4.4. Alternatively, the result can be proved by Eremenko's method in this case.

*Proof of Theorem 10.7 if $f$ has infinitely many poles.* For $R > 0$ we put

$$B(R) := \{z \in \mathbb{C}\colon |z| \geq R\} \cup \{\infty\}.$$

For each pole $p$ there exists $r(p) \in (0,1)$ and $R(p) > 0$ such that $f(D(p, r(p))) \supset B(R(p))$. Let now $p_0$ be pole and define a sequence $(p_n)_{n \geq 0}$ of poles recursively by choosing $|p_n| > R(p_{n-1}) + 1$. With $r_n := r(p_n)$ and $R_n := R(p_n)$ we then have

$$\overline{D}(p_n, r_n) \subset B(R_{n-1}) \subset f(D(p_{n-1}, r_{n-1})) \subset f^n(D(p_0, r_0)).$$

With $K_0 := \overline{D}(p_0, r_0)$ and $K_n := f^{-n}(\overline{D}(p_n, r_n)) \cap K_0$ we find that the sets $K_n$ are non-empty, compact, and satisfy

$$K_{n+1} \subset K_n \subset \overline{D}(p_0, r_0) \quad \text{and} \quad f^n(K) = \overline{D}(p_n, r_n).$$

Thus

$$K := \bigcap_{n=0}^{\infty} K_n \neq \emptyset.$$

Choosing the sequence $(p_n)$ such that $|p_n| \to \infty$ we have $K \subset I(f)$. $\square$

The proof allows us to construct points in $I(f)$ that escape to $\infty$ as fast as we like. We only have to choose the sequence $(p_n)$ such that it tends to infinity fast. So there is not really an analogue of the fast escaping set for meromorphic functions with infinitely many poles.

Domínguez [109, Theorems I and J] showed that Theorems 4.3 and 4.6 also hold for meromorphic functions; that is, we have $J(f) = \partial I(f)$ and $I(f) \cap J(f) \neq \emptyset$.

In contrast to entire functions, $I(f)$ and $J(f)$ may be totally disconnected for meromorphic functions. For example, if $0 < \lambda < 1$ and $f(z) = \lambda \tan z$, then $J(f) = \overline{I(f)}$ is a Cantor set.

Another difference between entire and meromorphic functions concerns the Hausdorff dimension of escaping sets. Here we only refer to the papers [8, 9, 57, 69, 99, 184], which show that the results for meromorphic functions are quite different from those for entire functions mentioned in §9.1.

On the other hand, as in the proof of Theorem 10.7, results and arguments for functions with finitely many poles are often similar to those for entire functions; see, e.g., Rippon and Stallard [238]. More generally, this is the case for functions with direct tracts; see, e.g., the paper by Bergweiler Rippon and Stallard [72] already discussed in §4.3 and that by Xuan and Zheng [288]. These two papers also study *escape rates within a tract*, and in particular a set $A(f, D)$ of points which escape fast in a direct tract $D$ of a meromorphic function $f$.

The iteration theory of functions meromorphic in the plane has been generalised to the study of functions meromorphic outside a small set, with various notions of "small"; see, e.g., [21, 86, 110]. For a function meromorphic in $\widehat{\mathbb{C}} \setminus E$ and $e \in E$, Baker, Domínguez and Herring [21, §4] considered the set

$$I(f, e) = \left\{ z \in \widehat{\mathbb{C}} \setminus E \colon f^n(z) \to e \right\}$$

and showed that under suitable hypothesis we again have $I(f, e) \cap J(f) \neq \emptyset$ and $J(f) = \partial I(f, e)$. Domínguez, Montes de Oca and Sienra [110, §6] studied properties of the (larger) set of points for which the $\omega$-limit set is in $O^-(E)$.

We mentioned holomorphic self-maps of the punctured plane already in §5.2. Martí-Pete [179] showed that for a holomorphic map $f \colon \mathbb{C} \setminus \{0\} \to \mathbb{C} \setminus \{0\}$ with essential singularities at 0 and $\infty$ not only both sets $I(f, 0)$ and $I(f, \infty)$ are non-empty, but there are also points with $\omega$-limit set $\{0, \infty\}$ that "oscillate" between 0 and $\infty$ in any prescribed way; see [179] for the precise statement.

# List of open questions

Here we collect all numbered open questions presented throughout this survey.

**Question 5.12.** Let $f(z) := \lambda e^z$, with $\lambda \in \mathbb{C} \setminus \{0\}$ such that $0 \in I(f)$. Which growth rate of the sequence $(f^n(0))$ ensures that $f$ is not ergodic? Which growth rate ensures that (5.18) holds?

**Question 5.15.** Let $D$ be a direct tract of $f$ and let $(a_n)$ be a sequence of positive numbers that tends to $\infty$. Does there exist $z \in D$ with $f^n(z) \in D$ for all $n \in \mathbb{N}$ and $|f^n(z)| \leq a_n$ for all large $n$?

**Question 5.17.** Let $D$ be a direct tract of $f$ such that $\partial D$ contains a curve tending to infinity. Let $(a_n)$ be a sequence of positive numbers that tends to $\infty$. Does there exist $z \in D$ and $C > 0$ such that $f^n(z) \in D$ for all $n \in \mathbb{N}$ and $a_n \leq |f^n(z)| \leq Ca_n$ for all large $n$?

**Question 6.17** (Bishop, 2018)**.** How smooth can the boundary of a multiply connected Fatou component of a transcendental entire function be? Can it consist of $C^2$ curves? Or $C^\infty$ curves? Or analytic curves?

**Question 6.18.** Can the boundary of a multiply connected Fatou component of a transcendental entire function consist of circles?

**Question 6.19.** For which bounded, multiply connected domains $D$ does there exist a transcendental entire function $f$ such that $D$ is a Fatou component of $f$?

**Question 6.20.** Does there exist a transcendental entire function, a multiply connected wandering domain $U$ of $f$, and a bounded connected component $K$ of $\mathbb{C} \setminus U$ such that $K$ is not a singleton and $\text{int}(K) \cap J(f) = \emptyset$?

If so, can $K$ be chosen such that $\text{int}(K) \neq \emptyset$?

**Question 6.29.** How sparse can the set of singularities be for an entire function with an invariant Baker domain?

**Question 6.30** (Fleischmann, 2008)**.** Does there exist a transcendental entire function $f$ with an invariant Baker domain for which there exist $K > 1$ and a sequence $(r_n)$ tending to infinity such that $\text{sing}(f^{-1}) \cap \text{ann}(r_n, Kr_n) = \emptyset$ for all $n$?

**Question 6.32.** Can the conditions $|p_{n+1}/p_n| \to 1$ and $\text{dist}(p_n, U) = o(|p_n|)$ in Theorem 6.31 be improved?

**Question 6.33** (Baker, 2001)**.** Do there exist transcendental entire functions of arbitrarily slow growth for which all Fatou components are simply connected and at least one Fatou component is in the escaping set?

**Question 6.34.** Do there exist transcendental entire functions of arbitrarily slow growth for which all Fatou components are simply connected and escaping?

**Question 6.37** (Martí-Pete, Rempe and Waterman, 2025)**.** Let $D$ be a bounded simply connected wandering domain of a transcendental entire function, and let $W$ be the unbounded connected component of $\mathbb{C} \setminus \overline{D}$. Is it true that $\partial W = \partial D$?

**Question 6.38** (Rippon and Stallard, 2011)**.** Let $f$ be a transcendental entire function and let $U$ be an escaping Fatou component of $f$. Does $\partial U$ contain escaping points?

**Question 6.40** (Rippon and Stallard, 2018)**.** Let $f$ be a transcendental entire function and let $U$ be a Baker domain of $f$. Do we have $\partial U \cap I(f) \neq \emptyset$?

**Question 6.43.** Is the hypothesis that $f|_U \colon U \to U$ has finite degree necessary in Theorem 6.42?

**Question 6.47** (Bishop, 2014)**.** Let $U$ be a simply connected escaping wandering domain of a transcendental entire function. Let $\varphi \colon \partial \mathbb{D} \to U$ be biholomorphic. Does the set of all $\xi \in \partial \mathbb{D}$ for which $\varphi(\xi) \notin I(f)$ have capacity zero?

**Question 7.1.** Is there a transcendental entire function $f$ such that $I(f)$ is a $G_{\delta\sigma}$ set? If yes, can $f$ be chosen such that $f \in \mathcal{B}$, or even $f \in \mathcal{S}$?

**Question 7.7.** What growth conditions on $f \in \mathcal{B}$ ensure that $f$ has the strong Eremenko property (that is, every point of $I(f)$ can be connected to infinity by a curve in $I(f)$)?

**Question 7.8.** Does every transcendental entire function of finite order have the strong Eremenko property?

**Question 7.10.** What growth conditions on $f \in \mathcal{B}$ ensure that $I(f)$ contains an arc?

**Question 7.12.** Let $f \in \mathcal{B}$ have the strong Eremenko property. Is $f$ necessarily criniferous?

**Question 7.14** (Pardo-Simón and Rempe, 2023)**.** Is there a disjoint-type entire function $f$ such that every point in $I(f)$ is contained in an unbounded connected set on which $f^n \to \infty$ uniformly, but such that $f$ is not criniferous?

**Question 7.15.** Suppose that $f \in \mathcal{B}$ is criniferous and that $I(f)$ contains no critical points. Let $C$ be a path-connected component of $I(f)$. Is $C$ the continuous injective image of an open or half-open interval?

**Question 7.17.** Suppose that $f \in \mathcal{B}$ is criniferous, and that $I(f)$ contains no critical points. If $\gamma \subset I(f)$ is an arc, does $f^n|_\gamma \to \infty$ uniformly?

**Question 7.18.** Does there exist a transcendental entire function $f$ with a multiply connected Fatou component such that $I(f)$ does not contain a curve to $\infty$?

**Question 7.19** (Viana, 1988)**.** Are the hairs in the escaping sets of the exponential family analytic?

**Question 7.22.** Let $K \subset \mathbb{C}$ be a planar continuum. When does there exist a transcendental entire function $f$ such that $K \subset I(f) \cap J(f)$, and such that $f^n(K) \cap f^m(K) = \emptyset$ for $n \neq m$?

**Question 7.42.** Let $f$ be a transcendental entire function of order $\rho(f) < 1/2$.

(a) Is $I(f)$ necessarily connected?

(b) Is $I(f)$ necessarily a spider's web?

(c) Is there necessarily a bounded simply connected domain $U$ with $U \cap J(f) \neq \emptyset$ and such that $f^n|_{\partial U} \to \infty$ uniformly?

**Question 7.51.** Is $I(f)$ a Jordan spider's web when

(a) $f$ is Fatou's function (2.3)?

(b) $f = \cosh$?

(c) $f = \cos + \cosh$?

(d) $f$ is a Poincaré function with spider's web escaping set, as in Theorem 7.39?

**Question 7.54** (Rempe, 2011)**.** Let $a \in \mathbb{C}$ be such that all periodic cycles of the exponential map $f(z) = ae^z$ are repelling. Is $I(f)$ connected?

**Question 7.55** (Rippon and Stallard, 2019)**.** Does there exist a transcendental entire function $f$ such that $I(f)$ is connected and $\mathbb{C} \setminus I(f)$ contains an unbounded *closed* connected set?

**Question 7.58** (Rippon and Stallard, 2019)**.** Let $f$ be a transcendental entire function. For each of the sets $I(f)$ and $A(f)$, is it the case that it is either connected or it has uncountably many connected components?

**Question 7.61** (Rippon and Stallard, 2019)**.** In Theorem 7.60, can we replace the "unbounded connected $F_\sigma$ sets" in $A_R(f)$ by "unbounded connected closed sets" in $A_R(f)$ or even by "connected components" of $A_R(f)$?

**Question 7.62** (Rippon and Stallard, 2019)**.** In Theorem 7.60, can we replace the "dense set" of values of $R \in [R_0, \infty)$ by "all" $R \in [R_0, \infty)$?

**Question 7.63** (Rippon and Stallard, 2019)**.** Does there exist a transcendental entire function $f$ and an admissible radius $R$ for $f$ such that $A_R(f)$ is connected, but $A_R(f)$ is not a spider's web?

**Question 7.67** (Martí-Pete, Rempe and Waterman, 2025)**.** Does Eremenko's conjecture hold for all functions in $\mathcal{B}$?

**Question 7.68.** Does Eremenko's conjecture hold for all functions in $\mathcal{S}$?

**Question 7.69** (Martí-Pete, Rempe and Waterman, 2025)**.** Is Eremenko's conjecture true for functions of finite order?

**Question 8.12.** Let $f, g \in \mathcal{B}$ with $f \underset{\text{qc}}{\overset{\infty}{\sim}} g$. Is it possible to construct a global generalised Böttcher map $\theta \colon \mathbb{C} \to \mathbb{C}$ explicitly, without using holomorphic motions or parameter arguments?

**Question 8.13.** Let $f, g \in \mathcal{B}$ with $f \underset{\text{qc}}{\overset{\infty}{\sim}} g$. Does there always exist a generalised Böttcher map $\theta$ for $f$ and $g$ so that the conjugacy relation (8.5) holds also on an open neighbourhood of $J_{\geq \rho}(g)$?

**Question 8.22** (Rempe, 2009)**.** Let $f \in \mathcal{B}$. Is there a proof of the absence of $f$-invariant line fields on $I(f)$ that uses only dynamical properties of the function $f$, and avoids arguments involving holomorphic motions in parameter space?

**Question 8.29.** Is there a geometrically finite entire function $f$ such that $J(f)$ contains no asymptotic values of $f$, but $f$ is not docile?

Is there an example of such a function for which the degree of critical points in $J(f)$ is unbounded, but $f$ is docile?

**Question 8.30.** Are there a transcendental entire function $f \in \mathcal{B}$ and a disjoint-type function $g \underset{\text{qc}}{\overset{\infty}{\sim}} f$ for which there exists a continuous surjection $\theta \colon J(g) \cup \{\infty\} \to J(f) \cup \{\infty\}$ such that

(a) $P(f) = \mathbb{C}$;

(b) $\theta(\infty) = \infty$ and $\theta^{-1}(\infty) = \{\infty\}$;

(c) $f \circ \theta = \theta \circ g$ on $J(g)$;

(d) $\theta \colon I(g) \to I(f)$ is a homeomorphism;

(e) for each component $C$ of $J(g)$, the map $\theta \colon C \cup \{\infty\} \to \theta(C) \cup \{\infty\}$ is a homeomorphism;

(f) there is $\rho > 0$ such that $\theta$ agrees with a generalised Böttcher map on $J_{\geq \rho}(g)$?

**Question 9.3** (Albrecht and Bishop, 2020)**.** Given $d \in [1, 2]$, does there exist $f \in \mathcal{S}$ such that $\dim_{\text{H}} I(f) = d$?

**Question 9.5.** Is the estimate (9.2) best possible also for $f \in \mathcal{S}$?

**Question 9.8.** Is the hypothesis that $f \in \mathcal{B}$ necessary in Theorem 9.4? In other words, does (9.1) imply (9.2) for all transcendental entire functions $f$?

**Question 9.9.** Is the hypothesis that $f \in \mathcal{B}$ necessary in Theorem 9.6? In other words, do we have $\dim_{\text{H}} I(f) = 2$ for every transcendental entire function $f$ of finite order?

**Question 9.10.** Is there a growth rate such that $\dim_{\text{H}} I(f) = 2$, or at least $\dim_{\text{H}} I(f) > 1$, for every transcendental entire function $f$ growing slower than this rate?

**Question 9.14.** Do we have $\dim_P I(f) = 2$ for every transcendental entire function $f$?

**Question 9.17.** Let $f$ be of the form (9.8) and suppose that there is a sector $S$ such that $|f(z)| \to 0$ as $|z| \to \infty$, $z \in S$. Does this imply that $\operatorname{meas} I(f) = 0$? Or at least that $\operatorname{meas} A(f) = 0$?

**Question 9.18.** Do we have $\operatorname{meas} I(f) > 0$ for every transcendental entire function $f$ of order less than $1/2$?

**Question 9.19.** Is there a growth rate such that $\operatorname{meas} I(f) > 0$ for every transcendental entire function $f$ growing slower than this rate?

**Question 9.20.** Let $(\lambda_n)$ be an increasing sequence of nonnegative integers and let
$$f(z) = \sum_{n=0}^{\infty} a_n z^{\lambda_n}$$
be an entire function. Are there are conditions on the $\lambda_n$ that ensure that $\dim_H I(f) = 2$ or even $\operatorname{meas} I(f) > 0$?

**Question 9.21.** Let $f$ be a random entire function in the sense of Littlewood and Offord. Do we have $\dim_H I(f) = 2$ or even $\operatorname{meas} I(f) > 0$ almost surely?

**Question 9.26.** Is condition (9.13) best possible in Theorem 9.25?

**Question 10.1.** Let $f$ be a quasiregular map of transcendental type. Let $R > 0$ and let $D$ be a connected component of $\{x \colon |f(x)| > R\}$. Does there exist $x \in D \cap I(f)$ such that $f^n(x) \in D$ for all $n \in \mathbb{N}$?

**Question 10.2** (Bergweiler, Fletcher and Nicks, 2014)**.** Do we have $J(f) = \partial A(f)$ for every quasiregular map $f$ of transcendental type?

**Question 10.3.** Is a hollow quasi-Fatou component of a quasiregular map of transcendental type necessarily bounded?

**Question 10.4.** Let $f$ be a Zorich map and $\lambda \in \mathbb{R} \setminus \{0\}$. Do we have $\operatorname{meas} I(\lambda f) = 0$?

**Question 10.5** (Bergweiler and Tsantaris, 2025)**.** What is the sharp lower bound for the Hausdorff dimension of the escaping set for quasiregular maps of transcendental type?

**Question 10.6.** Is there a growth rate such that for every quasiregular map $f \colon \mathbb{R}^d \to \mathbb{R}^d$ of transcendental type growing slower than this rate we have $\dim_H I(f) = d$, or even $\operatorname{meas} I(f) > 0$?

# List of Symbols

# Index

Walter Bergweiler
Mathematisches Seminar
Christian-Albrechts-Universität zu Kiel
Heinrich-Hecht-Platz 6
24098 Kiel
Germany
Email address: `bergweiler@math.uni-kiel.de`

Lasse Rempe
Department of Mathematics
The University of Manchester
Manchester M13 9PL
UK
Email address: `lasse.rempe@manchester.ac.uk`